\documentclass{article}
\usepackage{amsmath}
\usepackage{amssymb}
\usepackage{amsthm}
\usepackage{cite}
\usepackage[arrow,curve,matrix]{xy}
\begin{document}
%%%%%%%%%%%%%%%%%%%%%%%%%%%%%%%%%%%%%%%%%%%%%%%%%%%%%%%%%%%%%%%%%%%%%%%
%%%%%%%%%%%%%%%%%%%%%%%%%     Macros      %%%%%%%%%%%%%%%%%%%%%%%%%%%%%
%%%%%%%%%%%%%%%%%%%%%%%%%%%%%%%%%%%%%%%%%%%%%%%%%%%%%%%%%%%%%%%%%%%%%%%
\def\e#1\e{\begin{equation}#1\end{equation}}
\def\ea#1\ea{\begin{align}#1\end{align}}
\def\eq#1{{\rm(\ref{#1})}}
\theoremstyle{plain}% default
\newtheorem{thm}{Theorem}[section]
\newtheorem{lem}[thm]{Lemma}
\newtheorem{prop}[thm]{Proposition}
\newtheorem{cor}[thm]{Corollary}
\newtheorem{prob}[thm]{Problem}
\newtheorem{quest}[thm]{Question}
\newtheorem{conj}[thm]{Conjecture}
\theoremstyle{definition}
\newtheorem{dfn}[thm]{Definition}
\newtheorem{ex}[thm]{Example}
\newtheorem{ass}[thm]{Assumption}
\newtheorem{rem}[thm]{Remark}
\newtheorem{cond}[thm]{Condition}
\def\ind{{\rm ind}}
\def\red{{\rm red}}
\def\na{{\rm na}}
\def\stk{{\rm stk}}
\def\rk{\mathop{\rm rk}\nolimits}
\def\vi{{\rm vi}}
\def\td{\mathop{\rm td}\nolimits}
\def\dim{\mathop{\rm dim}\nolimits}
\def\Ker{\mathop{\rm Ker}}
\def\cha{\mathop{\rm char}}
\def\Re{\mathop{\rm Re}}
\def\Im{\mathop{\rm Im}}
\def\Gr{\mathop{\rm Gr}\nolimits}
\def\GL{\mathop{\rm GL}\nolimits}
\def\Spec{\mathop{\rm Spec}\nolimits}
\def\Sch{\mathop{\rm Sch}\nolimits}
\def\Sta{\mathop{\rm Sta\kern .05em}\nolimits}
\def\Stab{\mathop{\rm Stab}\nolimits}
\def\Var{\mathop{\rm Var}\nolimits}
\def\Quot{\mathop{\rm Quot}\nolimits}
\def\coh{\mathop{\rm coh}}
\def\qcoh{\mathop{\rm qcoh}}
\def\vect{\mathop{\rm vect}}
\def\ch{\mathop{\rm ch}\nolimits}
\def\Hom{\mathop{\rm Hom}\nolimits}
\def\Iso{\mathop{\rm Iso}\nolimits}
\def\Aut{\mathop{\rm Aut}}
\def\End{\mathop{\rm End}}
\def\Sub{\mathop{\rm Sub}\nolimits}
\def\Mor{\mathop{\rm Mor}\nolimits}
\def\CF{\mathop{\rm CF}\nolimits}
\def\CFi{\mathop{\rm CF}\nolimits^{\rm ind}}
\def\LCF{\mathop{\rm LCF}\nolimits}
\def\dLCF{{\dot{\rm LCF}}\kern-.1em\mathop{}\nolimits}
\def\dLCFi{{\dot{\rm LCF}}\kern-.1em\mathop{}\nolimits^{\rm ind}}
\def\SF{\mathop{\rm SF}\nolimits}
\def\SFi{\mathop{\rm SF}\nolimits^{\rm ind}}
\def\SFa{\mathop{\rm SF}\nolimits_{\rm al}}
\def\SFai{\mathop{\rm SF}\nolimits_{\rm al}^{\rm ind}}
\def\uSF{\mathop{\smash{\underline{\rm SF\!}\,}}\nolimits}
\def\uSFi{\mathop{\smash{\underline{\rm SF\!}\,}}\nolimits^{\rm ind}}
\def\oSF{\mathop{\bar{\rm SF}}\nolimits}
\def\oSFa{\mathop{\bar{\rm SF}}\nolimits_{\rm al}}
\def\oSFai{{\ts\bar{\rm SF}{}_{\rm al}^{\rm ind}}}
\def\uoSF{\mathop{\bar{\underline{\rm SF\!}\,}}\nolimits}
\def\uoSFa{\mathop{\bar{\underline{\rm SF\!}\,}}\nolimits_{\rm al}}
\def\uoSFi{\mathop{\bar{\underline{\rm SF\!}\,}}\nolimits_{\rm al}^{\rm
ind}}
\def\LSF{\mathop{\rm LSF}\nolimits}
\def\uLSF{\mathop{\smash{\underline{\rm LSF\!}\,}}\nolimits}
\def\dLSF{{\dot{\rm LSF}}\kern-.1em\mathop{}\nolimits}
\def\dLSFa{{\dot{\rm LSF}}\kern-.1em\mathop{}\nolimits_{\rm al}}
\def\dLSFai{{\dot{\rm LSF}}\kern-.1em\mathop{}\nolimits^{\rm ind}_{\rm
al}}
\def\ESF{\mathop{\rm ESF}\nolimits}
\def\Ht{{\mathcal H}{}^{\rm to}}
\def\bHt{\bar{\mathcal H}{}^{\rm to}}
\def\Hp{{\mathcal H}{}^{\rm pa}}
\def\bHp{\bar{\mathcal H}{}^{\rm pa}}
\def\Lt{{\mathcal L}{}^{\rm to}}
\def\bLt{\bar{\mathcal L}{}^{\rm to}}
\def\Lp{{\mathcal L}{}^{\rm pa}}
\def\bLp{\bar{\mathcal L}{}^{\rm pa}}
\def\Oss{\mathop{\rm Obj\kern .1em}\nolimits_{\rm ss}}
\def\Osi{\mathop{\rm Obj\kern .1em}\nolimits_{\rm si}}
\def\Ost{\mathop{\rm Obj\kern .1em}\nolimits_{\rm st}}
\def\Ext{\mathop{\rm Ext}\nolimits}
\def\id{\mathop{\rm id}\nolimits}
\def\Obj{\mathop{\rm Obj\kern .1em}\nolimits}
\def\fObj{\mathop{\mathfrak{Obj}\kern .05em}\nolimits}
\def\modA{\text{\rm mod-$A$}}
\def\modKQ{\text{\rm mod-$\K Q$}}
\def\modKQI{\text{\rm mod-$\K Q/I$}}
\def\nilKQ{\text{\rm nil-$\K Q$}}
\def\nilKQI{\text{\rm nil-$\K Q/I$}}
\def\modCQ{\text{\rm mod-$\C Q$}}
\def\nilCQ{\text{\rm nil-$\C Q$}}
\def\bs{\boldsymbol}
\def\ge{\geqslant}
\def\le{\leqslant\nobreak}
\def\pr{{\mathop{\preceq}\nolimits}}
\def\npr{{\mathop{\npreceq}\nolimits}}
\def\tl{\trianglelefteq\nobreak}
\def\ntl{\ntrianglelefteq}
\def\ps{\precsim\nobreak}
\def\ls{{\mathop{\lesssim\kern .05em}\nolimits}}
\def\ld{\lessdot}
\def\bF{{\mathbin{\mathbb F}}}
\def\N{{\mathbin{\mathbb N}}}
\def\R{{\mathbin{\mathbb R}}}
\def\Z{{\mathbin{\mathbb Z}}}
\def\Q{{\mathbin{\mathbb Q}}}
\def\C{{\mathbin{\mathbb C}}}
\def\CP{{\mathbin{\mathbb{CP}}}}
\def\K{{\mathbin{\mathbb K\kern .05em}}}
\def\KP{{\mathbin{\mathbb{KP}}}}
\def\cC{{\mathbin{\mathcal C}}}
\def\A{{\mathbin{\mathcal A}}}
\def\B{{\mathbin{\mathcal B}}}
\def\F{{\mathbin{\mathcal F}}}
\def\G{{\mathbin{\mathcal G}}}
\def\H{{\mathbin{\mathcal H}}}
\def\L{{\mathbin{\mathcal L}}}
\def\M{{\mathcal M}}
\def\O{{\mathbin{\mathcal O}}}
\def\P{{\mathbin{\mathcal P}}}
\def\T{{\mathbin{\mathcal T}}}
\def\U{{\mathbin{\mathcal U}}}
\def\fD{{\mathbin{\mathfrak D}}}
\def\fE{{\mathbin{\mathfrak E}}}
\def\fF{{\mathbin{\mathfrak F}}}
\def\fG{{\mathbin{\mathfrak G}}}
\def\fH{{\mathbin{\mathfrak H}}}
\def\fM{{\mathbin{\mathfrak M}}}
\def\fR{{\mathbin{\mathfrak R}}}
\def\fS{{\mathbin{\mathfrak S}}}
\def\fT{{\mathbin{\mathfrak T}}}
\def\fU{{\mathbin{\mathfrak U\kern .05em}}}
\def\fV{{\mathbin{\mathfrak V}}}
\def\sIp{{\smash{\sst(I,\pr)}}}
\def\sJp{{\smash{\sst(J,\pr)}}}
\def\sKp{{\smash{\sst(K,\pr)}}}
\def\sIt{{\smash{\sst(I,\tl)}}}
\def\sJt{{\smash{\sst(J,\tl)}}}
\def\sKt{{\smash{\sst(K,\tl)}}}
\def\sIl{{\smash{\sst(I,\ls)}}}
\def\sJl{{\smash{\sst(J,\ls)}}}
\def\sKl{{\smash{\sst(K,\ls)}}}
\def\dss{\de_{\rm ss}}
\def\dsi{\de_{\rm si}}
\def\dst{\de_{\rm st}}
\def\dssb{\de_{\smash{\rm ss}}^{\,\rm b}}
\def\dsib{\de_{\smash{\rm si}}^{\,\rm b}}
\def\dstb{\de_{\smash{\rm st}}^{\,\rm b}}
\def\bdss{\bar\de_{\rm ss}}
\def\bdsi{\bar\de_{\rm si}}
\def\bdst{\bar\de_{\rm st}}
\def\bdssb{\bar\de_{\smash{\rm ss}}^{\,\rm b}}
\def\bdsib{\bar\de_{\smash{\rm si}}^{\,\rm b}}
\def\bdstb{\bar\de_{\smash{\rm st}}^{\,\rm b}}
\def\Iss{I_{\rm ss}}
\def\Isi{I_{\rm si}}
\def\Ist{I_{\rm st}}
\def\Issb{I_{\smash{\rm ss}}^{\,\rm b}}
\def\Isib{I_{\smash{\rm si}}^{\,\rm b}}
\def\Istb{I_{\smash{\rm st}}^{\,\rm b}}
\def\Jsib{J_{\smash{\rm si}}^{\,\rm b}}
\def\Jstb{J_{\smash{\rm st}}^{\,\rm b}}
\def\Mss{{\mathcal M}_{\rm ss}}
\def\Msi{{\mathcal M}_{\rm si}}
\def\Mst{{\mathcal M}_{\rm st}}
\def\Mssb{{\mathcal M}_{\rm ss}^{\,\rm b}}
\def\Msib{{\mathcal M}_{\rm si}^{\,\rm b}}
\def\Mstb{{\mathcal M}_{\rm st}^{\,\rm b}}
\def\al{\alpha}
\def\be{\beta}
\def\ga{\gamma}
\def\de{\delta}
\def\bde{\bar\delta}
\def\io{\iota}
\def\ep{\epsilon}
\def\bep{\bar\epsilon}
\def\la{\lambda}
\def\ka{\kappa}
\def\th{\theta}
\def\ze{\zeta}
\def\up{\upsilon}
\def\vp{\varphi}
\def\si{\sigma}
\def\om{\omega}
\def\De{\Delta}
\def\La{\Lambda}
\def\Om{\Omega}
\def\Ga{\Gamma}
\def\Si{\Sigma}
\def\Th{\Theta}
\def\Up{\Upsilon}
\def\pd{\partial}
\def\ts{\textstyle}
\def\sst{\scriptscriptstyle}
\def\w{\wedge}
\def\sm{\setminus}
\def\bu{\bullet}
\def\op{\oplus}
\def\ot{\otimes}
\def\bigop{\bigoplus}
\def\bigot{\bigotimes}
\def\iy{\infty}
\def\ra{\rightarrow}
\def\ab{\allowbreak}
\def\longra{\longrightarrow}
\def\t{\times}
\def\ci{\circ}
\def\el{{\mathbin{\ell\kern .08em}}}
\def\ti{\tilde}
\def\ha{{\ts\frac{1}{2}}}
\def\md#1{\vert #1 \vert}
\def\bmd#1{\big\vert #1 \big\vert}
%%%%%%%%%%%%%%%%%%%%%%%%%%%%%%%%%%%%%%%%%%%%%%%%%%%%%%%%%%%%%%%%%%%%%%%
%%%%%%%%%%%%%%%%%%%%%     Text of paper    %%%%%%%%%%%%%%%%%%%%%%%%%%%%
%%%%%%%%%%%%%%%%%%%%%%%%%%%%%%%%%%%%%%%%%%%%%%%%%%%%%%%%%%%%%%%%%%%%%%%
\title{Configurations in abelian categories. IV. \\
Invariants and changing stability conditions}
\author{Dominic Joyce}
\date{}
\maketitle

\begin{abstract}
This is the last in a series on {\it configurations\/} in an abelian
category $\A$. Given a finite poset $(I,\pr)$, an
$(I,\pr)$-configuration $(\si,\io,\pi)$ is a finite collection of
objects $\si(J)$ and morphisms $\io(J,K)$ or
$\pi(J,K):\si(J)\ra\si(K)$ in $\A$ satisfying some axioms, where
$J,K$ are subsets of $I$. Configurations describe how an object $X$
in $\A$ decomposes into subobjects.

The first paper defined configurations and studied moduli spaces of
configurations in $\A$, using Artin stacks. It showed well-behaved
moduli stacks $\fObj_\A,\fM(I,\pr)_\A$ of objects and configurations
in $\A$ exist when $\A$ is the abelian category $\coh(P)$ of
coherent sheaves on a projective scheme $P$, or $\modKQ$ of
representations of a quiver $Q$. The second studied algebras of {\it
constructible functions\/} and {\it stack functions} on~$\fObj_\A$.

The third introduced {\it stability conditions\/} $(\tau,T,\le)$ on
$\A$, and showed the moduli space $\Oss^\al(\tau)$ of
$\tau$-semistable objects in class $\al$ is a constructible subset
in $\fObj_\A$, so its characteristic function $\dss^\al(\tau)$ is a
constructible function. It formed algebras $\Hp_\tau,\Ht_\tau,
\bHp_\tau,\bHt_\tau$ of constructible and stack functions on
$\fObj_\A$, and proved many identities in them.

In this paper, if $(\tau,T,\le)$ and $(\ti\tau,\ti T,\le)$ are
stability conditions on $\A$ we write $\dss^\al(\ti\tau)$ in terms
of the $\dss^\be(\tau)$, and deduce the algebras
$\Hp_\tau,\ldots,\bHt_\tau$ are independent of $(\tau,T,\le)$. We
study {\it invariants} $\Iss^\al(\tau)$ or $\Iss(I,\pr,\ka,\tau)$
`counting' $\tau$-semistable objects or configurations in $\A$,
which satisfy additive and multiplicative identities. We compute
them completely when $\A=\modKQ$ or $\A=\coh(P)$ for $P$ a smooth
curve. We also find invariants with special properties when
$\A=\coh(P)$ for $P$ a smooth surface with $K_P^{-1}$ nef, or a
Calabi--Yau 3-fold.
\end{abstract}

\section{Introduction}
\label{an1}

This is the fourth in a series of papers \cite{Joyc5,Joyc6,Joyc7} on
{\it configurations}. Given an abelian category $\A$ and a finite
partially ordered set (poset) $(I,\pr)$, we define an $(I,\pr)$-{\it
configuration} $(\si,\io,\pi)$ in $\A$ to be a collection of objects
$\si(J)$ and morphisms $\io(J,K)$ or $\pi(J,K):\si(J)\ra\si(K)$ in
$\A$ satisfying certain axioms, for~$J,K\subseteq I$.

The first paper \cite{Joyc5} defined configurations, developed their
basic properties, and studied moduli spaces of configurations in
$\A$, using the theory of Artin stacks. It proved well-behaved
moduli stacks $\fObj_\A,\fM(I,\pr)_\A$ of objects and configurations
exist when $\A$ is the abelian category $\coh(P)$ of coherent
sheaves on a projective $\K$-scheme $P$, or $\modKQ$ of
representations of a quiver $Q$. The second \cite{Joyc6} studied
algebras of {\it constructible functions} $\CF(\fObj_\A)$ and {\it
stack functions} $\SF(\fObj_\A)$ on $\fObj_\A$, motivated by {\it
Ringel--Hall algebras}.

The third paper \cite{Joyc7} studied ({\it weak\/}) {\it stability
conditions} $(\tau,T,\le)$ on $\A$, which include {\it slope
stability} on $\modKQ$ and {\it Gieseker stability} on $\coh(P)$. If
$(\tau,T,\le)$ is {\it permissible} then the moduli space
$\Oss^\al(\tau)$ of $\tau$-semistable objects $X$ in $\A$ with
$[X]=\al$ in $K(\A)$ is a {\it constructible set\/} in $\fObj_\A$,
so its characteristic function $\dss^\al(\tau)$ is a {\it
constructible function}. We used this to define interesting algebras
$\Hp_\tau, \Ht_\tau,\bHp_\tau,\bHt_\tau$ and Lie algebras $\Lp_\tau,
\Lt_\tau,\bLp_\tau,\bLt_\tau$ in $\CF(\fObj_\A)$ and
$\SF(\fObj_\A)$, and prove many identities in them.

The first goal of this paper is to understand {\it how moduli spaces
$\Oss^\al(\tau)$ transform as we change the} ({\it weak\/}) {\it
stability condition} $(\tau,T,\le)$ to $(\ti\tau,\ti T,\le)$. We
express this as the following identity in $\CF(\fObj_\A)$, which is
equation \eq{an5eq2} below:
\e
\begin{gathered}
\dss^\al(\ti\tau)=\sum_{\substack{\text{$\A$-data
$(\{1,\ldots,n\},\le,\ka):$}\\
\text{$\ka(\{1,\ldots,n\})=\al$}}} \!\!\!\!\!\!\!
\begin{aligned}[t]
S(\{1,&\ldots,n\},\le,\ka,\tau,\ti\tau)\cdot\\
&\dss^{\ka(1)}(\tau)*\dss^{\ka(2)}(\tau)*\cdots*
\dss^{\ka(n)}(\tau).
\end{aligned}
\end{gathered}
\label{an1eq}
\e
Here $S(*,\tau,\ti\tau)$ are {\it explicit combinatorial
coefficients\/} equal to 1,0 or $-1$, and `$*$' is the associative,
noncommutative multiplication on $\CF(\fObj_\A)$ studied
in~\cite{Joyc6}.

Roughly speaking, \eq{an1eq} characterizes whether $X\in\A$ is
$\ti\tau$-semistable in terms of $\tau$-semistability, via an {\it
inclusion-exclusion process} upon {\it filtrations}
$0=A_0\subset\cdots\subset A_n=X$ with $\tau$-semistable factors
$S_i=A_i/A_{i-1}$. Writing $\ka(i)=[S_i]$ in $K(\A)$, the
coefficient $S(\cdots)$ depends on the orderings of $\tau\ci\ka(i)$
and $\ti\tau\ci\ka(i)$ for $i=1,\ldots,n$ in the total orders
$(T,\le)$ and $(\ti T,\le)$, and this determines whether a
filtration is included, if $S(\cdots)=1$, or excluded,
if~$S(\cdots)=-1$.

We say that $(\ti\tau,\ti T,\le)$ {\it dominates} $(\tau,T,\le)$ if
$\tau(\al)\le\tau(\be)$ implies $\ti\tau(\al)\le\ti\tau(\be)$ for
all $\al,\be\in C(\A)$. In this case \eq{an1eq} follows easily from
the facts that each $X\in\A$ has a unique {\it Harder--Narasimhan
filtration} $0=A_0\subset\cdots\subset A_n=X$ with $S_i=A_i/A_{i-1}$
$\tau$-semistable and $\tau([S_1])>\cdots>\tau([S_n])$, and then $X$
is $\ti\tau$-semistable if and only if $\ti\tau([S_i])=\ti\tau(\al)$
for all $i$. For the general case we go via a weak stability
condition $(\hat\tau,\hat T,\le)$ dominating both $(\tau,T,\le)$ and
$(\ti\tau,\ti T,\le)$. For $\A=\coh(P)$ equation \eq{an1eq} may have
infinitely many nonzero terms, and {\it converges} in a suitable
sense.

The second goal of the paper is to study {\it systems of invariants}
of $\A$ and $(\tau,T,\le)$ which `count' $\tau$-semistable objects
or configurations in $\A$. Obviously there are many ways of doing
this, so we need to decide what are the most interesting, or useful,
ways to define invariants. The point of view we take is that the
invariants are interesting if they satisfy {\it natural identities},
and the more identities the better. Such identities are powerful
tools for calculating the invariants in examples, as we shall see.

We obtain our invariants $\Iss^\al(\tau)$ by applying some invariant
$\Up$ of constructible sets in Artin stacks to the moduli spaces
$\Oss^\al(\tau)$. If $\Up$ takes values in a vector space $\La$ and
$\Up(S\cup T)=\Up(S)+\Up(T)$ when $S,T$ are constructible sets with
$S\cap T=\emptyset$, then constructible function identities such as
\eq{an1eq} and many more in \cite{Joyc7} translate to {\it additive
identities} on the invariants, such as transformation laws under
change of (weak) stability condition. This is our basic assumption,
and holds for Euler characteristics, virtual Poincar\'e polynomials,
and so on.

If also $\La$ is a commutative algebra and $\Up$ has multiplicative
properties such as $\Up(S\t T)=\Up(S)\Up(T)$ or $\Up([X/G])=\Up(X)
\Up(G)^{-1}$ then we can derive extra {\it multiplicative
identities} on the $\Iss^\al(\tau)$. Usually these multiplicative
identities require extra conditions on the groups $\Ext^i(X,Y)$ for
$X,Y\in\A$ and $i>1$, and can be interpreted in terms of morphisms
from a stack (Lie) algebra in $\SF(\fObj_\A)$ to some explicit
algebra, as in \cite[\S 6]{Joyc6}. These assumptions on
$\Ext^i(X,Y)$ mean that our invariants have good properties in
special cases which we focus on, namely, when $\A=\modKQ$, or
$\A=\coh(P)$ for $P$ a smooth curve, or $P$ a smooth surface with
$K_P^{-1}$ nef, or $P$ a Calabi--Yau 3-fold.

Here is an overview of the paper. Section \ref{an2} gives background
material on Artin stacks, constructible functions and stack
functions from \cite{Joyc3,Joyc4}, and \S\ref{an3} briefly reviews
the first three papers \cite{Joyc5,Joyc6,Joyc7}. Section \ref{an4}
defines and studies the {\it transformation coefficients\/}
$S(*,\tau,\ti\tau)$ appearing in \eq{an1eq}, and related
coefficients $T,U(*,\tau,\ti\tau)$; this part of the paper is wholly
{\it combinatorial}.

In \S\ref{an5} we prove \eq{an1eq}, its stack function analogue, and
transformation laws for other families of constructible and stack
functions $\dss,\bdss(I,\pr,\ka,\tau)$ and $\ep^\al,\bep^\al(\tau)$.
When $\A=\modKQ$ equation \eq{an1eq} has only finitely many terms,
but when $\A=\coh(P)$ for $P$ a projective $\K$-scheme it might have
infinitely many nonzero terms, and holds with an appropriate notion
of convergence. We show there are only finitely many nonzero terms
in \eq{an1eq} if $P$ is a smooth surface.

Section \ref{an6} studies some families of invariants
$\Iss(I,\pr,\ka,\tau)^{\sst\La}$, $\Iss^\al(\tau)^{\sst\La}$,
$J^\al(\tau)^{\sst\La^\ci}$, $J^\al(\tau)^{\sst\Om}$ taking values
in $\Q$-algebras $\La,\La^\ci,\Om$ which `count' $\tau$-semistable
objects or configurations in $\A$. We determine their transformation
laws under change of stability condition, and additive and
multiplicative identities they satisfy under conditions on
$\Ext^i(X,Y)$ for $i>1$ and~$X,Y\in\A$.

We compute the invariants completely when $\A=\modKQ$ or
$\A=\coh(P)$ for $P$ a smooth curve, recovering results of Reineke
and Harder--Narasimhan--Atiyah--Bott. We also find invariants with
special multiplicative transformation laws when $\A=\coh(P)$ for $P$
a smooth surface with $K_P^{-1}$ nef or a Calabi--Yau 3-fold. For
surfaces $P$ with $c_1(P)=0$, such as $K3$ surfaces, we define
invariants $\bar J^\al(\tau)^{\sst\La}$ which are independent of the
choice of Gieseker stability condition $(\tau,T,\le)$ on $P$. We
discuss the connection of our invariants with {\it Donaldson
invariants} of surfaces and {\it Donaldson--Thomas invariants} of
Calabi--Yau 3-folds, and make some conjectures on the existence of
invariants combining good features of the various sets of
invariants.

Finally, \S\ref{an7} suggests problems for future research:
extending the whole programme to triangulated categories, combining
the invariants in generating functions, with applications to Mirror
Symmetry, and use of esoteric kinds of stacks to weaken the
assumptions we need on~$\Ext^i(X,Y)$.

A sequel \cite{Joyc8} explains how to encode some of the invariants
we study into {\it holomorphic generating functions\/} on the
complex manifold of stability conditions. These satisfy an
interesting p.d.e., that can be interpreted as the flatness of a
connection. This will be discussed in~\S\ref{an7}.
\medskip

\noindent{\it Acknowledgements.} I would like to thank Tom
Bridgeland, Alastair King, Frances Kirwan, Andrew Kresch, Ivan
Smith, Bal\'azs Szendr{\accent"7D o}i, Richard Thomas and Burt
Totaro for useful conversations. I was supported by an EPSRC
Advanced Research Fellowship whilst writing this paper.

\section{Background material}
\label{an2}

We begin with some background material on Artin stacks,
constructible functions, `stack functions', and motivic invariants.
Sections \ref{an21}--\ref{an23} are drawn from \cite{Joyc3,Joyc4},
and \S\ref{an24} is new.

\subsection{Artin $\K$-stacks and constructible functions}
\label{an21}

Let $\K$ be an algebraically closed field. There are four main
classes of `spaces' over $\K$ used in algebraic geometry, in
increasing order of generality:
\begin{equation*}
\text{$\K$-varieties}\subset \text{$\K$-schemes}\subset
\text{algebraic $\K$-spaces}\subset \text{algebraic $\K$-stacks}.
\end{equation*}

{\it Algebraic stacks} (also known as Artin stacks) were introduced
by Artin, generalizing {\it Deligne--Mumford stacks}. For a good
introduction to algebraic stacks see G\'omez \cite{Gome}, and for a
thorough treatment see Laumon and Moret-Bailly \cite{LaMo}. We make
the convention that all algebraic $\K$-stacks in this paper are {\it
locally of finite type}, and $\K$-substacks are {\it locally
closed}.

We define the set of $\K$-{\it points} of a stack.

\begin{dfn} Let $\fF$ be a $\K$-stack. Write $\fF(\K)$ for the set of
2-isomorphism classes $[x]$ of 1-morphisms $x:\Spec\K\ra\fF$.
Elements of $\fF(\K)$ are called $\K$-{\it points}, or {\it
geometric points}, of $\fF$. If $\phi:\fF\ra\fG$ is a 1-morphism
then composition with $\phi$ induces a map of
sets~$\phi_*:\fF(\K)\ra\fG(\K)$.

For a 1-morphism $x:\Spec\K\ra\fF$, the {\it stabilizer group}
$\Iso_\K(x)$ is the group of 2-morphisms $x\ra x$. When $\fF$ is an
algebraic $\K$-stack, $\Iso_\K(x)$ is an {\it algebraic $\K$-group}.
We say that $\fF$ {\it has affine geometric stabilizers} if
$\Iso_\K(x)$ is an affine algebraic $\K$-group for all
1-morphisms~$x:\Spec\K\ra\fF$.

As an algebraic $\K$-group up to isomorphism, $\Iso_\K(x)$ depends
only on the isomorphism class $[x]\in\fF(\K)$ of $x$ in
$\Hom(\Spec\K,\fF)$. If $\phi:\fF\ra\fG$ is a 1-morphism,
composition induces a morphism of algebraic $\K$-groups
$\phi_*:\Iso_\K([x])\ra\Iso_\K\bigr(\phi_*([x])\bigr)$,
for~$[x]\in\fF(\K)$.
\label{an2def1}
\end{dfn}

The theory of {\it constructible functions} on $\K$-stacks was
developed in~\cite{Joyc3}.

\begin{dfn} Let $\fF$ be an algebraic $\K$-stack. We call
$C\subseteq\fF(\K)$ {\it constructible} if $C=\bigcup_{i\in I}
\fF_i(\K)$, where $\{\fF_i:i\in I\}$ is a finite collection of
finite type algebraic $\K$-substacks $\fF_i$ of $\fF$. We call
$S\subseteq\fF(\K)$ {\it locally constructible} if $S\cap C$ is
constructible for all constructible~$C\subseteq\fF(\K)$.

A function $f:\fF(\K)\ra\Q$ is called {\it constructible} if
$f(\fF(\K))$ is finite and $f^{-1}(c)$ is a constructible set in
$\fF(\K)$ for each $c\in f(\fF(\K))\sm\{0\}$. A function
$f:\fF(\K)\ra\Q$ is called {\it locally constructible} if
$f\cdot\de_C$ is constructible for all constructible
$C\subseteq\fF(\K)$, where $\de_C$ is the characteristic function of
$C$. Write $\CF(\fF)$ and $\LCF(\fF)$ for the $\Q$-vector spaces of
$\Q$-valued constructible and locally constructible functions
on~$\fF$.
\label{an2def2}
\end{dfn}

Following \cite[Def.s~4.8, 5.1 \& 5.5]{Joyc3} we define {\it
pushforwards} and {\it pullbacks} of constructible functions along
1-morphisms. We need~$\cha\K=0$.

\begin{dfn} Let $\K$ have {\it characteristic zero} and $\fF$ be
an algebraic $\K$-stack with affine geometric stabilizers and
$C\subseteq\fF(\K)$ be constructible. Then \cite[Def.~4.8]{Joyc3}
defines the {\it na\"\i ve Euler characteristic} $\chi^\na(C)$ of
$C$. It is called {\it na\"\i ve} as it takes no account of
stabilizer groups. For $f\in\CF(\fF)$, define
\begin{equation*}
\chi^\na(\fF,f)=\ts\sum_{c\in f(\fF(\K))\sm\{0\}}c\,\chi^\na
\bigl(f^{-1}(c)\bigr)\quad\text{in $\Q$.}
\end{equation*}

Let $\fF,\fG$ be algebraic $\K$-stacks with affine geometric
stabilizers, and $\phi:\fF\ra\fG$ a representable 1-morphism. Then
for any $x\in\fF(\K)$ we have an injective morphism
$\phi_*:\Iso_\K(x)\!\ra\!\Iso_\K\bigl(\phi_*(x)\bigr)$ of affine
algebraic $\K$-groups. Define $m_\phi:\fF(\K)\!\ra\!\Z$ by
$m_\phi(x)\!=\!\chi\bigl(\Iso_\K(\phi_*(x))/\phi_*(\Iso_\K(x))
\bigr)$. For $f$ in $\CF(\fF)$, define $\CF^\stk(\phi)f:
\fG(\K)\!\ra\!\Q$ by $\CF^\stk(\phi)f(y)\!=\!\chi^\na\bigl(\fF,
m_\phi\cdot f\cdot \de_{\smash{\phi_*^{-1}(y)}}\bigr)$ for $y$ in
$\fG(\K)$, where $\de_{\smash{\phi_*^{-1}(y)}}$ is the
characteristic function of $\phi_*^{-1}(\{y\})\subseteq\fG(\K)$.
Then $\CF^\stk(\phi):\CF(\fF) \ra\CF(\fG)$ is a $\Q$-linear map
called the {\it stack pushforward}.

Let $\th:\fF\ra\fG$ be a finite type 1-morphism. If $C\subseteq
\fG(\K)$ is constructible then so is
$\th_*^{-1}(C)\subseteq\fF(\K)$. It follows that if $f\in\CF(\fG)$
then $f\ci\th_*$ lies in $\CF(\fF)$. Define the {\it pullback\/}
$\th^*:\CF(\fG)\ra\CF(\fF)$ by $\th^*(f)= f\ci\th_*$. It is a linear
map.
\label{an2def3}
\end{dfn}

Here \cite[Th.s~5.4, 5.6 \& Def.~5.5]{Joyc3} are some properties of
these.

\begin{thm} Let\/ $\fE,\fF,\fG,\fH$ be algebraic $\K$-stacks with
affine geometric stabilizers, and\/ $\be:\fF\ra\fG$, $\ga:\fG\ra\fH$
be $1$-morphisms. Then
\ea
\CF^\stk(\ga\ci\be)&=\CF^\stk(\ga)\ci\CF^\stk(\be):\CF(\fF)\ra\CF(\fH),
\label{an2eq1}\\
(\ga\ci\be)^*&=\be^*\ci\ga^*:\CF(\fH)\ra\CF(\fF),
\label{an2eq2}
\ea
supposing $\be,\ga$ representable in \eq{an2eq1}, and of finite type
in \eq{an2eq2}. If
\e
\begin{gathered}
\xymatrix{
\fE \ar[r]_\eta \ar[d]^\th & \fG \ar[d]_\psi \\
\fF \ar[r]^\phi & \fH }
\end{gathered}
\quad
\begin{gathered}
\text{is a Cartesian square with}\\
\text{$\eta,\phi$ representable and}\\
\text{$\th,\psi$ of finite type, then}\\
\text{the following commutes:}
\end{gathered}
\quad
\begin{gathered}
\xymatrix@C=35pt{
\CF(\fE) \ar[r]_{\CF^\stk(\eta)} & \CF(\fG) \\
\CF(\fF) \ar[r]^{\CF^\stk(\phi)} \ar[u]_{\th^*} & \CF(\fH).
\ar[u]^{\psi^*} }
\end{gathered}
\label{an2eq3}
\e
\label{an2thm1}
\end{thm}

As discussed in \cite[\S 3.3]{Joyc3} for the $\K$-scheme case,
equation \eq{an2eq1} is {\it false} for algebraically closed fields
$\K$ of characteristic $p>0$. In \cite[\S 5.3]{Joyc3} we extend
Definition \ref{an2def3} and Theorem \ref{an2thm1} to {\it locally
constructible functions} $\LCF(\fF)$. The main differences are in
which 1-morphisms must be of {\it finite type}.

\subsection{Stack functions}
\label{an22}

{\it Stack functions} are a universal generalization of
constructible functions introduced in \cite{Joyc4}. Here
\cite[Def.~3.1]{Joyc4} is the basic definition. Throughout $\K$ is
algebraically closed of arbitrary characteristic, except when we
specify~$\cha\K=0$.

\begin{dfn} Let $\fF$ be an algebraic $\K$-stack with affine
geometric stabilizers. Consider pairs $(\fR,\rho)$, where $\fR$ is a
finite type algebraic $\K$-stack with affine geometric stabilizers
and $\rho:\fR\ra\fF$ is a representable 1-morphism. We call two
pairs $(\fR,\rho)$, $(\fR',\rho')$ {\it equivalent\/} if there
exists a 1-isomorphism $\io:\fR\ra \fR'$ such that $\rho' \ci\io$
and $\rho$ are 2-isomorphic 1-morphisms $\fR\ra\fF$. Write
$[(\fR,\rho)]$ for the equivalence class of $(\fR,\rho)$. If
$(\fR,\rho)$ is such a pair and $\fS$ is a closed $\K$-substack of
$\fR$ then $(\fS,\rho\vert_\fS)$,
$(\fR\sm\fS,\rho\vert_{\fR\sm\fS})$ are pairs of the same kind.
Define $\SF(\fF)$ to be the $\Q$-vector space generated by
equivalence classes $[(\fR,\rho)]$ as above, with for each closed
$\K$-substack $\fS$ of $\fR$ a relation
\begin{equation*}
[(\fR,\rho)]=[(\fS,\rho\vert_\fS)]+[(\fR\sm\fS,\rho\vert_{\fR\sm\fS})].
\end{equation*}
Define $\uSF(\fF)$ in the same way, but without requiring the
1-morphisms $\rho$ to be representable.
Then~$\SF(\fF)\subseteq\uSF(\fF)$.
\label{an2def4}
\end{dfn}

In \cite[Def.~3.2]{Joyc4} we relate $\CF(\fF)$ and~$\SF(\fF)$.

\begin{dfn} Let $\fF$ be an algebraic $\K$-stack with affine
geometric stabilizers and $C\subseteq\fF(\K)$ be constructible. Then
$C=\coprod_{i=1}^n\fR_i(\K)$, for $\fR_1,\ldots,\fR_n$ finite type
$\K$-substacks of $\fF$. Let $\rho_i:\fR_i\ra\fF$ be the inclusion
1-morphism. Then $[(\fR_i,\rho_i)]\in\SF(\fF)$. Define
$\bde_C=\ts\sum_{i=1}^n[(\fR_i,\rho_i)]\in\SF(\fF)$. We think of
this stack function as the analogue of the characteristic function
$\de_C\in\CF(\fF)$ of $C$. Define a $\Q$-linear map
$\io_\fF:\CF(\fF)\ra\SF(\fF)$ by $\io_\fF(f)=\ts\sum_{0\ne c\in
f(\fF(\K))}c\cdot\bde_{f^{-1}(c)}$. For $\K$ of characteristic zero,
define a $\Q$-linear map $\pi_\fF^\stk:\SF(\fF)\ra\CF(\fF)$ by
\begin{equation*}
\pi_\fF^\stk\bigl(\ts\sum_{i=1}^nc_i[(\fR_i,\rho_i)]\bigr)=
\ts\sum_{i=1}^nc_i\CF^\stk(\rho_i)1_{\fR_i},
\end{equation*}
where $1_{\fR_i}$ is the function 1 in $\CF(\fR_i)$. Then
\cite[Prop.~3.3]{Joyc4} shows $\pi_\fF^\stk\ci\io_\fF$ is the
identity on $\CF(\fF)$. Thus, $\io_\fF$ is injective and
$\pi_\fF^\stk$ is surjective. In general $\io_\fF$ is far from being
surjective, and $\SF(\fF)$ is much larger than~$\CF(\fF)$.
\label{an2def5}
\end{dfn}

All the operations of constructible functions in \S\ref{an21} extend
to stack functions.

\begin{dfn} Define {\it multiplication} `$\,\cdot\,$' on $\uSF(\fF)$ by
\begin{equation*}
[(\fR,\rho)]\cdot[(\fS,\si)]=[(\fR\t_{\rho,\fF,\si}\fS,\rho\ci\pi_\fR)].
\end{equation*}
This extends to a $\Q$-bilinear product $\uSF(\fF)\t\uSF(\fF)\ra
\uSF(\fF)$ which is commutative and associative, and $\SF(\fF)$ is
closed under `$\,\cdot\,$'. Let $\phi:\fF\!\ra\!\fG$ be a 1-morphism
of algebraic $\K$-stacks with affine geometric stabilizers. Define
the {\it pushforward\/} $\phi_*:\uSF(\fF)\!\ra\!\uSF(\fG)$~by
\begin{equation*}
\phi_*:\ts\sum_{i=1}^mc_i[(\fR_i,\rho_i)]\longmapsto
\ts\sum_{i=1}^mc_i[(\fR_i,\phi\ci\rho_i)].
\end{equation*}
If $\phi$ is representable then $\phi_*$ maps $\SF(\fF)\!\ra\!
\SF(\fG)$. For $\phi$ of finite type, define {\it pullbacks}
$\phi^*:\SF(\fG)\!\ra\!\SF(\fF)$,
$\phi^*:\uSF(\fG)\!\ra\!\uSF(\fF)$~by
\begin{equation*}
\phi^*:\ts\sum_{i=1}^mc_i[(\fR_i,\rho_i)]\longmapsto
\ts\sum_{i=1}^mc_i[(\fR_i\t_{\rho_i,\fG,\phi}\fF,\pi_\fF)].
\end{equation*}
The {\it tensor product\/} $\ot\!:\!\SF(\fF)\!\t\!\SF(\fG)\!\ra
\!\SF(\fF\!\t\!\fG)$ or
$\uSF(\fF)\!\t\!\uSF(\fG)\!\ra\!\uSF(\fF\!\t\!\fG)$~is
\begin{equation*}
\bigl(\ts\sum_{i=1}^mc_i[(\fR_i,\rho_i)]\bigr)\!\ot\!
\bigl(\ts\sum_{j=1}^nd_j[(\fS_j,\si_j)]\bigr)\!=\!\ts
\sum_{i,j}c_id_j[(\fR_i\!\t\!\fS_j,\rho_i\!\t\!\si_j)].
\end{equation*}
\label{an2def6}
\end{dfn}

Here \cite[Th.~3.5]{Joyc4} is the analogue of Theorem~\ref{an2thm1}.

\begin{thm} Let\/ $\fE,\fF,\fG,\fH$ be algebraic $\K$-stacks with
affine geometric stabilizers, and\/ $\be:\fF\ra\fG$, $\ga:\fG\ra\fH$
be $1$-morphisms. Then
\e
\begin{aligned}
(\ga\!\ci\!\be)_*\!&=\!\ga_*\!\ci\!\be_*:\uSF(\fF)\!\ra\!\uSF(\fH),&
(\ga\!\ci\!\be)_*\!&=\!\ga_*\!\ci\!\be_*:\SF(\fF)\!\ra\!\SF(\fH),\\
(\ga\!\ci\!\be)^*\!&=\!\be^*\!\ci\!\ga^*:\uSF(\fH)\!\ra\!\uSF(\fF),&
(\ga\!\ci\!\be)^*\!&\!=\!\be^*\!\ci\!\ga^*:\SF(\fH)\!\ra\!\SF(\fF),
\end{aligned}
\label{an2eq4}
\e
for $\be,\ga$ representable in the second equation, and of finite
type in the third and fourth. If\/ $f,g\in\uSF(\fG)$ and\/ $\be$ is
finite type then $\be^*(f\cdot g)=\be^*(f)\cdot\be^*(g)$. If
\e
\begin{gathered}
\xymatrix{
\fE \ar[r]_\eta \ar[d]^{\,\th} & \fG \ar[d]_{\psi\,} \\
\fF \ar[r]^\phi & \fH }
\end{gathered}
\quad
\begin{gathered}
\text{is a Cartesian square with}\\
\text{$\eta,\phi$ representable and}\\
\text{$\th,\psi$ of finite type, then}\\
\text{the following commutes:}
\end{gathered}
\quad
\begin{gathered}
\xymatrix@C=35pt{
\SF(\fE) \ar[r]_{\eta_*} & \SF(\fG) \\
\SF(\fF) \ar[r]^{\phi_*} \ar[u]_{\,\th^*} & \SF(\fH).
\ar[u]^{\psi^*\,} }
\end{gathered}
\label{an2eq5}
\e
The same applies for $\uSF(\fE),\ldots,\uSF(\fH)$, without requiring
$\eta,\phi$ representable.
\label{an2thm2}
\end{thm}

In \cite[Prop.~3.7 \& Th.~3.8]{Joyc4} we relate pushforwards and
pullbacks of stack and constructible functions
using~$\io_\fF,\pi_\fF^\stk$.

\begin{thm} Let\/ $\K$ have characteristic zero, $\fF,\fG$ be
algebraic $\K$-stacks with affine geometric stabilizers, and\/
$\phi:\fF\ra\fG$ be a $1$-morphism. Then
\begin{itemize}
\setlength{\itemsep}{0pt}
\setlength{\parsep}{0pt}
\item[{\rm(a)}] $\phi^*\!\ci\!\io_\fG\!=\!\io_\fF\!\ci\!\phi^*:
\CF(\fG)\!\ra\!\SF(\fF)$ if\/ $\phi$ is of finite type;
\item[{\rm(b)}] $\pi^\stk_\fG\ci\phi_*=\CF^\stk(\phi)\ci\pi_\fF^\stk:
\SF(\fF)\ra\CF(\fG)$ if\/ $\phi$ is representable; and
\item[{\rm(c)}] $\pi^\stk_\fF\ci\phi^*=\phi^*\ci\pi_\fG^\stk:
\SF(\fG)\ra\CF(\fF)$ if\/ $\phi$ is of finite type.
\end{itemize}
\label{an2thm3}
\end{thm}

In \cite[\S 5.2]{Joyc4} we define {\it projections}
$\Pi^\vi_n:\SF(\fF)\ra\SF(\fF)$ and $\uSF(\fF)\ra\uSF(\fF)$ which
project to stack functions whose stabilizer groups have `virtual
rank'~$n$.

In \cite[\S 3]{Joyc4} we define {\it local stack functions}
$\LSF,\uLSF(\fF)$, the analogue of locally constructible functions.
Analogues of Definitions \ref{an2def5}--\ref{an2def6} and Theorems
\ref{an2thm2}--\ref{an2thm3} hold for $\LSF,\uLSF(\fF)$, with
differences in which 1-morphisms are required to be of finite type.

\subsection{Motivic invariants of Artin stacks}
\label{an23}

In \cite[\S 4]{Joyc4} we extend {\it motivic} invariants of
quasiprojective $\K$-varieties to Artin stacks. We need the
following data,~\cite[Assumptions 4.1 \& 6.1]{Joyc4}.

\begin{ass} Suppose $\La$ is a commutative $\Q$-algebra with
identity 1, and
\begin{equation*}
\Up:\{\text{isomorphism classes $[X]$ of quasiprojective
$\K$-varieties $X$}\}\longra\La
\end{equation*}
a map for $\K$ an algebraically closed field, satisfying the
following conditions:
\begin{itemize}
\setlength{\itemsep}{0pt} \setlength{\parsep}{0pt}
\item[(i)] If $Y\subseteq X$ is a closed subvariety then
$\Up([X])=\Up([X\sm Y])+\Up([Y])$;
\item[(ii)] If $X,Y$ are quasiprojective $\K$-varieties then
$\Up([X\!\t\!Y])\!=\!\Up([X])\Up([Y])$;
\item[(iii)] Write $\el=\Up([\K])$ in $\La$, regarding $\K$ as a
$\K$-variety, the affine line (not the point $\Spec\K$). Then $\el$
and $\el^k-1$ for $k=1,2,\ldots$ are invertible in~$\La$.
\end{itemize}
Suppose $\La^\ci$ is a $\Q$-subalgebra of $\La$ containing the image
of $\Up$ and the elements $\el^{-1}$ and
$(\el^k+\el^{k-1}+\cdots+1)^{-1}$ for $k=1,2,\ldots$, but {\it
not\/} containing $(\el-1)^{-1}$. Let $\Om$ be a commutative
$\Q$-algebra, and $\pi:\La^\ci\ra\Om$ a surjective $\Q$-algebra
morphism, such that $\pi(\el)=1$. Define
\begin{equation*}
\Th:\{\text{isomorphism classes $[X]$ of quasiprojective
$\K$-varieties $X$}\}\longra\Om
\end{equation*}
by $\Th=\pi\ci\Up$. Then~$\Th([\K])=1$.
\label{an2ass1}
\end{ass}

We chose the notation `$\el$' as in motivic integration $[\K]$ is
called the {\it Tate motive} and written $\mathbb L$. We have
$\Up\bigl([\GL(m,\K)]\bigr)=\el^{m(m-1)/2}\prod_{k=1}^m(\el^k-1)$,
so (iii) ensures $\Up([\GL(m,\K)])$ is invertible in $\La$ for all
$m\ge 1$. Here \cite[Ex.s 4.3 \& 6.3]{Joyc4} is an example of
suitable $\La,\Up,\ldots$; more are given in~\cite[\S 4.1 \& \S
6.1]{Joyc4}.

\begin{ex} Let $\K$ be an algebraically closed field. Define
$\La=\Q(z)$, the algebra of rational functions in $z$ with
coefficients in $\Q$. For any quasiprojective $\K$-variety $X$, let
$\Up([X])=P(X;z)$ be the {\it virtual Poincar\'e polynomial\/} of
$X$. This has a complicated definition in \cite[Ex.~4.3]{Joyc4}
which we do not repeat, involving Deligne's weight filtration when
$\cha\K=0$ and the action of the Frobenius on $l$-adic cohomology
when $\cha\K>0$. If $X$ is smooth and projective then $P(X;z)$ is
the ordinary Poincar\'e polynomial $\sum_{k=0}^{2\dim X}b^k(X)z^k$,
where $b^k(X)$ is the $k^{\rm th}$ Betti number in $l$-adic
cohomology, for $l$ coprime to $\cha\K$. Also~$\el=P(\K;z)=z^2$.

Let $\La^\ci$ be the subalgebra of $P(z)/Q(z)$ in $\La$ for which
$z\pm 1$ do not divide $Q(z)$. Here are two possibilities for
$\Om,\pi$. Assumption \ref{an2ass1} holds in each case.
\begin{itemize}
\setlength{\itemsep}{0pt}
\setlength{\parsep}{0pt}
\item[(a)] Set $\Om=\Q$ and $\pi:f(z)\mapsto f(-1)$. Then
$\Th([X])=\pi\ci\Up([X])$ is the {\it Euler characteristic} of~$X$.
\item[(b)] Set $\Om=\Q$ and $\pi:f(z)\mapsto f(1)$. Then
$\Th([X])=\pi\ci\Up([X])$ is the {\it sum of the virtual Betti
numbers} of~$X$.
\end{itemize}
\label{an2ex1}
\end{ex}

We need a few facts about {\it algebraic $\K$-groups}. A good
reference is Borel \cite{Bore}. Following Borel, we define a
$\K$-{\it variety} to be a $\K$-scheme which is reduced, separated,
and of finite type, but {\it not\/} necessarily irreducible. An
algebraic $\K$-group is then a $\K$-variety $G$ with identity $1\in
G$, multiplication $\mu:G\t G\ra G$ and inverse $i:G\ra G$ (as
morphisms of $\K$-varieties) satisfying the usual group axioms. We
call $G$ {\it affine} if it is an affine $\K$-variety. {\it
Special\/} $\K$-groups are studied by Serre and Grothendieck
in~\cite[\S 1, \S 5]{Chev}.

\begin{dfn} An algebraic $\K$-group $G$ is called {\it special\/} if
every principal $G$-bundle is locally trivial. Properties of special
$\K$-groups can be found in \cite[\S\S 1.4, 1.5 \& 5.5]{Chev} and
\cite[\S 2.1]{Joyc4}. In \cite[Lem.~4.6]{Joyc4} we show that if
Assumption \ref{an2ass1} holds and $G$ is special then $\Up([G])$ is
invertible in~$\La$.
\label{an2def7}
\end{dfn}

In \cite[Th.~4.9]{Joyc4} we extend $\Up$ to Artin stacks, using
Definition~\ref{an2def7}.

\begin{thm} Let Assumption \ref{an2ass1} hold. Then there exists a
unique morphism of\/ $\Q$-algebras $\Up':\uSF(\Spec\K)\ra\La$ such
that if\/ $G$ is a special algebraic $\K$-group acting on a
quasiprojective $\K$-variety $X$ then~$\Up'\bigl(\bigl[[X/G]
\bigr]\bigr)=\Up([X])/\Up([G])$.
\label{an2thm4}
\end{thm}

Thus, if $\fR$ is a finite type algebraic $\K$-stack with affine
geometric stabilizers the theorem defines $\Up'([\fR])\in\La$.
Taking $\La,\Up$ as in Example \ref{an2ex1} yields the {\it virtual
Poincar\'e function} $P(\fR;z)=\Up'([\fR])$ of $\fR$, a natural
extension of virtual Poincar\'e polynomials to stacks. In \cite[\S
6]{Joyc4} we overcome the restriction that $\Up([G])^{-1}$ exists
for all special $\K$-groups $G$ by defining a finer extension of
$\Up$ to stacks that keeps track of maximal tori of stabilizer
groups, and allows $\Up=\chi$. This can then be used with $\Th,\Om$
in Assumption~\ref{an2ass1}.

In \cite[\S 4--\S 6]{Joyc4} we define other classes of stack
functions $\uSF,\uoSF,\oSF(\fF,\Up,\La),\ab \uoSF,\ab
\oSF(\fF,\Up,\La^\ci),\uoSF,\oSF(\fF,\Th,\Om)$ `twisted' by the
motivic invariants $\Up,\Th$ of Assumption \ref{an2ass1}; the basic
facts are explained in \cite[\S 2.5]{Joyc6}. All the material of
\S\ref{an22} applies to these spaces, except that
$\pi_\fF^\stk,\Pi^\vi_n$ are not always defined. For the purposes of
this paper the differences between these spaces are unimportant, so
we shall not explain them.

\subsection{Essential stack functions and convergent sums}
\label{an24}

Motivated by ideas in Behrend and Dhillon \cite{BeDh}, we extend the
theory of \S\ref{an22}--\S\ref{an23} to include certain local stack
functions, and convergent infinite series. This is new material, not
contained in \cite{Joyc4}. It will be applied in \S\ref{an63}. We
develop it only for $\SF,\LSF(\fF)$, but the extensions to
$\uSF,\uLSF(\fF)$ are obvious.

\begin{dfn} Let $\K$ be algebraically closed and $\fF$ an
algebraic $\K$-stack with affine stabilizers. As in
\cite[p.~98-9]{LaMo} a $\K$-stack $\fR$ has a {\it dimension}
$\dim\fR$ in $\Z\cup\{-\iy,\iy\}$, with $\dim[X/G]=\dim X-\dim G$
for a global quotient. For $m\in\Z$ define $\SF,\LSF(\fF)_m$ to be
the subspaces of $\SF,\LSF(\fF)$ spanned by $[(\fR,\rho)]$ with
$\dim\fR\le m$, where for $\LSF(\fF)_m$ we allow infinite sums
$\sum_{i\in I}c_i[(\fR_i,\rho_i)]$ with $\dim\fR_i\le m$ for all
$i\in I$. Then $\SF(\fF)_m\subseteq\SF(\fF)_n$ if $m\le n$, and
\e
\SF(\fF)=\bigcup_{\!\!m\in\Z\!\!}\SF(\fF)_m,\;
\bigcap_{\!\!m\in\Z\!\!}\LSF(\fF)_m=\{0\},\;
\SF(\fF)_m=\SF(\fF)\cap\LSF(\fF)_m.
\label{an2eq6}
\e
\label{an2def8}
\end{dfn}

Behrend and Dhillon \cite[Def.~2.2]{BeDh} define an algebraic
$\K$-stack $\fR$ to be {\it essentially of finite type} if
$\fR=\coprod_{n\ge 1}\fR_n$ for finite type $\K$-substacks $\fR_n$
with $\dim\fR_n\ra -\iy$ as $n\ra\iy$. This motivates our next
definition, as $[(\fR,\rho)]\in\LSF(\fF)$ lies in $\ESF(\fF)$ if and
only if $\fR$ is essentially of finite type.

\begin{dfn} For $\fF$ as above define $\ESF(\fF)$ to be the subspace
of $f\in\LSF(\fF)$ such that for each $m\in\Z$ we may write $f=g+h$
for $g\in\SF(\fF)$ and $h\in\LSF(\fF)_m$. Then
$\SF(\fF)\subseteq\ESF(\fF)\subseteq\LSF(\fF)$. Elements of
$\ESF(\fF)$ will be called {\it essential stack functions}. Write
$\ESF(\fF)_m=\ESF(\fF)\cap\LSF(\fF)_m$.
\label{an2def9}
\end{dfn}

Here are two notions of convergence of infinite sums
in~$\ESF,\LSF(\fF)$.

\begin{dfn} Let $\fF$ be as above. A possibly infinite sum
$\sum_{i\in I}f_i$ with $f_i\in\LSF(\fF)$ is called {\it
convergent\/} if for all finite type $\K$-substacks $\fG$ in $\fF$,
the restriction of $f_i$ to $\fG$ is nonzero for only finitely many
$i\in I$. Write $f_i=\sum_{a\in A_i}c_i^a[(\fR_i^a,\rho_i^a)]$ for
each $i\in I$, where $[(\fR_i^a,\rho_i^a)]$ is supported on the
support of $f_i$ for each $a\in A_i$. One can then show
$f=\sum_{i\in I}\sum_{a\in A_i}c_i^a[(\fR_i^a, \rho_i^a)]$ is a
well-defined element of $\LSF(\fF)$. We call $f$ the {\it limit\/}
of $\sum_{i\in I}f_i$, and write $\sum_{i\in I}f_i=f$. Note that
$\LSF(\fF)_m$ for $m\in\Z$ are {\it closed under limits}. The same
notion of convergence and limits also works for infinite sums
in~$\LCF(\fF)$.

A sum $\sum_{i\in I}f_i$ with $f_i\in\ESF(\fF)$ is called {\it
strongly convergent\/} if it is convergent, and for all $m\in\Z$ we
have $f_i\in\ESF(\fF)_m$ for all but finitely many $i\in I$. Write
$f=\sum_{i\in I}f_i$ as above, and let $m\in\Z$. Then
$f_i\in\ESF(\fF)_m$ for all $i\in I\sm J$, for finite $J\subseteq
I$. As $f_j\in\ESF(\fF)$ for $j\in J$ we have $f_j=g_j+h_j$ with
$g_j\in\SF(\fF)$ and $h_j\in\LSF(\fF)_m$. Define $g=\sum_{j\in
J}g_j$ and $h=\sum_{j\in J}h_j+\sum_{i\in I\sm J}f_i$. Then $g$ lies
in $\SF(\fF)$ as $g_j$ does and $J$ is finite, and $h$ lies in
$\LSF(\fF)_m$ as $h_j$ for $j\in J$ and $f_i$ for $i\in I\sm J$ do
and $\LSF(\fF)_m$ is closed under limits. Therefore $f\in\ESF(\fF)$,
and $\ESF(\fF)$ is {\it closed under strongly convergent limits}.
\label{an2def10}
\end{dfn}

We modify the first part of Assumption~\ref{an2ass1}.

\begin{ass} Let $\La$ be a commutative $\Q$-algebra and
$\La_m\subset\La$ for $m\in\Z$ a vector subspace, such that
$\La_m\subseteq\La_n$ when $m\le n$ and $\La_m\cdot\La_n
\subseteq\La_{m+n}$ for all $m,n$, with $1\in\La_0$, and
$\La=\bigcup_{m\in\Z}\La_m$, $\bigcap_{m\in\Z}\La_m=\{0\}$. Suppose
\begin{equation*}
\Up:\{\text{isomorphism classes $[X]$ of quasiprojective
$\K$-varieties $X$}\}\longra\La
\end{equation*}
is a map for $\K$ an algebraically closed field, with
$\Up([X])\in\La_{\dim X}$, satisfying:
\begin{itemize}
\setlength{\itemsep}{0pt}
\setlength{\parsep}{0pt}
\item[(i)] If $Y\subseteq X$ is a closed subvariety then
$\Up([X])=\Up([X\sm Y])+\Up([Y])$;
\item[(ii)] If $X,Y$ are quasiprojective $\K$-varieties then
$\Up([X\!\t\!Y])\!=\!\Up([X])\Up([Y])$;
\item[(iii)] Write $\el=\Up([\K])$ in $\La$. Then $\el$ is
invertible in $\La$, with $\el^{-1}\in\La_{-1}$.
\item[(iv)] Suppose we are given elements $\la_m\in\La/\La_m$ for
$m\in\Z$, such that $\la_m+\La_n=\la_n$ whenever $m<n$, using the
inclusion $\La_m\subset\La_n$. Then there exists $\la\in\La$ with
$\la+\La_m=\la_m$ for all $m\in\Z$. This $\la$ is unique
as~$\bigcap_{m\in\Z}\La_m=\{0\}$.
\end{itemize}
\label{an2ass2}
\end{ass}

Here are two examples of suitable $\Up,\La$, the first modifying
Example~\ref{an2ex1}.

\begin{ex} Let $\K$ be an algebraically closed field. Define
$\La$ to be the $\Q$-algebra of $\Q$-Laurent series of the form
$\sum_{k=-\iy}^nc_kz^k$ for $c_k\in\Q$ and $n\in\Z$, that is, power
series in $z^k$ where $k\in\Z$ is bounded above but not necessarily
below. For $m\in\Z$ define $\La_m$ to be the vector subspace of
series $\sum_{k=-\iy}^{2m}c_kz^k$ involving powers of $z$ bounded
above by $2m$. For any quasiprojective $\K$-variety $X$, let
$\Up([X])=P(X;z)$ be the {\it virtual Poincar\'e polynomial\/} of
$X$, as in Example \ref{an2ex1}. Then Assumption \ref{an2ass2}
holds, with~$\el=z^2$.
\label{an2ex2}
\end{ex}

\begin{ex} Let $\K$ be an algebraically closed field. As in Craw
\cite[\S 2.3]{Craw} and Behrend and Dhillon \cite[\S 2.1]{BeDh} we
define the Grothendieck ring $K_0(\Var_\K)$ of the category of
$\K$-varieties $\Var_\K$. Setting $\el=[\K]\in K_0(\Var_\K)$ we form
the ring of fractions $K_0(\Var_\K)[\el^{-1}]$ by inverting $\el$.
For $m\in\Z$ define $K_0(\Var_\K)[\el^{-1}]_m$ to be the subspace
generated by elements $\el^{-n}[X]$ for $n\ge 0$ and $\K$-varieties
$X$ with $\dim X\le m+n$. This defines a {\it filtration} of
$K_0(\Var_\K)[\el^{-1}]$. Write $\hat K_0(\Var_\K)$ to be the {\it
completion} of $K_0(\Var_\K)[\el^{-1}]$ with respect to this
filtration. It is naturally filtered by subspaces $\hat
K_0(\Var_\K)_m$ for~$m\in\Z$.

Define a $\Q$-algebra $\La^{\rm uni}=\hat K_0(\Var_\K)\ot_\Z\Q$, and
set $\La^{\rm uni}_m=\hat K_0(\Var_\K)_m\ot_\Z\Q$ for $m\in\Z$.
Define $\Up^{\rm uni}$ to be the natural map taking a $\K$-variety
$X$ to its class $[X]$ in $K_0(\Var_\K)$ projected to $\La$. Then
Assumption \ref{an2ass2} holds for $\Up^{\rm uni},\La^{\rm uni}$,
and they are {\it universal\/} in that any $\Up,\La$ satisfying
Assumption \ref{an2ass2} factor via a filtered algebra morphism
$\La^{\rm uni}\ra\La$. Note that the ring $\hat K_0(\Var_\K)$ is
used all the time in the subject of {\it motivic integration}, as in
Craw~\cite{Craw}.
\label{an2ex3}
\end{ex}

There is a natural notion of {\it convergence} of infinite sums
in~$\La$.

\begin{dfn} Let Assumption \ref{an2ass2} hold. A possibly infinite sum
$\sum_{i\in I}\la_i$ with $\la_i\in\La$ is called {\it convergent\/}
if for all $m\in\Z$ we have $\la_i\in\La_m$ for all but finitely
many $i\in I$. For $m\in\Z$ define $\la_m\in\La/\La_m$ by
$\la_m=\sum_{j\in J_m}\la_j+\La_m$, where $J_m$ is the finite subset
of $j\in I$ with $\la_j\notin\La_m$. Then
$\la_m\!+\!\La_n\!=\!\la_n$ whenever $m<n$, so Assumption
\ref{an2ass2}(iv) gives a unique $\la\!\in\!\La$ such that
$\la\!+\!\La_m\!=\!\la_m$ for all $m\!\in\!\Z$. We call $\la$ the
{\it limit\/} of $\sum_{i\in I}\la_i$, and write~$\sum_{i\in
I}\la_i\!=\!\la$.
\label{an2def11}
\end{dfn}

For $k\ge 1$, the usual geometric series proof shows $\sum_{n\ge
1}\el^{-nk}$ converges in $\La$ and its limit is $(\el^k-1)^{-1}$.
Therefore Assumption \ref{an2ass1}(i)--(iii) hold, which implies
Theorem \ref{an2thm4}. We can also strengthen it: as
$\el^{-1}\in\La_{-1}$ and $(\el^k-1)^{-1}\in\La_{-k}$ we deduce from
\cite[Lem.s 4.5 \& 4.6]{Joyc4} that $\Up([G])^{-1}\in\La_{-\dim G}$
for all special $\K$-groups $G$. Thus in $\Up'([[X/G]])=
\Up([X])/\Up([G])$ we have $\Up([X])\in\La_{\dim X}$ and
$\Up([G])^{-1}\in \La_{-\dim G}$, so $\Up'([[X/G]])\in\La_{\dim
X-\dim G}=\La_{\dim[X/G]}$. As any finite type $\fR$ with affine
stabilizers is a disjoint union of $[X/G]$ we deduce:

\begin{thm} Let Assumption \ref{an2ass2} hold. Then there exists a
unique morphism of\/ $\Q$-algebras $\Up':\uSF(\Spec\K)\ra\La$ such
that if\/ $G$ is a special algebraic $\K$-group acting on a
$\K$-variety $X$ then $\Up'\bigl(\bigl[[X/G]\bigr]
\bigr)=\Up([X])/\Up([G])$, and if\/ $\fR$ is a finite type algebraic
$\K$-stack with affine geometric stabilizers
then~$\Up'\bigl([\fR]\bigr)\in\La_{\dim\fR}$.
\label{an2thm5}
\end{thm}

The next definition was suggested to the author by Behrend and
Dhillon's definition \cite[\S 2.2]{BeDh} of the {\it motive} of an
essentially of finite type $\K$-stack.

\begin{dfn} Let Assumption \ref{an2ass2} hold, $\Up'$ be as in Theorem
\ref{an2thm5}, $\fF$ be as above, and $\Pi:\fF\ra\Spec\K$ be the
projection 1-morphism. Then $\Pi_*$ maps $\SF(\fF)\ra
\uSF(\Spec\K)$, so $\Up'\ci\Pi_*:\SF(\fF)\ra\La$, with
$\Up'\ci\Pi_*:[(\fR,\rho)]\mapsto\Up'([\fR])$. Since $\SF(\fF)_m$ is
spanned by $[(\fR,\rho)]$ with $\dim\fR\le m$, so that
$\Up'\bigl([\fR]\bigr)\in\La_m$ by Theorem \ref{an2thm5}, we have
$\Up'\ci\Pi_*:\SF(\fF)_m\ra\La_m$ for~$m\in\Z$.

Let $f\in\ESF(\fF)$. For $m\in\Z$ write $f=g_m+h_m$ for
$g_m\in\SF(\fF)$ and $h_m\in\LSF(\fF)_m$. Set
$\la_m=\Up'\ci\Pi_*(g_m)+\La_m$ in $\La/\La_m$. If $g_m',h_m'$ are
alternative choices of $g_m,h_m$ then $g_m+h_m=g_m'+h_m'$, so
$g_m-g_m'=h_m'-h_m$, which lies in $\SF(\fF)\cap\LSF(\fF)_m=
\SF(\fF)_m$ by \eq{an2eq6}. Thus $\Up'\ci\Pi_*(g_m-g_m')\in\La_m$,
so $\Up'\ci\Pi_* (g_m)+\La_m=\Up'\ci\Pi_*(g_m')+\La_m$, and $\la_m$
is independent of choices.

If $m<n$ then we may define $\la_n$ using $g_m,h_m$ instead of
$g_n,h_n$, giving $\la_m+\La_n=\la_n$. Thus Assumption
\ref{an2ass2}(iv) gives a unique $\la\in\La$ with $\la+\La_m=\la_m$
for all $m\in\Z$. Define $\Pi_\La(f)=\la$. This gives a $\Q$-linear
map $\Pi_\La:\ESF(\fF)\ra\La$. If $f\in\SF(\fF)$ we may take $g_m=f$
and $h_m=0$ for all $m\in\Z$, giving $\la_m=\Up'\ci\Pi_*(f)+\La_m$,
so $\la=\Up'\ci\Pi_*(f)$ by uniqueness. Thus $\Pi_\La\!=\!
\Up'\ci\Pi_*$ on $\SF(\fF)$. It is easy to show $\Pi_\La$ maps
$\ESF(\fF)_m\!\ra\!\La_m$ for~$m\!\in\!\Z$.
\label{an2def12}
\end{dfn}

This $\Pi_\La$ commutes with (strongly) convergent limits in
$\ESF(\fF)$ and~$\La$.

\begin{prop} Let Assumption \ref{an2ass2} hold, $\fF$ be an algebraic
$\K$-stack with affine geometric stabilizers, and\/ $\sum_{i\in
I}f_i$ be a strongly convergent sum in $\ESF(\fF)$ with limit\/ $f$.
Then $\sum_{i\in I}\Pi_\La(f_i)$ is convergent in $\La$ with
limit\/~$\Pi_\La(f)$.
\label{an2prop}
\end{prop}

\begin{proof} Let $m\in\Z$. As $\sum_{i\in I}f_i$ is strongly
convergent we have $f_i\in\ESF(\fF)_m$ for all $i\in I\sm J$, where
$J\subseteq I$ is finite. But $\Pi_\La$ maps $\ESF(\fF)_m\ra\La_m$,
so $\Pi_\La(f_i)\in\La_m$ for all $i\in I\sm J$ with $J$ finite, and
thus $\sum_{i\in I}\Pi_\La(f_i)$ converges in $\La$. Let the limits
be $f=\sum_{i\in I}f_i$ in $\ESF(\fF)$ and $\la=\sum_{i\in
I}\Pi_\La(f_i)$ in $\La$. For each $j\in J$ write $f_j=g_j+h_j$ for
$g_j\in\SF(\fF)$ and $h_j\in\LSF(\fF)_m$. Define $g=\sum_{j\in
J}g_j$ and $h=\sum_{j\in J}h_j+\sum_{i\in I\sm J}f_i$. Then $f=g+h$
with $g\in\SF(\fF)$ and $h\in\LSF(\fF)_m$. Therefore using
Definitions \ref{an2def10}, \ref{an2def11} and \ref{an2def12} we
have
\begin{equation*}
\ts\Pi_\La(f)\!+\!\La_m\!=\!\Pi_\La(g)\!+\!\La_m\!=\!\sum_{j\in
J}\Pi_\La(g_j)\!+\!\La_m\!=\!
\sum_{j\in J}\Pi_\La(f_j)\!+\!\La_m\!=\!\la\!+\!\La_m,
\end{equation*}
for all $m\in\Z$. As $\bigcap_{m\in\Z}\La_m=\{0\}$ this
forces~$\Pi_\La(f)=\la$.
\end{proof}

\section{Background on configurations from \cite{Joyc5,Joyc6,Joyc7}}
\label{an3}

We now recall in \S\ref{an31}--\S\ref{an32} the main definitions and
results from \cite{Joyc5} on $(I,\pr)$-configurations and their
moduli stacks that we will need later, in \S\ref{an33} some facts
about algebras of constructible and stack functions from
\cite{Joyc6}, and in \S\ref{an34}--\S\ref{an35} some material on
(weak) stability conditions from~\cite{Joyc7}.

\subsection{Basic definitions}
\label{an31}

Here is some notation for {\it finite posets}, taken from
\cite[Def.s~3.2, 4.1 \& 6.1]{Joyc5}.

\begin{dfn} A {\it finite partially ordered set\/} or {\it
finite poset\/} $(I,\pr)$ is a finite set $I$ with a partial
order $I$. Define $J\subseteq I$ to be an {\it f-set\/} if
$i\in I$ and $h,j\in J$ and $h\pr i\pr j$ implies $i\in J$.
Define $\F_\sIp$ to be the set of f-sets of $I$. Define
$\G_\sIp$ to be the subset of $(J,K)\in\F_\sIp\t\F_\sIp$
such that $J\subseteq K$, and if $j\in J$ and $k\in K$
with $k\pr j$, then $k\in J$. Define $\H_\sIp$ to be the
subset of $(J,K)\in\F_\sIp\t\F_\sIp$ such that
$K\subseteq J$, and if $j\in J$ and $k\in K$ with
$k\pr j$, then~$j\in K$.

Let $I$ be a finite set and $\pr,\tl$ partial orders on $I$ such
that if $i\pr j$ then $i\tl j$ for $i,j\in I$. Then we say that
$\tl$ {\it dominates} $\pr$. A partial order $\tl$ on $I$ is called
a {\it total order} if $i\tl j$ or $j\tl i$ for all $i,j\in I$. Then
$(I,\tl)$ is canonically isomorphic to $(\{1,\ldots,n\},\le)$
for~$n=\md{I}$.
\label{an3def1}
\end{dfn}

We define $(I,\pr)$-{\it configurations},~\cite[Def.~4.1]{Joyc5}.

\begin{dfn} Let $(I,\pr)$ be a finite poset, and use the
notation of Definition \ref{an3def1}. Define an $(I,\pr)$-{\it
configuration} $(\si,\io,\pi)$ in an abelian category $\A$ to be
maps $\si:\F_\sIp\ra\Obj(\A)$, $\io:\G_\sIp\ra\Mor(\A)$, and
$\pi:\H_\sIp\ra\Mor(\A)$, where
\begin{itemize}
\setlength{\itemsep}{0pt}
\setlength{\parsep}{0pt}
\item[(i)] $\si(J)$ is an object in $\A$ for $J\in\F_\sIp$,
with~$\si(\emptyset)=0$.
\item[(ii)] $\io(J,K):\si(J)\!\ra\!\si(K)$ is injective
for $(J,K)\!\in\!\G_\sIp$, and~$\io(J,J)\!=\!\id_{\si(J)}$.
\item[(iii)] $\pi(J,K)\!:\!\si(J)\!\ra\!\si(K)$ is surjective
for $(J,K)\!\in\!\H_\sIp$, and~$\pi(J,J)\!=\!\id_{\si(J)}$.
\end{itemize}
These should satisfy the conditions:
\begin{itemize}
\setlength{\itemsep}{0pt}
\setlength{\parsep}{0pt}
\item[(A)] Let $(J,K)\in\G_\sIp$ and set $L=K\sm J$. Then the
following is exact in~$\A$:
\e
\xymatrix@C=40pt{ 0 \ar[r] &\si(J) \ar[r]^{\io(J,K)} &\si(K)
\ar[r]^{\pi(K,L)} &\si(L) \ar[r] & 0. }
\label{an3eq1}
\e
\item[(B)] If $(J,K)\in\G_\sIp$ and $(K,L)\in\G_\sIp$
then~$\io(J,L)=\io(K,L)\ci\io(J,K)$.
\item[(C)] If $(J,K)\in\H_\sIp$ and $(K,L)\in\H_\sIp$
then~$\pi(J,L)=\pi(K,L)\ci\pi(J,K)$.
\item[(D)] If $(J,K)\in\G_\sIp$ and $(K,L)\in\H_\sIp$ then
\begin{equation*}
\pi(K,L)\ci\io(J,K)=\io(J\cap L,L)\ci\pi(J,J\cap L).
\end{equation*}
\end{itemize}

A {\it morphism} $\al:(\si,\io,\pi)\ra(\si',\io',\pi')$ of
$(I,\pr)$-configurations in $\A$ is a collection of morphisms
$\al(J):\si(J)\ra\si'(J)$ for each $J\in\F_\sIp$ satisfying
\begin{align*}
\al(K)\ci\io(J,K)&=\io'(J,K)\ci\al(J)&&
\text{for all $(J,K)\in\G_\sIp$, and}\\
\al(K)\ci\pi(J,K)&=\pi'(J,K)\ci\al(J)&&
\text{for all $(J,K)\in\H_\sIp$.}
\end{align*}
It is an {\it isomorphism} if $\al(J)$ is an isomorphism for
all~$J\in\F_\sIp$.
\label{an3def2}
\end{dfn}

In \cite[Prop.~4.7]{Joyc5} we relate the classes $[\si(J)]$
in~$K_0(\A)$.

\begin{prop} Let\/ $(\si,\io,\pi)$ be an $(I,\pr)$-configuration
in an abelian category $\A$. Then there exists a unique map
$\ka:I\ra K_0(\A)$ such that\/ $[\si(J)]=\sum_{j\in J}\ka(j)$
in $K_0(\A)$ for all f-sets~$J\subseteq I$.
\label{an3prop1}
\end{prop}

Here \cite[Def.s~5.1, 5.2]{Joyc5} are two ways to construct new
configurations.

\begin{dfn} Let $(I,\pr)$ be a finite poset and $J\in\F_\sIp$. Then
$(J,\pr)$ is also a finite poset, and $\F_\sJp,\G_\sJp,\H_\sJp\!
\subseteq\!\F_\sIp,\G_\sIp,\H_\sIp$. Let $(\si,\io,\pi)$ be an
$(I,\pr)$-configuration in an abelian category $\A$. Define the
$(J,\pr)$-{\it subconfiguration} $(\si',\io',\pi')$ of
$(\si,\io,\pi)$ by $\si'\!=\!\si\vert_{\F_\sJp}$,
$\io'\!=\!\io\vert_{\G_\sJp}$ and~$\pi'\!=\!\pi\vert_{\H_\sJp}$.

Let $(I,\pr),(K,\tl)$ be finite posets, and $\phi:I\!\ra\!K$ be
surjective with $i\pr j$ implies $\phi(i)\!\tl \!\phi(j)$. Using
$\phi^{-1}$ to pull subsets of $K$ back to $I$ maps
$\F_\sKt,\G_\sKt,\ab\H_\sKt\!\ra\!\F_\sIp,\G_\sIp,\H_\sIp$. Let
$(\si,\io,\pi)$ be an $(I,\pr)$-configuration in an abelian category
$\A$. Define the {\it quotient\/ $(K,\tl)$-configuration}
$(\ti\si,\ti\io,\ti\pi)$ by
$\ti\si(A)\!=\!\si\bigl(\phi^{-1}(A)\bigr)$ for $A\!\in\!\F_\sKt$,
$\ti\io(A,B)\!=\!\io\bigl(\phi^{-1}(A),\ab \phi^{-1}(B)\bigr)$ for
$(A,B)\!\in\!\G_\sKt$, and
$\ti\pi(A,B)\!=\!\pi\bigl(\phi^{-1}(A),\phi^{-1}(B)\bigr)$
for~$(A,B)\!\in\!\H_\sKt$.
\label{an3def3}
\end{dfn}

\subsection{Moduli stacks of configurations}
\label{an32}

Here \cite[Assumptions 7.1 \& 8.1]{Joyc5} is the data we require.

\begin{ass} Let $\K$ be an algebraically closed field and $\A$
a $\K$-linear noetherian abelian category with $\Ext^i(X,Y)$
finite-dimensional $\K$-vector spaces for all $X,Y\in\A$ and $i\ge
0$. Let $K(\A)$ be the quotient of the Grothendieck group $K_0(\A)$
by some fixed subgroup. Suppose that if $X\in\A$ with $[X]=0$ in
$K(\A)$ then~$X\cong 0$.

To define moduli stacks of objects or configurations in $\A$, we
need some {\it extra data}, to tell us about algebraic families of
objects and morphisms in $\A$, parametrized by a base scheme $U$. We
encode this extra data as a {\it stack in exact categories\/}
$\fF_\A$ on the {\it category of\/ $\K$-schemes\/} $\Sch_\K$, made
into a {\it site\/} with the {\it \'etale topology}. The
$\K,\A,K(\A),\fF_\A$ must satisfy some complex additional conditions
\cite[Assumptions 7.1 \& 8.1]{Joyc5}, which we do not give.
\label{an3ass}
\end{ass}

In \cite[\S 9--\S 10]{Joyc5} we define the data $\A,K(\A),\fF_\A$ in
some large classes of examples, and prove Assumption \ref{an3ass}
holds in each case. Note that \cite{Joyc5,Joyc6} did not assume $\A$
{\it noetherian}, but as in \cite{Joyc7} we need this to make
$\tau$-semistability well-behaved. All the examples of \cite[\S
9--\S 10]{Joyc5} have $\A$ noetherian.

To apply the constructible functions material of \S\ref{an21} we
need the ground field $\K$ to have {\it characteristic zero}, but
the stack functions of \S\ref{an22}--\S\ref{an24} work for $\K$ of
{\it arbitrary characteristic}. As we develop the two strands in
parallel, in this section for brevity we make the convention that
$\cha\K=0$ in the parts dealing with $\CF,\LCF(\fObj_\A)$ and
$\cha\K$ is arbitrary otherwise.

\begin{dfn} We work in the situation of Assumption \ref{an3ass}.
Define
\e
C(\A)=\bigl\{[X]\in K(\A):X\in\A,\;\> X\not\cong 0\bigr\}
\subset K(\A).
\label{an3eq2}
\e
That is, $C(\A)$ is the collection of classes in $K(\A)$ of {\it
nonzero objects} $X\in\A$. Note that $C(\A)$ is {\it closed under
addition}, as $[X\op Y]=[X]+[Y]$. Note also that $0\notin C(\A)$, as by
Assumption \ref{an3ass} if $X\not\cong 0$ then $[X]\ne 0$ in~$K(\A)$.

In \cite{Joyc5,Joyc6} we worked mostly with $\bar
C(\A)=C(\A)\cup\{0\}$, the collection of classes in $K(\A)$ of all
objects $X\in\A$. But here and in \cite{Joyc7} we find $C(\A)$ more
useful, as stability conditions will be defined only on nonzero
objects. We think of $C(\A)$ as the `positive cone' and $\bar C(\A)$
as the `closed positive cone' in~$K(\A)$.

Define a set of $\A$-{\it data} to be a triple $(I,\pr,\ka)$
such that $(I,\pr)$ is a finite poset and $\ka:I\ra C(\A)$ a
map. We {\it extend\/ $\ka$ to the set of subsets of\/} $I$ by
defining $\ka(J)=\sum_{j\in J}\ka(j)$. Then $\ka(J)\in C(\A)$
for all $\emptyset\ne J\subseteq I$, as $C(\A)$ is closed
under addition. Define an $(I,\pr,\ka)$-{\it configuration}
to be an $(I,\pr)$-configuration $(\si,\io,\pi)$ in $\A$ with
$[\si(\{i\})]=\ka(i)$ in $K(\A)$ for all $i\in I$. Then $[\si(J)]
=\ka(J)$ for all $J\in\F_\sIp$, by Proposition~\ref{an3prop1}.
\label{an3def4}
\end{dfn}

In the situation above, we define the following $\K$-stacks
\cite[Def.s 7.2 \& 7.4]{Joyc5}:
\begin{itemize}
\setlength{\itemsep}{0pt}
\setlength{\parsep}{0pt}
\item The {\it moduli stacks} $\fObj_\A$ of {\it objects in} $\A$,
and $\fObj^\al_\A$ of {\it objects in $\A$ with class $\al$ in}
$K(\A)$, for each $\al\in\bar C(\A)$. They are algebraic
$\K$-stacks. The underlying geometric spaces
$\fObj_\A(\K),\fObj_\A^\al(\K)$ are the sets of isomorphism classes
of objects $X$ in $\A$, with $[X]=\al$ for~$\fObj_\A^\al(\K)$.
\item The {\it moduli stacks\/} $\fM(I,\pr)_\A$ of $(I,\pr)$-{\it
configurations} and $\fM(I,\pr,\ka)_\A$ of $(I,\pr,\ka)$-{\it
configurations in} $\A$, for all finite posets $(I,\pr)$ and
$\ka:I\ra\bar C(\A)$. They are algebraic $\K$-stacks. Write
$\M(I,\pr)_\A,\M(I,\pr,\ka)_\A$ for the underlying geometric spaces
$\fM(I,\pr)_\A(\K),\fM(I,\pr,\ka)_\A(\K)$. Then $\M(I,\pr)_\A$,
$\M(I,\pr,\ka)_\A$ are the {\it sets of isomorphism classes of\/
$(I,\pr)$- and\/ $(I,\pr,\ka)$-configurations in} $\A$,
by~\cite[Prop.~7.6]{Joyc5}.
\end{itemize}

In \cite[Def.~7.7 \& Prop.~7.8]{Joyc5} we define 1-morphisms of
$\K$-stacks, as follows:
\begin{itemize}
\setlength{\itemsep}{0pt}
\setlength{\parsep}{0pt}
\item For $(I,\pr)$ a finite poset, $\ka:I\ra\bar C(\A)$ and
$J\in\F_\sIp$, we define $\bs\si(J):\fM(I,\pr)_\A\ra\fObj_\A$
or $\bs\si(J):\fM(I,\pr,\ka)_\A\ra\fObj_\A^{\ka(J)}$. The
induced maps $\bs\si(J)_*:\M(I,\pr)_\A\ra\fObj_\A(\K)$ or
$\M(I,\pr,\ka)_\A\ra\fObj_\A^{\ka(J)}(\K)$ act
by~$\bs\si(J)_*:[(\si,\io,\pi)]\mapsto[\si(J)]$.
\item For $(I,\pr)$ a finite poset, $\ka:I\ra\bar C(\A)$ and
$J\in\F_\sIp$, we define the $(J,\pr)$-{\it subconfiguration
$1$-morphism} $S(I,\pr,J):\fM(I,\pr,\ka)_\A\ra
\fM(J,\pr,\ka\vert_J)_\A$. Then $S(I,\pr,J)_*$ takes
$[(\si,\io,\pi)]\!\mapsto\![(\si',\io',\pi')]$, for $(\si,\io,\pi)$
an $(I,\pr,\ka)$-configuration in $\A$, and $(\si',\io',\pi')$ its
$(J,\pr)$-subconfiguration.
\item Let $(I,\pr)$, $(K,\tl)$ be finite posets, $\ka:I\!\ra\!\bar C(\A)$,
and $\phi:I\!\ra\!K$ be surjective with $i\!\pr\!j$ implies $\phi(i)
\!\tl\!\phi(j)$ for $i,j\in I$. Define $\mu:K\!\ra\!\bar C(\A)$ by
$\mu(k)\!=\!\ka(\phi^{-1}(k))$. The {\it quotient\/
$(K,\tl)$-configuration\/ $1$-morphism} is
$Q(I,\pr,K,\tl,\phi):\fM(I,\pr,\ka)_\A\!\ra\!\fM(K,\tl,\mu)_\A$.
Then $Q(I,\pr,K,\tl,\phi)_*$ takes
$[(\si,\io,\pi)]\!\mapsto\![(\ti\si,\ti\io,\ti\pi)]$, where
$(\si,\io,\pi)$ is an $(I,\pr,\ka)$-configuration in $\A$, and
$(\ti\si,\ti\io,\ti\pi)$ its quotient $(K,\tl)$-configuration.
\end{itemize}

\subsection{Algebras of constructible and stack functions}
\label{an33}

Next we summarize parts of \cite{Joyc6}, which define and study
associative multiplications $*$ on $\CF(\fObj_\A)$ and
$\SF(\fObj_\A)$, based on {\it Ringel--Hall algebras}.

\begin{dfn} Let Assumption \ref{an3ass} hold with $\K$ of
characteristic zero. Write $\de_{[0]}\in\CF(\fObj_\A)$ for the
characteristic function of $[0]\in\fObj_\A(\K)$. Following
\cite[Def.~4.1]{Joyc6}, using the diagrams of 1-morphisms and
pullbacks, pushforwards
\begin{equation*}
\text{
\begin{footnotesize}
$\displaystyle
\begin{gathered}
\xymatrix@C=54pt@R=10pt{ \fObj_\A\t\fObj_\A & \fM(\{1,2\},\le)_\A
\ar[l]_{\bs\si(\{1\})\t\bs\si(\{2\})} \ar[r]^{\bs\si(\{1,2\})}
& \fObj_\A,\\
\CF\bigl(\fObj_\A\bigr)\t\CF\bigl(\fObj_\A\bigr)
\!\!\!\!\!\!\!\!\!\!\!\!\!\!\!
\ar@<.7ex>[dr]^(0.6){\qquad\qquad(\bs\si(\{1\}))^*\cdot(\bs\si(\{2\}))^*}
\ar[d]_{\ot}
\\
\CF\bigl(\fObj_\A\!\t\!\fObj_\A\bigr)
\ar[r]^(.47){(\bs\si(\{1\})\t\bs\si(\{2\}))^*} &
\CF\bigl(\fM(\{1,2\},\le)_\A\bigr)
\ar[r]^(.56){\CF^\stk(\bs\si(\{1,2\}))} & \CF\bigl(\fObj_\A\bigr), }
\end{gathered}
$
\end{footnotesize}}
\end{equation*}
define a bilinear operation $*:\CF(\fObj_\A)\t\CF(\fObj_\A)
\ra\CF(\fObj_\A)$ by
\e
f*g=\CF^\stk(\bs\si(\{1,2\}))\bigl[\bs\si(\{1\})^*(f)\cdot
\bs\si(\{2\})^*(g)\bigr].
\label{an3eq3}
\e
Then \cite[Th.~4.3]{Joyc6} shows $*$ is {\it associative}, and
$\CF(\fObj_\A)$ is a $\Q$-{\it algebra}, with identity $\de_{[0]}$
and multiplication~$*$.

This extends to {\it locally constructible functions}, \cite[\S
4.2]{Joyc6}. Write $\dLCF(\fObj_\A)$ for the subspace of
$f\in\LCF(\fObj_\A)$ supported on $\coprod_{\al\in S}
\fObj_\A^\al(\K)$ for $S\subset\bar C(\A)$ a finite subset. Then $*$
in \eq{an3eq3} is well-defined on $\dLCF(\fObj_\A)$, and makes
$\dLCF(\fObj_\A)$ into a $\Q$-algebra containing $\CF(\fObj_\A)$ as
a subalgebra.

Following \cite[Def.~4.8]{Joyc6}, write $\CFi,\dLCFi(\fObj_\A)$ for
the vector subspaces of $f$ in $\CF,\dLCF(\fObj_\A)$ {\it supported
on indecomposables}, that is, $f\bigl([X]\bigr)\ne 0$ implies
$0\not\cong X$ is indecomposable. Define bilinear brackets $[\,,\,]$
on $\CF,\dLCF(\fObj_\A)$ by $[f,g]=f*g-g*f$. Then
\cite[Th.~4.9]{Joyc6} shows $\CFi,\dLCFi(\fObj_\A)$ are closed under
$[\,,\,],$ and so are $\Q$-Lie algebras.
\label{an3def5}
\end{dfn}

In \cite[\S 5]{Joyc6} we extend much of the above to {\it stack
functions}, as in \S\ref{an22}. Here are a few of the basic
definitions and results.

\begin{dfn} Suppose Assumption \ref{an3ass} holds. Define
$\dLSF(\fObj_\A)$ to be the subspace of $f\in\LSF(\fObj_\A)$
supported on $\coprod_{\al\in S}\fObj_\A^\al(\K)$ for $S\subset\bar
C(\A)$ finite. By analogy with \eq{an3eq3}, using
$\fM(\{1,2\},\le)_\A$ define \cite[Def.~5.1 \& \S 5.4]{Joyc6}
bilinear operations $*:\SF(\fObj_\A)\t \SF(\fObj_\A)\ra
\SF(\fObj_\A)$ and $*:\dLSF(\fObj_\A)\t\dLSF(\fObj_\A)\ra
\dLSF(\fObj_\A)$~by
\begin{equation*}
f*g=\bs\si(\{1,2\})_*\bigl[(\bs\si(\{1\})\t\bs\si(\{2\}))^*(f\ot
g)\bigr].
\end{equation*}
Write $\bde_{[0]}\!\in\!\SF(\fObj_\A)$ for $\bde_C$ in Definition
\ref{an2def5} with $C\!=\!\{[0]\}$. Then \cite[Th.~5.2]{Joyc6} shows
$\SF,\dLSF(\fObj_\A)$ are $\Q$-algebras with associative
multiplication $*$ and identity $\bde_{[0]}$. When $\K$ has
characteristic zero there are $\Q$-algebra morphisms
\e
\pi_{\fObj_\A}^\stk:\SF(\fObj_\A)\ra\CF(\fObj_\A),\;\>
\pi_{\fObj_\A}^\stk:\dLSF(\fObj_\A)\ra\dLCF(\fObj_\A).
\label{an3eq4}
\e
As in \cite[Def.~5.5]{Joyc6} write $\SFa,\dLSFa(\fObj_\A)$ for the
subspaces of $\SF,\dLSF(\fObj_\A)$ spanned by $[(\fR,\rho)]$ such
that for all $r\in\fR(\K)$ with $\rho_*(r)=[X]$, the $\K$-subgroup
$\rho_*\bigl(\Iso_\K(r)\bigr)$ in $\Aut(X)$ is the $\K$-group of
invertible elements in a $\K$-subalgebra of $\End(X)$. Then
$\io_{\fObj_\A}$ in Definition \ref{an2def5} maps
$\CF(\fObj_\A)\!\ra\!\SFa(\fObj_\A)$ and
$\dLCF(\fObj_\A)\!\ra\!\dLSFa(\fObj_\A)$, and
\cite[Prop.~5.6]{Joyc6} shows $\SFa,\dLSFa(\fObj_\A)$ are closed
under $*$ and so are $\Q$-subalgebras.
\label{an3def6}
\end{dfn}

\begin{dfn} Suppose Assumption \ref{an3ass} holds. Following
\cite[Def.~5.13]{Joyc6}, define $\SFai,\dLSFai(\fObj_\A)$ to be the
subspaces of $f\in\SFa,\dLSFa(\fObj_\A)$ with $\Pi^\vi_1(f)=f$,
where $\Pi^\vi_1$ is the operator of \cite[\S 5.2]{Joyc4},
interpreted as projecting to stack functions `supported on virtual
indecomposables'. Write $[f,g]=f*g-g*f$ for
$f,g\in\SFa,\dLSFa(\fObj_\A)$. As $*$ is associative $[\,,\,]$
satisfies the {\it Jacobi identity}, and makes
$\SFa,\dLSFa(\fObj_\A)$ into $\Q$-{\it Lie algebras}. Then
\cite[Th.~5.17]{Joyc6} shows $\SFai,\dLSFai(\fObj_\A)$ are closed
under $[\,,\,]$, and are {\it Lie subalgebras}. When $\cha\K=0$,
\eq{an3eq4} restricts to Lie algebra morphisms
\begin{equation*}
\pi_{\fObj_\A}^\stk\!:\!\SFai(\fObj_\A)\!\ra\!\CFi(\fObj_\A),\;
\pi_{\fObj_\A}^\stk\!:\!\dLSFai(\fObj_\A)\!\ra\!\dLCFi(\fObj_\A).
\end{equation*}
\label{an3def7}
\end{dfn}

All this also works for other stack function spaces on $\fObj_\A$,
in particular for $\oSF(\fObj_\A,\Up,\La),\oSF(\fObj_\A,\Up,
\La^\ci)$ and $\oSF(\fObj_\A,\Th,\Om)$, giving algebras
$\oSF,\oSFa(\fObj_\A,*,*)$ and Lie algebras $\oSFai(\fObj_\A,*,*)$.
In \cite[\S 6]{Joyc6} under extra conditions on $\A$, we define
({\it Lie}) {\it algebra morphisms} from $\SFa,\SFai(\fObj_\A,*,*)$
to explicit algebras $A(\A,\La,\chi), \ldots,C(\A,\Om,\chi)$, which
will be important in~\S\ref{an6}.

\subsection{(Weak) stability conditions}
\label{an34}

We now summarize the material of \cite[\S 4]{Joyc7}, beginning
with~\cite[Def.s 4.1--4.3]{Joyc7}.

\begin{dfn} Let $\A$ be an abelian category, $K(\A)$ be
the quotient of $K_0(\A)$ by some fixed subgroup, and $C(\A)$ as in
\eq{an3eq2}. Suppose $(T,\le)$ is a totally ordered set, and
$\tau:C(\A)\ra T$ a map. We call $(\tau,T,\le)$ a {\it stability
condition} on $\A$ if whenever $\al,\be,\ga\in C(\A)$ with
$\be=\al+\ga$ then either $\tau(\al)\!<\!\tau(\be) \!<\!\tau(\ga)$,
or $\tau(\al)\!>\!\tau(\be)\!>\!\tau(\ga)$, or
$\tau(\al)\!=\!\tau(\be)\!=\!\tau(\ga)$. We call $(\tau,T,\le)$ a
{\it weak stability condition} on $\A$ if whenever $\al,\be, \ga\in
C(\A)$ with $\be=\al+\ga$ then either $\tau(\al)\!\le\!
\tau(\be)\!\le\!\tau(\ga)$, or $\tau(\al)\!\ge\!\tau(\be)\!\ge
\!\tau(\ga)$. The alternative $\tau(\al)\!\le\!\tau(\be)\!\le\!
\tau(\ga)$ or $\tau(\al)\!\ge\!\tau(\be)\!\ge\!\tau(\ga)$ is called
the {\it weak seesaw inequality}.
\label{an3def8}
\end{dfn}

We use many ordered sets in the paper: finite posets $(I,\pr),
(J,\ls),(K,\tl)$ for $(I,\pr)$-configurations, and now total
orders $(T,\le)$ for stability conditions. As the number of
order symbols is limited, we will always use `$\le$' for the
total order, so that $(\tau,T,\le)$, $(\ti\tau,\ti T,\le)$
may denote two different stability conditions, with two {\it
different\/} total orders on $T,\ti T$ both denoted by~`$\le$'.

\begin{dfn} Let $(\tau,T,\le)$ be a weak stability condition on
$\A,K(\A)$ as above. Then we say that a nonzero object $X$ in $\A$ is
\begin{itemize}
\setlength{\itemsep}{0pt}
\setlength{\parsep}{0pt}
\item[(i)] $\tau$-{\it semistable} if for all $S\subset X$ with
$S\not\cong 0,X$ we have $\tau([S])\le\tau([X/S])$;
\item[(ii)] $\tau$-{\it stable} if for all $S\subset X$ with
$S\not\cong 0,X$ we have $\tau([S])<\tau([X/S])$; and
\item[(iii)] $\tau$-{\it unstable} if it is not $\tau$-semistable.
\end{itemize}
\label{an3def9}
\end{dfn}

\begin{dfn} Let $(\tau,T,\le)$ be a weak stability condition on
$\A,K(\A)$. We say $\A$ is $\tau$-{\it artinian} if there exist
no infinite chains of subobjects $\cdots\!\subset\!A_2\!\subset\!
A_1\!\subset\!X$ in $\A$ with $A_{n+1}\!\ne\!A_n$ and
$\tau([A_{n+1}])\!\ge\!\tau([A_n/A_{n+1}])$ for all~$n$.
\label{an3def10}
\end{dfn}

Here is \cite[Th.~4.4]{Joyc7}, based on Rudakov \cite[Th.~2]{Ruda}.
We call $0\!=\!A_0\!\subset\!\cdots\!\subset\!A_n\!=\!X$ in Theorem
\ref{an3thm} the {\it Harder--Narasimhan filtration} of~$X$.

\begin{thm} Let\/ $(\tau,T,\le)$ be a weak stability condition
on an abelian category $\A$. Suppose $\A$ is noetherian and\/
$\tau$-artinian. Then each\/ $X\in\A$ admits a unique filtration
$0\!=\!A_0\!\subset\!\cdots\!\subset\!A_n\!=\!X$ for $n\ge 0$, such
that\/ $S_k\!=\!A_k/A_{k-1}$ is $\tau$-semistable for
$k=1,\ldots,n$, and\/~$\tau([S_1])>\tau([S_2])>\cdots>\tau([S_n])$.
\label{an3thm}
\end{thm}

We define some notation,~\cite[Def.s 4.6, 8.1, 4.7 \& 4.10]{Joyc7}.

\begin{dfn} Let Assumption \ref{an3ass} hold and $(\tau,T,\le)$
be a weak stability condition on $\A$. Then $\fObj_\A^\al$ is an
algebraic $\K$-stack for $\al\in C(\A)$, with $\fObj_\A^\al(\K)$
the set of isomorphism classes of $X\in\A$ with class $\al$ in
$K(\A)$. Define
\begin{align*}
\Oss^\al(\tau)&=\bigl\{[X]\in\fObj_\A^\al(\K):\text{$X$ is
$\tau$-semistable}\bigr\}\subset\fObj_\A(\K),\\
\Ost^\al(\tau)&=\bigl\{[X]\in\fObj_\A^\al(\K): \text{$X$ is
$\tau$-stable}\bigr\}.
\end{align*}
For $\A$-data $(I,\pr,\ka)$, define $\Mss(I,\pr,\ka,\tau)_\A=\bigl\{
[(\si,\io,\pi)]\in\M(I,\pr,\ka)_\A:\si(\{i\})$ is $\tau$-semistable
for all $i\in I\bigr\}$. Then $\Oss^\al,\Ost^\al(\tau)$ and
$\Mss(I,\pr,\ka,\tau)_\A$ are open sets in the natural topologies,
and so are {\it locally constructible sets} in the stacks
$\fObj_\A,\fM(I,\pr,\ka)_\A$. Write $\dss^\al,\dst^\al(\tau):
\fObj_\A(\K)\ra\{0,1\}$ and $\dss(I,\pr,\ka,\tau):\M(I,\pr,\ka)_\A
\ra\{0,1\}$ for their characteristic functions. Define local stack
functions $\bdss^\al(\tau)=\bde_{\Oss^\al(\tau)}$ and
$\bdss(I,\pr,\ka,\tau)=\bde_{\Mss(I,\pr,\ka,\tau)_\A}$, using the
local generalization \cite[Def.~3.10]{Joyc4} of Definition
\ref{an2def5}. Then
\e
\begin{aligned}
\dss^\al,\dst^\al(\tau)&\in\dLCF(\fObj_\A),&
\dss(I,\pr,\ka,\tau)&\in\LCF\bigl(\fM(I,\pr,\ka)_\A\bigr),\\
\bdss^\al(\tau)&\in\dLSFa(\fObj_\A),&
\bdss(I,\pr,\ka,\tau)&\in\LSF\bigl(\fM(I,\pr,\ka)_\A\bigr).
\end{aligned}
\label{an3eq5}
\e
\label{an3def11}
\end{dfn}

\begin{dfn} Let Assumption \ref{an3ass} hold and $(\tau,T,\le)$ be
a weak stability condition on $\A$. We call $(\tau,T,\le)$ {\it
permissible} if:
\begin{itemize}
\setlength{\itemsep}{0pt}
\setlength{\parsep}{0pt}
\item[(i)] $\A$ is $\tau$-artinian, in the sense of Definition
\ref{an3def10}, and
\item[(ii)] $\Oss^\al(\tau)$ is a constructible subset in
$\fObj_\A^\al$ for all~$\al\in C(\A)$.
\end{itemize}
Examples of (weak) stability conditions on $\A=\modKQ$ and
$\A=\coh(P)$ are given in~\cite[\S 4.3--\S 4.4]{Joyc7}.
\label{an3def12}
\end{dfn}

\begin{dfn} Let $(\tau,T,\le)$ and $(\ti\tau,\ti T,\le)$ be
weak stability conditions on an abelian category $\A$, with
the same $K(\A)$. We say $(\ti\tau,\ti T,\le)$ {\it dominates}
$(\tau,T,\le)$ if $\tau(\al)\le\tau(\be)$ implies $\ti\tau(\al)
\le\ti\tau(\be)$ for all~$\al,\be\in C(\A)$.
\label{an3def13}
\end{dfn}

If $(\tau,T,\le)$ is permissible then \cite[Th.~4.8]{Joyc7} implies
$\Mss(I,\pr,\ka,\tau)_\A$ is constructible. Together with Definition
\ref{an3def12}(ii) this gives, following \eq{an3eq5},
\e
\begin{aligned}
\dss^\al,\dst^\al(\tau)&\in\CF(\fObj_\A),&
\dss(I,\pr,\ka,\tau)&\in\CF\bigl(\fM(I,\pr,\ka)_\A\bigr),\\
\bdss^\al(\tau)&\in\SFa(\fObj_\A),&
\bdss(I,\pr,\ka,\tau)&\in\SF\bigl(\fM(I,\pr,\ka)_\A\bigr).
\end{aligned}
\label{an3eq6}
\e
Now \cite{Joyc6} also studies other locally constructible sets:
$\Osi^\al(\tau)$ of $\tau$-semistable indecomposable objects in
class $\al$, and $\Msi,\Mst,\Mssb,\Msib,\Mstb(I,\pr,\ka,\tau)_\A$ of
(best) $(I,\pr,\ka)$-configurations $(\si,\io,\pi)$ whose smallest
objects $\si(\{i\})$ for $i\in I$ are $\tau$-(semi)stable
(indecomposable). Write $\dsi^\al(\tau)$,
$\dsi,\dst,\dssb,\dsib,\dstb(I,\pr,\ka,\tau)$ for their
characteristic functions. We also define local stack function
versions $\bdsi^\al(\tau)$ and $\bdsi,\ldots,\bdstb
(I,\pr,\ka,\tau)$. When $(\tau,T,\le)$ is permissible these sets and
functions are all constructible, and the local stack functions are
stack functions.

In \cite[\S 5--\S 8]{Joyc7} we prove many identities relating these
functions under pushforwards. For example, if $(\tau,T,\le)$ is a
permissible stability condition then
\e
\begin{aligned}
\dst^\al(\tau)=
\sum_{\substack{\text{iso. classes}\\ \text{of finite sets $I$}}}
\frac{1}{\md{I}!}\,\cdot &\sum_{\begin{subarray}{l}
\text{$\pr,\ka$: $(I,\pr,\ka)$ is $\A$-data,}\\
\text{$\ka(I)=\al$, $\tau\ci\ka\equiv\tau(\al)$} \end{subarray}
\!\!\!\!\!\!\!\!\!\!\!\!\!\!\!\!\!\!\!\!\!\!\!\!\!\!\!}
\CF^\stk(\bs\si(I))\dss(I,\pr,\ka,\tau)\cdot\\
&
\raisebox{-6pt}{\begin{Large}$\displaystyle\Bigl[$\end{Large}}
\sum_{\text{p.o.s $\ls$ on $I$ dominating $\pr$}}
\!\!\!\!\!\!\!\!\!\!\!\!\!\!\! n(I,\pr,\ls)\,N(I,\ls)
\raisebox{-6pt}{\begin{Large}$\displaystyle\Bigr]$\end{Large}},
\end{aligned}
\label{an3eq7}
\e
with only finitely many nonzero terms in the sum. This follows from
\cite[Th.s 6.3 \& 6.12]{Joyc7}, with integers $n(I,\pr,\ls),
N(I,\ls)$ defined in~\cite[Def.s 6.1 \& 6.9]{Joyc7}.

\subsection{Algebras $\Ht_\tau,\Hp_\tau,\bHt_\tau,\bHp_\tau$
and elements $\ep^\al(\tau),\bep^\al(\tau)$}
\label{an35}

In \cite[\S 7--\S 8]{Joyc7} we study {\it subalgebras} $\Hp_\tau,
\Ht_\tau,\bHp_\tau,\bHt_\tau$ in~$\CF,\SFa(\fObj_\A)$.

\begin{dfn} Let Assumption \ref{an3ass} hold, and $(\tau,T,\le)$
be a permissible weak stability condition on $\A$. As in \cite[Def.s
7.1 \& 8.4]{Joyc7}, define $\Hp_\tau,\Ht_\tau$ when $\cha\K=0$ and
$\bHp_\tau,\bHt_\tau$ for all $\K$ by
\ea
\Hp_\tau&=\bigl\langle\CF^\stk(\bs\si(I))\dss(I,\pr,\ka,\tau):
\text{$(I,\pr,\ka)$ is $\A$-data}\bigr\rangle{}_\Q
\subseteq\CF(\fObj_\A),
\label{an3eq8}\\
\Ht_\tau&=\bigl\langle\de_{[0]},\dss^{\al_1}(\tau)*\cdots*
\dss^{\al_n}(\tau):\al_1,\ldots,\al_n\in C(\A)\bigr\rangle{}_\Q
\subseteq\CF(\fObj_\A),
\label{an3eq9}\\
\bHp_\tau&=\bigl\langle\bs\si(I)_*\,\bdss(I,\pr,\ka,
\tau):\text{$(I,\pr,\ka)$ is $\A$-data}\bigr\rangle{}_\Q
\subseteq\SFa(\fObj_\A),
\label{an3eq10}\\
\bHt_\tau&=\bigl\langle\bde_{[0]},
\bdss^{\al_1}(\tau)*\cdots*\bdss^{\al_n}(\tau):
\al_1,\ldots,\al_n\in C(\A)\bigr\rangle{}_\Q
\subseteq\SFa(\fObj_\A).
\label{an3eq11}
\ea
Here $\langle\cdots\rangle_\Q$ is the set of all finite $\Q$-linear
combinations of the elements~`$\,\cdots$'.

To relate $\Hp_\tau,\bHp_\tau$ and $\Ht_\tau,\bHt_\tau$, let
$(\{1,\ldots,n\},\le,\ka)$ be $\A$-data. Then
\e
\begin{split}
&\dss^{\ka(1)}(\tau)*\!\cdots\!*\dss^{\ka(n)}(\tau)
=\CF^\stk(\bs\si(\{1,\ldots,n\}))\dss(\{1,\ldots,n\},\le,\ka,\tau),\\
&\bdss^{\ka(1)}(\tau)*\!\cdots\!*\bdss^{\ka(n)}(\tau)
=\bs\si(\{1,\ldots,n\})_*\,\bdss(\{1,\ldots,n\},\le,\ka,\tau).
\end{split}
\label{an3eq12}
\e
Thus $\Hp_\tau,\bHp_\tau$ are the spans of $\CF^\stk(\bs\si(I))
\dss(I,\pr,\ka,\tau),\bs\si(I)_*\,\bdss(I,\pr,\ka,\tau)$ for
$\A$-data $(I,\pr,\ka)$ with $\pr$ a {\it partial order}, and
$\Ht_\tau,\bHt_\tau$ the spans with $\pr$ a {\it total order}. Hence
$\Ht_\tau\!\subseteq\!\Hp_\tau$ and $\bHt_\tau\!\subseteq\!
\bHp_\tau$. By \cite[Prop.~7.2 \& Th.~8.5]{Joyc7},
$\Hp_\tau,\Ht_\tau$ are subalgebras of $\CF(\fObj_\A)$ and
$\bHp_\tau,\bHt_\tau$ subalgebras of $\SFa(\fObj_\A)$, and
$\pi^\stk_{\fObj_\A}$ induces surjective morphisms
$\bHp_\tau\ra\Hp_\tau$ and~$\bHt_\tau \ra\Ht_\tau$.
\label{an3def14}
\end{dfn}

Equation \eq{an3eq7} implies that $\dst^\al(\tau)\in\Hp_\tau$. Other
identities in \cite[\S 5--\S 8]{Joyc7} imply that the five families
of functions $\CF^\stk(\bs\si(I))\dsi,\dst,\dssb,\dsib,\dstb
(*,\tau)$ lie in $\Hp_\tau$ and are alternative spanning sets, so
that the identities yield {\it basis change formulae} in $\Hp_\tau$,
and similarly for $\bHp_\tau$. One moral is that $\Hp_\tau,
\bHp_\tau$ contain information on both $\tau$-stability and
$\tau$-semistability, and so are good tools for studying invariants
counting $\tau$-stable and $\tau$-semistable objects, whereas the
smaller algebras $\Ht_\tau,\bHt_\tau$ only really contain
information about $\tau$-semistability.

In \cite[Def.s 7.6 \& 8.1]{Joyc7} we define interesting
elements~$\ep^\al(\tau),\bep^\al(\tau)$.

\begin{dfn} Let Assumption \ref{an3ass} hold, and $(\tau,T,\le)$
be a permissible weak stability condition on $\A$. For $\al\in
C(\A)$ define $\ep^\al(\tau)$ in $\CF(\fObj_\A)$ for $\cha\K=0$, and
$\bep^\al(\tau)$ in $\SFa(\fObj_\A)$ for all $\K$, by
\ea
\ep^\al(\tau)&=
\!\!\!\!\!\!\!\!\!\!\!\!\!\!\!\!
\sum_{\substack{\text{$\A$-data $(\{1,\ldots,n\},\le,\ka):$}\\
\text{$\ka(\{1,\ldots,n\})=\al$, $\tau\ci\ka\equiv\tau(\al)$}}}
\!\!\!\!\!\!\!
\frac{(-1)^{n-1}}{n}\,\,\dss^{\ka(1)}(\tau)*\dss^{\ka(2)}(\tau)*
\cdots*\dss^{\ka(n)}(\tau),
\label{an3eq13}\\
\bep^\al(\tau)&=
\!\!\!\!\!\!\!\!\!\!\!\!\!\!\!\!
\sum_{\substack{\text{$\A$-data $(\{1,\ldots,n\},\le,\ka):$}\\
\text{$\ka(\{1,\ldots,n\})=\al$, $\tau\ci\ka\equiv\tau(\al)$}}}
\!\!\!\!\!\!\!
\frac{(-1)^{n-1}}{n}\,\,\bdss^{\ka(1)}(\tau)*\bdss^{\ka(2)}(\tau)*
\cdots*\bdss^{\ka(n)}(\tau).
\label{an3eq14}
\ea
Then $\pi^\stk_{\fObj_\A}(\bep^\al(\tau))=\ep^\al(\tau)$, and
\cite[Th.s 7.8 \& 8.7]{Joyc7} show $\ep^\al(\tau)\in \CFi(\fObj_\A)$
and $\bep^\al(\tau)\in\SFai(\fObj_\A)$. In \cite[Th.s 7.7 \&
8.2]{Joyc7} we {\it invert\/} \eq{an3eq13} and \eq{an3eq14}, giving
\ea
\dss^\al(\tau)= \!\!\!\!\!\!\!\!\!\!\!\!\!\!\!\!
\sum_{\substack{\text{$\A$-data $(\{1,\ldots,n\},\le,\ka):$}\\
\text{$\ka(\{1,\ldots,n\})=\al$, $\tau\ci\ka\equiv\tau(\al)$}}}
\!\!\!\! \frac{1}{n!}\,\,\ep^{\ka(1)}(\tau)*\ep^{\ka(2)}(\tau)*
\cdots*\ep^{\ka(n)}(\tau),
\label{an3eq15}\\
\bdss^\al(\tau)=
\!\!\!\!\!\!\!\!\!\!\!\!\!\!\!\!
\sum_{\substack{\text{$\A$-data $(\{1,\ldots,n\},\le,\ka):$}\\
\text{$\ka(\{1,\ldots,n\})=\al$, $\tau\ci\ka\equiv\tau(\al)$}}}
\!\!\!\! \frac{1}{n!}\,\,\bep^{\ka(1)}(\tau)*\bep^{\ka(2)}(\tau)*
\cdots*\bep^{\ka(n)}(\tau).
\label{an3eq16}
\ea
There are {\it only finitely many nonzero terms} in
\eq{an3eq13}--\eq{an3eq16}. We have
\e
\ep^\al(\tau)\bigl([X])=
\begin{cases} 1, & \text{$X$ is $\tau$-stable,} \\
\text{in $\Q$,} & \text{$X$ is strictly $\tau$-semistable and
indecomposable,} \\
0, & \text{$X$ is $\tau$-unstable or decomposable,} \end{cases}
\label{an3eq17}
\e
so $\ep^\al(\tau)$ interpolates between $\dss^\al,\dst^\al(\tau)$.
In \cite[eq.s (95) \& (123)]{Joyc7} we show that
\ea
\Ht_\tau&=\bigl\langle\de_{[0]},\ep^{\al_1}(\tau)*\cdots*
\ep^{\al_n}(\tau):\al_1,\ldots,\al_n\in C(\A)\bigr\rangle{}_\Q,
\label{an3eq18}\\
\bHt_\tau&=\bigl\langle\bde_{[0]},\bep^{\al_1}(\tau)*\cdots*
\bep^{\al_n}(\tau):\al_1,\ldots,\al_n\in C(\A)\bigr\rangle{}_\Q.
\label{an3eq19}
\ea
Thus $\ep^\al,\bep^\al(\tau)$ are {\it alternative generators} for
$\Ht_\tau,\bHt_\tau$ in~$\CFi,\SFai(\fObj_\A)$.
\label{an3def15}
\end{dfn}

In \cite[Def.s 7.1 \& 8.9]{Joyc7} we define {\it Lie
algebras}~$\Lp_\tau,\Lt_\tau,\bLp_\tau,\bLt_\tau$.

\begin{dfn} Let Assumption \ref{an3ass} hold, and $(\tau,T,\le)$
be a permissible weak stability condition on $\A$. Define
$\Lp_\tau\!=\!\Hp_\tau\!\cap\!\CFi(\fObj_\A)$,
$\Lt_\tau\!=\!\Ht_\tau\!\cap\!\CFi(\fObj_\A)$ and
$\bLp_\tau\!=\!\bHp_\tau\!\cap\!\SFai(\fObj_\A)$. Then
$\Lp_\tau,\Lt_\tau,\bLp_\tau$ are {\it Lie subalgebras} of
$\CFi,\SFai(\fObj_\A)$. We also define $\bLt_\tau$ to be the Lie
subalgebra of $\SFai(\fObj_\A)$ generated by the $\bep^\al(\tau)$
for all $\al\in C(\A)$. Then $\Lt_\tau\subseteq\Lp_\tau$,
$\bLt_\tau\subseteq\bLp_\tau$, and $\pi^\stk_{\fObj_\A}$ induces
surjective morphisms $\bLp_\tau\ra\Lp_\tau$ and~$\bLt_\tau
\ra\Lt_\tau$.

In \cite[Prop.~7.5 \& Def.~8.9]{Joyc7} we show $\Lp_\tau,\bLp_\tau$
are spanned by elements $\CF^\stk(\bs\si(I))\dsib(I,\pr,\ka,\tau)$,
$\bs\si(I)_*\,\bdsib(I,\pr,\ka,\tau)$ for connected $(I,\pr)$, and
deduce that $\Lp_\tau,\bLp_\tau$ generate $\Hp_\tau,\bHp_\tau$ as
algebras. This induces surjective morphisms $\Phi^{\rm
pa}_\tau:U(\Lp_\tau)\!\ra\!\Hp_\tau$, $\bar\Phi^{\rm pa}_\tau:
U(\bLp_\tau)\!\ra\!\bHp_\tau$, where $U(\Lp_\tau),U(\bLp_\tau)$ are
the {\it universal enveloping algebras} of $\Lp_\tau,\bLp_\tau$.
Moreover $\Phi^{\rm pa}_\tau$ is an isomorphism.

In \cite[Cor.~7.9]{Joyc7} we show $\Lt_\tau$ is the Lie subalgebra
of $\CFi(\fObj_\A)$ generated by the $\ep^\al(\tau)$ for $\al\in
C(\A)$. So from \eq{an3eq18}--\eq{an3eq19} we see
$\Lt_\tau,\bLt_\tau$ generate $\Ht_\tau,\bHt_\tau$ as algebras,
giving surjective algebra morphisms $\Phi^{\rm
to}_\tau:U(\Lt_\tau)\ra\Ht_\tau$ and $\bar\Phi^{\rm
to}_\tau:U(\bLt_\tau)\ra\bHt_\tau$. Again, $\Phi^{\rm to}_\tau$ is
an isomorphism.
\label{an3def16}
\end{dfn}

\section{Transformation coefficients $S,T,U(*,\tau,\ti\tau)$}
\label{an4}

Let $\A$ satisfy Assumption \ref{an3ass} and $(\tau,T,\le),
(\ti\tau,\ti T,\le)$ be permissible weak stability conditions
on $\A$. In \S\ref{an5} we shall prove {\it transformation
laws} from $(\tau,T,\le)$ to $(\ti\tau,\ti T,\le)$ of the
form, for all $\al\in C(\A)$ and $\A$-data~$(K,\tl,\mu)$,
\begin{gather}
\begin{gathered}
\sum_{\substack{\text{$\A$-data $(\{1,\ldots,n\},\le,\ka):$}\\
\text{$\ka(\{1,\ldots,n\})=\al$}}}
\!\!\!\!\!\!\!\!\!\!\!\!
\begin{aligned}[t]
S(\{1,&\ldots,n\},\le,\ka,\tau,\ti\tau)\cdot\\
&\dss^{\ka(1)}(\tau)*\dss^{\ka(2)}(\tau)*\cdots*
\dss^{\ka(n)}(\tau)=\dss^\al(\ti\tau),
\end{aligned}
\end{gathered}
\label{an4eq1}\\
\begin{aligned}
\sum_{\substack{\text{iso.}\\ \text{classes}\\
\text{of finite }\\ \text{sets $I$}}} \frac{1}{\md{I}!}\cdot
\sum_{\substack{
\text{$\pr,\ka,\phi$: $(I,\pr,\ka)$ is $\A$-data,}\\
\text{$(I,\pr,K,\phi)$ is dominant,}\\
\text{$\tl=\P(I,\pr,K,\phi)$,}\\
\text{$\ka(\phi^{-1}(k))=\mu(k)$ for $k\in K$}}}\,\,\,
\begin{aligned}[t]
&T(I,\pr,\ka,K,\phi,\tau,\ti\tau)\cdot\\
&\CF^\stk\bigl(Q(I,\pr,K,\tl,\phi)\bigr)\\
&\quad\dss(I,\pr,\ka,\tau)=\dss(K,\tl,\mu,\ti\tau),
\end{aligned}
\end{aligned}
\label{an4eq2}\\
\begin{gathered}
\sum_{\substack{\text{$\A$-data $(\{1,\ldots,n\},\le,\ka):$}\\
\text{$\ka(\{1,\ldots,n\})=\al$}}}
\!\!\!\!\!\!\!\!\!\!\!\!\!
\begin{aligned}[t]
U(\{1,&\ldots,n\},\le,\ka,\tau,\ti\tau)\cdot\\
&\ep^{\ka(1)}(\tau)*\ep^{\ka(2)}(\tau)*\cdots*
\ep^{\ka(n)}(\tau)=\ep^\al(\ti\tau),
\end{aligned}
\end{gathered}
\label{an4eq3}
\end{gather}
and analogues of these equations for stack functions
$\bdss^\al(\ti\tau),\bdss(K,\tl,\mu,\ti\tau),\bep^\al(\ti\tau)$.
Here $S,T,U(*,\tau,\ti\tau)$ are explicit {\it transformation
coefficients} in $\Q$, and $*$ is as in \S\ref{an33}. This
section will define and study the $S,T,U(*,\tau,\ti\tau)$, in
preparation for \S\ref{an5}. Its definitions and results are
all {\it combinatorial\/} in nature.

Suppose Assumption \ref{an3ass} holds, and $(\tau,T,\le)$ is a weak
stability condition on $\A$. Then from \S\ref{an3} the following
hold:

\begin{cond}(i) $K(\A)$ is an abelian group.
\begin{itemize}
\setlength{\itemsep}{0pt}
\setlength{\parsep}{0pt}
\item[(ii)] $C(\A)$ is a subset of $K(\A)$, closed under addition
and not containing zero.
\item[(iii)] a set of $\A$-{\it data} $(I,\pr,\ka)$ is by definition
a finite poset $(I,\pr)$ and a map $\ka:I\ra C(\A)$. For
$J\subseteq I$ we define~$\ka(J)=\sum_{j\in J}\ka(j)$.
\item[(iv)] $(T,\le)$ is a total order.
\item[(v)] $\tau:C(\A)\ra T$ is a map such that if $\al,\be,\ga\in C(\A)$
with $\be=\al+\ga$ then $\tau(\al)\!\le\!\tau(\be)\!\le\!\tau(\ga)$
or~$\tau(\al)\!\ge\!\tau(\be)\!\ge\!\tau(\ga)$.
\end{itemize}
More generally, if we have more than one weak stability
condition we shall say that `Condition \ref{an4cond} holds for
$(\ti\tau,\ti T,\le),(\hat\tau,\hat T,\le),\ldots$' if
$(\ti T,\le),(\hat T,\le),\ldots$ are total orders, and (v)
holds for $\ti\tau:C(\A)\ra\ti T$, $\hat\tau:C(\A)\ra\hat T,\ldots$.
\label{an4cond}
\end{cond}

These are the only properties of $\A,K(\A),(\tau,T,\le)$ that
we will use in this section. We shall not use Assumption
\ref{an3ass}, or suppose $(\tau,T,\le)$ is permissible.

\subsection{Basic definitions and main results}
\label{an41}

We begin by defining {\it transformation coefficients}
$S,T,U(*,\tau,\ti\tau)$.

\begin{dfn} Suppose Condition \ref{an4cond} holds for $(\tau,T,\le)$
and $(\ti\tau,\ti T,\le)$, and let $(\{1,\ldots,n\},\le,\ka)$ be
$\A$-data. If for all $i=1,\ldots,n-1$ we have either
\begin{itemize}
\setlength{\itemsep}{0pt}
\setlength{\parsep}{0pt}
\item[{\rm(a)}] $\tau\ci\ka(i)\le\tau\ci\ka(i+1)$ and
$\ti\tau\ci\ka(\{1,\ldots,i\})>\ti\tau\ci\ka(\{i+1,\ldots,n\})$ or
\item[{\rm(b)}] $\tau\ci\ka(i)>\tau\ci\ka(i+1)$ and~$\ti\tau
\ci\ka(\{1,\ldots,i\})\le\ti\tau\ci\ka(\{i+1,\ldots,n\})$,
\end{itemize}
then define $S(\{1,\ldots,n\},\le,\ka,\tau,\ti\tau)=(-1)^r$,
where $r$ is the number of $i=1,\ldots,n-1$ satisfying (a).
Otherwise define~$S(\{1,\ldots,n\},\le,\ka,\tau,\ti\tau)=0$.

If $(I,\pr,\ka)$ is $\A$-data with $\pr$ a {\it total order},
there is a unique bijection $\phi:\{1,\ldots,n\}\ra I$ with
$n=\md{I}$ and $\phi_*(\le)=\pr$, and $(\{1,\ldots,n\},\le,
\ka\ci\phi)$ is $\A$-data. Define~$S(I,\pr,\ka,\tau,\ti\tau)
=S(\{1,\ldots,n\},\le,\ka\ci\phi,\tau,\ti\tau)$.
\label{an4def1}
\end{dfn}

\begin{dfn} Let $(I,\pr)$ be a finite poset, $K$ a
finite set, and $\phi:I\ra K$ a surjective map. We
call $(I,\pr,K,\phi)$ {\it dominant\/} if there exists
a partial order $\tl$ on $K$ such that $(\phi^{-1}(\{k\}),
\pr)$ is a total order for all $k\in K$, and if
$i,j\in I$ with $\phi(i)\ne\phi(j)$ then $i\pr j$ if
and only if $\phi(i)\tl\phi(j)$. Then $\tl$ is determined
uniquely by $(I,\pr,K,\phi)$, and we write $\tl=\P(I,\pr,K,\phi)$.
Note that $i\pr j$ implies $\phi(i)\tl\phi(j)$ for $i,j\in I$. We
use the notation `dominant' as if $\ls$ is a partial order on
$I$ dominating $\pr$ with $i\ls j$ implies $\phi(i)\tl\phi(j)$,
then $\ls=\pr$. That is, the partial order $\pr$ is as strong as
it can be, given~$I,K,\phi,\tl$.

Now suppose Condition \ref{an4cond} holds for $(\tau,T,\le)$
and $(\ti\tau,\ti T,\le)$, and let $(I,\pr,\ka)$ be $\A$-data
and $(I,\pr,K,\phi)$ be dominant. Define
\e
T(I,\pr,\ka,K,\phi,\tau,\ti\tau)=\ts\prod_{k\in K}S\bigl(\phi^{-1}
(\{k\}),\pr\vert_{\phi^{-1}(\{k\})},\ka\vert_{\phi^{-1}(\{k\})},\tau,
\ti\tau\bigr).
\label{an4eq4}
\e
\label{an4def2}
\end{dfn}

\begin{dfn} Suppose Condition \ref{an4cond} holds for $(\tau,T,\le)$
and $(\ti\tau,\ti T,\le)$, and let $(\{1,\ldots,n\},\le,\ka)$ be
$\A$-data. Define
\e
\begin{gathered}
U(\{1,\ldots,n\},\le,\ka,\tau,\ti\tau)=
\sum_{1\le l\le m\le n}\\
\sum_{\substack{
\text{surjective $\psi:\{1,\ldots,n\}\!\ra\!\{1,\ldots,m\}$}\\
\text{and\/ $\xi:\{1,\ldots,m\}\!\ra\!\{1,\ldots,l\}$:}\\
\text{$i\!\le\!j$ implies $\psi(i)\!\le\!\psi(j)$,
$i\!\le\!j$ implies $\xi(i)\!\le\!\xi(j)$.}\\
\text{Define $\la:\{1,\ldots,m\}\ra C(\A)$ by $\la(b)=\ka(\psi^{-1}(b))$.}\\
\text{Define $\mu:\{1,\ldots,l\}\ra C(\A)$ by $\mu(a)=\la(\xi^{-1}(a))$.}\\
\text{Then $\tau\ci\ka\equiv\tau\ci\la\ci\mu:I\ra T$ and\/
$\ti\tau\ci\mu\equiv\ti\tau(\al)$}}}
\!\!\!\!\!\!\!\!
\begin{aligned}[t]
\prod_{a=1}^lS(\xi^{-1}(\{a\}),\le,\la,\tau,\ti\tau)\cdot&\\
\frac{(-1)^{l-1}}{l}\cdot\prod_{b=1}^m\frac{1}{\md{\psi^{-1}(b)}!}\,&.
\end{aligned}
\end{gathered}
\label{an4eq5}
\e

If $(I,\pr,\ka)$ is $\A$-data with $\pr$ a {\it total order},
there is a unique bijection $\phi:\{1,\ldots,n\}\ra I$ with
$n=\md{I}$ and $\phi_*(\le)=\pr$, and $(\{1,\ldots,n\},\le,
\ka\ci\phi)$ is $\A$-data. Define~$U(I,\pr,\ka,\tau,\ti\tau)
=U(\{1,\ldots,n\},\le,\ka\ci\phi,\tau,\ti\tau)$.
\label{an4def3}
\end{dfn}

Then $S,T(*,\tau,\ti\tau)$ are 1,0 or $-1$, and $U(*,\tau,\ti\tau)$
lies in $\Q$. Here is our main result on the properties of the
$S(*,\tau,\ti\tau)$. It is easy to see \eq{an4eq6} and \eq{an4eq7}
are necessary if \eq{an4eq1} is to hold: \eq{an4eq6} means \eq{an4eq1}
reduces to $\dss^\al(\tau)=\dss^\al(\tau)$ when $(\ti\tau,\ti T,\le)=
(\tau,T,\le)$, and \eq{an4eq7} is the condition for transforming from
$(\tau,T,\le)$ to $(\hat\tau,\hat T,\le)$ and then from $(\hat\tau,\hat T,
\le)$ to $(\ti\tau,\ti T,\le)$ using \eq{an4eq1} to give the same
answer as transforming from $(\tau,T,\le)$ directly to~$(\ti\tau,
\ti T,\le)$.

\begin{thm} Let Condition \ref{an4cond} hold for
$(\tau,T,\le),(\hat\tau,\hat T,\le)$ and\/ $(\ti\tau,\ti T,\le)$.
Suppose $(\{1,\ldots,n\},\le,\ka)$ is $\A$-data. Then
\begin{gather}
S(\{1,\ldots,n\},\le,\ka,\tau,\tau)
=\begin{cases} 1, & \text{$n=1$,} \\
0, & \text{otherwise,} \end{cases}
\label{an4eq6}
\\
\begin{aligned}
\sum_{m=1}^n\,
\sum_{\substack{\text{$\psi:\{1,\ldots,n\}\ra\{1,\ldots,m\}:$}\\
\text{$\psi$ is surjective,}\\
\text{$1\!\le\!i\!\le\!j\!\le\!n$ implies $\psi(i)\!\le\!\psi(j)$,}\\
\text{define $\la:\{1,\ldots,m\}\!\ra\!C(\A)$}\\
\text{by $\la(k)=\ka(\psi^{-1}\!(k))$}}}\,
\begin{aligned}[t]
&S(\{1,\ldots,m\},\le,\la,\hat\tau,\ti\tau)\cdot\\
&\prod_{k=1}^mS\bigl(\psi^{-1}(\{k\}),\le,
\ka\vert_{\psi^{-1}(\{k\})},\tau,\hat\tau)\\
&\qquad =S(\{1,\ldots,n\},\le,\ka,\tau,\ti\tau).
\end{aligned}
\end{aligned}
\label{an4eq7}
\end{gather}
\label{an4thm1}
\end{thm}

We also give a criterion for when $S(\{1,\ldots,n\},\le,
\ka,\tau,\ti\tau)\ne 0$ based on the minimum and maximum
values of the~$\tau\ci\ka(i)$.

\begin{thm} Suppose Condition \ref{an4cond} holds for
$(\tau,T,\le)$ and\/ $(\ti\tau,\ti T,\le)$, and\/
$(\{1,\ldots,n\},\le,\ka)$ is $\A$-data for $n\ge 1$ with\/
$S(\{1,\ldots,n\},\le,\ka,\tau,\ti\tau)\ne 0$. Then
there exist\/ $k,l=1,\ldots,n$ such that\/ $\tau\ci\ka(k)
\le\tau\ci\ka(i)\le\tau\ci\ka(l)$ for all\/ $i=1,\ldots,n$,
and\/~$\ti\tau\ci\ka(k)\ge\ti\tau\ci\ka(\{1,\ldots,n\})
\ge\ti\tau\ci\ka(l)$.
\label{an4thm2}
\end{thm}

Here are the analogues of Theorem \ref{an4thm1} for the
$T,U(*,\tau,\ti\tau)$. Again, it is easy to see
\eq{an4eq8}--\eq{an4eq9} and \eq{an4eq10}--\eq{an4eq11} are
necessary if \eq{an4eq2} and \eq{an4eq3} are to hold.

\begin{thm} Let Condition \ref{an4cond} hold for
$(\tau,T,\le),(\hat\tau,\hat T,\le)$ and\/ $(\ti\tau,\ti T,\le)$.
Suppose $(I,\pr,\ka)$ is $\A$-data, and\/ $\phi:I\ra K$ is
surjective with\/ $(I,\pr,K,\phi)$ dominant. Then
\begin{gather}
T(I,\pr,\ka,K,\phi,\tau,\tau)
=\begin{cases} 1, & \text{$\phi$ is a bijection,} \\
0, & \text{otherwise,} \end{cases}
\label{an4eq8}
\\
\begin{aligned}
\sum_{\substack{\text{iso.}\\ \text{classes}\\
\text{of finite }\\ \text{sets $J$}}}
\!\!\! \frac{1}{\md{J}!}\cdot \!\!\!
\sum_{\substack{
\text{$\psi:I\!\ra\!J$, $\xi:J\!\ra\!K$ surjective:
$\phi\!=\!\xi\!\ci\!\psi$,}\\
\text{$(I,\pr,J,\psi)$ is dominant, $\ls\!=\!\P(I,\pr,J,\psi)$,}\\
\text{define $\la:J\!\ra\!K(\A)$ by $\la(j)\!=\!\ka(\psi^{-1}\!(j))$}}}
\,
\begin{aligned}[t]
&T(I,\pr,\ka,J,\psi,\tau,\hat\tau)\cdot\\[-3pt]
&T(J,\ls,\la,K,\xi,\hat\tau,\ti\tau)\!=\\[-3pt]
&\;\> T(I,\pr,\ka,K,\phi,\tau,\ti\tau).
\end{aligned}
\end{aligned}
\label{an4eq9}
\end{gather}
\label{an4thm3}
\end{thm}

This follows immediately from Theorem \ref{an4thm1}: to get
\eq{an4eq8} and \eq{an4eq9} we take the product over $k\in K$
of equations \eq{an4eq6} and \eq{an4eq7} with $(\{1,\ldots,n\},
\le,\ka)$ replaced by $(\phi^{-1}(\{k\}),\pr\vert_{\phi^{-1}
(\{k\})},\ka\vert_{\phi^{-1}(\{k\})})$, and use \eq{an4eq4} and
some simple combinatorics. We leave the details to the reader.

\begin{thm} Let Condition \ref{an4cond} hold for
$(\tau,T,\le),(\hat\tau,\hat T,\le)$ and\/ $(\ti\tau,\ti T,\le)$.
Suppose $(\{1,\ldots,n\},\le,\ka)$ is $\A$-data. Then
\ea
U(\{1,\ldots,n\},\le,\ka,\tau,\tau)
=\begin{cases} 1, & \text{$n=1$,} \\
0, & \text{otherwise,} \end{cases}
\label{an4eq10}
\\
\begin{aligned}
\sum_{m=1}^n\,
\sum_{\substack{\text{$\psi:\{1,\ldots,n\}\ra\{1,\ldots,m\}:$}\\
\text{$\psi$ is surjective,}\\
\text{$i\!\le\!j$ implies $\psi(i)\!\le\!\psi(j)$,}\\
\text{define $\la:\{1,\ldots,m\}\!\ra\!C(\A)$}\\
\text{by $\la(k)=\ka(\psi^{-1}\!(k))$}}}\,
\begin{aligned}[t]
&U(\{1,\ldots,m\},\le,\la,\hat\tau,\ti\tau)\cdot\\
&\ts\prod_{k=1}^mU\bigl(\psi^{-1}(\{k\}),\le,
\ka\vert_{\psi^{-1}(\{k\})},\tau,\hat\tau)\\
&\qquad =U(\{1,\ldots,n\},\le,\ka,\tau,\ti\tau).
\end{aligned}
\end{aligned}
\label{an4eq11}
\ea
\label{an4thm4}
\end{thm}

Theorem \ref{an4thm4} follows directly from Theorem \ref{an4thm1}
and the following proposition, which is a combinatorial consequence
of the proof in \cite[Th.~7.7]{Joyc7} that \eq{an3eq15} is the
inverse of~\eq{an3eq13}.

\begin{prop} Let\/ $1\le l\le n$ and\/ $\phi:\{1,\ldots,n\}\ra
\{1,\ldots,l\}$ be surjective with\/ $1\le i\le j\le n$ implies
$\phi(i)\le\phi(j)$. Then
\begin{gather*}
\sum_{m=l}^n
\sum_{\substack{
\text{$\psi:\{1,\ldots,n\}\!\ra\!\{1,\ldots,m\}$}\\
\text{and $\xi:\{1,\ldots,m\}\!\ra\!\{1,\ldots,l\}$}\\
\text{surjective with\/ $\phi=\xi\ci\psi$:}\\
\text{$i\!\le\!j$ implies $\psi(i)\!\le\!\psi(j)$,}\\
\text{$i\!\le\!j$ implies $\xi(i)\!\le\!\xi(j)$}}}
\begin{aligned}[t]
\prod_{a=1}^l\frac{1}{\md{\xi^{-1}(a)}!}
&\cdot\prod_{b=1}^m\frac{(-1)^{\md{\psi^{-1}(b)}-1}}{\md{\psi^{-1}(b)}}\\
=\,&\begin{cases} 1, & \text{$l=n$,}
\\ 0, & \text{otherwise,}\end{cases}
\end{aligned}
\\
\sum_{m=l}^n
\sum_{\substack{
\text{$\psi:\{1,\ldots,n\}\!\ra\!\{1,\ldots,m\}$}\\
\text{and $\xi:\{1,\ldots,m\}\!\ra\!\{1,\ldots,l\}$}\\
\text{surjective with\/ $\phi=\xi\ci\psi$:}\\
\text{$i\!\le\!j$ implies $\psi(i)\!\le\!\psi(j)$,}\\
\text{$i\!\le\!j$ implies $\xi(i)\!\le\!\xi(j)$}}}
\begin{aligned}[t]
\prod_{a=1}^l\frac{(-1)^{\md{\xi^{-1}(a)}-1}}{\md{\xi^{-1}(a)}}
&\cdot\prod_{b=1}^m\frac{1}{\md{\psi^{-1}(b)}!}\\
=\,&\begin{cases} 1, & \text{$l=n$,}
\\ 0, & \text{otherwise.}\end{cases}
\end{aligned}
\end{gather*}
\label{an4prop1}
\end{prop}

In the remainder of the section we prove Theorems \ref{an4thm1}
and~\ref{an4thm2}.

\subsection{Proof of Theorem \ref{an4thm2} and equation \eq{an4eq6}}
\label{an42}

In the situation of Theorem \ref{an4thm2}, as $S(\{1,\ldots,
n\},\le,\ka,\tau,\ti\tau)\ne 0$ one of Definition \ref{an4def1}(a)
or (b) holds for each $i=1,\ldots,n-1$. If $n=1$ the result is
trivial, so suppose $n>1$. Let $k=1,\ldots,n$ be the unique value
such that $\tau\ci\ka(k)$ is minimal amongst all $\tau\ci\ka(i)$ in
the total order $(T,\le)$, and $k$ is least with this condition.
Then $\tau\ci\ka(k)\le\tau\ci\ka(i)$ for all $i=1,\ldots,n$. To
prove $\ti\tau\ci\ka(k)\ge\ti\tau\ci\ka(\{1,\ldots,n\})$, we divide
into three cases (i) $k=1$, (ii) $k=n$, and (iii)~$1<k<n$.

In case (i) we have $\tau\ci\ka(1)\le\tau\ci\ka(2)$, which
excludes (b) for $i=1$. Hence (a) holds, giving $\ti\tau\ci
\ka(\{1\})>\ti\tau\ci\ka(\{2,\ldots,n\})$, so $\ti\tau\ci\ka(1)
\ge\ti\tau\ci\ka(\{1,\ldots,n\})$ by Condition \ref{an4cond}(v),
as we want. In (ii) we have $\tau\ci\ka(n-1)>\tau\ci\ka(n)$,
since $k=n$ is least with the minimal value of $\tau\ci\ka(i)$.
This excludes (a) for $i=n-1$, so (b) holds, giving $\ti\tau\ci
\ka(\{1,\ldots,n-1\})\le\ti\tau\ci\ka(\{n\})$. Condition
\ref{an4cond}(v) then gives $\ti\tau\ci\ka(n)\ge\ti\tau\ci
\ka(\{1,\ldots,n\})$, as we want.

Case (iii) implies that $\tau\ci\ka(k-1)>\tau\ci\ka(k)\le\tau
\ci\ka(k+1)$, since $k$ is least with the minimal value of
$\tau\ci\ka(i)$. Therefore (b) holds for $i=k-1$ and (a) holds for
$i=k$, giving $\ti\tau\ci\ka(\{1,\ldots,k-1\})
\le\ti\tau\ci\ka(\{1,\ldots,n\})<\ti\tau\ci\ka(\{1,\ldots,k\})$ by
Condition \ref{an4cond}(v). As $\ti\tau\ci\ka(\{1,\ldots,k-1\})
<\ti\tau\ci\ka(\{1,\ldots,k\})$ Condition \ref{an4cond}(v) gives
$\ti\tau\ci\ka(\{1,\ldots,k\})\le\ti\tau\ci\ka(k)$, and so $\ti\tau
(k)\ge\ti\tau\ci\ka(\{1,\ldots,n\})$, as we want. For the second
part, let $l=1,\ldots,n$ be the unique value such that $\tau\ci
\ka(l)$ is maximal amongst all $\tau\ci\ka(i)$ in the total order
$(T,\le)$, and $l$ is greatest with this condition, and argue in the
same way. This proves Theorem~\ref{an4thm2}.

The first line of \eq{an4eq6} is immediate from Definition
\ref{an4def1}. For the second, suppose for a contradiction
that $n>1$ and $S(\{1,\ldots,n\},\le,\ka,\tau,\tau)\ne 0$.
Then Theorem \ref{an4thm2} gives $k,l$ with $\tau\ci\ka(k)
\le\tau\ci\ka(i)\le\tau\ci\ka(l)$ for all $i$ and $\tau\ci
\ka(k)\ge\tau\ci\ka(\{1,\ldots,n\})\ge\tau\ci\ka(l)$. Thus
$\tau\ci\ka(k)=\tau\ci\ka(l)$ and all $\tau\ci\ka(i)$ are
equal. Condition \ref{an4cond}(v) and induction on $i,j$
then implies that $\tau\ci\ka(\{i,\ldots,j\})$ are equal
for all $1\le i\le j\le n$. So neither of Definition
\ref{an4def1}(a),(b) apply for any $i$, as the strict
inequalities do not hold. This proves~\eq{an4eq6}.

\subsection{An alternative formula for $S(*,\tau,\ti\tau)$}
\label{an43}

We will need the following notation.

\begin{dfn} Let Condition \ref{an4cond} hold, and $(I,\pr,\ka)$
be $\A$-data. We say
\begin{itemize}
\setlength{\itemsep}{0pt}
\setlength{\parsep}{0pt}
\item[(i)] $(I,\pr,\ka)$ is $\tau$-{\it semistable} if $\tau(\ka(J))
\!\le\!\tau(\ka(I\sm J))$ for all $(I,\pr)$ s-sets~$J\!\ne\!\emptyset,I$.
\item[(ii)] $(I,\pr,\ka)$ is $\tau$-{\it reversing} if $i\npr j$
implies $\tau\ci\ka(i)<\tau\ci\ka(j)$ for~$i,j\in I$.
\end{itemize}
\label{an4def4}
\end{dfn}

Part (i) is based on Definition \ref{an3def9}, replacing subobjects
by s-sets. If $(\si,\io,\pi)$ is an $(I,\pr,\ka)$-configuration with
$\si(I)$ $\tau$-semistable, then Definition \ref{an3def9} implies
$(I,\pr,\ka)$ is $\tau$-semistable. In (ii), if $(I,\pr,\ka)$ is
$\tau$-reversing then $(I,\pr)$ is a {\it total order}, and
$i\mapsto\tau\ci\ka(i)$ {\it reverses} the order of~$(I,\pr)$.

\begin{prop} Suppose Condition \ref{an4cond} holds for
$(\tau,T,\le)$ and\/ $(\ti\tau,\ti T,\le),$ and\/
$(\{1,\ldots,n\},\le,\ka)$ is $\A$-data. Then
\e
\begin{gathered}
S(\{1,\ldots,n\},\le,\ka,\tau,\ti\tau)=\\
\sum_{1\le b\le a\le n}
\!\!\!\sum_{\substack{
\text{surjective $\al:\{1,\ldots,n\}\!\ra\!\{1,\ldots,a\}$,
$\be:\{1,\ldots,a\}\!\ra\!\{1,\ldots,b\}$:}\\
\text{$1\!\le\!i\!\le\!j\!\le\!n$ implies $\al(i)\!\le\!\al(j)$,
$1\!\le\!i\!\le\!j\!\le\!a$ implies $\be(i)\!\le\!\be(j)$,}\\
\text{$(\al^{-1}(j),\le,\ka)$ $\tau$-reversing, $j=1,\ldots,a$.
Define $\nu:\{1,\ldots,b\}\!\ra\!C(\A)$}\\
\text{by $\nu(k)=\ka((\be\ci\al)^{-1}(k))$. Then
$(\{1,\ldots,b\},\le,\nu)$ is $\ti\tau$-semistable}}}
\!\!\!\!\!\!\!\!\!\!\!\!\!\!\!\!\!\!\!\!\!\!\!\!\!\!\!
(-1)^{a-b}.
\end{gathered}
\label{an4eq12}
\e
\label{an4prop2}
\end{prop}

\begin{proof} The $\tau$-reversing, $\ti\tau$-semistable
conditions in \eq{an4eq12} may be rewritten:
\begin{itemize}
\setlength{\itemsep}{0pt}
\setlength{\parsep}{0pt}
\item $1\le i<n$ and $\al(i)=\al(i+1)$ implies
$\tau\ci\ka(i)>\tau\ci\ka(i+1)$,
\item $1\!\le\!i\!<\!n$ and $\be\ci\al(i)\!\ne\!\be\ci\al(i\!+\!1)$
implies $\ti\tau\ci\ka(\{1,\ldots,i\})\!\le\!\ti\tau\ci
\ka(\{i\!+\!1,\ldots,n\})$.
\end{itemize}
Suppose first that each $1\le i<n$ satisfies Definition
\ref{an4def1}(a) or (b), and $a,b,\al,\be$ are as in
\eq{an4eq12}. If $i$ satisfies (a) then $\al(i)\!\ne\!\al(i+1)$,
so $\al(i+1)\!=\!\al(i)+1$, and $\be\ci\al(i+1)\!=\!\be\ci\al(i)$.
If $i$ satisfies (b) there are three possibilities:
\begin{itemize}
\setlength{\itemsep}{0pt}
\setlength{\parsep}{0pt}
\item[(i)] $\al(i+1)=\al(i)$ and $\be\ci\al(i+1)=\be\ci\al(i)$,
\item[(ii)] $\al(i+1)=\al(i)+1$ and $\be\ci\al(i+1)=\be\ci\al(i)$, or
\item[(iii)] $\al(i+1)=\al(i)+1$ and $\be\ci\al(i+1)=\be\ci\al(i)+1$.
\end{itemize}
Let $r$ be the number of $i=1,\ldots,n-1$ satisfying (a), and
$s_1,s_2,s_3$ be the numbers of $i$ satisfying (b) and (i),
(ii) or (iii). Then~$r+s_1+s_2+s_3=n-1$.

It is not difficult to see that for all of the $n-1-r$ values of
$i$ satisfying (b), any of (i)--(iii) is possible independently,
and these choices of (i)--(iii) determine $a,b,\be,\al$ uniquely.
Hence there are $3^{n-1-r}$ possible quadruples $a,b,\al,\be$ in
\eq{an4eq12}. Furthermore, as $a-1$ is the number of $i$ with
$\al(i+1)=\al(i)+1$ and $b-1$ the number of $i$ with
$\be\ci\al(i+1)=\be\ci\al(i)+1$, we see that $a=1+r+s_2+s_3$
and $b=1+s_3$, so that $(-1)^{a-b}=(-1)^r\cdot(-1)^{s_2}$.
Therefore the sum over all $3^{n-1-r}$ possibilities of
$(-1)^r\cdot(-1)^{s_2}$ equals $(-1)^r\cdot(1-1+1)^{n-1-r}=(-1)^r$,
since $(1-1+1)^{n-1-r}$ is the sum over all $(n-1-r)$-tuples of
choices of (i), (ii) or (iii) of the product $(-1)^{s_2}$
of 1 for each choice of (i), $-1$ for each (ii) and 1 for
each (iii). Hence both sides of \eq{an4eq12} are~$(-1)^r$.

Now suppose some $i$ does not satisfy Definition \ref{an4def1}(a)
or (b). Then either
\begin{itemize}
\setlength{\itemsep}{0pt}
\setlength{\parsep}{0pt}
\item[(c)] $\tau\ci\ka(i)\le\tau\ci\ka(i+1)$ and $\ti\tau
\ci\ka(\{1,\ldots,i\})\le\ti\tau\ci\ka(\{i+1,\ldots,n\})$, or
\item[(d)] $\tau\ci\ka(i)>\tau\ci\ka(i+1)$ and~$\ti\tau
\ci\ka(\{1,\ldots,i\})>\ti\tau\ci\ka(\{i+1,\ldots,n\})$.
\end{itemize}
We shall show both sides of \eq{an4eq12} are zero in each case.

Suppose $i$ satisfies (c), and let $a,b,\al,\be$ be as in
\eq{an4eq12}. Then $\tau\ci\ka(i)\le\tau\ci\ka(i+1)$ implies
that $\al(i)\ne\al(i+1)$. Let $j=\al(i)$, so that
$j+1=\al(i+1)$, and suppose $\be(j)=\be(j+1)$. Define
$b'=b+1$ and $\be':\{1,\ldots,a\}\ra\{1,\ldots,b+1\}$ by
$\be'(k)=\be(j)$ for $1\le k\le j$, and $\be'(k)=\be(k)+1$
for $j+1\le k\le a$. Then $a,b',\al,\be'$ satisfy the
conditions in \eq{an4eq12}, with~$\be'(j)\ne\be'(j+1)$.

This establishes a 1-1 correspondence between quadruples
$(a,b,\al,\be)$ in \eq{an4eq12} with $\be\ci\al(i)=\be\ci
\al(i+1)$, and $(a,b',\al,\be')$ in \eq{an4eq12} with
$\be'\ci\al(i)\ne\be'\ci\al(i+1)$. Each such pair
contributes $(-1)^{a-b}+(-1)^{a-b'}=0$ to \eq{an4eq12}, as
$b'=b+1$, so both sides are zero. In the same way, if $i$
satisfies (d) then $\be\ci\al(i)=\be\ci\al(i+1)$ for
any $a,b,\al,\be$ in \eq{an4eq12}. We construct a 1-1
correspondence between quadruples $(a,b,\al,\be)$ in
\eq{an4eq12} with $\al(i)=\al(i+1)$ and $(a',b,\al',\be')$
with $\al(i)\ne\al(i+1)$, where $a'=a+1$. The contribution
of each pair to \eq{an4eq12} is zero, so both sides are zero.
This completes the proof.
\end{proof}

\subsection{Proof of equation \eq{an4eq7}}
\label{an44}

We begin with two propositions.

\begin{prop} Let Condition \ref{an4cond} hold, $(\{1,\ldots,b\},\le,
\mu)$ be $\A$-data, and\/ $\ga:\{1,\ldots,b\}\!\ra\!\{1,\ldots,c\}$
be surjective with\/ $i\!\le\!j$ implies $\ga(i)\!\le\!\ga(j)$. Then
there exist unique $m=c,\ldots,b$ and surjective $\phi:\{1,
\ldots,b\}\!\ra\!\{1,\ldots,m\}$ and\/ $\chi:\{1,\ldots,m\}\!\ra\!
\{1,\ldots,c\}$ with\/ $\ga=\chi\ci\phi$ and\/ $i\!\le\!j$ implies
$\phi(i)\!\le\!\phi(j)$ and\/ $\chi(i)\!\le\!\chi(j),$ such that
if\/ $\la:\{1,\ldots,m\}\ra C(\A)$ is given by $\la(i)=\mu\bigl(
\phi^{-1}(i)\bigr)$ then $(\phi^{-1}(i),\le,\mu)$ is
$\tau$-semistable for all\/ $i=1,\ldots,m,$ and\/
$(\chi^{-1}(j),\le,\la)$ is $\tau$-reversing for
all\/~$j=1,\ldots,c$.
\label{an4prop3}
\end{prop}

\begin{proof} First consider the case $c=1$. By induction
choose a unique sequence $k_0,\ldots,k_m$ with $0=k_0<k_1<
\cdots<k_m=b$, as follows. Set $k_0=0$. Having chosen $k_i$,
if $k_i=b$ then put $m=i$ and finish. Otherwise, let $k_{i+1}$
be largest such that $k_i<k_{i+1}\le b$ and $(\{k_i+1,k_i+2,
\ldots,k_{i+1}\},\le,\mu)$ is $\tau$-semistable. As $(\{k_i+1\},
\le,\mu)$ is $\tau$-semistable, this is the maximum of a
nonempty finite set, and $k_{i+1}$ is well-defined. Now
define $\phi:\{1,\ldots,b\}\ra\{1,\ldots,m\}$ by $\phi(j)=i$
if $k_{i-1}<j\le k_i$. Clearly $\phi$ is surjective with
$1\!\le\!i\!\le\!j\!\le\!b$ implies $\phi(i)\!\le\!\phi(j)$.
Define $\la$ as above, and~$\chi\equiv 1$.

By definition $\bigl(\phi^{-1}(i),\le,\mu\bigr)=(\{k_{i-1}+1,
\ldots,k_i\},\le,\mu)$ is $\tau$-{\it semistable} for all $i$,
as we want. We shall show $(\{1,\ldots,m\},\le,\la)$ is
$\tau$-{\it reversing}. Suppose for a contradiction that
$1\le i<m$ and $\tau\ci\la(i)\le\tau\ci\la(i+1)$. This implies that
\e
\text{
\begin{small}
$\tau\!\ci\!\mu(\{k_{i-1}\!+\!1,\ldots,k_i\})\!\le\!
\tau\!\ci\!\mu(\{k_{i-1}\!+\!1,\ldots,k_{i+1}\})\!\le\!
\tau\!\ci\!\mu(\{k_i\!+\!1,\ldots,k_{i+1}\}),$
\end{small}
}\!\!
\label{an4eq13}
\e
by Condition \ref{an4cond}(v). Now $(\{k_{i-1}+1,\ldots,k_{i+1}\},\le,
\mu)$ is not $\tau$-semistable by choice of $k_i$, so by Condition
\ref{an4cond}(v) there exists $k_{i-1}<j<k_{i+1}$ such that
\e
\text{
\begin{small}
$\tau\!\ci\!\mu(\{k_{i-1}\!+\!1,\ldots,j\})\!\ge\!
\tau\!\ci\!\mu(\{k_{i-1}\!+\!1,\ldots,k_{i+1}\})\!\ge\!
\tau\!\ci\!\mu(\{j\!+\!1,\ldots,k_{i+1}\}),$
\end{small}
}
\label{an4eq14}
\e
with at least one of these inequalities {\it strict}. Divide into cases
(i) $j=k_i$;
\begin{itemize}
\setlength{\itemsep}{0pt} \setlength{\parsep}{0pt}
\item[(ii)] $k_{i-1}\!<\!j\!<\!k_i$ and $\tau\ci\mu(\{k_{i-1}+1,
\ldots,j\})\!>\!\tau\ci\mu(\{k_{i-1}+1,\ldots,k_{i+1}\})$;
\item[(iii)] $k_{i-1}\!<\!j\!<\!k_i$ and $\tau\ci\mu(\{k_{i-1}+1,
\ldots,k_{i+1}\})\!>\!\tau\ci\mu(\{j+1,\ldots,k_{i+1}\})$;
\item[(iv)] $k_i\!<\!j\!<\!k_{i+1}$ and $\tau\ci\mu(\{k_{i-1}+1,
\ldots,j\})\!>\!\tau\ci\mu(\{k_{i-1}+1,\ldots,k_{i+1}\})$;
\item[(v)] $k_i\!<\!j\!<\!k_{i+1}$ and $\tau\ci\mu(\{k_{i-1}+1,
\ldots,k_{i+1}\})\!>\!\tau\ci\mu(\{j+1,\ldots,k_{i+1}\})$.
\end{itemize}

In case (i) \eq{an4eq13} and \eq{an4eq14} are contradictory, as
\eq{an4eq14} has a strict inequality. In (ii) \eq{an4eq13} gives
$\tau\ci\mu(\{k_{i-1}\!+\!1,\ldots,j\})\!>\!\tau\ci\mu(\{k_{i-1}\!+\!1,
\ldots,k_i\})$, implying $\tau\ci\mu(\{k_{i-1}\!+\!1,\ldots,j\})\!>
\!\tau\ci\mu(\{j\!+\!1,\ldots,k_i\})$ by Condition \ref{an4cond}(v),
which contradicts $(\{k_{i-1}\!+\!1,\ldots,k_i\},\le,\mu)$
$\tau$-semistable. For (iii), as
$\tau\ci\mu(\{j\!+\!1,\ldots,k_{i+1}\})$ lies between
$\tau\ci\mu(\{j\!+\!1,\ldots,k_i\})$ and $\tau\ci\mu(\{k_i
\!+\!1,\ldots,k_{i+1}\})$ by Condition \ref{an4cond}(v) we have {\it
either}
$\tau\ci\mu(\{k_{i-1}\!+\!1,\ldots,k_{i+1}\})\!>\!\tau\ci\mu(
\{j\!+\!1,\ldots,k_i\})$, giving $\tau\ci\mu(\{k_{i-1}\!+\!1,
\ldots,j\})>\!\tau\ci\mu(\{j\!+\!1,\ldots,k_i\})$ by \eq{an4eq14}
contradicting $(\{k_{i-1}\!+\!1,\ldots,k_i\},\le,\mu)$
$\tau$-semistable, {\it or}
$\tau\ci\mu(\{k_{i-1}\!+\!1,\ldots,k_{i+1}\})\!>\!\tau\ci\mu(
\{k_i\!+\!1,\ldots,k_{i+1}\})$, contradicting \eq{an4eq13}.
Similarly (iv),(v) give contradictions, using
$(\{k_i\!+\!1,\ldots,k_{i+1}\}, \le,\mu)$ $\tau$-semistable.
Therefore $(\{1,\ldots,m\},\le,\la)$ is $\tau$-reversing.

This proves existence of $m,\phi,\chi,\la$. For uniqueness,
suppose $m',\phi',\chi',\la'$ also satisfy the conditions.
Then $\phi'^{-1}(\{1,\ldots,i\})\!=\!\{1,\ldots,k_i'\}$ for
unique $0=k_0'<k_1'<\cdots<k_m'=b$. As $(\{1,\ldots,k_1'\},\le,\mu)$
is $\tau$-semistable we have $k_1'\le k_1$ by choice of $k_1$.
Suppose $k_1'<k_1$. Then $\tau\ci\mu(\{1,\ldots,k_1\})\ge\tau\ci\mu
(\{1,\ldots,k_1'\})$ by Condition \ref{an4cond}(v), as $1\le k_1'
<k_1$ and $(\{1,\ldots,k_1\},\le,\mu)$ is $\tau$-semistable. Also
\begin{equation*}
\tau\ci\mu(\{1,\ldots,k_1'\})\!>\!\tau\ci\mu(\{k_1'+1,\ldots,k_2'\})
\!>\!\cdots\!>\!\tau\ci\mu(\{k_{m-1}'+1,\ldots,k_m'\})
\end{equation*}
as $(\{1,\ldots,m'\},\le,\la')$ is $\tau$-reversing. Now
$k_i'<k_1\le k_{i+1}'$ for some $1\le i<m'$, and using the
inequalities above and Condition \ref{an4cond}(v) gives
$\tau\ci\mu(\{1,\ldots,k_1\})>\tau\ci\mu(\{1,\ldots,k_i'\})$,
so that
\begin{equation*}
\tau\ci\mu(\{k_i'+1,\ldots,k_1\})>\tau\ci\mu(\{1,\ldots,k_i'\})
>\tau\ci\mu(\{k_i'+1,\ldots,k_{i+1}'\}),
\end{equation*}
contradicting both $k_1=k_{i+1}'$, and $\tau$-semistability of
$\bigl(\{k_i'+1,\ldots,k_{i+1}'\},\le,\mu\bigr)$ if $k_1<k_{i+1}'$
by Condition \ref{an4cond}(v). Therefore $k_1'=k_1$. Extending this
argument shows $k_i'=k_i$ for $i=1,2,\ldots$ and $m=m'$, proving
uniqueness of $m,\phi,\chi$, and the proposition for~$c=1$.

For $c>1$ we apply the $c=1$ case to $(\ga^{-1}(\{j\}),\le,\mu)$ for
$j=1,\ldots,c$. This gives unique $m_j$, $\phi_j:\ga^{-1}(\{j\})\ra
\{1,\ldots,m_j\}$ and $\la_j$ such that $(\phi_j^{-1}(i),\le,\mu)$
is $\tau$-semistable for $i=1,\ldots,m_j,$ and $(\{1,\ldots,m_j\},
\le,\la_j)$ is $\tau$-reversing for $j=1,\ldots,c$. Define $m\!=\!m_1
\!+\!\cdots\!+\!m_c$ and $\phi:\{1,\ldots,b\}\!\ra\!\{1,\ldots,m\}$
by $\phi(i)\!=\!m_1\!+\!\cdots\!+\!m_{j-1}\!+\!\phi_j(i)$ when
$\ga(i)\!=\!j$, and $\chi:\{1,\ldots,m\}\!\ra\!\{1,\ldots,c\}$
by $\chi(i)\!=\!j$ when $m_1\!+\!\cdots\!+\!m_{j-1}\!<\!i\!\le\!
m_1\!+\!\cdots\!+\!m_j$. Then $m,\phi,\chi$ satisfy the proposition.
Uniqueness of $m,\phi,\chi$ follows easily from that of the~$m_j,\phi_j$.
\end{proof}

\begin{prop} Let\/ $1\le d\le b$ and\/ $\ep:\{1,\ldots,b\}\ra
\{1,\ldots,d\}$ be surjective with\/ $1\le i\le j\le b$ implies
$\ep(i)\le\ep(j)$. Then
\e
\sum_{c=d}^b \!\!\!
\sum_{\substack{
\text{surjective $\ga:\{1,\ldots,b\}\!\ra\!\{1,\ldots,c\}$ and}\\
\text{$\de:\{1,\ldots,c\}\!\ra\!\{1,\ldots,d\}$ with\/
$\ep=\de\ci\ga$:}\\
\text{$i\!\le\!j$ implies $\ga(i)\!\le\!\ga(j)$,
$i\!\le\!j$ implies $\de(i)\!\le\!\de(j)$}}}
\!\!\!\!\!\!\!\!\!\!\!\!\!\!\!\!\!\!\!\!
(-1)^{b-c}=\begin{cases} 1, & \text{$d=b$,}
\\ 0, & \text{otherwise.}\end{cases}
\label{an4eq15}
\e
\label{an4prop4}
\end{prop}

\begin{proof} The first line of \eq{an4eq15} is immediate as if
$d=b$ the only possibility is $c=b$ and $\ep=\ga= \de=\id$. So
suppose $d<b$. Define $S=\{i=1,\ldots,b-1:\ep(i+1)=\ep(i)\}$. Then
$\md{S}=b-d$. For $\ga,\de$ as in \eq{an4eq15}, define $T=\{i\in
S:\ga(i+1)=\ga(i)\}$. Then $\md{T}=c-d$, and $T$ determines $\ga$ by
$\ga(1)=1$, $\ga(i+1)=\ga(i)$ if $i\in T$ and $\ga(i+1)=\ga(i)+1$ if
$i\in\{1,\ldots,b-1\}\sm T$. Also $c,\ep,\ga$ determine $\de$ by
$\ep=\de\ci\ga$. This establishes a 1-1 correspondence between
choices of $c,\ga,\de$ in \eq{an4eq15} and subsets $T\subseteq S$
with $\md{T}=c-d$. Hence the l.h.s.\ of \eq{an4eq15} becomes
\begin{equation*}
\sum_{\text{subsets $T\subseteq S$}}(-1)^{\md{S\sm T}}=
\sum_{j=0}^{b-d}\binom{b-d}{j}(-1)^{b-d-j}=(1-1)^{b-d}=0,
\end{equation*}
as there are $\bigl(\begin{smallmatrix} b-d \\ j \end{smallmatrix}
\bigr)$ subsets $T\subseteq S$ with~$\md{T}=j$.
\end{proof}

To complete the proof of Theorem \ref{an4thm1} we use sets and maps
as follows:
\begin{equation*}
\text{
\begin{small}
$\displaystyle
\xymatrix@C=11pt{
\{1,\ldots,n\} \ar[r]^\al \ar@/_1.1pc/[rrr]^\psi
& \{1,\ldots,a\} \ar[r]^\be
& \{1,\ldots,b\} \ar[r]^\phi \ar@<2ex>@/^.4pc/[rr]_\ga
& \{1,\ldots,m\} \ar[r]^\chi \ar@/_1.1pc/[rr]^\ep
& \{1,\ldots,c\} \ar[r]^\de
& \{1,\ldots,d\},
}$
\end{small}
}
\end{equation*}
and maps $\ka:\{1,\ldots,n\}\!\ra\!C(\A)$, $\la:\{1,\ldots,m\}\!\ra\!C(\A)$,
$\mu:\{1,\ldots,b\}\!\ra\!C(\A)$ and $\nu:\{1,\ldots,d\}\!\ra\!C(\A)$. We
deduce \eq{an4eq7} from Propositions \ref{an4prop2}--\ref{an4prop4} by:
\begin{gather}
\sum_{m=1}^n\,
\sum_{\substack{\text{$\psi:\{1,\ldots,n\}\ra\{1,\ldots,m\}:$}\\
\text{$\psi$ is surjective,}\\
\text{$1\!\le\!i\!\le\!j\!\le\!n$ implies $\psi(i)\!\le\!\psi(j)$,}\\
\text{define $\la:\{1,\ldots,m\}\!\ra\!C(\A)$}\\
\text{by $\la(k)=\ka(\psi^{-1}\!(k))$}}}\,
\begin{aligned}[t]
&S(\{1,\ldots,m\},\le,\la,\hat\tau,\ti\tau)\cdot\\
&\prod_{k=1}^mS\bigl(\psi^{-1}(\{k\}),\le,
\ka\vert_{\psi^{-1}(\{k\})},\tau,\hat\tau)=
\end{aligned}
\nonumber
\allowdisplaybreaks
\\
\sum_{\substack{\text{$m=1,\ldots,n$ and}\\
\text{$\psi:\{1,\ldots,n\}\ra\{1,\ldots,m\}:$}\\
\text{$\psi$ is surjective,}\\
\text{$i\!\le\!j$ implies $\psi(i)\!\le\!\psi(j)$,}\\
\text{define $\la:\{1,\ldots,m\}\!\ra\!C(\A)$}\\
\text{by $\la(k)=\ka(\psi^{-1}\!(k))$}}}\,
\raisebox{-12pt}{\begin{Large}$\displaystyle\Biggl[$\end{Large}}
\sum_{1\le d\le c\le m}
\!\!\!\!\!\!\!\!\!\!\!\!\!
\sum_{\substack{
\text{surjective $\chi:\{1,\ldots,m\}\!\ra\!\{1,\ldots,c\}$}\\
\text{and $\de:\{1,\ldots,c\}\!\ra\!\{1,\ldots,d\}$:}\\
\text{$i\!\le\!j$ implies $\chi(i)\!\le\!\chi(j)$,
$i\!\le\!j$ implies $\de(i)\!\le\!\de(j)$,}\\
\text{$(\chi^{-1}(j),\le,\la)$ $\hat\tau$-reversing, $1\le j\le c$.}\\
\text{Let $\nu:\{1,\ldots,d\}\!\ra\!C(\A)$ be
$\nu:k\!\mapsto\!\la((\de\ci\chi)^{-1}(k))$.}\\
\text{Then $(\{1,\ldots,d\},\le,\nu)$ is $\ti\tau$-semistable}}}
\!\!\!\!\!\!\!\!\!\!\!\!\!\!\!\!\!\!\!
(-1)^{c-d}
\raisebox{-12pt}{\begin{Large}$\displaystyle\Biggr]$\end{Large}}
\cdot
\nonumber
\allowdisplaybreaks
\\
\raisebox{-12pt}{\begin{Large}$\displaystyle\Biggl[$\end{Large}}
\sum_{m\le b\le a\le n}\,\,\,
\sum_{\substack{
\text{surjective $\al:\{1,\ldots,n\}\!\ra\!\{1,\ldots,a\}$,
$\be:\{1,\ldots,a\}\!\ra\!\{1,\ldots,b\}$}\\
\text{and $\phi:\{1,\ldots,b\}\!\ra\!\{1,\ldots,m\}$ with
$\psi=\phi\ci\be\ci\al$:}\\
\text{$i\!\le\!j$ implies $\al(i)\!\le\!\al(j)$,
$i\!\le\!j$ implies $\be(i)\!\le\!\be(j)$,}\\
\text{$i\!\le\!j$ implies $\phi(i)\!\le\!\phi(j)$,
$(\al^{-1}(j),\le,\ka)$ $\tau$-reversing, $1\!\le\!j\!\le\!a$.}\\
\text{Let $\mu:\{1,\ldots,b\}\!\ra\!C(\A)$ be
$\mu:k\!\mapsto\!\ka((\be\ci\al)^{-1}(k))$.}\\
\text{Then $(\phi^{-1}(i),\le,\mu)$ is $\hat\tau$-semistable,
$1\!\le\!i\!\le\!m$}}}
\!\!\!\!\!\!\!\!\!\!\!\!\!\!\!\!\!\!\!\!
(-1)^{a-b}
\raisebox{-12pt}{\begin{Large}$\displaystyle\Biggr]$\end{Large}}=
\nonumber
\allowdisplaybreaks
\\
\sum_{d\le b\le a\le n}\,
\sum_{\substack{
\text{surjective $\al:\{1,\ldots,n\}\!\ra\!\{1,\ldots,a\}$,}\\
\text{$\be:\{1,\ldots,a\}\!\ra\!\{1,\ldots,b\}$ and
$\ep:\{1,\ldots,b\}\!\ra\!\{1,\ldots,d\}$:}\\
\text{$i\!\le\!j$ implies $\al(i)\!\le\!\al(j)$,
$i\!\le\!j$ implies $\be(i)\!\le\!\be(j)$,}\\
\text{$i\!\le\!j$ implies $\ep(i)\!\le\!\ep(j)$,
$(\al^{-1}(j),\le,\ka)$ $\tau$-reversing, $1\!\le\!j\!\le\!a$.}\\
\text{Let $\nu:\{1,\ldots,d\}\!\ra\!C(\A)$ be
$\nu:k\!\mapsto\!\ka((\ep\ci\be\ci\al)^{-1}(k))$.}\\
\text{Then $(\{1,\ldots,d\},\le,\nu)$ is $\ti\tau$-semistable}}}
\!\!\!\!\!\!\!\!\!\!\!\!\!\!\!\!\!
(-1)^{a-d}\cdot
\nonumber
\allowdisplaybreaks
\\
\raisebox{-12pt}{\begin{Large}$\displaystyle\Biggl[$\end{Large}}
\sum_{d\le c\le b}\,\,
\sum_{\substack{
\text{surjective $\ga:\{1,\ldots,b\}\!\ra\!\{1,\ldots,c\}$ and}\\
\text{$\de:\{1,\ldots,c\}\!\ra\!\{1,\ldots,d\}$ with $\ep=\de\ci\ga$:}\\
\text{$i\!\le\!j$ implies $\ga(i)\!\le\!\ga(j)$,
$i\!\le\!j$ implies $\de(i)\!\le\!\de(j)$}}}
\!\!\!\!\!\!\!\!\!\!\!\!\!\!\!\!
(-1)^{b-c}
\raisebox{-12pt}{\begin{Large}$\displaystyle\Biggr]$\end{Large}}
\label{an4eq16}
\\
=S(\{1,\ldots,n\},\le,\ka,\tau,\ti\tau).
\nonumber
\end{gather}

Here the term $[\cdots]$ on the second line is $S(\{1,\ldots,m\},
\le,\la,\hat\tau,\ti\tau)$ by \eq{an4eq12}, and $[\cdots]$ on the
third line is $\prod_{k=1}^mS\bigl(\psi^{-1}(\{k\}),\le,\ka
\vert_{\psi^{-1}(\{k\})},\tau,\hat\tau)$ by \eq{an4eq12}, where
we have relabelled and combined the product over $k=1,\ldots,m$
of sums over $a_k$, $\al_k:\psi^{-1}(\{k\})\!\ra\!\{1,\ldots,a_k\}$
and $b_k$, $\be_k:\{1,\ldots,a_k\}\!\ra\!\{1,\ldots,b_k\}$ into one
large sum over $a,b,\al,\be$ and $\phi$, where $a\!=\!a_1\!+\!
\cdots\!+\!a_m$ and~$b\!=\!b_1\!+\!\cdots\!+\!b_m$.

To deduce the fourth and fifth lines of \eq{an4eq16} from
the second and third, we set $\ga=\chi\ci\phi$. Then the
conditions $i\!\le\!j$ implies $\chi(i)\!\le\!\chi(j)$ and
$(\chi^{-1}(j),\le,\la)$ $\tau$-reversing for $j=1,\ldots,c$
in the second line, and $i\!\le\!j$ implies $\phi(i)\!\le\!
\phi(j)$, $(\phi^{-1}(i),\le,\mu)$ $\hat\tau$-semistable for
$i=1,\ldots,m$ in the third line, are the hypotheses of
Proposition \ref{an4prop3} with $\hat\tau$ in place of $\tau$.
Thus Proposition \ref{an4prop3} shows that for all choices of
$b,c,\ga$ and $\mu$ there are unique choices of $m,\phi,\chi$
satisfying the conditions involving $\hat\tau$. So we drop the
sums over $m,\phi,\chi,\psi$ and the $\hat\tau$ conditions,
and also rearrange the sums.

To deduce the sixth and last line of \eq{an4eq16}, note
that the term $[\cdots]$ on the fifth line is 1 if $d=b$
and 0 otherwise, by Proposition \ref{an4prop4}. If $d=b$ then
$\ep$ is the identity, and the fourth line of \eq{an4eq16}
reduces to $S(\{1,\ldots,n\},\le,\ka,\tau,\ti\tau)$ by
\eq{an4eq12}. This completes the proofs of \eq{an4eq7}
and Theorem~\ref{an4thm1}.

\section{Transforming between stability conditions}
\label{an5}

We now prove the transformation laws \eq{an4eq1}--\eq{an4eq3}
from $(\tau,T,\le)$ to $(\ti\tau,\ti T,\le)$, and their stack
function analogues. For most of the section we do not suppose
$(\tau,T,\le),(\ti\tau,\ti T,\le)$ are {\it permissible}.
Therefore our equations are identities in $\dLCF,\dLSFa(\fObj_\A)$
rather than $\CF,\SFa(\fObj_\A)$, which may have {\it infinitely
many nonzero terms}, and must be interpreted using a notion
of {\it convergence}.

This is partly for greater generality, but mostly because even when
$(\tau,T,\le),\ab (\ti\tau,\ti T,\le)$ are permissible, to prove
\eq{an4eq1}--\eq{an4eq3} we may need to go via an intermediate weak
stability condition $(\hat\tau,\hat T,\le)$ which is not
permissible, which happens when $\A=\coh(P)$. It is far from obvious
that \eq{an4eq1}--\eq{an4eq3} have only finitely many nonzero terms
even for permissible $(\tau,T,\le),(\ti\tau, \ti T,\le)$. We prove
this for $\A\!=\!\coh(P)$ when $P$ is a smooth surface and
$\tau,\ti\tau$ Gieseker stability conditions.

The (locally) constructible functions material below needs the
ground field $\K$ to have {\it characteristic zero}, but the (local)
stack function versions work for $\K$ of {\it arbitrary
characteristic}. As we develop the two strands in parallel, for
brevity we make the convention that $\cha\K=0$ in the parts dealing
with $\CF,\dLCF(\fObj_\A)$ and $\cha\K$ is arbitrary otherwise, and
will mostly not remark on it.

\subsection{Basic definitions and main results}
\label{an51}

We will need the following finiteness conditions.

\begin{dfn} Let Assumption \ref{an3ass} hold and $(\tau,T,\le),
(\ti\tau,\ti T,\le)$ be weak stability conditions on $\A$. We say
{\it the change from $(\tau,T,\le)$ to $(\ti\tau,\ti T,\le)$ is
locally finite} if for all constructible $U\subseteq\fObj_\A(\K)$,
there are only finitely many sets of $\A$-data $(\{1,\ldots,n\},\le,
\ka)$ for which $S(\{1,\ldots,n\},\le,\ka,\tau,\ti\tau)\!\ne\!0$ and
\e
U\cap\bs\si(\{1,\ldots,n\})_*\bigl(\Mss(\{1,\ldots,n\},\le,\ka,
\tau)_\A\bigr)\ne\emptyset.
\label{an5eq1}
\e
We say {\it the change from $(\tau,T,\le)$ to $(\ti\tau,\ti T,\le)$
is globally finite} if this holds for $U=\fObj_\A^\al(\K)$ (which
is not constructible, in general) for all $\al\in C(\A)$. Since
any constructible $U\subseteq\fObj_\A(\K)$ is contained in a
finite union of $\fObj_\A^\al(\K)$, globally finite implies
locally finite.

When we say {\it the changes between $(\tau,T,\le)$ and\/}
$(\ti\tau,\ti T,\le)$, we mean both the change from $(\tau,T,\le)$
to $(\ti\tau,\ti T,\le)$ and the change from $(\ti\tau,\ti T,\le)$
to~$(\tau,T,\le)$.
\label{an5def1}
\end{dfn}

Here is the main result of this section.

\begin{thm} Let Assumption \ref{an3ass} hold and\/ $(\tau,T,\le),
(\ti\tau,\ti T,\le),(\hat\tau,\hat T,\le)$ be weak stability
conditions on $\A$ with\/ $(\hat\tau,\hat T,\le)$ dominating
$(\tau,T,\le),(\ti\tau,\ti T,\le)$, and suppose the change from
$(\hat\tau,\hat T,\le)$ to $(\ti\tau,\ti T,\le)$ is locally finite.
Then the change from $(\tau,T,\le)$ to $(\ti\tau,\ti T,\le)$ is
locally finite.

Let\/ $\al\in C(\A)$ and\/ $(K,\tl,\mu)$ be $\A$-data. Then
\begin{gather}
\begin{gathered}
\sum_{\substack{\text{$\A$-data $(\{1,\ldots,n\},\le,\ka):$}\\
\text{$\ka(\{1,\ldots,n\})=\al$}}} \!\!\!\!\!\!\!\!\!
\begin{aligned}[t]
S(\{1,&\ldots,n\},\le,\ka,\tau,\ti\tau)\cdot\\
&\dss^{\ka(1)}(\tau)*\dss^{\ka(2)}(\tau)*\cdots*
\dss^{\ka(n)}(\tau)=\dss^\al(\ti\tau),
\end{aligned}
\end{gathered}
\label{an5eq2}
\allowdisplaybreaks
\\
\begin{gathered}
\sum_{\substack{\text{$\A$-data $(\{1,\ldots,n\},\le,\ka):$}\\
\text{$\ka(\{1,\ldots,n\})=\al$}}} \!\!\!\!\!\!\!\!\!
\begin{aligned}[t]
S(\{1,&\ldots,n\},\le,\ka,\tau,\ti\tau)\cdot\\
&\bdss^{\ka(1)}(\tau)*\bdss^{\ka(2)}(\tau)*\cdots*
\bdss^{\ka(n)}(\tau)=\bdss^\al(\ti\tau),
\end{aligned}
\end{gathered}
\label{an5eq3}
\allowdisplaybreaks
\\
\begin{aligned}
\sum_{\substack{\text{iso.}\\ \text{classes}\\
\text{of finite }\\ \text{sets $I$}}} \frac{1}{\md{I}!}\cdot
\sum_{\substack{
\text{$\pr,\ka,\phi$: $(I,\pr,\ka)$ is $\A$-data,}\\
\text{$(I,\pr,K,\phi)$ is dominant,}\\
\text{$\tl=\P(I,\pr,K,\phi)$,}\\
\text{$\ka(\phi^{-1}(k))=\mu(k)$ for $k\in K$}}}\,\,\,
\begin{aligned}[t]
&T(I,\pr,\ka,K,\phi,\tau,\ti\tau)\cdot\\
&\CF^\stk\bigl(Q(I,\pr,K,\tl,\phi)\bigr)\\
&\quad\dss(I,\pr,\ka,\tau)=\dss(K,\tl,\mu,\ti\tau),
\end{aligned}
\end{aligned}
\label{an5eq4}
\allowdisplaybreaks
\\
\begin{aligned}
\sum_{\substack{\text{iso.}\\ \text{classes}\\
\text{of finite }\\ \text{sets $I$}}} \frac{1}{\md{I}!}\cdot
\sum_{\substack{
\text{$\pr,\ka,\phi$: $(I,\pr,\ka)$ is $\A$-data,}\\
\text{$(I,\pr,K,\phi)$ is dominant,}\\
\text{$\tl=\P(I,\pr,K,\phi)$,}\\
\text{$\ka(\phi^{-1}(k))=\mu(k)$ for $k\in K$}}}\,\,\,
\begin{aligned}[t]
&T(I,\pr,\ka,K,\phi,\tau,\ti\tau)\cdot\\
&Q(I,\pr,K,\tl,\phi)_*\\
&\quad\bdss(I,\pr,\ka,\tau)=\bdss(K,\tl,\mu,\ti\tau).
\end{aligned}
\end{aligned}
\label{an5eq5}
\end{gather}

If also $(\tau,T,\le),(\ti\tau,\ti T,\le)$ are permissible then
\begin{gather}
\begin{gathered}
\sum_{\substack{\text{$\A$-data $(\{1,\ldots,n\},\le,\ka):$}\\
\text{$\ka(\{1,\ldots,n\})=\al$}}} \!\!\!\!\!\!\!\!\!\!\!\!\!\!\!
\begin{aligned}[t]
U(\{1,&\ldots,n\},\le,\ka,\tau,\ti\tau)\cdot\\
&\ep^{\ka(1)}(\tau)*\ep^{\ka(2)}(\tau)*\cdots*
\ep^{\ka(n)}(\tau)=\ep^\al(\ti\tau),
\end{aligned}
\end{gathered}
\label{an5eq6}
\allowdisplaybreaks
\\
\begin{gathered}
\sum_{\substack{\text{$\A$-data $(\{1,\ldots,n\},\le,\ka):$}\\
\text{$\ka(\{1,\ldots,n\})=\al$}}} \!\!\!\!\!\!\!\!\!\!\!\!\!\!\!
\begin{aligned}[t]
U(\{1,&\ldots,n\},\le,\ka,\tau,\ti\tau)\cdot\\
&\bep^{\ka(1)}(\tau)*\bep^{\ka(2)}(\tau)*\cdots*
\bep^{\ka(n)}(\tau)=\bep^\al(\ti\tau).
\end{aligned}
\end{gathered}
\label{an5eq7}
\end{gather}
Here \eq{an5eq2}--\eq{an5eq7} are equations in
$\dLCF,\dLSFa(\fObj_\A),\LCF,\LSF(\fM(K,\tl,\mu)_\A),\ab
\CF,\SFa(\fObj_\A)$ respectively, with\/ $\cha\K=0$ in \eq{an5eq2},
\eq{an5eq4}, \eq{an5eq6}. There may be infinitely many nonzero terms
on the left hand sides of\/ \eq{an5eq2}--\eq{an5eq7}, but they hold
as convergent sums in $\LCF,\LSF(\cdots)$ in the sense of
Definition~\ref{an2def10}.

If the change from $(\tau,T,\le)$ to $(\ti\tau,\ti T,\le)$ is
globally finite then there are only finitely many nonzero terms in
\eq{an5eq2}--\eq{an5eq7}, and they hold as finite sums.
\label{an5thm1}
\end{thm}

We assume $(\tau,T,\le),(\ti\tau,\ti T,\le)$ {\it permissible} in
\eq{an5eq6}--\eq{an5eq7}, as we defined $\ep^\al(\tau)$,
$\bep^\al(\tau)$ this way in \cite[\S 7--\S 8]{Joyc7}. But
\eq{an5eq6}--\eq{an5eq7} also hold in $\dLCF,\dLSFa(\fObj_\A)$ if
$(\tau,T,\le),(\ti\tau,\ti T,\le)$ are stability conditions and $\A$
is $\tau$- and $\ti\tau$-artinian.

Suppose \eq{an5eq6} holds as a {\it finite} sum. Then $\ep^\al
(\ti\tau)$ lies in $\Ht_\tau$ and $\CFi(\fObj_\A)$, so by Definition
\ref{an3def16} it lies in the Lie algebra $\Lt_\tau$, which is
generated by the $\ep^\be(\tau)$. That is, $\ep^\al(\ti\tau)$ is a
$\Q$-linear combination of multiple commutators such~as
\begin{equation*}
[\,[\cdots[\,[\ep^{\ka(1)}(\tau),\ep^{\ka(2)}(\tau)],
\ep^{\ka(3)}(\tau)],\cdots],\ep^{\ka(n)}(\tau)].
\end{equation*}
Thus it is natural to ask whether we can rewrite \eq{an5eq6} as a
sum of commutators of $\ep^{\ka(i)}(\tau)$ rather than a sum of
products, so that it becomes an identity in the Lie algebra
$\CFi(\fObj_\A)$ rather than the algebra $\CF(\fObj_\A)$. The next
theorem shows that we can. First we introduce some notation.

\begin{dfn} For $I$ a finite set, write $\Q[I]$ for the
associative, noncommutative $\Q$-algebra generated by
$i\in I$. That is, elements of $\Q[I]$ are finite linear
combinations of {\it words\/} $i_1i_2\cdots i_m$ for
$m\ge 0$ and $i_1,\ldots,i_m\in I$, with multiplication
given by juxtaposition. If $J\subseteq I$ and $\pr$ is a
total order on $J$, write
\begin{equation*}
\prod_{\text{$j\in J$ in order $\pr$}}\!\!\!\!\!\!\!\!
j\,\,\,\,=\!j_1j_2\cdots j_n
\;\>\text{if $J\!=\!\{j_1,\ldots,j_n\}$ with
$j_1\pr j_2\pr\cdots\pr j_n$.}
\end{equation*}
We can regard $\Q[I]$ as a Lie algebra, with commutator
$[i_1\cdots i_m,j_1\cdots j_n]=i_1\cdots i_mj_1\cdots j_n
-j_1\cdots j_ni_1\cdots i_m$. Define $\L[I]$ to be the
Lie subalgebra of $\Q[I]$ generated by $i\in I$. Then
$\L[I]$ is the free $\Q$-Lie algebra generated by $I$,
and $\Q[I]$ is naturally isomorphic to~$U(\L[I])$.
\label{an5def2}
\end{dfn}

\begin{thm} Let Condition \ref{an4cond} hold for $(\tau,T,\le)$
and\/ $(\ti\tau,\ti T,\le)$, and $I$ be a finite set and\/
$\ka:I\ra C(\A)$. Then the expression in~$\Q[I]$
\e
\sum_{\text{total orders $\pr$ on $I$}}U(I,\pr,\ka,\tau,\ti\tau)
\cdot\prod_{\text{$i\in I$ in order $\pr$}}\!\!\!\!\! i
\qquad\text{lies in $\L[I]$.}
\label{an5eq8}
\e
\label{an5thm2}
\end{thm}

The author has a purely combinatorial proof of Theorem \ref{an5thm2},
but it is long and complicated, so we instead give a shorter
constructible functions proof. By an easy combinatorial argument
\eq{an5eq6}--\eq{an5eq7} may
be rewritten
\begin{gather}
\begin{aligned}
&\ep^\al(\ti\tau)=\\
&\sum_{\substack{\text{iso classes}\\
\text{of finite}\\ \text{sets $I$}}}
\frac{1}{\md{I}!}\!
\sum_{\substack{\ka:I\ra C(\A):\\ \ka(I)=\al}}\!
\raisebox{-12pt}{\begin{Large}$\displaystyle\Biggl[$\end{Large}}
\!\!
\sum_{\substack{\text{total orders $\pr$ on $I$.}\\
\text{Write $I=\{i_1,\ldots,i_n\}$,}\\
\text{$i_1\pr i_2\pr\cdots\pr i_n$}}}\,
\begin{aligned}[t]
&U(I,\pr,\ka,\tau,\ti\tau)\cdot\\
&\ep^{\ka(i_1)}(\tau)\!*\!\cdots\!*\!\ep^{\ka(i_n)}(\tau)
\end{aligned}
\raisebox{-12pt}{\begin{Large}$\displaystyle\Biggr]$\end{Large}},
\end{aligned}
\label{an5eq9}
\allowdisplaybreaks
\\
\begin{aligned}
&\bep^\al(\ti\tau)=\\
&\sum_{\substack{\text{iso classes}\\
\text{of finite}\\ \text{sets $I$}}}
\frac{1}{\md{I}!}\!
\sum_{\substack{\ka:I\ra C(\A):\\ \ka(I)=\al}}\!
\raisebox{-12pt}{\begin{Large}$\displaystyle\Biggl[$\end{Large}}
\!\!
\sum_{\substack{\text{total orders $\pr$ on $I$.}\\
\text{Write $I=\{i_1,\ldots,i_n\}$,}\\
\text{$i_1\pr i_2\pr\cdots\pr i_n$}}}\,
\begin{aligned}[t]
&U(I,\pr,\ka,\tau,\ti\tau)\cdot\\
&\bep^{\ka(i_1)}(\tau)\!*\!\cdots\!*\!\bep^{\ka(i_n)}(\tau)
\end{aligned}
\raisebox{-12pt}{\begin{Large}$\displaystyle\Biggr]$\end{Large}}.
\end{aligned}
\label{an5eq10}
\end{gather}
By Theorem \ref{an5thm2}, the terms $[\cdots]$ are linear
combinations of commutators of $\ep^{\ka(i)},\bep^{\ka(i)}$
for $i\in I$, so \eq{an5eq9}--\eq{an5eq10} are identities
in the Lie algebras $\CFi(\fObj_\A)$ and $\SFai(\fObj_\A)$.
We can now deduce:

\begin{cor} Let Assumption \ref{an3ass} hold and\/ $(\tau,T,\le),
(\ti\tau,\ti T,\le),(\hat\tau,\hat T,\le)$ be weak stability
conditions on $\A$ with\/ $(\tau,T,\le),(\ti\tau,\ti T,\le)$
permissible and\/ $(\hat\tau,\hat T,\le)$ dominating $(\tau,T,\le),
(\ti\tau,\ti T,\le)$. Suppose the changes from $(\hat\tau,\hat T,\le)$
to $(\tau,T,\le),\ab
(\ti\tau,\ti T,\le)$ are locally finite, and the changes between
$(\tau,T,\le)$ and\/ $(\ti\tau,\ti T,\le)$ are globally finite.
Then in the notation of\/ \S\ref{an35} we have
\e
\begin{aligned}
\Hp_\tau&=\Hp_{\ti\tau},& \Ht_\tau&=\Ht_{\ti\tau},&
\bHp_\tau&=\bHp_{\ti\tau}, & \bHt_\tau&=\bHt_{\ti\tau},\\
\Lp_\tau&=\Lp_{\ti\tau},& \Lt_\tau&=\Lt_{\ti\tau},&
\bLp_\tau&=\bLp_{\ti\tau}, & \bLt_\tau&=\bLt_{\ti\tau}.
\end{aligned}
\label{an5eq11}
\e
\label{an5cor1}
\end{cor}

\begin{proof} Theorem \ref{an5thm1} implies
\eq{an5eq2}--\eq{an5eq7} hold as finite sums in $\CF,\SFa(\fObj_\A)$
and $\CF,\SF(\fM(K,\tl,\mu)_\A)$. Thus
$\dss^\al(\ti\tau)\!\in\!\Ht_\tau$,
$\bdss^\al(\ti\tau)\!\in\!\bHt_\tau$, and applying
$\CF^\stk(\bs\si(K)),\bs\si(K)_*$ to \eq{an5eq4}--\eq{an5eq5} gives
$\CF^\stk(\bs\si(K))\dss(K,\tl,\mu,\ti\tau)\!\in\!\Hp_\tau$ and
$\bs\si(K)_*\,\bdss(K,\tl,\mu,\ti\tau)\!\in\!\bHp_\tau$. Hence
$\Hp_{\ti\tau}\subseteq\Hp_\tau,\ldots,
\bHt_{\ti\tau}\subseteq\bHt_\tau$. Exchanging $\tau,\ti\tau$ proves
the top line of \eq{an5eq11}. The first three equations of the
bottom line follow from Definition \ref{an3def16}. For the last,
\eq{an5eq10} and Theorem \ref{an5thm2} imply
$\bep^\al(\ti\tau)\!\in\!\bLt_\tau$, so
$\bLt_{\ti\tau}\!\subseteq\!\bLt_\tau$, and exchanging
$\tau,\ti\tau$ gives the result.
\end{proof}

In \cite[Ex.s 10.5--10.9]{Joyc5} we define data $\A,K(\A),\fF_\A$
satisfying Assumption \ref{an3ass} with $\A=\modKQ$ or $\nilKQ$ for
$Q=(Q_0,Q_1,b,e)$ a {\it quiver}, and $\A=\modKQI$ or $\nilKQI$ for
$(Q,I)$ a {\it quiver with relations}, and $\A=\modA$ for $A$ a {\it
finite-dimensional\/ $\K$-algebra}. The last case $\A=\modA$ also
has an associated quiver, the {\it Ext-quiver} $Q$ of $A$, and
$K(\A)\cong\Z^{Q_0}$ in each case. All weak stability conditions on
these $\A$ are permissible on by \cite[Cor.~4.13]{Joyc7}, and as for
each $\al\in C(\A)$ there are only finitely many $\be,\ga\in C(\A)$
with $\al=\be+\ga$, the changes between any two weak stability
conditions are globally finite.

For any such $\A$ we have the trivial stability condition
$(0,\{0\},\le)$, mapping $\al\mapsto 0$ for all $\al\in C(\A)$.
This dominates all weak stability conditions on $\A$. So applying
Theorem \ref{an5thm1} with $(\hat\tau,\hat T,\le)=(0,\{0\},\le)$
and Corollary \ref{an5cor1} gives:

\begin{cor} Let\/ $\A,K(\A),\fF_\A$ be as in one of the quiver
examples {\rm\cite[Ex.s 10.5--10.9]{Joyc5}}, and let\/
$(\tau,T,\le),(\ti\tau,\ti T,\le)$ be any two weak stability
conditions on $\A$, such as the slope functions of\/ {\rm
\cite[Ex.~4.14]{Joyc7}}. Then $(\tau,T,\le),(\ti\tau,\ti T,\le)$ are
permissible, and the changes between them are globally finite, and\/
\eq{an5eq2}--\eq{an5eq7} hold as finite sums in $\CF,\SFa(\fObj_\A)$
or $\CF,\SF(\fM(K,\tl,\mu)_\A)$, and\/ \eq{an5eq11} holds. Thus
$\Hp_\tau,\ldots,\bLt_\tau$ are independent of the weak stability
condition~$(\tau,T,\le)$.
\label{an5cor2}
\end{cor}

In \cite[\S 4.4]{Joyc7} we studied weak stability conditions on
$\A=\coh(P)$, the abelian category of {\it coherent sheaves} on a
projective $\K$-scheme $P$, with data $K(\coh(P)),\fF_{\coh(P)}$ as
in \cite[Ex.~9.1 or 9.2]{Joyc5}. With $m=\dim P$, we defined
\cite[Ex.s 4.16--4.18]{Joyc7} a permissible stability condition
$(\ga,G_m,\le)$ from {\it Gieseker stability}, a permissible weak
stability condition $(\mu,M_m,\le)$ from $\mu$-{\it stability}, and
a non-permissible weak stability condition $(\de,D_m,\le)$ from {\it
purity} and the {\it torsion filtration}, which dominates
$(\ga,G_m,\le)$ and~$(\mu,M_m,\le)$.

\begin{prop} In the situation of\/ {\rm\cite[\S 4.4]{Joyc7}}
with\/ $\A=\coh(P)$, the changes from $(\de,D_m,\le)$ to
$(\ga,G_m,\le),(\mu,M_m,\le)$ are locally finite.
\label{an5prop1}
\end{prop}

\begin{proof} We first verify Definition \ref{an5def1} with
$(\tau,T,\le)=(\de,D_m,\le)$, $(\ti\tau,\ti T,\le)=(\ga,G_m,\le)$
and $U=\{[X]\}$ with $[X]=\al\in C(\A)$. If $(\{1,\ldots,n\},\le,
\ka)$ is $\A$-data with $\ka(\{1,\ldots,n\})=\al$ then the argument
of \eq{an5eq18} below shows that
\begin{equation*}
S(\{1,\ldots,n\},\le,\ka,\de,\ga)=
\begin{cases} (-1)^{n-1}, & \begin{aligned}
&\text{$\ga\ci\ka(\{1,\ldots,i\})\!>\!\ga\ci\ka(\{i\!+\!1,
\ldots,n\})$}\\
&\text{for all $1\le i<n$, and $\de\ci\ka\equiv\de(\al)$,}
\end{aligned} \\
0, &\text{otherwise.} \end{cases}
\end{equation*}
Suppose $S(\{1,\ldots,n\},\le,\ka,\de,\ga)\ne 0$ and
$[(\si,\io,\pi)]\in\Mss(\{1,\ldots,n\},\le,\ka,\de)_\A$ with
$\si(\{1,\ldots,n\})=X$. As the $\si(\{i\})$ are $\de$-semistable
and $\de\equiv\de(\al)=\de([X])$, by \cite[Ex.~4.18]{Joyc7} the
$\si(\{i\})$ are all {\it pure} with $\dim\si(\{i\})=\dim X$.

By induction we see that $\si(\{i+1,\ldots,n\})$ is pure of
dimension $\dim X$ for all $1\le i<n$. But we have an exact sequence
\begin{equation*}
\xymatrix@C=13pt{
0 \ar[r] & \si(\{1,\ldots,i\})
\ar[rrrr]
\ar@<.3ex>@{}[rrrr]^(0.56){\io(\{1,\ldots,i\},\{1,\ldots,n\})}
&&&& X \ar[rrrr]
\ar@<.3ex>@{}[rrrr]^(0.37){\pi(\{1,\ldots,n\},\{i+1,\ldots,n\})}
&&&& \si(\{i\!+\!1,\ldots,n\}) \ar[r] & 0,
}
\end{equation*}
so $\si(\{i+1,\ldots,n\})$ is a {\it quotient sheaf\/} of $X$.
Also, $\ga\ci\ka(\{1,\ldots,i\})>\ga\ci\ka(\{i+1,\ldots,n\})$
implies $\ga(\al)\ge\ga\ci\ka(\{i+1,\ldots,n\})$ by the weak seesaw
inequality, and as $\deg\ga(\al)=\deg\ga\ci\ka(\{i+1,\ldots,n\})$
this implies that $\hat\mu\bigl(\si(\{i+1,\ldots,n\})\bigr)\le
\hat\mu(X)$, where $\hat\mu$ is the {\it slope} of
\cite[Def.~1.6.8]{HuLe}, basically the second coefficient in
the $\ga$ polynomial.

Now Huybrechts and Lehn \cite[Lem.~1.7.9 \& Rem.~1.7.10]{HuLe} prove
that if $X\in\coh(P)$ if fixed and $\mu_0\in\R$, the family of all
pure quotient sheaves $Y$ of $X$ with $\dim Y=\dim X$ and
$\hat\mu([Y])\le\mu_0$ is constructible, and hence realizes only
finitely many classes in $C(\A)$. Applying this with $\mu_0=
\hat\mu(X)$ shows there are only finitely many possibilities for
$\ka(\{i+1,\ldots,n\})$ in $C(\A)$. It easily follows that there are
only finitely many possibilities for $n,\ka$, as we want.

This extends to arbitrary constructible $U\subseteq\fObj_\A(\K)$
using a families version of \cite[Lem.~1.7.9 \& Rem.~1.7.10]{HuLe},
so the change from $(\de,D_m,\le)$ to $(\ga,G_m,\le)$ is locally
finite. The proof for $(\mu,M_m,\le)$ is the same, since
$\mu(\al)\ge\mu\ci\ka(\{i+1,\ldots,n\})$ and
$\deg\mu(\al)=\deg\mu\ci\ka(\{i+1,\ldots,n\})$ also
imply~$\hat\mu\bigl(\si(\{i+1,\ldots,n\})\bigr)\le\hat\mu(X)$.
\end{proof}

In Theorem \ref{an5thm1} choose $(\hat\tau,\hat
T,\le)=(\de,D_m,\le)$ and $(\tau,T,\le),(\ti\tau,\ti T,\le)$
arbitrary from \cite[\S 4.4]{Joyc7}. Then $(\hat\tau,\hat T,\le)$
dominates $(\tau,T,\le),(\ti\tau,\ti T,\le)$, and the change from
$(\hat\tau,\hat T,\le)$ to $(\ti\tau,\ti T,\le)$ is locally finite.
So Theorem \ref{an5thm1} implies:

\begin{cor} Let\/ $P$ be a projective $\K$-scheme and\/
$\A=\coh(P)$, $K(\coh(P))$, $\fF_{\coh(P)}$ be as in
{\rm\cite[Ex.~9.1 or 9.2]{Joyc5}}. Suppose $(\tau,T,\le)$ and\/
$(\ti\tau,\ti T,\le)$ are any two weak stability conditions on
$\coh(P)$ from {\rm\cite[\S 4.4]{Joyc7}}, which may be defined using
different ample line bundles $E,\ti E$ on $P$. Then
\eq{an5eq2}--\eq{an5eq5} hold for $\coh(P)$, as infinite convergent
sums. If\/ $(\tau,T,\le),(\ti\tau,\ti T,\le)$ are chosen from
{\rm\cite[Ex.s 4.16 or 4.17]{Joyc7}}, so that they are permissible,
then \eq{an5eq6}--\eq{an5eq7} also hold as infinite convergent sums.
\label{an5cor3}
\end{cor}

For {\it surfaces} changes between Gieseker stability conditions are
globally finite.

\begin{thm} Let\/ $P$ be a smooth projective surface, and\/
$\A\!=\!\coh(P),\ab K(\A),\ab\fF_\A$ as in
{\rm\cite[Ex.~9.1]{Joyc5}}. Suppose $(\tau,T,\le),(\ti\tau,\ti
T,\le)$ are any two permissible weak stability conditions on $\A$
from {\rm\cite[Ex.s 4.16--4.17]{Joyc7}}, which may be defined using
different ample line bundles $E,\ti E$. Then the changes between
$(\tau,T,\le)$ and\/ $(\ti\tau,\ti T,\le)$ are globally finite in
the sense of Definition \ref{an5def1}. So \eq{an5eq2}--\eq{an5eq7}
hold as finite sums in $\CF,\SFa(\fObj_\A)$ and\/
$\CF,\SF(\fM(K,\tl,\mu)_\A),$ and\/ \eq{an5eq11} holds.
\label{an5thm3}
\end{thm}

The author can prove some partial results on global finiteness of
changes between permissible weak stability conditions on
$\A=\coh(P)$ when $\dim P\ge 3$, but we will not give them as they
are complicated and inelegant. The main conclusion is this. Suppose
$(\tau,T,\le),(\ti\tau,\ti T,\le)$ are any two permissible weak
stability conditions on $\A$ from \cite[\S 4.4]{Joyc7}. Then we can
construct a finite sequence $(\tau,T,\le)=(\tau_0,T_0,\le),
(\tau_1,T_1,\le),\ldots,(\tau_n,T_n,\le)=(\ti\tau,\ti T,\le)$ of
permissible weak stability conditions on $\A$, such that the changes
between $(\tau_{i-1},T_{i-1},\le)$ and $(\tau_i,T_i,\le)$ are
globally finite for $i=1,\ldots,n$. This implies that \eq{an5eq11}
holds, so $\Hp_\tau,\ldots,\bLt_\tau$ are independent of choice
of~$(\tau,T,\le)$.

We show the {\it locally finite} condition in Theorem
\ref{an5thm1} is necessary.

\begin{ex} Take $P$ to be the projective space $\KP^2$, with
$\A=\coh(\KP^2)$ and $K(\A),\fF_\A$ as in \cite[Ex.~9.1]{Joyc5}.
Define $(\tau,T,\le)$ to be the trivial stability condition
$(0,\{0\},\le)$ on $\A$ and $(\ti\tau,\ti T,\le)$ the Gieseker
stability condition $(\ga,G_2,\le)$ of \cite[Ex.~4.16]{Joyc7},
defined using $\O_{\KP^2}(1)$. Suppose $x_1,\ldots,x_n$ are distinct
points in $\KP^2(\K)$. Then there is an exact sequence
\e
0\ra Y_{x_1,\ldots,x_n}\ra\O_{\KP^2}(-1)\op\O_{\KP^2}(0)\ra
\O_{\KP^2}(-1)\op\ts\bigop_{i=1}^n\O_{x_i}\ra 0,
\label{an5eq12}
\e
where $\O_{x_i}$ is the structure sheaf of $x_i$. Define
$X=\O_{\KP^2}(-1)\op\O_{\KP^2}(0)$ in $\A$, $\al=[X]$ in $K(\A)$,
and $\ka:\{1,2\}\ra C(\A)$ by $\ka(1)=[Y_{x_1,\ldots,x_n}]=
[\O_{\KP^2}(0)]-n[\O_{x_1}]$ and $\ka(2)=[\O_{\KP^2}(-1)]+
n[\O_{x_1}]$, so that~$\ka(\{1,2\})=\al$.

Calculation shows $X,Y_{x_1,\ldots,x_n}$ have Hilbert
polynomials $p_X(t)=t^2+2t+1$ and $p_{Y_{x_1,\ldots,x_n}}=
\frac{1}{2}t^2+\frac{3}{2}t+1-n$. Thus $\ga\ci\ka(1)=
t^2+3t+2-2n$ and $\ga\ci\ka(2)=t^2+t+2n$, so that
$\ga\ci\ka(1)>\ga\ci\ka(\{2\})$ and $S(\{1,2\},\le,\ka,0,\ga)=-1$.
As all sheaves are $0$-semistable, \eq{an5eq12} implies
$\bdss^{\ka(1)}(0)*\bdss^{\ka(2)}(0)\ne 0$ over~$[X]$.

Hence for each $n\ge 0$ we have distinct $(\{1,2\},\le,\ka)$ giving
nonzero terms over $[X]$ in \eq{an5eq3}. Therefore \eq{an5eq3} does
not make sense, even with the notion of convergence in Definition
\ref{an2def10}. This works because the change from $(0,\{0\},\le)$
to $(\ga,G_2,\le)$ is not locally finite.
\label{an5ex}
\end{ex}

The rest of the section proves Theorems \ref{an5thm1}, \ref{an5thm2}
and~\ref{an5thm3}.

\subsection{Transforming to a dominant weak stability condition}
\label{an52}

Let Assumption \ref{an3ass} hold and $(\tau,T,\le),(\ti\tau,\ti
T,\le)$ be weak stability conditions on $\A$ with $(\ti\tau,\ti T,
\le)$ {\it dominating} $(\tau,T,\le)$ in the sense of Definition
\ref{an3def13}. Suppose $(\{1,\ldots,n\},\le,\ka)$ is $\A$-data with
$S(\{1,\ldots,n\},\le,\ka,\tau,\ti\tau)\ne 0$. Then Theorem
\ref{an4thm2} gives $k,l$ with $\tau\ci\ka(k)\le\tau\ci\ka(l)$ and
$\ti\tau\ci \ka(k)\ge\ti\tau\ci\ka(l)$. But as $(\ti\tau,\ti T,\le)$
dominates $(\tau,T,\le)$, $\tau\ci\ka(k)\le\tau\ci\ka(l)$ implies
$\ti\tau\ci \ka(k)\le\ti\tau\ci\ka(l)$, and so
$\ti\tau\ci\ka(k)=\ti\tau\ci \ka(l)$. Therefore all
$\ti\tau\ci\ka(i)$ are equal, by Theorem~\ref{an4thm2}.

Using this and Definition \ref{an4def1} we see that if
$(\ti\tau,\ti T,\le)$ dominates $(\tau,T,\le)$ and
$(\{1,\ldots,n\},\le,\ka)$ is $\A$-data with $\ka(\{1,
\ldots,n\})=\al$ then
\e
S(\{1,\ldots,n\},\le,\ka,\tau,\ti\tau)=
\begin{cases} 1, & \tau\ci\ka(1)\!>\!\cdots\!>\!\tau\ci\ka(n),
\;\ti\tau\ci\ka\!\equiv\!\ti\tau(\al), \\
0, &\text{otherwise.} \end{cases}
\label{an5eq13}
\e
In our next theorem, \eq{an5eq13} shows \eq{an5eq14}--\eq{an5eq15}
are special cases of~\eq{an5eq2}--\eq{an5eq3}.

\begin{thm} Let Assumption \ref{an3ass} hold, $(\tau,T,\le),
(\ti\tau,\ti T,\le)$ be weak stability conditions on $\A$ with\/
$(\ti\tau,\ti T,\le)$ dominating $(\tau,T,\le)$, and\/ $\al\in
C(\A)$. Then
\begin{gather}
\sum_{\substack{\text{$\A$-data $(\{1,\ldots,n\},\le,\ka)$:
$\ti\tau\ci\ka\equiv\ti\tau(\al)$,}\\
\text{$\tau\ci\ka(1)\!>\!\cdots\!>\!\tau\ci\ka(n)$,
$\ka(\{1,\ldots,n\})\!=\!\al$}}}\!\!\!\!\!\!\!\!\!\!\!\!\!\!\!\!\!\!
\dss^{\ka(1)}(\tau)*\dss^{\ka(2)}(\tau)*\cdots*
\dss^{\ka(n)}(\tau)=\dss^\al(\ti\tau),
\label{an5eq14}\\
\sum_{\substack{\text{$\A$-data $(\{1,\ldots,n\},\le,\ka)$:
$\ti\tau\ci\ka\equiv\ti\tau(\al)$,}\\
\text{$\tau\ci\ka(1)\!>\!\cdots\!>\!\tau\ci\ka(n)$,
$\ka(\{1,\ldots,n\})\!=\!\al$}}}\!\!\!\!\!\!\!\!\!\!\!\!\!\!\!\!\!\!
\bdss^{\ka(1)}(\tau)*\bdss^{\ka(2)}(\tau)*\cdots*
\bdss^{\ka(n)}(\tau)=\bdss^\al(\ti\tau).
\label{an5eq15}
\end{gather}
These are potentially infinite sums, converging as in
Definition~\ref{an2def10}.
\label{an5thm4}
\end{thm}

\begin{proof} Consider $\A$-data $(\{1,\ldots,n\},\le,\ka)$ with
$\tau\ci\ka(1)\!>\!\cdots\!>\!\tau\ci\ka(n)$ and $\dss^{\ka(1)}
(\tau)\!*\!\cdots\!*\!\dss^{\ka(n)}(\tau)\!\ne\!0$ at
$[X]\!\in\!\fObj_\A(\K)$. By \eq{an3eq12}
$\CF^\stk(\bs\si(\{1,\ldots,n\}))\ab
\dss(\{1,\ldots,n\},\le,\ka,\tau)\!\ne\!0$ at $[X]$, so we have
$[(\si,\io,\pi)]\!\in\!\Mss(\{1,\ldots,n\},\le,\ka,\tau)_\A$ with
$\si(\{1,\ldots,n\})\!=\!X$. From \cite[Cor.~4.4]{Joyc5}, such
$[(\si,\io,\pi)]$ are in 1-1 correspondence with filtrations
$0=A_0\subset\cdots\subset A_n=X$ with $S_i=A_i/A_{i-1}\cong
\si(\{i\})$ for $i=1,\ldots,n$. Thus $S_i$ is $\tau$-semistable and
$\tau([S_1])>\cdots>\tau([S_n])$. Theorem \ref{an3thm} now shows
that $0=A_0\subset\cdots\subset A_n=X$ is the unique {\it
Harder--Narasimhan filtration} of $X$. So $n,\ka$ and $[(\si,
\io,\pi)]$ are unique for fixed~$[X]$.

Now $\Iso_\K\bigl([(\si,\io,\pi)]\bigr)\!\cong\!\Aut((\si,\io,\pi))
\!\cong\!\Aut(A_0\!\subset\!\cdots\!\subset\!A_n\!=\!X)$ and
$\Iso_\K\bigl([X]\bigr)\ab \!\cong\!\Aut(X)$. But
$\Aut(A_0\!\subset\! \cdots\!\subset\!A_n\!=\!X)\!=\!\Aut(X)$ as $X$
determines $A_0\!\subset\!\cdots\!\subset\!A_n\!=\!X$ by Theorem
\ref{an3thm}. Thus $\bs\si(\{1,\ldots,n\})_*:\Iso_\K\bigl(
[(\si,\io,\pi)]\bigr)\ra\Iso_\K\bigl([X]\bigr)$ is an isomorphism of
$\K$-groups, and $m_{\bs\si(\{1,\ldots,n\})}\bigl([(\si,\io,\pi)]
\bigr)\!=\!1$ in Definition \ref{an2def3}. Equation \eq{an3eq12}
then gives $\dss^{\ka(1)}(\tau)*\cdots*\dss^{\ka(n)}(\tau)=1$
at~$[X]$.

As $(\ti\tau,\ti T,\le)$ dominates $(\tau,T,\le)$ we have
$\ti\tau\ci\ka(1)\!\ge\!\cdots\!\ge\!\ti\tau\ci\ka(n)$, and the
$\si(\{i\})$ are $\ti\tau$-semistable by \cite[Lem.~4.11]{Joyc7}. It
easily follows that $X$ is $\ti\tau$-semistable if and only if
$\ti\tau\ci\ka\equiv\ti\tau(\al)$. Thus, if $X$ is
$\ti\tau$-semistable and $[X]=\al$ there exist unique $n,\ka$ in
\eq{an5eq14} with $\dss^{\ka(1)}(\tau)*\cdots*\dss^{\ka(n)}(\tau)=1$
at $[X]$ and all other terms zero at $[X]$, and otherwise all terms
in \eq{an5eq14} are zero at $[X]$. That is, the l.h.s.\ of
\eq{an5eq14} is 1 at $[X]$ if $[X]\in\Oss^\al(\ti\tau)$, and 0
otherwise. This proves~\eq{an5eq14}.

Consider the map on $\fObj_\A(\K)$ taking $[X]$ to $(n,\ka)$
constructed from the Harder--Narasimhan filtration as above. In the
natural topology on $\fObj_\A(\K)$, it is easy to see $n$ is an
upper semicontinuous function of $[X]$, and $\ka$ is locally
constant on the subset of $\fObj_\A(\K)$ with fixed $n$. Thus the
map $[X]\mapsto(n,\ka)$ is {\it locally constructible}. Hence if
$\fG\subset\fObj_\A$ is a finite type substack then $[X]\mapsto
(n,\ka)$ takes only finitely many values on $\fG(\K)$, and the
restriction of \eq{an5eq14} to $\fG$ has only finitely many nonzero
terms, so \eq{an5eq14} converges.

Next we extend all this to \eq{an5eq15}. The argument above shows
that $[X]\in\fObj_\A(\K)$ lies in $\Oss^\al(\ti\tau)$ if and
only if $[X]=\bs\si(\{1,\ldots,n\})_*\bigl([(\si,\io,\pi)]\bigr)$
for some $[(\si,\io,\pi)]\in\Mss(\{1,\ldots,n\},\le,\ka,\tau)_\A$
and $n,\ka$ as in \eq{an5eq15}, which are then unique. This gives
a disjoint union
\begin{equation*}
\Oss^\al(\ti\tau)=
\!\!\!\!\!\!\!\!\!\!\!\!\!\!\!
\coprod_{\substack{\text{$\A$-data $(\{1,\ldots,n\},\le,\ka)$:
$\ti\tau\ci\ka\equiv\ti\tau(\al)$,}\\
\text{$\tau\ci\ka(1)\!>\!\cdots\!>\!\tau\ci\ka(n)$,
$\ka(\{1,\ldots,n\})\!=\!\al$}}}
\!\!\!\!\!\!\!\!\!\!\!\!\!\!\!
\bs\si(\{1,\ldots,n\})_*\bigl(\Mss(\{1,\ldots,n\},\le,\ka,\tau)_\A\bigr),
\end{equation*}
where for any finite type substack $\fG\subseteq\fObj_\A$, only
finitely many sets on the r.h.s.\ intersect $\fG(\K)$. This
translates into an identity in $\dLSF(\fObj_\A)$:
\e
\bdss^\al(\ti\tau)= \!\!\!\!\!\!\!\!\!\!\!\!
\sum_{\substack{\text{$\A$-data $(\{1,\ldots,n\},\le,\ka)$:
$\ti\tau\ci\ka\equiv\ti\tau(\al)$,}\\
\text{$\tau\ci\ka(1)\!>\!\cdots\!>\!\tau\ci\ka(n)$,
$\ka(\{1,\ldots,n\})\!=\!\al$}}}
\!\!\!\!\!\!\!\!\!\!\!\!
\bde_{\bs\si(\{1,\ldots,n\})_*(\Mss(\{1,\ldots,n\},\le,\ka,\tau)_\A)},
\label{an5eq16}
\e
which is a potentially infinite convergent sum.

Now for $n,\ka$ as in \eq{an5eq15}, the argument above shows that
$\bs\si(\{1,\ldots,n\})_*:\Mss(\{1,\ldots,n\},\le,\ka,\tau)_\A\ra
\bs\si(\{1,\ldots,n\})_*\bigl(\Mss(\{1,\ldots,n\},\le,\ka,\tau)_\A
\bigr)$ is a bijection, inducing isomorphisms of stabilizer groups.
Here each side is the set of $\K$-points in a locally closed
$\K$-substack of $\fM(\{1,\ldots,n\},\le,\ka)_\A$ and $\fObj_\A$, and
one can strengthen the argument to show that $\bs\si(\{1,\ldots,n\})$
induces a 1-isomorphism of these substacks. Therefore
\begin{equation*}
\bs\si(\{1,\ldots,n\})_*\bigl(\bde_{\Mss(\{1,\ldots,n\},
\le,\ka,\tau)_\A}\bigr)=\bde_{\bs\si(\{1,\ldots,n\})_*
(\Mss(\{1,\ldots,n\},\le,\ka,\tau)_\A)}.
\end{equation*}
Combining this with Definition \ref{an3def11} and \eq{an3eq12} gives
\e
\bdss^{\ka(1)}(\tau)*\cdots*\bdss^{\ka(n)}(\tau)=
\bde_{\bs\si(\{1,\ldots,n\})_*(\Mss(\{1,\ldots,n\},\le,\ka,\tau)_\A)}.
\label{an5eq17}
\e
Equation \eq{an5eq15}, and its convergence, then
follow from \eq{an5eq16} and~\eq{an5eq17}.
\end{proof}

\subsection{Inverting equations \eq{an5eq14}--\eq{an5eq15}}
\label{an53}

In the situation of \S\ref{an52}, if $(\ti\tau,\ti T,\le)$
dominates $(\tau,T,\le)$ and $(\{1,\ldots,m\},\le,\la)$ is
$\A$-data with $\la(\{1,\ldots,m\})=\al$ then a similar
argument to \eq{an5eq13} shows that
\ea
\!\! S(\{1,\ldots,m\},\le,\la,\ti\tau,\tau)\!=\!
\begin{cases} (-1)^{m\!-\!1}, &\begin{aligned}
&\ts\text{$\tau\!\ci\!\la(\{1,\ldots,i\})\!>\!\tau\!\ci\!
\la(\{i\!+\!1,\ldots,m\})$}\\
&\ts\text{for $i\!=\!1,\ldots,m-1$, and
$\ti\tau\!\ci\la\!\equiv\!\ti\tau(\al)$,}
\end{aligned} \\
0, &\text{otherwise.} \end{cases}
\nonumber\\[-16pt]
\label{an5eq18}
\ea
Equation \eq{an5eq18} shows \eq{an5eq19}--\eq{an5eq20} are
special cases of~\eq{an5eq2}--\eq{an5eq3}.

\begin{thm} Let Assumption \ref{an3ass} hold, $(\tau,T,\le),
(\ti\tau,\ti T,\le)$ be weak stability conditions on $\A$ with\/
$(\ti\tau,\ti T,\le)$ dominating $(\tau,T,\le)$ and the change
from $(\ti\tau,\ti T,\le)$ to $(\tau,T,\le)$ locally finite,
and\/ $\al\in C(\A)$. Then
\begin{gather}
\sum_{\substack{\text{$\A$-data $(\{1,\ldots,m\},\le,\la)$:
$\ti\tau\ci\la\equiv\ti\tau(\al)$,}\\
\text{$\tau\!\ci\!\la(\{1,\ldots,i\})\!>\!\tau\!\ci\!
\la(\{i\!+\!1,\ldots,m\})$}\\
\text{for $1\!\le\!i\!<\!m$, $\la(\{1,\ldots,m\})\!=\!\al$}}
\!\!\!\!\!\!\!\!\!\!\!\!}
\!\!\!\!\!\!\!\!\!\!\!\!\!\!\!\!\!\!\!\!\!\!\!
(-1)^{m-1}\dss^{\la(1)}(\ti\tau)*\dss^{\la(2)}(\ti\tau)*\cdots*
\dss^{\la(m)}(\ti\tau)=\dss^\al(\tau),
\label{an5eq19}\\
\sum_{\substack{\text{$\A$-data $(\{1,\ldots,m\},\le,\la)$:
$\ti\tau\ci\la\equiv\ti\tau(\al)$,}\\
\text{$\tau\!\ci\!\la(\{1,\ldots,i\})\!>\!\tau\!\ci\!
\la(\{i\!+\!1,\ldots,m\})$}\\
\text{for $1\!\le\!i\!<\!m$, $\la(\{1,\ldots,m\})\!=\!\al$}}
\!\!\!\!\!\!\!\!\!\!\!\!}
\!\!\!\!\!\!\!\!\!\!\!\!\!\!\!\!\!\!\!\!\!\!\!
(-1)^{m-1}\bdss^{\la(1)}(\ti\tau)*\bdss^{\la(2)}(\ti\tau)*\cdots*
\bdss^{\la(m)}(\ti\tau)=\bdss^\al(\tau).
\label{an5eq20}
\end{gather}
These are potentially infinite sums, converging as in
Definition~\ref{an2def10}.
\label{an5thm5}
\end{thm}

\begin{proof} As the change from $(\ti\tau,\ti T,\le)$ to
$(\tau,T,\le)$ is locally finite, \eq{an5eq18} implies that
\eq{an5eq19}--\eq{an5eq20} converge. Equation \eq{an5eq19} follows
from
\begin{gather}
\sum_{\substack{\text{$\A$-data $(\{1,\ldots,m\},\le,\la)$:
$\ti\tau\ci\la\equiv\ti\tau(\al)$,}\\
\text{$\tau\!\ci\!\la(\{1,\ldots,i\})\!>\!\tau\!\ci\!
\la(\{i\!+\!1,\ldots,m\})$}\\
\text{for $1\!\le\!i\!<\!m$, $\la(\{1,\ldots,m\})\!=\!\al$}}
\!\!\!\!\!\!\!\!\!\!\!\!}
\!\!\!\!\!\!\!\!\!\!\!\!\!\!\!\!\!\!\!\!\!\!\!
(-1)^{m-1}\dss^{\la(1)}(\ti\tau)*\cdots*\dss^{\la(m)}(\ti\tau)
=\sum_{\substack{\text{$\A$-data}\\
\text{$(\{1,\ldots,m\},\le,\la):$}\\
\text{$\la(\{1,\ldots,m\})=\al$}}}
\nonumber
\\
S(\{1,\ldots,m\},\le,\la,\ti\tau,\tau)\cdot\!\!\!\!\!\!\!\!\!\!
\sum_{\substack{\text{$\A$-data}\\
\text{$(\{1,\ldots,n_i\},\le,\ka_i)$,}\\
\text{$i=1,\ldots,m:$}\\
\text{$\ka_i(\{1,\ldots,n_i\})=\la_i$,}\\
\text{all $i=1,\ldots,m$}}}\!\!\!\!\!\!\!\!\!\!\!\!
\begin{aligned}[t]
\ts\prod_{i=1}^m &S(\{1,\ldots,n_i\},\le,\ka_i,\tau,\ti\tau)\cdot \\
&\bigl(\dss^{\ka_1(1)}(\tau)*\!\cdots\!*\dss^{\ka_1(n_1)}
(\tau)\bigr)*\cdots\\
&*\bigl(\dss^{\ka_m(1)}(\tau)*\cdots*\dss^{\ka_m(n_m)}(\tau)\bigr)\!=
\end{aligned}
\nonumber
\allowdisplaybreaks
\\
\sum_{\substack{\text{$\A$-data $(\{1,\ldots,n\},\le,\ka):$}\\
\text{$\ka(\{1,\ldots,n\})=\al$}}}
\dss^{\ka(1)}(\tau)*\cdots*\dss^{\ka(n)}(\tau)\cdot
\nonumber
\\
\raisebox{-12pt}{\begin{Large}$\displaystyle\Biggl[$\end{Large}}
\sum_{m=1}^n\,
\sum_{\substack{\text{$\psi:\{1,\ldots,n\}\ra\{1,\ldots,m\}:$}\\
\text{$\psi$ is surjective,}\\
\text{$1\!\le\!i\!\le\!j\!\le\!n$ implies $\psi(i)\!\le\!\psi(j)$,}\\
\text{define $\la:\{1,\ldots,m\}\!\ra\!C(\A)$}\\
\text{by $\la(k)=\ka(\psi^{-1}\!(k))$}}}\,
\begin{aligned}[t]
&S(\{1,\ldots,m\},\le,\la,\ti\tau,\tau)\cdot\\
&\prod_{k=1}^mS\bigl(\psi^{-1}(\{k\}),\le,
\ka\vert_{\psi^{-1}(\{k\})},\tau,\ti\tau\bigr)
\end{aligned}
\raisebox{-12pt}{\begin{Large}$\displaystyle\Biggr]$\end{Large}}
\nonumber
\\[-16pt]
\qquad\qquad\qquad\qquad=\dss^\al(\tau).
\label{an5eq21}
\end{gather}

Here in the first step we use \eq{an5eq13} and \eq{an5eq18} and
substitute \eq{an5eq14} in for $\dss^{\la(i)}(\ti\tau)$, and in the
second we replace the sums over $m,\la$ and $n_i,\ka_i$ for $1\le
i\le m$ by sums over $n,\ka$ and $m,\psi,\la$, where
$n\!=\!n_1\!+\!\cdots \!+\!n_m$, and $\ka:\{1,\ldots,n\}\ra C(\A)$
is given by $\ka(j)= \ka_i(k)$ if
$j\!=\!m_1\!+\!\cdots\!+\!m_{i-1}\!+\!k$, $1\!\le\!k\! \le\!m_i$,
and $\psi:\{1,\ldots,n\}\!\ra\!\{1,\ldots,m\}$ is given by
$\psi(j)\!=\!i$ if $m_1\!+\!\cdots\!+\!m_{i-1}\!<\!j\!\le\!m_1
\!+\!\cdots\!+\!m_i$, and use associativity of $*$. For the third
and final step we note that $[\cdots]$ in the fourth line of
\eq{an5eq21} is $S(\{1,\ldots,n\},\le,\ka,\tau,\tau)$ by
\eq{an4eq7}, but this is 1 if $n=1$ and 0 otherwise by \eq{an4eq6}.
So the only contribution is $n=1$ and $\ka(1)=\al$,
giving~$\dss^\al(\tau)$.

As \eq{an5eq21} involves infinite sums, we must also check our
rearrangements of these are valid. Let $\fG\subseteq\fObj_\A$ be a
finite type $\K$-substack. We shall show that the restriction of
\eq{an5eq21} to $\fG$ has only finitely many nonzero terms at every
stage, so our argument above just rearranges finite sums over $\fG$,
and \eq{an5eq21} holds as a convergent sum. Since the change from
$(\ti\tau,\ti T,\le)$ to $(\tau,T,\le)$ is locally finite and
$\fG(\K)$ is constructible, by Definition \ref{an5def1} there are
only finitely many $(\{1,\ldots,m\},\le,\la)$ with
$S(\{1,\ldots,m\},\le,\la,\ti\tau,\tau)\ne 0$ and
\e
\fG(\K)\cap\bs\si(\{1,\ldots,m\})_*\bigl(\Mss(\{1,\ldots,m\},\le,
\la,\ti\tau)_\A\bigr)\ne\emptyset.
\label{an5eq22}
\e
Thus there are finitely many nonzero terms over $\fG$ in the first
step of~\eq{an5eq21}.

The second and third steps of \eq{an5eq21} are equivalent as they
are a relabelling. Suppose $n,\ka,m,\psi,\la$ give a nonzero term
over $\fG$ in the third step. Then $S(\{1,\ldots,m\},\le,\la,\ti
\tau,\tau)\ne 0$ and $S(\psi^{-1}(\{k\}),\le,\ka\vert_{
\psi^{-1}(\{k\})},\tau,\ti\tau)\ne 0$ for $k=1,\ldots,m$, and
\e
\fG(\K)\cap\bs\si(\{1,\ldots,n\})_*\bigl(\Mss(\{1,\ldots,n\},\le,
\ka,\tau)_\A\bigr)\ne\emptyset.
\label{an5eq23}
\e
Let $[(\si,\io,\pi)]\in\Mss(\{1,\ldots,n\},\le,\ka,\tau)_\A$. Then
$\si(\{i\})$ is $\tau$-semistable for all $i$, so that $\si(\{i\})$
is $\ti\tau$-semistable by \cite[Lem.~4.11]{Joyc7} as $(\ti\tau,\ti
T,\le)$ dominates $(\tau,T,\le)$. Also,
$S(\psi^{-1}(\{k\}),\le,\ka\vert_{\psi^{-1}(\{k\})},\tau,
\ti\tau)\ne 0$ implies that $\ti\tau\ci\ka$ is constant on
$\psi^{-1}(\{k\})$ by \eq{an5eq13}. It follows that
$\si(\psi^{-1}(\{k\}))$ is $\ti\tau$-semistable for $k=1,\ldots,m$,
which implies that $Q(\{1,\ldots,n\},
\le,\{1,\ldots,m\},\le,\psi)_*$ maps $\Mss(\{1,\ldots,n\},\ab
\le,\ab \ka,\ab \tau)_\A\ra\Mss(\{1,\ldots,m\},\le,\la,\ti\tau)_\A$.

As $\bs\si(\{1,\ldots,n\})=\bs\si(\{1,\ldots,m\})\ci
Q(\{1,\ldots,n\},\le,\{1,\ldots,m\},\le,\psi)$, this and
\eq{an5eq23} imply \eq{an5eq22}. Thus from the previous
part, there are only finitely many possibilities for $m,\la$
giving nonzero terms in the second and third steps of
\eq{an5eq21}. Fix one such choice. Then any
$(\{1,\ldots,n_i\},\le,\ka_i)$ for $i=1,\ldots,m$ giving a
nonzero term in the second step comes from the (unique)
$\tau$-Harder--Narasimhan filtration of some point $[X_i]$ in
\e
\bs\si(\{i\})_*\bigl[\bs\si(\{1,\ldots,m\})_*^{-1}(\fG(\K))
\cap\Mss(\{1,\ldots,m\},\le,\la,\ti\tau)_\A\bigr].
\label{an5eq24}
\e

Since $\bs\si(\{1,\ldots,m\})$ here is of finite type by
\cite[Th.~8.4(b)]{Joyc5} and $\fG(\K)$ is constructible,
$\bs\si(\{1,\ldots,m\})_*^{-1}(\fG(\K))$ is constructible, and so
\eq{an5eq24} is constructible. As in the proof of Theorem
\ref{an5thm1} the map $[X_i]\mapsto(n_i,\ka_i)$ is locally
constructible, and so takes only finitely many values on
\eq{an5eq24}. Therefore there are only finitely many choices for
$n_i,\ka_i$ giving nonzero terms, and the restriction of
\eq{an5eq21} to $\fG$ has only finitely many nonzero terms at each
step. This proves \eq{an5eq19}. The proof of \eq{an5eq20} is the
same, substituting \eq{an5eq15} in for~$\bdss^{\la(i)}(\ti\tau)$.
\end{proof}

\subsection{Proof of Theorem \ref{an5thm1}}
\label{an54}

Let $U\subseteq\fObj_\A(\K)$ be constructible, and suppose
$(\{1,\ldots,n\},\le,\ka)$ is $\A$-data for which $S(\{1,\ldots,n\},
\le,\ka,\tau,\ti\tau)\!\ne\!0$ and \eq{an5eq1} holds. Then by
\eq{an4eq7} there exist $\A$-data $(\{1,\ldots,m\},\le,\la)$
and surjective $\psi:\{1,\ldots,n\}\ra\{1,\ldots,m\}$ with $i\le j$
implies $\psi(i)\le\psi(j)$ and $\la(k)=\ka(\psi^{-1}(k))$ for
$k=1,\ldots,m$, such that $S(\{1,\ldots,m\},\le,\la,\hat\tau,
\ti\tau)\ne 0$ and $S\bigl(\psi^{-1}(\{k\}),\le,\ka\vert_{\psi^{-1}
(\{k\})},\tau,\hat\tau)\ne 0$ for all~$k$.

Since $(\hat\tau,\hat T,\le)$ dominates $(\tau,T,\le)$, the argument
of Theorem \ref{an5thm5} shows that $Q(\{1,\ldots,n\},\le,\{1,
\ldots,m\},\le,\psi)_*$ maps $\Mss(\{1,\ldots,n\},\le,\ka,\tau)_\A
\ra\ab \Mss(\{1,\ab \ldots,\ab m\},\le,\la,\hat\tau)_\A$, and
$U\cap\bs\si
(\{1,\ldots,m\})_*\bigl(\Mss(\{1,\ldots,m\},\le,\la,\hat\tau)_\A
\bigr)\ne\emptyset$. As the change from $(\hat\tau,\hat T,\le)$ to
$(\ti\tau,\ti T,\le)$ is locally finite, by Definition \ref{an5def1}
there are only finitely many possibilities for $m,\la$. The argument
of Theorem \ref{an5thm5} then shows there are only finitely many
possibilities for $n,\ka$. Hence the change from $(\tau,T,\le)$ to
$(\ti\tau,\ti T,\le)$ is locally finite, as we want. This implies
\eq{an5eq2}--\eq{an5eq3} are convergent in the sense of
Definition~\ref{an2def10}.

To prove \eq{an5eq2}--\eq{an5eq3} we now substitute
\eq{an5eq14}--\eq{an5eq15} with $\hat\tau,\la(i)$ in place of
$\ti\tau,\al$ for $i=1,\ldots,m$ into \eq{an5eq19}--\eq{an5eq20}
with $\hat\tau,\ti\tau$ in place of $\ti\tau,\tau$ respectively. We
then rearrange these equations using \eq{an4eq7} exactly as for
\eq{an5eq21}, but without using \eq{an4eq6}. The argument of Theorem
\ref{an5thm5} shows that for finite type $\fG\subseteq\fObj_\A$
there are only finitely many nonzero terms over $\fG$ at each stage,
so the rearrangements are valid and \eq{an5eq2}--\eq{an5eq3} hold.

Next we prove \eq{an5eq4}--\eq{an5eq5}. Let $I,\pr,\ka,\phi$ be as
in \eq{an5eq4}, and set $n_k=\md{\phi^{-1}(\{k\})}$ for $k\in K$. By
Definition \ref{an4def2} $(\phi^{-1}(\{k\}),\pr)$ is a total order,
so there exists a unique bijection $\psi_k:\{1,
\ldots,n_k\}\ra\phi^{-1}(\{k\})$ with $\psi_k^*(\pr)=\le$. Since
$\dss(K,\tl,\mu,\ti\tau)=\prod_{k\in K}\bs\si(\{k\})^*
\bigl(\dss^{\mu(k)}(\ti\tau)\bigr)$, taking the product over $k\in
K$ of $\bs\si(\{k\})^*$ applied to \eq{an5eq2} with $\mu(k)$ in
place of $\al$ and using \eq{an3eq12} gives
\begin{equation}
\text{
\begin{small}
$\displaystyle
\begin{gathered}
\dss(K,\tl,\mu,\ti\tau)\!=\!\prod_{k\in K}\!\!\!
\sum_{\substack{\text{$\A$-data}\\
\text{$(\{1,\ldots,n_k\},\le,\ka_k):$}\\
\text{$\ka_k(\{1,\ldots,n_k\})=\mu(k)$}}}
\!\!\!\!\!\!\!\!\!\!\!\!\!\!\!
\begin{aligned}[t]
S(\{1,&\ldots,n_k\},\le,\ka_k,\tau,\ti\tau)\cdot\\
&\bs\si(\{k\})^*\ci\CF^\stk(\bs\si(\{1,\ldots,n_k\}))\\
&\dss(\{1,\ldots,n_k\},\le,\ka_k,\tau)
\end{aligned}
\\
=\!\!\!\sum_{\substack{\text{iso.}\\ \text{classes}\\
\text{of finite }\\ \text{sets $I$}}} \frac{1}{\md{I}!}\cdot\!\!\!\!
\sum_{\substack{
\text{$\pr,\ka,\phi$: $(I,\pr,\ka)$ is $\A$-data,}\\
\text{$(I,\pr,K,\phi)$ is dominant,}\\
\text{$\tl=\P(I,\pr,K,\phi)$,}\\
\text{$\ka(\phi^{-1}(k))=\mu(k)$ for $k\in K$}}}\!\!\!\!\!\!\!\!\!\!\!\!\!
\begin{aligned}[t]
\ts\prod_{k\in K}&S(\phi^{-1}(\{k\}),\pr,\ka,\tau,\ti\tau)\cdot\\
&\ts\prod_{k\in K}\bigl[\bs\si(\{k\})^*\ci
\CF^\stk(\bs\si(\phi^{-1}(\{k\})))\\
&\dss(\phi^{-1}(\{k\}),\pr,\ka,\tau)\bigr].
\end{aligned}
\end{gathered}
$
\end{small}
}
\!\!\!\!\!\!\!\!\!\!\!\!\!\!\!
\label{an5eq25}
\end{equation}

Here in the second line of \eq{an5eq25} we note that as above, to
each choice of $I,\pr,\ka,\phi$ in \eq{an5eq4} we can associate
$n_k,\psi_k$ as above and $\ka_k=\ka\ci\psi_k$ as in the first line
of \eq{an5eq25}. Conversely, the data $n_k,\ka_k$ for $k\in K$
determine $I,\pr,\ka,\phi$ uniquely up to isomorphism. However, the
sum in \eq{an5eq25} is not over isomorphism classes of quadruples
$(I,\pr,\ka,\phi)$, but rather over isomorphism classes of $I$,
followed by a sum over all $\pr,\ka,\phi$. It is easy to see that
$(I,\pr,\ka,\phi)$ and $(I,\pr',\ka',\phi')$ are isomorphic if and
only if $(\pr',\ka',\phi')=\pi_*(\pr,\ka,\phi)$ for some permutation
$\pi:I\ra I$, and different $\pi$ give different
$(\pr',\ka',\phi')$. So we include the factor $1/\md{I}!$ to cancel
the number $\md{I}!$ of permutations of~$I$.

As \eq{an5eq25} involves infinite sums, we must also consider the
convergence issues, and whether our rearrangements of sums are
valid. Let $\fG\subseteq\fM(K,\tl,\mu)_\A$ be a finite type
$\K$-substack. Then $\fG(\K)$ is constructible, so
$\bs\si(\{k\})_*(\fG(\K))$ is constructible in $\fObj_\A(\K)$ for
each $k\in K$, so by the convergence of \eq{an5eq2}--\eq{an5eq3}
there are only finitely many choices of $n_k,\ka_k$ giving nonzero
terms over $\fG$ in \eq{an5eq25}, and so only finitely many choices
up to isomorphism of $I,\pr,\ka,\phi$. This proves \eq{an5eq25}, as
a convergent infinite sum in~$\LCF(\fM(K,\tl,\mu)_\A)$.

Now for $I,\pr,\ka,\phi$ as in \eq{an5eq4}, by applying the proof of
\cite[Th.~7.10]{Joyc5} $\md{K}$ times, we can show that the
following is a Cartesian square:
\e
\begin{gathered}
\xymatrix@C=100pt@R=12pt{
\fM(I,\pr,\ka)_\A \ar@<-4ex>[d]^{\prod_{k\in K}S(I,\pr,\phi^{-1}(\{k\}))}
\ar[r]_{Q(I,\pr,K,\tl,\phi)} &
\fM(K,\tl,\mu)_\A \ar@<4ex>[d]_{\prod_{k\in K}\bs\si(\{k\})} \\
\prod_{k\in K}\fM(\phi^{-1}(\{k\}),\pr,\ka)_\A
\ar[r]^(0.55){\prod_{k\in K}\bs\si(\phi^{-1}(\{k\}))} &
\prod_{k\in K}\fObj_\A^{\mu(k)}.
}
\end{gathered}
\label{an5eq26}
\e
By \cite[Th.~8.4]{Joyc5} the rows of \eq{an5eq26} are representable
and the right 1-morphism is finite type, so the left 1-morphism is
too as \eq{an5eq26} is Cartesian. Applying Theorem \ref{an2thm1}
then shows that
\begin{equation*}
\text{
\begin{small}
$\displaystyle
\begin{aligned}
&\ts\prod_{k\in K}\bigl[\bs\si(\{k\})^*\ci\CF^\stk\bigl(\bs\si(\phi^{-1}
(\{k\}))\bigr)\dss(\phi^{-1}(\{k\}),\pr,\ka,\tau)\bigr]\\
&=\!\ts\bigl(\prod_{k\in K}\bs\si(\{k\})\bigr)^*
\ci\CF^\stk\bigl(\prod_{k\in K}\bs\si(\phi^{-1}(\{k\}))\bigr)
\bigl[\prod_{k\in K}\dss(\phi^{-1}(\{k\}),\pr,\ka,\tau)\bigr]\\
&=\!\ts\CF^\stk\bigl(Q(I,\pr,K,\tl,\phi)\bigr)\!\ci\!
\bigl(\prod\limits_{k\in K\!\!\!\!}S(I,\pr,\phi^{-1}(\{k\}))\bigr)^*\bigl[
\ts\prod\limits_{k\in K\!\!\!\!}\dss(\phi^{-1}(\{k\}),\pr,\ka,\tau)\bigr]
\\[-4pt]
&=\!\CF^\stk\bigl(Q(I,\pr,K,\tl,\phi)\bigr)\dss(I,\pr,\ka,\tau).
\end{aligned}
$
\end{small}
}
\end{equation*}

Substituting this into the last line of \eq{an5eq25} and using
\eq{an4eq4} then proves \eq{an5eq4}, and its convergence. For
\eq{an5eq5} we use the same argument with \eq{an5eq3} instead of
\eq{an5eq2}, $(\cdots)_*$ instead of $\CF^\stk(\cdots)$, Theorem
\ref{an2thm2} instead of Theorem \ref{an2thm1}, and replacing
products of functions such as $\prod_{k\in K}
\bs\si(\{k\})^*\bigl(\dss^{\mu(k)}(\ti\tau)\bigr)$ by expressions
like~$\bigl(\prod_{k\in K}\bs\si(\{k\})\bigr)^* \bigl(\bigot_{k\in
K}\bdss^{\mu(k)}(\ti\tau)\bigr)$.

To prove \eq{an5eq6}--\eq{an5eq7}, let $(\tau,T,\le)$ and
$(\ti\tau,\ti T,\le)$ be permissible. Substituting
\eq{an4eq5} into the left hand side of \eq{an5eq6} and
rewriting gives
\begin{gather*}
\sum_{\substack{\text{$\A$-data}\\
\text{$(\{1,\ldots,m\},\le,\la):$}\\
\text{$\la(\{1,\ldots,m\})=\al$}}}\,
\sum_{1\le l\le m}\frac{(-1)^{l-1}}{l}\cdot
\!\!\!\!\!\!\!\!\!\!\!\!\!\!\!\!\!\!\!\!
\sum_{\substack{\text{surjective}\\
\text{$\xi:\{1,\ldots,m\}\!\ra\!\{1,\ldots,l\}$:}\\
\text{$i\!\le\!j$ implies $\xi(i)\!\le\!\xi(j)$.}\\
\text{Define $\mu:\{1,\ldots,l\}\ra C(\A)$ by}\\
\text{$\mu(a)\!=\!\la(\xi^{-1}(a))$.
Then $\ti\tau\ci\mu\!\equiv\!\ti\tau(\al)$}}}
\!\!\!\!\!\!\!\!\!\!\!\!\!\!\!\!\!\!\!\!\!\!\!\!
\begin{aligned}[t]
\ts\prod_{a=1}^lS(\xi^{-1}(\{a\}),\le,\la,\tau,\ti\tau)\cdot\\
\dss^{\la(1)}(\tau)*\cdots*\dss^{\la(m)}(\tau)
\end{aligned}
\\
=\!\!\!\!\!
\sum_{\substack{\text{$\A$-data $(\{1,\ldots,l\},\le,\mu):$}\\
\text{$\mu(\{1,\ldots,l\})=\al$, $\ti\tau\ci\mu\equiv\ti\tau(\al)$}}}
\!\!\!\!\!\!\!\!\!\!
\frac{(-1)^{l-1}}{l}\cdot
\dss^{\mu(1)}(\ti\tau)*\cdots*\dss^{\mu(l)}(\ti\tau)=\ep^\al(\ti\tau),
\end{gather*}
using \eq{an3eq15} at the first step, \eq{an5eq2} at the second, and
\eq{an3eq13} at the third. These rearrangements are valid, and the
equation holds as an infinite convergent sum, since \eq{an5eq2} is
convergent and \eq{an3eq13}, \eq{an3eq15} have only finitely many
nonzero terms. This proves \eq{an5eq6} and its convergence, and
\eq{an5eq7} is the same, using \eq{an3eq16}, \eq{an5eq3} and
\eq{an3eq14}. The final part of Theorem \ref{an5thm1} is immediate.

\subsection{Proof of Theorem \ref{an5thm2}}
\label{an55}

Let Condition \ref{an4cond} hold for $(\tau,T,\le),(\ti\tau,\ti
T,\le)$. When $\md{I}=1$ equation \eq{an5eq8} is trivial. Suppose by
induction that $n>1$ and \eq{an5eq8} holds whenever $\md{I}<n$, and
let $\ka:I\ra C(\A)$ with $\md{I}=n$. We shall prove \eq{an5eq8}
for~$I,\ka$.

As in \cite[Ex.~7.10]{Joyc7}, define a quiver $Q=(Q_0,Q_1,b,e)$ to
have vertices $I$ and an edge $\smash{\mathop{\bu}\limits^{\sst
i}\ra\mathop{\bu}\limits^{\sst j}}$ for all $i,j\in I$. Consider the
abelian category $\nilCQ$ of nilpotent $\C$-representations of $Q$,
with data $K(\nilCQ),\fF_\nilCQ$ satisfying Assumption \ref{an3ass}
as in \cite[Ex.~10.6]{Joyc5}. Then $K(\nilCQ)=\Z^I$ and
$C(\nilCQ)=\N^I\sm\{0\}$. Write elements of $\Z^I$ as $\ga:I\ra\Z$.
For $i\in I$ define $e_i\in C(\nilCQ)$ by $e_i(j)=1$ if $j=i$ and
$e_i(j)=0$ otherwise. Then $\sum_{i\in I}e_i=1$ in $C(\nilCQ)$.
Define $\tau':C(\nilCQ)\ra T$ and $\ti\tau':C(\nilCQ)\ra\ti T$ by
$\tau'(\ga)=\tau\bigl(\sum_{i\in I}\ga(i)\ka(i)\bigr)$ and
$\ti\tau'(\ga)=\ti\tau\bigl(\sum_{i\in I}\ga(i)\ka(i)\bigr)$. Then
Condition \ref{an4cond}(iv),(v) imply $(\tau',T,\le),(\ti\tau',\ti
T,\le)$ are weak stability conditions on $\nilCQ$, which are
permissible by~\cite[Cor.~4.13]{Joyc7}.

Apply \eq{an5eq9} in $\nilCQ$ with $\al=1\in C(\nilCQ)$ and
$\tau',\ti\tau',J,\la$ in place of $\tau,\ti\tau,I,\ka$. By
\cite[Prop.~7.11(a)]{Joyc7}, if $J$ is finite and $\la:J\ra
C(\nilCQ)$ with $\la(J)=1$ then $\md{J}\le\md{I}=n$, and if
$\md{J}=n$ then there is a unique bijection $\imath:J\ra I$ with
$\la(j)=e_{\imath(j)}\in C(\A)$ for all $j\in J$. The number of
$\la:J\ra C(\nilCQ)$ is the number $\md{I}!$ of bijections $\imath$,
which cancels the factor $1/\md{I}!$ in \eq{an5eq9}. Taking $J=I$
and $\la:i\mapsto e_i$, the definition of $\tau',\ti\tau'$ implies
that $U(I,\pr,\la,\tau',\ti\tau')=U(I,\pr,\ka,\tau,\ti\tau)$. Thus
we may rewrite \eq{an5eq9} as
\begin{gather}
\begin{gathered}
\sum_{\substack{\text{total orders $\pr$ on $I$.}\\
\text{Write $I=\{i_1,\ldots,i_n\}$,}\\
\text{$i_1\pr i_2\pr\cdots\pr i_n$}}}\!\!\!\!
U(I,\pr,\ka,\tau,\ti\tau)\cdot\ep^{e_{i_1}}(\tau')*\!
\cdots\!*\ep^{e_{i_n}}(\tau')=\ep^1(\ti\tau')-
\end{gathered}
\label{an5eq27}
\\
\sum_{\substack{\text{iso classes of}\\
\text{finite sets $J$}\\ \text{with $\md{J}<n$}}}
\frac{1}{\md{J}!}\!
\sum_{\substack{\la:J\ra C(\nilCQ):\\ \la(J)=1}}\!
\raisebox{-12pt}{\begin{Large}$\displaystyle\Biggl[$\end{Large}}
\!\!
\sum_{\substack{\text{total orders $\pr$ on $J$.}\\
\text{Write $J=\{j_1,\ldots,j_m\}$,}\\
\text{$j_1\pr j_2\pr\cdots\pr j_m$}}}\,
\begin{aligned}[t]
&U(J,\pr,\la,\tau',\ti\tau')\cdot\\
&\ep^{\la(j_1)}(\tau')\!*\!\cdots\!*\!\ep^{\la(j_m)}(\tau')
\end{aligned}
\raisebox{-12pt}{\begin{Large}$\displaystyle\Biggr]$\end{Large}}.
\nonumber
\end{gather}

Now $\ep^1(\ti\tau'),\ep^{\la(j_a)}(\tau')$ lie in
$\CFi(\fObj_\nilCQ)$ by \cite[Th.~7.8]{Joyc7}. By induction, as
$\md{J}<n$ the term $[\cdots]$ is a sum of commutators of
$\ep^{\la(j_a)}(\tau')$ and so lies in $\CFi(\fObj_\nilCQ)$. So
every term on the right of \eq{an5eq27} lies in
$\CFi(\fObj_\nilCQ)$, and thus so does the left hand side. However,
\cite[Prop.~7.11(c)]{Joyc7} implies the subalgebra $A_{I,\tau'}$ of
$\CF(\fObj_\nilCQ)$ generated by the $\ep^{e_i}(\tau')$ for $i\in I$
is {\it freely generated}, so that $A_{I,\tau'}\cong\Q[I]$. We can
then deduce from \cite[Prop.~4.14]{Joyc6} that the Lie algebra
$A_{I,\tau'}\cap\CFi(\fObj_\nilCQ)$ is freely generated by the
$\ep^{e_i}(\tau')$ for $i\in I$, so that $A_{I,\tau'}\cap
\CFi(\fObj_\nilCQ)\cong\L[I]$. Equation \eq{an5eq8} for $I,\ka$
follows, as the l.h.s.\ of \eq{an5eq27} lies in
$A_{I,\tau'}\cap\CFi(\fObj_\nilCQ)$. This completes the proof.

\subsection{Proof of Theorem \ref{an5thm3}}
\label{an56}

Let $P$ be a smooth projective $\K$-scheme of dimension $m$ and
define $\A=\coh(P),K(\A),\fF_\A$ satisfying Assumption \ref{an3ass}
as in \cite[Ex.~9.1]{Joyc5}. We allow $m$ arbitrary for our first
two propositions, and then restrict to $m=2$. As in
\cite[Ex.~9.1]{Joyc5} the {\it Chern character} gives a natural
injective group homomorphism
\e
\ch:K(\coh(P))\ra\ts\bigop_{i=1}^mH^{2i}(P,\Z).
\label{an5eq28}
\e
Here $H^*(P,\Z)$ makes sense when $\K=\C$. For general $\K$ we can
instead use the {\it Chow ring} $A(P)$, but we will not worry about
this. Write $\ch_i(\al)$ for the component of $\ch(\al)$ in
$H^{2i}(P,\Z)$. These may be written in terms of the Chern classes
$c_i(X)$ and the rank $\rk(X)$ as on \cite[p.~432]{Hart} by
\e
\ch_0([X])\!=\!\rk(X),\;\> \ch_1([X])\!=\!c_1(X),\;\>
\ch_2([X])\!=\!\ts\frac{1}{2}\bigl(c_1(X)^2\!-\!2c_2(X)\bigr),
\label{an5eq29}
\e
and so on. Now let $E$ be an ample line bundle on $P$. Then by the
Riemann--Roch Theorem \cite[Th.~A.4.1]{Hart} the {\it Hilbert
polynomial\/} $p_X$ of $X$ w.r.t.\ $E$ is
\e
p_X(n)=\chi(X\ot E^n)=\int_P(1+c_1(E))^n\cdot\ch([X])\cdot\td(P).
\label{an5eq30}
\e
Here $\int_P(\cdots)$ means the image of the component in
$H^{2m}(P,\Z)$ in $(\cdots)$ under the natural homomorphism
$H^{2m}(P,\Z)\ra\Z$, and $\td(P)$ is the {\it Todd class} of $P$,
given in terms of Chern classes as on \cite[p.~432]{Hart} by
\e
\td(P)=\ts 1+\frac{1}{2}\,c_1(P)+\frac{1}{12}\bigl(c_1(P)^2
+c_2(P)\bigr)+\frac{1}{24}\,c_1(P)c_2(P)+\cdots.
\label{an5eq31}
\e

For $0\not\cong X\in\A$, write $\dim X$ for the dimension of the
support of $X$ in $P$. This depends only on the class $[X]\in
C(\coh(P))$, and so defines a map $\dim:C(\coh(P))\ra\N$. Note that
for all the weak stability conditions $(\tau,T,\le)$ on $\A$ defined
in \cite[\S 4.4]{Joyc7}, $\tau(\al)$ is a polynomial of degree
$\dim\al$ for all~$\al\in C(\A)$.

\begin{prop} Let\/ $(\tau,T,\le),(\ti\tau,\ti T,\le)$ be
weak stability conditions on $\A$ from {\rm\cite[\S 4.4]{Joyc7}}
and\/ $(\{1,\ldots,n\},\le,\ka)$ be $\A$-data with\/
$\ka(\{1,\ldots,n\})\!=\!\al\in C(\A)$ and\/
$S(\{1,\ldots,n\},\le,\ka,\tau,\ti\tau)\!\ne\!0$.
Then~$\dim\ci\ka\equiv\dim\al$.
\label{an5prop2}
\end{prop}

\begin{proof} Apply \eq{an4eq7} with $(\hat\tau,\hat T,\le)=
(\de,D_m,\le)$ from \cite[Ex.~4.18]{Joyc7}. Then there exist
$m,\psi,\la$ in \eq{an4eq7} with $S\bigl(\psi^{-1}(\{k\}),\le,
\ka\vert_{\psi^{-1}(\{k\})},\tau,\de)\ne 0$ for $k=1,\ldots,m$ and
$S(\{1,\ldots,m\},\le,\la,\de,\ti\tau)\ne 0$. As $(\de,D_m,\le)$
dominates $(\tau,T,\le),(\ti\tau,\ti T,\le)$, equations \eq{an5eq13}
and \eq{an5eq18} then imply that $\de\ci\la\equiv\de(\al)$ and
$\de\ci\ka\vert_{\psi^{-1}(\{k\})} \equiv\de\ci\la(k)$ for
$k=1,\ldots,m$. Together these give $\de\ci\ka\equiv\de(\al)$. But
$\de(\be)=t^{\dim\be}$ for $\be\in C(\A)$,
so~$\dim\ci\ka\equiv\dim\al$.
\end{proof}

Here is a finiteness result for $\A$-data of dimension zero or one.

\begin{prop} Suppose $(\tau,T,\le),(\ti\tau,\ti T,\le)$ are weak
stability conditions on $\A$ from {\rm\cite[Ex.~4.16 or
4.17]{Joyc7}}. Let\/ $\al\in C(\A)$ with\/ $\dim\al=0$ or $1$. Then
there exist at most finitely many sets of\/ $\A$-data
$(\{1,\ldots,n\},\le,\ka)$ with\/ $\ka(\{1,\ldots,n\})=\al$
and\/~$S(\{1,\ldots,n\},\le,\ka,\tau,\ti\tau)\ne 0$.
\label{an5prop3}
\end{prop}

\begin{proof} If $\dim\al=0$ and $n,\ka$ are as above then
$\dim\ci\ka\equiv 0$, so that $\tau\ci\ka\equiv\ti\tau\ci\ka\equiv
t^0$ in all the examples of \cite[\S 4.4]{Joyc7}. Thus the strict
inequalities in Definition \ref{an4def1}(a),(b) cannot hold, and
$S(\{1,\ldots,n\},\le,\ka,\tau,\ti\tau)=0$ unless $n=1$ and
$\ka(1)=\al$. Therefore there is only one set of $\A$-data.

Let $\dim\al=1$, and suppose for simplicity that $P$ is {\it
connected}. If it is not, we may apply the argument below to each
connected component of $P$. Identify $H^{2m}(P,\Z)\cong\Z$. Let
$(\tau,T,\le),(\ti\tau,\ti T,\le)$ be defined using ample line
bundles $E,\ti E$ on $P$ respectively. Then using
\eq{an5eq30}--\eq{an5eq31} we find that
\e
\textstyle \tau(\al)=t+\frac{\text{\begin{footnotesize}
$\ch_m(\al)\!+\!\frac{1}{2}_{\phantom{d\!\!\!}}
c_1(P)\ch_{m-1}(\al)$\end{footnotesize}}}{\text{\begin{footnotesize}
$c_1(E)\ch_{m-1}(\al)$\end{footnotesize}}\phantom{{}^{l^l}}}\,,
\;\>\ti\tau(\al)=t+\frac{\text{\begin{footnotesize}
$\ch_m(\al)\!+\!\frac{1}{2}_{\phantom{d\!\!\!}}
c_1(P)\ch_{m-1}(\al)$\end{footnotesize}}}{\text{\begin{footnotesize}
$c_1(\ti E)\ch_{m-1}(\al)$\end{footnotesize}}\phantom{{}^{l^l}}}\,,
\label{an5eq32}
\e
whichever of \cite[Ex.~4.16]{Joyc7} or \cite[Ex.~4.17]{Joyc7} are
used to define~$(\tau,T,\le),(\ti\tau,\ti T,\le)$.

Suppose $(\{1,\ldots,n\},\le,\ka)$ is $\A$-data with
$S(\{1,\ldots,n\},\le,\ka,\tau,\ti\tau)\!\ne\!0$ and
$\ka(\{1,\ldots,n\})\!=\!\al$. Then $\dim\ci\ka\equiv 1$ by
Proposition \ref{an5prop2}, forcing $\ch_j(\ka(i))=0$ for $j<m-1$
and $c_1(E')\ch_{m-1}(\ka(i))>0$ for all $i$ and ample line bundles
$E'$ on $P$. Since $\sum_{i=1}^n\ch_{m-1}(\ka(i))=\ch_{m-1}(\al)$ it
is not difficult to see there are only finitely many possible
choices for $n$ and $\ch_{m-1}(\ka(i))$, $i=1,\ldots,n$.

Let $k,l$ be as in Theorem \ref{an4thm2}. Using \eq{an5eq32} and
$\ti\tau\ci\ka(k)\!\ge\!\ti\tau(\al)$, $\tau\ci\ka(k)
\!\le\!\tau\ci\ka(i)$, $1\!\le\!c_1(\ti E)\ch_{m-1}(\ka(k))\!\le\!
c_1(\ti E)\ch_{m-1}(\al)$, $1\!\le\!c_1(E)\ch_{m-1}(\ka(k))$ and
$c_1(E)\ch_{m-1}(\ka(i))\!\le\!c_1(E)\ch_{m-1}(\al)$ we deduce the
first inequality of
\e
\begin{split}
c_1(E)\ch_{m-1}(\al)\cdot
\min\bigl(\ch_m(\al)+\ts\frac{1}{2}\,c_1(P)\ch_{m-1}(\al),0\bigr)&\le\\
\ch_m(\ka(i))+\ts\frac{1}{2}\,c_1(P)\ch_{m-1}(\ka(i))&\le\\
c_1(E)\ch_{m-1}(\al)\cdot
\max\bigl(\ch_m(\al)+\ts\frac{1}{2}\,c_1(P)\ch_{m-1}(\al),0\bigr)&
\end{split}
\label{an5eq33}
\e
for all $i=1,\ldots,n$. The second follows in the same way from
$\ti\tau(\al)\ge\ti\tau\ci\ka(l)$ and $\tau\ci\ka(i)\le\tau\ci
\ka(l)$. For each of the finitely many choices for $n$ and
$\ch_{m-1}(\ka(i))$, equation \eq{an5eq33} shows there are only
finitely many possibilities for $\ch_m(\ka(i))$ in
$H^{2m}(P,\Z)\cong\Z$. So there are finitely many choices for
$\ch(\ka(i))$, and thus for $\ka(i)$ as \eq{an5eq28} is injective.
Hence there are only finitely many possibilities for~$n,\ka$.
\end{proof}

Now suppose $P$ is a smooth 2-dimensional $\K$-variety, that is, a
$\K$-{\it surface}. Following \cite[p.~71]{HuLe}, define the {\it
discriminant\/} of $X\in\coh(P)$ to be
\e
\begin{split}
\De(X)&=2\rk(X)c_2(X)-\bigl(\rk(X)-1\bigr)c_1(X)^2\\
&=\ch_1([X])^2-2\ch_0([X])\ch_2([X])
\end{split}
\label{an5eq34}
\e
in $H^4(P,\Z)\cong\Z$, where the second line follows from
\eq{an5eq29}. More generally, if $\al\in C(\coh(P))$ we write
$\De(\al)=\ch_1(\al)^2-2\ch_0(\al)\ch_2(\al)$ in $\Z$. We are
interested in discriminants because of the following important
inequality due to Bogomolov \cite[Th.~3.4.1]{HuLe} when $\cha\K=0$,
and Langer \cite[Th.~3.3]{Lang} when~$\cha\K>0$.

\begin{thm} Let\/ $P$ be a $\K$-surface, $(\tau,T,\le)$ a weak
stability condition from {\rm\cite[Ex.~4.16 or 4.17]{Joyc7}} on
$\A=\coh(P)$ defined using an ample line bundle $E$, and\/ $X\in\A$
be $\tau$-semistable. Then $\De(X)\ge C\rk(X)^2(\rk(X)\!-\!1)^2$,
where $C=0$ if\/ $\cha\K=0$, and\/ $C<0$ depends only on $P,\cha\K$
and\/ $E$ if\/~$\cha\K>0$.
\label{an5thm6}
\end{thm}

Drawing on ideas from Huybrechts and Lehn \cite[\S 4.C]{HuLe}, we
prove a finiteness result for dimension two $\A$-data
$(\{1,\ldots,n\},\le,\ka)$. It is related to Yoshioka \cite[\S
2]{Yosh2}, who proves that for $P$ a smooth projective surface over
$\C$ and $\al\in C(\coh(P))$ fixed the ample cone of $P$ is divided
into finitely many {\it chambers}, and if $(\ga,G_2,\le)$ is
Gieseker stability w.r.t.\ an ample line bundle $E$ then
$\Oss^\al(\ga)$ depends only on which chamber $c_1(E)$ lies in.

\begin{thm} Let\/ $P$ be a smooth projective surface, and\/
$\A\!=\!\coh(P),K(\A)$ as in {\rm\cite[Ex.~9.1]{Joyc5}}. Suppose
$(\tau,T,\le),(\ti\tau,\ti T,\le)$ are permissible weak stability
conditions on $\A$ from {\rm\cite[Ex.~4.16 or 4.17]{Joyc7}}, defined
using ample line bundles $E,\ti E$. Let\/ $\al\in C(\A)$ with\/
$\dim\al=2$. Then there exist at most finitely many sets of\/
$\A$-data $(\{1,\ldots,n\},\le,\ka)$ with\/ $\ka(\{1,\ldots,n\})=
\al$, $S(\{1,\ldots,n\},\le,\ka,\tau,\ti\tau)\ne 0$, and\/
$\De(\ka(i))\ge C\rk(\ka(i))^2(\rk(\ka(i))\!-\!1)^2$ for
$i=1,\ldots,n$ and\/ $C$ as above, and if there are any such\/
$n,\ka$ then~$\De(\al)\ge C\rk(\al)^2(\rk(\al)\!-\!1)^2$.
\label{an5thm7}
\end{thm}

\begin{proof} Identifying $H^4(P,\Z)\cong\Z$, for $\be\in C(\A)$
with $\dim\be=2$ define
\e
\mu(\be)=c_1(\be)c_1(E)/\rk(\be) \quad\text{and}\quad
\ti\mu(\be)=c_1(\be)c_1(\ti E)/\rk(\be)\quad\text{in $\Q$.}
\label{an5eq35}
\e
For $s\in[0,1]$ define $\mu_s(\be)=(1-s)\mu(\be)+s\ti\mu(\be)$ in
$\R$. By computing Hilbert polynomials w.r.t.\ $E,\ti E$ using
\eq{an5eq29}--\eq{an5eq31} and referring to \cite[Ex.s
4.16--4.17]{Joyc7}, we can show that if $\be,\ga\in C(\A)$ with
$\dim\be=2=\dim\ga$ then $\mu(\be)<\mu(\ga)$ implies
$\tau(\be)<\tau(\ga)$ and $\ti\mu(\be)<\ti\mu(\ga)$ implies
$\ti\tau(\be)<\ti\tau(\ga)$. So we shall use $\mu,\ti\mu$ as
substitutes for $\tau,\ti\tau$ in the proof, and $\mu_s$ for
$s\in[0,1]$ interpolate between them.

Let $n,\ka$ be as in the theorem. Then $\dim\ci\ka\equiv 2$ by
Proposition \ref{an5prop2}. If $n=1$ the only possibility is
$\ka(1)=\al$, so suppose $n>1$. For $s\in[0,1]$ and $1\le i<n$,
consider the three (not exhaustive) alternatives:
\begin{itemize}
\setlength{\itemsep}{0pt}
\setlength{\parsep}{0pt}
\item[(a)] $\mu_s\ci\ka(\{1,\ldots,i\})>\mu_s(\al)$ and
$\tau\ci\ka(i)\le\tau\ci\ka(i+1)$;
\item[(b)] $\mu_s\ci\ka(\{1,\ldots,i\})<\mu_s(\al)$ and
$\tau\ci\ka(i)>\tau\ci\ka(i+1)$; or
\item[(c)] $\mu_s\ci\ka(\{1,\ldots,i\})=\mu_s(\al)$.
\end{itemize}
The point of this is that as $\ti\mu(\be)<\ti\mu(\ga)$ implies
$\ti\tau(\be)<\ti\tau(\ga)$ from above, when $s=1$ parts (a),(b)
above imply Definition \ref{an4def1}(a),(b), and Definition
\ref{an4def1}(a),(b) imply (a),(b) or (c). Thus, $S(\{1,\ldots,n\},
\le,\ka,\tau,\ti\tau)\!\ne\!0$ and Definition \ref{an4def1} imply
that one of (a)--(c) above hold for all $1\le i<n$ when~$s=1$.

Similarly, when $s=0$ parts (a)--(c) above are related in the same
way to Definition \ref{an4def1}(a),(b) with $\tau$ in place of
$\ti\tau$. But $S(\{1,\ldots,n\},\le,\ka,\tau,\tau)=0$ by
\eq{an4eq6} as $n>1$. Therefore Definition \ref{an4def1} implies
that there exists $1\le i<n$ for which neither (a) nor (b) above
holds when $s=0$. As one of (a)--(c) holds for this $i$ when $s=1$,
the sign of $\mu_s\ci\ka(\{1,\ldots,i\})-\mu_s(\al)$ must change
over $[0,1]$, so there exist $s\in[0,1]$ and $1\le i<n$ such that
(c) holds. Choose such values $s_0,i_0$ of $s,i$ with $s_0$ as large
as possible. Then since one of (a)--(c) hold for all $1\le i<n$ at
$s=1$, it is easy to see that one of (a)--(c) hold for all $1\le
i<n$ at $s=s_0$, and (c) holds for $s=s_0$, $i=i_0$. Also
$s_0\in\Q$, as $\mu,\ti\mu$ take values in~$\Q$.

Define $\al',\al''\in C(\A)$ and $\A$-data
$(\{1,\ldots,n'\},\le,\ka')$, $(\{1,\ldots,n''\},\le,\ka'')$ by
$\al'=\ka(\{1,\ldots,i_0\})$, $\al''= \ka(\{i_0+1,\ldots,n\})$,
$n'=i_0$, $n''=n-i_0$, $\ka'(i)=\ka(i)$, and $\ka''(i)=\ka(i+i_0)$.
Then $\ka'(\{1,\ldots,n'\})=\al'$, $\ka''(\{1,\ldots,n''\})=\al''$.
As one of (a)--(c) hold for all $1\le i<n$ at $s=s_0$, and (c) holds
for $s=s_0$, $i=i_0$, it is not difficult to see that the analogues
of one of (a)--(c) for $\al'$ and $(\{1,\ldots,n'\},\le,\ka')$ hold
for $s=s_0$ and all $1\le i<n'$, and the analogues of one of
(a)--(c) for $\al''$ and $(\{1,\ldots,n''\},\le,\ka'')$ hold for
$s=s_0$ and all~$1\le i<n''$.

Set $\xi=\rk(\al')c_1(\al'')-\rk(\al'')c_1(\al')$ in $H^2(P,\Z)$.
Then \eq{an5eq29}, \eq{an5eq34} yield
\e
\frac{1}{\rk(\al)}\,\De(\al)=\frac{1}{\rk(\al')}\,\De(\al')+
\frac{1}{\rk(\al'')}\,\De(\al'')-\frac{\xi^2}{\rk(\al)
\rk(\al')\rk(\al'')}\,.
\label{an5eq36}
\e
As $s_0\in\Q\cap[0,1]$ we can write $s_0=p/(p+q)$ for integers
$p,q\ge 0$ with $p+q>0$. Since (c) holds for $s=s_0$, $i=i_0$ we
find that
\e
\bigl(q\,c_1(E)+p\,c_1(\ti E)\bigr)\xi=0.
\label{an5eq37}
\e
But $E^q\ot\ti E^p$ is ample, so the {\it Hodge Index Theorem}
\cite[Th.~V.1.9]{Hart} implies that $\xi^2\le 0$, and thus
\eq{an5eq36} gives
\e
\frac{1}{\rk(\al)}\,\De(\al)\ge \frac{1}{\rk(\al')}\,\De(\al')+
\frac{1}{\rk(\al'')}\,\De(\al'').
\label{an5eq38}
\e

Suppose now that we know that
\e
\De(\al')\ge C\rk(\al')^2(\rk(\al')\!-\!1)^2 \quad\text{and}\quad
\De(\al'')\ge C\rk(\al'')^2(\rk(\al'')\!-\!1)^2.
\label{an5eq39}
\e
Then from \eq{an5eq36} and $\rk(\al'),\rk(\al'')\le\rk(\al)$ we see
that
\begin{equation*}
2C\rk(\al)^4(\rk(\al)\!-\!1)^2-\rk(\al)^2\De(\al)\le\xi^2\le 0,
\end{equation*}
so there are only finitely many possibilities for the integer
$\xi^2$. Combining this with \eq{an5eq37} we see as in the proof of
\cite[Lem.~4.C.2]{HuLe} that $\xi$ lies in a bounded, and hence
finite, subset of~$H^2(P,\Z)$.

Since $1\le\rk(\al'),\rk(\al'')\le\rk(\al)$ there are only finitely
many choices for $\ch_0(\al')=\rk(\al')$ and $\ch_0(\al'')=\rk
(\al'')$. But as $c_1(\al')+c_1(\al'')=c_1(\al)$, $\rk(\al'),
\rk(\al'')$ and $\xi$ determine $\ch_1(\al')=c_1(\al')$ and
$\ch_1(\al'')=c_1(\al'')$, so there are only finitely many
possibilities for these. From \eq{an5eq34} and \eq{an5eq39} we see
that
\begin{align*}
\ch_2(\al')&\le \ch_1(\al')^2/2\ch_0(\al')-\ha
C\ch_0(\al')\bigl(\ch_0(\al')-1\bigr)^2,\\
\ch_2(\al'')&\le\ch_1(\al'')^2/2\ch_0(\al'')-\ha
C\ch_0(\al'')\bigl(\ch_0(\al'')-1\bigr)^2,
\end{align*}
so as $\ch_2(\al')+\ch_2(\al'')=\ch_2(\al)$ we see that once
$\ch_i(\al'),\ch_i(\al'')$ are fixed for $i=0,1$ there are only
finitely many choices for $\ch_2(\al'),\ch_2(\al'')$ in $\Z$. Thus,
given $\al,C$ there are only finitely many possibilities for
$\ch(\al'),\ch(\al'')$, and hence for $\al',\al''$ as \eq{an5eq28}
is injective.

We can now prove the following inductive hypothesis, for given~$l\ge
1$:
\begin{itemize}
\item[$(*\kern .05em{}_l)$] Suppose $\al\in C(\A)$ with $\dim\al=2$ and
$\rk(\al)\le l$, and $(\{1,\ldots,n\},\le,\ka)$ is $\A$-data with
$\ka(\{1,\ldots,n\})=\al$ and $\De(\ka(i))\ge C\rk(\ka(i))^2
(\rk(\ka(i))\!-\!1)^2$ for $i=1,\ldots,n$, such that for some
$s\in[0,1]$ and all $1\le i<n$ one of (a)--(c) above holds. Then
$\De(\al)\ge C\rk(\al)^2(\rk(\al)\!-\!1)^2$. Moreover, for fixed
$\al$ there are only finitely many possibilities for such~$n,\ka$.
\end{itemize}
When $l=1$ this is trivial as the only possibility is $n=1$ and
$\ka(1)=\al$. Suppose by induction that $(*\kern .05em{}_l)$ holds
for $l\ge 1$, and let $\al\in C(\A)$ with $\rk(\al)=l+1$, and
$n,\ka$ be as in $(*\kern .05em{}_{l+1})$. Construct $s_0\in[0,s]$
and $\al',n',\ka'$, $\al'',n'',\ka''$ as in the first part of the
proof. Then $\rk(\al'),\rk(\al'')\le l$, so $\al',n',\ka'$ and
$\al'',n'',\ka''$ satisfy the conditions of $(*\kern .05em{}_l)$,
which implies \eq{an5eq39} holds. Combining this with \eq{an5eq38}
and $\rk(\al)=\rk(\al')+ \rk(\al'')$ gives $\De(\al)\ge
C\rk(\al)^2(\rk(\al) \!-\!1)^2$, as we want.

By the first part there are only finitely many possibilities for
$\al',\al''$. For each of these, $(*\kern .05em{}_l)$ shows there
are only finitely many choices for $n',\ka'$ and $n'',\ka''$. But
these determine $n,\ka$, so there are only finitely many
possibilities for $n,\ka$, completing the inductive step. Thus
$(*\kern .05em{}_l)$ holds for all $l\ge 1$. The theorem now follows
from the first part of the proof, which showed that for $n,\ka$ as
in the theorem the hypotheses of $(*\kern .05em{}_l)$ hold
with~$s=1$.
\end{proof}

Definition \ref{an5def1}, Proposition \ref{an5prop3} and Theorems
\ref{an5thm6}--\ref{an5thm7} show the change from $(\tau,T,\le)$ to
$(\ti\tau,\ti T,\le)$ is {\it globally finite}, so the change from
$(\ti\tau,\ti T,\le)$ to $(\tau,T,\le)$ is too by symmetry. Theorem
\ref{an5thm3} now follows from Corollary~\ref{an5cor3}.

\section{Invariants}
\label{an6}

Given a permissible weak stability condition $(\tau,T,\le)$ on $\A$,
we now consider how best to define {\it systems of invariants}
$\Iss(\cdots)$ of $\A$ and $(\tau,T,\le)$ which `count'
$\tau$-semistable objects in class $\al\in K(\A)$, or more generally
`count' $(I,\pr,\ka)$-configurations $(\si,\io,\pi)$ with
$\si(\{i\})$ $\tau$-semistable for all $i\in I$. Obviously there are
many ways of doing this, so we need to decide what are the most
interesting, or useful, ways to define invariants. For instance, we
could ask that the invariants we choose can be calculated in
examples, or are important in other areas of mathematics.

The criterion we shall use to select interesting invariants is that
{\it they should satisfy natural identities}, and the more
identities the better. Of course, this is not unrelated to other
criteria; for instance, such identities are powerful tools for
calculating the invariants in examples, as we shall see, and may be
a reason for the invariants to be important in other areas.

We shall divide the identities we are interested in into {\it
additive identities}, for which the $\Iss(\cdots)$ should take
values in an abelian group or $\Q$-vector space, and {\it
multiplicative identities}, for which the $\Iss(\cdots)$ should take
values in a ring or $\Q$-algebra. Here is a very general way of
defining invariants satisfying useful additive identities.

\begin{dfn} Let $\K,\A,K(\A)$ satisfy Assumption \ref{an3ass}, and
suppose $\La$ is a $\Q$-vector space and $\rho:\SFa(\fObj_\A)\ra
\La$ a $\Q$-linear map. Let $(\tau,T,\le)$ be a permissible weak
stability condition on $\A$. For all $\A$-data $(I,\pr,\ka)$, define
{\it invariants} $\Iss(I,\pr,\ka,\tau)$ in $\La$
by~$\Iss(I,\pr,\ka,\tau)=
\rho\ci\bs\si(I)_*\,\bdss(I,\pr,\ka,\tau)$.
\label{an6def1}
\end{dfn}

\begin{rem} Here are some ways we might choose the linear map~$\rho$:
\begin{itemize}
\setlength{\itemsep}{0pt}
\setlength{\parsep}{0pt}
\item[(a)] When $\cha\K=0$ we could take $\rho=\rho'\ci
\pi^\stk_{\fObj_\A}$, for some linear $\rho':\CF(\fObj_\A)\ra\La$.
This would have the advantage of making the meaning of the
invariants -- what they are `counting' -- clearer. For example,
$\Ist^\al(\tau)$ in \eq{an6eq2} below would measure the moduli space
$\Ost^\al(\tau)$ of $\tau$-stable elements in class $\al$, rather
than the more obscure stack function $\bdst^\al(\tau)$ of \cite[\S
8]{Joyc7}. In the same way, we could take $\rho=\rho'\ci\bar\Pi^{
\Th,\Om}_{\fObj_\A}$ for some linear $\rho':\oSFa(\fObj_\A,\Th,\Om)
\ra\La$. Then as in the discussion before \cite[Th.~8.7]{Joyc7}, we
could interpret $\Ist^\al(\tau)$ as counting `virtual
$\tau$-stables' in class $\al$, etc.
\item[(b)] As in \S\ref{an61} we could take $\rho=\rho'\ci\Pi_*$,
where $\Pi:\fObj_\A\ra\Spec\K$ is the projection, $\Pi_*$ maps
$\SFa(\fObj_\A)\ra\uSF(\Spec\K)$ (not $\SF(\Spec\K)$ as $\Pi$ is not
representable), and $\rho':\uSF(\Spec\K)\ra\La$ is linear. Then
$\Iss(I,\pr,\ka,\tau)=\rho'([\Mss(I,\pr,\ka,\tau)])$, so
$\Iss(I,\pr,\ka,\tau)$ depends only on the moduli space
$\Mss(I,\pr,\ka,\tau)$ and not on its projection to~$\fObj_\A$.
\item[(c)] We could take $\La$ to be a $\Q$-algebra and $\rho:\SFa
(\fObj_\A)\ra\La$ to be an {\it algebra morphism}, such as those
constructed in \cite[\S 6]{Joyc6}. This would imply {\it
multiplicative identities} on the invariants in $\La$. In the same
way, we could restrict to invariants coming from $\bLp_\tau$ and
take $\rho:\SFai(\fObj_\A)\ra\La$ to be a {\it Lie algebra
morphism}, as in~\cite[\S 6.6]{Joyc6}.
\end{itemize}
\label{an6rem1}
\end{rem}

Applying $\rho\ci\bs\si(K)_*$ to \eq{an5eq5} in Theorem
\ref{an5thm1} gives:

\begin{thm} Let Assumption \ref{an3ass} hold, $(\tau,T,\le),
(\ti\tau,\ti T,\le)$ be permissible weak stability conditions on
$\A$, and\/ $\Iss(*)$ be as in Definition \ref{an6def1}. Suppose the
change from $(\tau,T,\le)$ to $(\ti\tau,\ti T,\le)$ is globally
finite, and there exists a weak stability condition $(\hat\tau,\hat
T,\le)$ on $\A$ dominating $(\tau,T,\le),(\ti\tau,\ti T,\le)$ with
the change from $(\hat\tau,\hat T,\le)$ to $(\ti\tau,\ti T,\le)$
locally finite. Then for all\/ $\A$-data $(K,\tl,\mu)$ the following
holds in $\La$, with only finitely many nonzero terms:
\e
\begin{aligned}
\sum_{\substack{\text{iso.}\\ \text{classes}\\
\text{of finite }\\ \text{sets $I$}}} \frac{1}{\md{I}!}\cdot
\sum_{\substack{
\text{$\pr,\ka,\phi$: $(I,\pr,\ka)$ is $\A$-data,}\\
\text{$(I,\pr,K,\phi)$ is dominant,}\\
\text{$\tl=\P(I,\pr,K,\phi)$,}\\
\text{$\ka(\phi^{-1}(k))=\mu(k)$ for $k\in K$}}}\,\,\,
\begin{aligned}[t]
T(I,\pr,\ka,K,\phi,\tau,\ti\tau)&\cdot\\
\Iss(I,\pr,\ka,\tau)&=\\
\Iss(K,\tl,\mu&,\ti\tau).
\end{aligned}
\end{aligned}
\label{an6eq1}
\e
\label{an6thm1}
\end{thm}

Knowing $\Iss(I,\pr,\ka,\tau)$ for all $(I,\pr,\ka)$ is equivalent
to knowing the restriction $\rho:\bHp_\tau\ra\La$. We know from
\S\ref{an5} that under mild conditions, if $(\tau,T,\le)$ and
$(\ti\tau,\ti T,\le)$ are permissible weak stability conditions on
$\A$ then $\bHp_\tau=\bHp_{\ti\tau}$, as in \eq{an5eq11}, so knowing
$\Iss(I,\pr,\ka,\tau)$ for all $(I,\pr,\ka)$ is equivalent to
knowing $\Iss(I,\pr,\ka,\ti\tau)$ for all $(I,\pr,\ka)$. Theorem
\ref{an6thm1} makes this equivalence explicit.

In the same way, subsystems of the $\Iss(*,\tau)$ are equivalent to
knowing the restrictions of $\rho$ to $\bHt_\tau,\bLp_\tau$ and
$\bLt_\tau$. For example, knowing $\Iss(I,\pr,\ka,\tau)$ for all
$(I,\pr,\ka)$ with $\pr$ a {\it total order} is equivalent to
knowing $\rho:\bHt_\tau\ra\La$. If $(K,\tl)$ is a total order then
all $(I,\pr)$ occurring in \eq{an6eq1} are total orders,
corresponding to the fact that $\bHt_\tau=\bHt_{\ti\tau}$. Thus, the
$\Iss(I,\pr,\ka,\tau)$ with $(I,\pr)$ a total order form a {\it
closed subsystem} of the system $\Iss(*,\tau)$ under change of
stability condition.

We can also define other families of invariants $\Isi,\Ist,\Issb,
\Isib,\Istb(I,\pr,\ka,\tau)$ in $\La$ by applying
$\rho\ci\bs\si(I)_*$ to the stack functions
$\bdsi,\bdst,\bdssb,\bdsib,\bdstb(I,\pr,\ka, \tau)$ of \cite[\S
8]{Joyc7}, and $\Iss^\al,\Isi^\al,\Ist^\al,J^\al (\tau)$ by applying
$\rho$ to $\bdss^\al,\bdsi^\al,\bdst^\al, \bep^\al(\tau)$. Then the
identities of \cite[\S 8]{Joyc7} yield many additive identities
between these families of invariants. For instance, the stack
function analogue of \eq{an3eq7} implies that
\e
\Ist^\al(\tau)\!=\!\!\!\!
\sum_{\substack{\text{iso. classes}\\ \text{of finite sets $I$}}}
\frac{1}{\md{I}!}\,\cdot\! \sum_{\begin{subarray}{l}
\text{$\pr,\ka$: $(I,\pr,\ka)$ is $\A$-data,}\\
\text{$\ka(I)=\al$, $\tau\ci\ka\equiv\tau(\al)$} \end{subarray}
\!\!\!\!\!\!\!\!\!\!\!\!\!\!\!\!\!\!\!\!\!\!\!\!\!\!\!\!\!\!
\!\!\!\!\!\!\!\!\! } \Iss(I,\pr,\ka,\tau)\cdot\!\!
\sum_{\text{p.o.s $\ls$ on $I$ dominating $\pr$}
\!\!\!\!\!\!\!\!\!\!\!\!\!\!\!\!\!\!\!\!\!\!\!\!\!\!\!\!\!\!\!\!
\!\!\!\!\!\!\!\!\!\!\!\!\!\!\!\!\!\!} n(I,\pr,\ls)\,N(I,\ls).
\label{an6eq2}
\e

Here $\Ist^\al(\tau)$ is an invariant which `counts' $\tau$-stable
objects in class $\al\in K(\A)$ (though see Remark \ref{an6rem1}(a)
on this point). One overall moral is that all the invariants can be
written in terms of the $\Iss(*,\tau)$, so we use $\Iss(*,\tau)$ for
preference. The other invariants will be useful for some purposes
though; for instance, from \cite[Def.~8.9]{Joyc7} we see that
$\rho:\bLp_\tau\ra\La$ is determined by the $\Isib(I,\pr,\ka,\tau)$
or $\Istb(I,\pr,\ka,\tau)$ for {\it connected\/} $(I,\pr)$, so these
$\Isib(*,\tau)$ or $\Istb(*,\tau)$ form a closed subsystem under
change of stability condition.

In the rest of the section we study examples in which $\La$ is a
$\Q$-algebra and the invariants satisfy additional multiplicative
identities. The obvious way to do this is to arrange that
$\rho:\SFa(\fObj_\A)\ra\La$ is an algebra morphism, or
$\rho:\SFai(\fObj_\A)\ra\La$ a Lie algebra morphism, as in the
constructions of \cite[\S 6]{Joyc6}. But in \cite[\S 5.1]{Joyc6} we
also define multilinear operations $P_\sIp$ on $\SFa(\fObj_\A)$, and
in some cases we can arrange for these to commute with operations
$P_\sIp$ on~$\La$.

Having multiplicative identities means that the full system of
invariants $\Iss(I,\pr,\ka,\tau)$ is wholly determined by a smaller
generating subset of invariants. For instance, if $\rho$ is an
algebra morphism then $\rho$ on $\bHp_\tau$ or $\bHt_\tau$ is
determined by its value on a set of generators for $\bHp_\tau$ or
$\bHt_\tau$, such as the $\bdss^\al(\tau)$ or $\bep^\al(\tau)$ for
$\bHt_\tau$. This reduces the amount of data needed to specify the
$\Iss(*,\tau)$, and so simplifies the problem of computing
invariants in examples. So where we can we will focus on the
invariants $\Iss^\al(\tau),J^\al(\tau)$ associated
to~$\bdss^\al(\tau),\bep^\al(\tau)$.

\subsection{Multiplicative relations from disconnected $(I,\pr)$}
\label{an61}

As in \cite[\S 7.2]{Joyc7}, if $(I,\pr)$ is a finite poset let
$\approx$ be the equivalence relation on $I$ generated by $i\approx
j$ if $i\pr j$ or $j\pr i$, and define the {\it connected
components} of $(I,\pr)$ to be the $\approx$-equivalence classes. We
shall now study invariants $\Iss(I,\pr,\ka,\tau)$ that are
multiplicative over disjoint unions of connected components.

\begin{dfn} Let Assumption \ref{an3ass} hold and $(\tau,T,\le)$ be
a permissible weak stability condition on $\A$. Suppose $\La$ is a
commutative $\Q$-algebra and $\rho':\uSF(\Spec\K)\ra\La$ an algebra
morphism. As in Remark \ref{an6rem1}(b), define invariants
$\Iss(I,\pr,\ka,\tau)$ in $\La$ for all $\A$-data $(I,\pr,\ka)$ by
\e
\Iss(I,\pr,\ka,\tau)=\rho'\ci\Pi_*\,\bdss(I,\pr,\ka,\tau)=
\rho'\bigl([\Mss(I,\pr,\ka,\tau)]\bigr),
\label{an6eq3}
\e
where $\Pi:\fObj_\A\ra\Spec\K$ is the projection
and~$\Pi_*:\SFa(\fObj_\A)\ra\uSF(\Spec\K)$.

There are many ways of constructing such $\rho'$. For example, we
can take $\La=\uSF(\Spec\K)$ and $\rho'=\id_{\sst\La}$. Or for any
$\Up,\La$ satisfying Assumption \ref{an2ass1}, such as Example
\ref{an2ex1} or others in \cite[\S 4.1]{Joyc4}, we can take $\rho'$
to be the morphism $\Up'$ of Theorem \ref{an2thm4}. Also, using the
notation of \cite[\S 4.4]{Joyc3}, if $\cha\K=0$ and $w$ is any {\it
allowable multiplicative weight function} on affine algebraic
$\K$-groups taking values in $\Q$, then $\La=\Q$ and
$\rho':[\fR]\mapsto\chi_w(\fR(\K))$ defines an algebra morphism.
Examples of such $w$ are given in~\cite[Prop.~4.16]{Joyc3}.
\label{an6def2}
\end{dfn}

These invariants satisfy a multiplicative identity. Note that
\eq{an6eq4} holds without any additional conditions on $\A$, in
contrast to the identities of \S\ref{an62}--\S\ref{an65}, which
require assumptions on the groups $\Ext^i(X,Y)$ for~$X,Y\in\A$.

\begin{prop} Let Assumption \ref{an3ass} hold and\/ $(\tau,T,\le)$
be a permissible weak stability condition on $\A$, and define
invariants $\Iss(I,\pr,\ka,\tau)$ as in Definition \ref{an6def2}.
Suppose $(J,\ls,\la)$, $(K,\tl,\mu)$ are $\A$-data with\/ $J\cap
K=\emptyset$. Define $\A$-data $(I,\pr,\ka)$ by $I=J\amalg K$,
$\ka\vert_J=\la$, $\ka\vert_K=\mu$, and\/ $i\pr i'$ for $i,i'\in I$
if either $i,i'\in J$ and\/ $i\ls i'$, or $i,i'\in K$ and\/ $i\tl
i'$. Then
\e
\Iss(I,\pr,\ka,\tau)=\Iss(J,\ls,\la,\tau)\cdot\Iss(K,\tl,\mu,\tau).
\label{an6eq4}
\e
Equation \eq{an6eq4} also holds with $\Iss$ replaced by
$\Isi,\Ist,\Issb,\Isib$ or $\Istb$ throughout.
\label{an6prop1}
\end{prop}

\begin{proof} As in \cite[\S 7.4]{Joyc5} we can prove we have a
1-isomorphism:
\begin{equation*}
S(I,\pr,J)\t S(I,\pr,K):\fM(I,\pr,\ka)_\A\longra\fM(J,\ls,\la)_\A
\t\fM(K,\tl,\mu)_\A.
\end{equation*}
This implies $\bigl(S(I,\pr,J)\t S(I,\pr,K)\bigr){}_*\bigl(
\Mss(I,\pr,\ka, \tau)_\A\bigr)=\Mss(J,\ls,\la,\tau)_\A\t\ab
\Mss(K,\tl,\mu,\tau)_\A$, so $[\Mss(I,\pr,\ka,\tau)_\A]\!=\!
[\Mss(J,\ls,\la,\tau)_\A]\cdot[\Mss(K,\tl,\mu,\tau)_\A]$ in the
algebra $\uSF(\Spec\K)$. Equation \eq{an6eq4} then follows from
\eq{an6eq3} and $\rho'$ an algebra morphism. The analogues for
$\Isi,\ldots,\Istb$ can be deduced from \eq{an6eq4} and the additive
identities relating $\Iss(*,\tau)$ and $\Isi,\ldots,\Istb(*,\tau)$
that follow from the additive identities on
$\bdss,\ldots,\bdstb(*,\tau)$ in~\cite[\S 8]{Joyc7}.
\end{proof}

We see from \eq{an6eq4} that all $\Iss(I,\pr,\ka,\tau)$ are
determined by the subset of $\Iss(I,\pr,\ka,\tau)$ with $(I,\pr)$
{\it connected}. Equivalently, they are determined by the subset of
$\Isib(I,\pr,\ka,\tau)$ with $(I,\pr)$ connected. But by
\cite[Def.~8.9]{Joyc7}, the Lie algebra $\bLp_\tau$ is spanned by
$\bdsib(I,\pr,\ka,\tau)$ with $(I,\pr)$ connected. Therefore
\eq{an6eq4} implies that $\rho=\rho'\ci\Pi_*:\bHp_\tau\ra\La$ is
determined by~$\rho:\bLp_\tau\ra\La$.

Here is a different way to define invariants $\Iss(*,\tau)$
satisfying~\eq{an6eq4}.

\begin{dfn} Let Assumption \ref{an3ass} hold, with $\K$ of
characteristic zero. Then each $X\in\A$ may be written $X\cong
X_1\op\cdots\op X_n$, for $X_1,\ldots,X_n$ {\it indecomposable} in
$\A$, and unique up to order and isomorphism. Identifying $X$ with
$X_1\op\cdots\op X_n$, define an algebraic $\K$-subgroup $T_X$ of
$\Aut(X)$ by
\begin{equation*}
T_X=\bigl\{\la_1\id_{X_1}+\la_2\id_{X_2}+\cdots+\la_n\id_{X_n}:
\la_1,\ldots,\la_n\in\K^\t\bigr\}\cong(\K^\t)^n.
\end{equation*}
Then $T_X$ is a maximal torus in $\Aut(X)$, so $\Aut(X)/T_X$ is a
quasiprojective $\K$-variety, and its Euler characteristic
$\chi(\Aut(X)/T_X)$ exists as in~\cite[\S 3.3]{Joyc3}.

Now $\End(X)$ is a finite-dimensional $\K$-algebra. Applying the
{\it Wedderburn structure theorem} \cite[\S 1.3]{Bens} gives an
algebra isomorphism $\End(X)/J(\End(X))\cong\ts
\bigop_{i=1}^r\End(\K^{n_i})$, where $J(\End(X))$ is the {\it
Jacobson radical}, a nilpotent, two-sided ideal in $\End(X)$, and
$n=n_1+\cdots+n_r$. This implies an isomorphism of varieties
$\Aut(X)/T_X\cong\K^l\t\prod_{i=1}^n\GL(n_i,\K)/(\K^\t)^{n_i}$,
where $l=\dim J(\End(X))$ and $(\K^\t)^{n_i}\subseteq\GL(n_i,\K)$ is
the subgroup of diagonal matrices. Elementary calculation then shows
that $\chi\bigl(\Aut(X)/T_X\bigr)=\prod_{i=1}^rn_i!$, which is {\it
nonzero}.

Define $\La=\Q[t]$, the $\Q$-algebra of rational polynomials in $t$.
Define a {\it weight function} $w:\fObj_\A(\K)\ra\La$ by
$w([X])=\chi\bigl(\Aut(X)/T_X\bigr){}^{-1}t^n$, for $n,T_X$ as
above. This is well-defined as $\chi(\Aut(X)/T_X)\ne 0$ from above,
and is {\it locally constructible} on $\fObj_\A$. One can prove that
this weight function has the following multiplicative property. Let
$\bu$ be the partial order on $\{1,2\}$ with $i\bu j$ only if\/
$i=j$, and $P_{\sst(\{1,2\},\bu)}$ be the bilinear operation on
$\CF(\fObj_\A)$ studied in \cite[\S 4.8]{Joyc6}. Then for all
$f,g\in\CF(\fObj_\A)$ we have
\e
\chi^\na\bigl(\fObj_\A,w\cdot P_{\sst(\{1,2\},\bu)}(f,g)\bigr)=
\chi^\na\bigl(\fObj_\A,w\cdot f\bigr)
\cdot\chi^\na\bigl(\fObj_\A,w\cdot g\bigr) \quad\text{in $\La$,}
\label{an6eq5}
\e
using the {\it na\"\i ve weighted Euler characteristic} of \cite[\S
4.1]{Joyc3}. The proof is elementary: we calculate the multiples of
$f([X])g([Y])$ contributing to each side at $[X\op Y]$ and show they
are the same.

Now let $(\tau,T,\le)$ be a permissible stability condition on $\A$
and define
\e
\begin{split}
\Iss(I,\pr,\ka,\tau)&=\chi^\na\bigl(\fObj_\A,w\cdot
\CF^\stk(\bs\si(I))\dss(I,\pr,\ka,\tau)\bigr)\\
&=\chi^\na\bigl(\fObj_\A,w\cdot\pi^\stk_{\fObj_\A}
\ci\bs\si(I)_*\,\bdss(I,\pr,\ka,\tau)\bigr)
\end{split}
\label{an6eq6}
\e
for all $(I,\pr,\ka)$. Then using $\CF^\stk\bigl(\bs\si(I)\bigr)
\dss(I,\pr,\ka,\tau)=P_\sIp\bigl(\dss^{\ka(i)}(\tau):i\in I\bigr)$,
\eq{an6eq5} and \cite[Th.~4.22]{Joyc6} we find these invariants
satisfy~\eq{an6eq4}.
\label{an6def3}
\end{dfn}

The $\Iss(*,\tau)$ of \eq{an6eq6} do not come from Definition
\ref{an6def2}. Remark \ref{an6rem1}(a) applies to them, so that
invariants such as $\Ist^\al(\tau)$ in \eq{an6eq2} have a clear
interpretation. In fact, combining \eq{an3eq7}, \eq{an6eq2} and
\eq{an6eq6} yields $\Ist^\al(\tau)=\chi^\na(\fObj_\A,
w\cdot\dst^\al(\tau))$. But any $\tau$-stable $X$ has
$\Aut(X)\cong\K^\t$, so $w\cong t$ on $\Ost^\al(\tau)$, giving
$\Ist^\al(\tau)=t\,\chi^\na(\Ost^\al(\tau))$, arguably the simplest,
most obvious way to `count' $\tau$-stables.

\subsection{When $\Ext^i(X,Y)=0$ for all $X,Y\in\A$ and $i>1$}
\label{an62}

Recall that a $\K$-linear abelian category $\A$ is called of {\it
finite type} if $\Ext^i(X,Y)$ is a finite-dimensional $\K$-vector
space for all $X,Y\in\A$ and $i\ge 0$, and $\Ext^i(X,Y)=0$ for $i\gg
0$. Then there is a unique biadditive map $\chi:K_0(\A)\t
K_0(\A)\ra\Z$ on the Grothendieck group $K_0(\A)$ known as the {\it
Euler form}, satisfying
\e
\chi\bigl([X],[Y]\bigr)=\ts\sum_{i\ge 0}(-1)^i\dim_\K\Ext^i(X,Y)
\quad\text{for all $X,Y\in\A$.}
\label{an6eq7}
\e
We shall suppose $K(\A)$ in Assumption \ref{an3ass} is chosen such
that $\chi$ factors through the projection $K_0(\A)\ra K(\A)$, and
so descends to $\chi:K(\A)\t K(\A)\ra\Z$. This holds for nearly all
the examples of \cite[\S 9--\S 10]{Joyc5}.

Now assume $\Ext^i(X,Y)=0$ for all $X,Y\in\A$ and $i>1$. Then
\eq{an6eq7} becomes
\e
\dim_\K\Hom(X,Y)-\dim_\K\Ext^1(X,Y)=\chi\bigl([X],[Y]\bigr)
\quad\text{for all $X,Y\in\A$.}
\label{an6eq8}
\e
This happens for $\A=\coh(P)$ in \cite[Ex.~9.1]{Joyc5} with $P$ a
smooth projective curve, and for $\A=\modKQ$ in
\cite[Ex.~10.5]{Joyc5}.

Supposing \eq{an6eq8} and Assumption \ref{an2ass1}, \cite[\S
6.2]{Joyc6} defined an {\it algebra morphism}
\e
\Phi^{\sst\La}\ci\Pi^{\Up,\La}_{\fObj_\A}:\SFa(\fObj_\A)\ra
A(\A,\La,\chi),
\label{an6eq9}
\e
where $A(\A,\La,\chi)$ is an explicit algebra depending only on
$K(\A),C(\A),\chi$ and $\La$. In fact we defined $\Phi^{\sst\La}$ on
the algebra $\uSF(\fObj_\A,\Up,\La)$, but composing with the
projection from $\SFa(\fObj_\A)$ gives an algebra morphism from
$\SFa(\fObj_\A)$. Furthermore, we showed \eq{an6eq9} intertwines the
multilinear operations $P_\sIp$ on $\SFa(\fObj_\A)$ of \cite[\S
5.1]{Joyc6} with operations $P_\sIp$ on~$A(\A,\La,\chi)$.

We shall define families of invariants $\Iss^\al(\tau)^{\sst\La},
J^\al(\tau)^{\sst\La^\ci,\sst\Om}$ and
$\Iss(I,\pr,\ka,\tau)^{\sst\La},\ab
\Jsib(I,\pr,\ka,\tau)^{\sst\La^\ci,\sst\Om}$ taking values in the
algebras $\La,\La^\ci,\Om$ of Assumption \ref{an2ass1}, which will
satisfy multiplicative identities when \eq{an6eq9} is an algebra
morphism. For later use we define them without extra assumptions
on~$\Ext^i(X,Y)$.

\begin{dfn} Let Assumptions \ref{an2ass1} and \ref{an3ass} hold
and $(\tau,T,\le)$ be a permissible weak stability condition. For
all $\al\in C(\A)$ and $\A$-data $(I,\pr,\ka)$ define
\e
\Iss^\al(\tau)^{\sst\La}\!=\!\Up'\ci\Pi_*\,\bdss^\al(\tau)
\;\>\text{and}\;\>
\Iss(I,\pr,\ka,\tau)^{\sst\La}\!=\!\Up'\ci\Pi_*\ci\bs\si(I)_*\,
\bdss(I,\pr,\ka,\tau)
\label{an6eq10}
\e
in $\La$, where $\Up':\uSF(\Spec\K)\ra\La$ is as in Theorem
\ref{an2thm4} and $\Pi:\fObj_\A\ra\Spec\K$ is the projection, so
that $\Pi_*$ maps $\SFa(\fObj_\A)\ra\uSF(\Spec\K)$. Now suppose
$(I,\pr)$ is {\it connected}. Then \cite[\S 8]{Joyc7} defines stack
functions $\bep^\al(\tau),\bdsib(I,\pr,\ka,\tau)$, which lie in
$\SFai(\fObj_\A)$ by \cite[Th.~8.7]{Joyc7}. Define
\e
\begin{split}
J^\al(\tau)^{\sst\La^\ci}&=(\el-1)\Up'\ci\Pi_*\,\bep^\al(\tau)
\qquad\text{and}\\
\Jsib(I,\pr,\ka,\tau)^{\sst\La^\ci}&=(\el-1)\Up'\ci\Pi_*\ci
\bs\si(I)_*\,\bdsib(I,\pr,\ka,\tau).
\end{split}
\label{an6eq11}
\e

Suppose $f\in\SFai(\fObj_\A)$, so that $\Pi^\vi_1(f)=f$. Then
$\Pi^\vi_1$ is also the identity on the projection
$\bar\Pi^{\Up,\La^\ci}_{\Spec\K}\ci\Pi_*(f)$ of $f$ to
$\uoSF(\Spec\K,\Up,\La^\ci)$, since $\Pi^\vi_1$ commutes with
$\bar\Pi^{\Up,\La^\ci}_{\Spec\K}$ and $\Pi_*$. Using the explicit
description \cite[Prop.~6.11]{Joyc4} of $\uoSF(\Spec\K,\Up,\La^\ci)$
we now see that $\bar\Pi^{\Up,\La^\ci}_{\Spec\K}\ci\Pi_*(f)=
\be[[\Spec\K/\K^\t]]$ for some $\be\in\La^\ci$. Now $\Up'$ factors
via $\bar\Pi^{\Up,\La^\ci}_{\Spec\K}$, and it easily follows that
$(\el-1)\Up'\ci\Pi_*(f)=\be$, so that $(\el-1)\Up'\ci\Pi_*$ maps
$\SFai(\fObj_\A)\ra\La^\ci\subset\La$. It follows that
$J^\al(\tau)^{\sst\La^\ci},\Jsib(I,\pr,\ka,\tau)^{\sst\La^\ci}$
actually lie in $\La^\ci$. Thus we may define
\e
J^\al(\tau)^{\sst\Om}=\pi\bigl(J^\al(\tau)^{\sst\La^\ci}\bigr)
\quad\text{and}\quad
\Jsib(I,\pr,\ka,\tau)^{\sst\Om}=\pi\bigl(\Jsib(I,\pr,\ka,
\tau)^{\sst\La^\ci}\bigr),
\label{an6eq12}
\e
where $\pi:\La^\ci\ra\Om$ is as in Assumption~\ref{an2ass1}.
\label{an6def4}
\end{dfn}

Note that Remark \ref{an6rem1}(b) applies, so that
\begin{equation*}
\Iss^\al(\tau)^{\sst\La}=\Up'\bigl(\bigl[\Oss^\al(\tau)\bigr]\bigr)
\;\>\text{and}\;\> \Iss(I,\pr,\ka,\tau)^{\sst\La}=\Up'\bigl(\bigl[
\Mss(I,\pr,\ka,\tau)_\A\bigr]\bigr).
\end{equation*}
When $(\tau,T,\le)$ is a stability condition we can also define
$\Jstb(I,\pr,\ka,\tau)^{\sst\La^\ci,\Om}$ in the same way, using
$\bdstb(I,\pr,\ka,\tau)\in\SFai(\fObj_\A)$ from \cite[\S 8]{Joyc7}.
With the obvious notation we have $\Jsib(I,\pr,
\ka,\tau)^{\sst\La^\ci}=(\el-1)\Isib(I,\pr,\ka,\tau)^{\sst\La}$.

\begin{thm} Let Assumptions \ref{an2ass1} and \ref{an3ass} hold
and\/ $\chi:K(\A)\t K(\A)\ra\Z$ be biadditive and satisfy
\eq{an6eq8}. Then for all\/ $\al\in C(\A)$ and\/ $\A$-data
$(I,\pr,\ka)$ the following hold, with only finitely many nonzero
terms:
\ea
J^\al(\tau)^{\sst\La^\ci}&=
\sum_{\begin{subarray}{l}\text{$\A$-data $(\{1,\ldots,n\},\le,\ka):$}\\
\text{$\ka(\{1,\ldots,n\})=\al$, $\tau\ci\ka\equiv\tau(\al)$}
\end{subarray}
\!\!\!\!\!\!\!\!\!\!\!\!\!\!\!\!\!\!\!\!\!\!\!\!\!\!\!\!\!\!\!\!
\!\!\!\!\!\!\!\!\!\!\!\!\!\!\!\!\!\!\!\!\!\!\! }
\el^{-\sum_{1\le i<j\le n}\chi(\ka(j),\ka(i))}
\frac{(-1)^{n-1}(\el-1)}{n}\,
\prod_{i=1}^n\Iss^{\ka(i)}(\tau)^{\sst\La},
\label{an6eq13}
\\
\Iss^\al(\tau)^{\sst\La}&=
\sum_{\begin{subarray}{l}\text{$\A$-data $(\{1,\ldots,n\},\le,\ka):$}\\
\text{$\ka(\{1,\ldots,n\})=\al$, $\tau\ci\ka\equiv\tau(\al)$}
\end{subarray}
\!\!\!\!\!\!\!\!\!\!\!\!\!\!\!\!\!\!\!\!\!\!\!\!\!\!\!\!\!\!\!\!
\!\!\!\!\!\!\!\!\!\!\!\!\!\!\!\!\!\!\!\!\!\!\! }
\el^{-\sum_{1\le i<j\le n}\chi(\ka(j),\ka(i))}\frac{(\el-1)^{-n}}{n!}\,
\prod_{i=1}^nJ^{\ka(i)}(\tau)^{\sst\La^\ci},
\label{an6eq14}\\
\text{and}&\qquad \Iss(I,\pr,\ka,\tau)^{\sst\La}=
\el^{-\sum_{i\ne j\in I:\,\,i\pr j}\chi(\ka(j),\ka(i))}\cdot
\ts\prod_{i\in I}\Iss^{\ka(i)}(\tau)^{\sst\La}.
\label{an6eq15}
\ea

Suppose $(\tau,T,\le),(\ti\tau,\ti T,\le)$ are permissible weak
stability conditions on $\A$, the change from $(\tau,T,\le)$ to
$(\ti\tau,\ti T,\le)$ is globally finite, and there exists a weak
stability condition $(\hat\tau,\hat T,\le)$ on $\A$ dominating
$(\tau,T,\le),(\ti\tau,\ti T,\le)$ with the change from
$(\hat\tau,\hat T,\le)$ to $(\ti\tau,\ti T,\le)$ locally finite.
Then for all\/ $\al\in C(\A)$ the following hold, with only finitely
many nonzero terms:
\ea
\begin{gathered}
\Iss^\al(\ti\tau)^{\sst\La}=
\sum_{\begin{subarray}{l}
\text{$\A$-data $(\{1,\ldots,n\},\le,\ka):$}\\
\text{$\ka(\{1,\ldots,n\})=\al$}\end{subarray}
\!\!\!\!\!\!\!\!\!\!\!\!\!\!\!\!\!\!\!\!\!\!\!\!
\!\!\!\!\!\!\!\!\!\!\!\!\!\!\!\!\! }
\begin{aligned}[t]
S(\{1,\ldots,n\},\le,\ka,\tau,\ti\tau)\cdot
\el^{-\sum_{1\le i<j\le n}\chi(\ka(j),\ka(i))}\cdot&\\
\ts\prod_{i=1}^n\Iss^{\ka(i)}(\tau)^{\sst\La}&,
\end{aligned}
\end{gathered}
\label{an6eq16}\\
\begin{gathered}
J^\al(\ti\tau)^{\sst\La^\ci}=
\sum_{\begin{subarray}{l}
\text{$\A$-data $(\{1,\ldots,n\},\le,\ka):$}\\
\text{$\ka(\{1,\ldots,n\})=\al$}\end{subarray}
\!\!\!\!\!\!\!\!\!\!\!\!\!\!\!\!\!\!\!\!\!\!\!\!
\!\!\!\!\!\!\!\!\!\!\!\!\!\!\!\!\! }
\begin{aligned}[t]
U(\{1,\ldots,n\},\le,\ka,\tau,\ti\tau)\cdot
\el^{-\sum_{1\le i<j\le n}\chi(\ka(j),\ka(i))}\cdot&\\
(\el-1)^{1-n}\ts\prod_{i=1}^nJ^{\ka(i)}(\tau)^{\sst\La^\ci}&.
\end{aligned}
\end{gathered}
\label{an6eq17}
\ea
\label{an6thm2}
\end{thm}

\begin{proof} The definition of $\Phi^{\sst\La}$ in \cite[\S
6.2]{Joyc6} gives $\Phi^{\sst\La}\ci\Pi^{\Up,\La}_{\fObj_\A}
\bigl(\bdss^\al(\tau)\bigr)\!=\!\Iss^\al(\tau)^{\sst\La}a^\al$,
$\Phi^{\sst\La}\ci\Pi^{\Up,\La}_{\fObj_\A}\bigl(\bs\si(I)_*\,
\bdss(I,\pr,\ka,\tau)\bigr)\!=\!\Iss(I,\pr,\ka,\tau)^{\sst\La}
a^{\ka(I)}$ and $\Phi^{\sst\La}\ci\Pi^{\Up,\La}_{\fObj_\A}\bigl(
\ep^\al(\tau)\bigr)\!=\!(\el\!-\!1)^{-1}J^\al(\tau)^{\sst\La^\ci}
a^\al$, where $a^\al$ for $\al\in\bar C(\A)$ are a $\La$-basis for
$A(\A,\La,\chi)$. Equations \eq{an6eq13}--\eq{an6eq14} and
\eq{an6eq16}--\eq{an6eq17} follow from \eq{an3eq14}, \eq{an3eq16},
\eq{an5eq3} and \eq{an5eq7} respectively, the definition
\cite[Def.~6.3]{Joyc6} of multiplication in $A(\A,\La,\chi)$, and
the fact \cite[Th.~6.4]{Joyc6} that \eq{an6eq9} is an algebra
morphism. The powers of $\el-1$ compensate for the factor $\el-1$ in
\eq{an6eq11}. There are only finitely many nonzero terms in each
equation as this holds for \eq{an3eq14}, \eq{an3eq16} \eq{an5eq3},
and \eq{an5eq7}. Equation \eq{an6eq15} follows from
$\bs\si(I)_*\,\bdss(I,\pr,\ka, \tau)=P_\sIp\bigl(\bdss^{\ka(i)}
(\tau):i\in I\bigr)$, the definition \cite[Def.~6.3]{Joyc6} of
$P_\sIp$ in $A(\A,\La,\chi)$, and the fact \cite[Th.~6.4]{Joyc6}
that \eq{an6eq9} intertwines the $P_\sIp$ on $\SFa(\fObj_\A)$
and~$A(\A,\La,\chi)$.
\end{proof}

Equations \eq{an6eq13}--\eq{an6eq14} show that when \eq{an6eq8}
holds the systems of invariants $\Iss^\al(\tau)^{\sst\La}$ and
$J^\al(\tau)^{\sst\La^\ci}$ for $\al\in C(\A)$ are equivalent. But
one can argue the $J^\al(\tau)^{\sst\La^\ci}$ are better, since they
take their values in the smaller algebra $\La^\ci$, and so it takes
less information to describe them. Equivalently, the
$\Iss^\al(\tau)^{\sst\La}$ satisfy natural identities, that the
right hand side of \eq{an6eq13} lies in $\La^\ci$ for all~$\al\in
C(\A)$.

Suppose now that $\chi$ is {\it symmetric}, that is, $\chi(\al,\be)=
\chi(\be,\al)$ for all $\al,\be\in K(\A)$. Then $A(\A,\La,\chi)$ is
a {\it commutative} algebra, as its multiplication relations are
$a^\al\star a^\be=\el^{-\chi(\be,\al)}a^{\al+\be}$. Now Theorem
\ref{an5thm2} shows that \eq{an5eq7} may be rewritten
\begin{equation*}
\bep^\al(\ti\tau)=\bep^\al(\tau)+\text{sum of multiple commutators
of two or more $\bep^{\ka(i)}(\tau)$.}
\end{equation*}
When we project this to the commutative algebra $A(\A,\La,\chi)$
under the algebra morphism \eq{an6eq9} the multiple commutators go
to zero, giving $\Phi^{\sst\La}\ci\Pi^{\Up,\La}_{\fObj_\A}
\bigl(\bep^\al(\ti\tau)\bigr)=\Phi^{\sst\La}\ci\Pi^{\Up,
\La}_{\fObj_\A}\bigl(\bep^\al(\tau)\bigr)$. So as $\Phi^{\sst\La}
\ci\Pi^{\Up,\La}_{\fObj_\A}\bigl(\ep^\al(\tau)\bigr)\!=\!(\el-1)^{-1}
J^\al(\tau)^{\sst\La^\ci}a^\al$ we deduce:

\begin{cor} In Theorem \ref{an6thm2}, if\/ $\chi$ is symmetric then
$J^\al(\ti\tau)^{\sst\La^\ci}=J^\al(\tau)^{\sst\La^\ci}$ and\/
$J^\al(\ti\tau)^{\sst\Om}=J^\al(\tau)^{\sst\Om}$ for all\/~$\al\in
C(\A)$.
\label{an6cor1}
\end{cor}

We compute the invariants of Definition \ref{an6def4} explicitly
when~$\A\!=\!\modKQ$.

\begin{ex} Let $Q=(Q_0,Q_1,b,e)$ be a {\it quiver}. That is, $Q_0$
is a finite set of {\it vertices}, $Q_1$ is a finite set of {\it
arrows}, and $b,e:Q_1\ra Q_0$ are maps giving the {\it beginning}
and {\it end\/} of each arrow. Fix an algebraically closed field
$\K$, and take $\A$ to be the abelian category $\modKQ$ of {\it
representations} of $Q$. Objects ${\bf X}=(X_v,\rho_a)\in\modKQ$
comprise of finite-dimensional $\K$-vector spaces $X_v$ for all
$v\in Q_0$ and linear maps $\rho_a:X_{b(a)}\ra X_{e(a)}$ for
all~$a\in Q_1$.

Define data $K(\A),\fF_\A$ satisfying Assumption \ref{an3ass} as in
\cite[Ex.~10.5]{Joyc5}. Then $K(\A)=\Z^{Q_0}$, with elements of
$K(\A)$ written as maps $Q_0\ra\Z$, and $[(X_v,\rho_a)]=\al\in
K(\A)$ if $\dim_\K X_v=\al(v)$ for all $v\in Q_0$. It is well-known
that $\Ext^m(X,Y)=0$ for all $X,Y\in\A$ and $m\ge 2$, and
\eq{an6eq8} holds with $\chi$ given by
\begin{equation*}
\chi(\al,\be)=\ts\sum_{v\in Q_0}\al(v)\be(v)-\sum_{a\in
Q_1}\al(b(a))\be(e(a)) \;\>\text{for $\al,\be\in\Z^{Q_0}$.}
\end{equation*}
As in the proof of \cite[Th.~10.11]{Joyc5}, for $\al\in C(\A)$ there
is a 1-isomorphism
\e
\fObj_\A^\al\cong\bigl[\K^{\sum_{a\in Q_1}\al({b(a)})\al(e(a))}/
\ts\prod_{v\in Q_0}\GL(\al(v),\K)\bigr].
\label{an6eq18}
\e

Now suppose Assumption \ref{an2ass1} holds with this $\K$. Then
Theorem \ref{an2thm4} defines an algebra morphism
$\Up':\uSF(\Spec\K)\ra\La$. Noting that $\prod_{v\in
Q_0}\GL(\al(v),\K)$ is a {\it special\/} algebraic $\K$-group and
using Theorem \ref{an2thm4}, $\Up([\K^n])=\el^n$ and a formula for
$\Up([\GL(m,\K)])$ in \cite[Lem.~4.5]{Joyc4}, we deduce from
\eq{an6eq18} that
\e
\Up'\bigl([\fObj_\A^\al]\bigr)=\frac{\ts\el^{\sum_{a\in Q_1}
\al(b(a))\al(e(a))-\sum_{v\in Q_0}\al(v)(\al(v)-1)/2}}{
\ts\prod_{v\in Q_0}\prod_{k=1}^{\al(v)}\bigl(\el^k-1\bigr)}\,.
\label{an6eq19}
\e

Let $(\tau,T,\le)$ be any weak stability condition on $\A$, such as
the slope stability conditions of \cite[Ex.~4.14]{Joyc5}. Then
$(\tau,T,\le)$ is permissible by \cite[Cor.~4.13]{Joyc5}. Define
invariants $\Iss^\al(\tau)^{\sst\La},\Iss(I,\pr,\ka,\tau)^{\sst\La}$
as in \eq{an6eq10}. One possibility for $(\tau,T,\le)$ is the
trivial stability condition $(0,\{0\},\le)$. Since every object is
0-semistable we have $\Oss^\al(0)=\fObj_\A^\al(\K)$, and thus
$\Iss^\al(\tau)^{\sst\La}=\Up'\bigl([\fObj_\A^\al]\bigr)$, which is
given by \eq{an6eq19}. Applying Theorem \ref{an6thm2} with
$(0,\{0\},\le)$ in place of $(\tau,T,\le),(\hat\tau,\hat T,\le)$ and
$(\tau,T,\le)$ in place of $(\ti\tau,\ti T,\le)$, and simplifying
$S(\{1,\ldots,n\},\le,\ka,0,\tau)$ using \eq{an5eq18}, from
\eq{an6eq16} and \eq{an6eq19} we deduce that for all $\al\in C(\A)$
we have
\begin{align}
\sum_{\begin{subarray}{l}\text{$\A$-data}\\
\text{$(\{1,\ldots,n\},\le,\ka):$}\\
\text{$\ka(\{1,\ldots,n\})=\al$,}\\
\text{$\tau\!\ci\!\ka(\{1,\ldots,i\})\!>\!\tau(\al)$,}\\
\text{$1\le i<n$}\end{subarray}\!\!\!\!\!\!\!\!\!\!\!\!}
\!\!\!\!\!\!\!\!\! &\begin{aligned}[t]
(-1)^{n-1}\cdot\,&\displaystyle \biggl[\,\prod_{1\le i<j\le
n}\!\!\!\el^{\sum_{a\in Q_1}\!\ka(j)(b(a))\ka(i)(e(a))-\sum_{v\in
Q_0}\!\ka(j)(v)\ka(i)(v)}\biggr]
\cdot\\
&\displaystyle \biggl[\,\prod_{i=1}^n\frac{\ts\el^{\sum_{a\in Q_1}\!
\ka(i)(b(a))\ka(i)(e(a))\!-\!\sum_{v\in Q_0}\!
\ka(i)(v)(\ka(i)(v)\!-\!1)/2}}{\ts\prod_{v\in
Q_0}\prod_{k=1}^{\ka(i)(v)}\bigl(\el^k-1\bigr)}\biggr]
\end{aligned}
\nonumber\\
& \qquad\qquad\qquad\qquad\qquad\qquad\qquad
=\Iss^\al(\tau)^{\sst\La}.
\label{an6eq20}
\end{align}
Combining this with \eq{an6eq13} and \eq{an6eq15} gives expressions
for $J^\al(\tau)^{\sst\La^\ci},\Iss(I,\pr,\ka,\tau)^{\sst \La}$. It
might be interesting rewrite this formula for
$J^\al(\tau)^{\sst\La^\ci}$ so that every term lies in $\La^\ci$.
Then projecting to $\Om$ would give a formula for
$J^\al(\tau)^{\sst\Om}$. This would, for instance, enable easier
calculation of Euler characteristics of quiver moduli spaces in the
case~$\Oss^\al(\tau)=\Ost^\al(\tau)$.
\label{an6ex1}
\end{ex}

We can relate Example \ref{an6ex1} to results of Reineke
\cite{Rein}. Take $\K=\C$, let $(\tau,T,\le)$ be a slope stability
condition on $\A=\modCQ$ as in \cite[Ex.~4.14]{Joyc7}, and let
$\al\in C(\A)$ be `coprime', which essentially means that
$\Ost^\al(\tau)=\Oss^\al(\tau)$, that is, $\tau$-stability and
$\tau$-semistability coincide in class $\al$. Regarding
$\Oss^\al(\tau)$ as an open $\C$-substack of $\fObj_\A$ we have a
1-isomorphism $\Oss^\al(\tau)\cong\Mss^\al(\tau)\t[\Spec\C/\C^\t]$,
where $\Mss^\al(\tau)$ is a {\it nonsingular complex projective
variety}. Here the factor $[\Spec\C/\C^\t]$ arises as
$\Iso_\C(X)\cong\C^\t$ for all $X\in\Ost^\al(\tau)=\Oss^\al(\tau)$,
and we can split the stabilizer groups off as a product with
$[\Spec\C/\C^\t]$. Since $\Up([\C^\t])=\el-1$, we see that
\e
\Up\bigl(\bigl[\Mss^\al(\tau)\bigr]\bigr)=(\el-1)
\Up'\bigl(\bigl[\Mss^\al(\tau)\t[\Spec\C/\C^\t]\bigr]\bigr)=
(\el-1)\Iss^\al(\tau)^{\sst\La}.
\label{an6eq21}
\e
Equations \eq{an6eq20} and \eq{an6eq21} thus give an explicit
expression for~$\Up([\Mss^\al(\tau)])$.

Now in \cite[Cor.~6.8]{Rein}, Reineke gives a formula for the
Poincar\'e polynomial $P\bigl(\Mss^\al(\tau);z\bigr)$, as in Example
\ref{an2ex1}. Putting $\el=z^2$ in \eq{an6eq20} and \eq{an6eq21},
careful comparison shows that our formula agrees with Reineke's. We
have reproved Reineke's result by quite different methods, and also
extended it: Reineke restricts to $\K=\C$, to slope functions, to
Poincar\'e polynomials, and to coprime $\al$, but Example
\ref{an6ex1} holds for algebraically closed $\K$, arbitrary weak
stability conditions, any motivic invariant $\Up$ satisfying
Assumption \ref{an2ass1} (such as virtual Hodge polynomials), and
all $\al\in C(\A)$. In particular, we can now interpret Reineke's
formula for non-coprime~$\al$.

\subsection{Counting vector bundles on smooth curves}
\label{an63}

In two seminal papers which have been elaborated on by many other
authors since, Harder and Narasimhan \cite{HaNa} and Atiyah and Bott
\cite{AtBo} found recursive formulae for the Poincar\'e polynomials
of moduli spaces of semistable vector bundles with fixed rank and
determinant over a smooth curve $P$ of genus $g$. We will now recast
this in our own framework, using the ideas of~\S\ref{an24}.

Let $\K$ be an algebraically closed field, $P$ a smooth projective
curve over $\K$ with genus $g$, and take $\A=\coh(P)$ with data
$K(\A),\fF_\A$ satisfying Assumption \ref{an3ass} as in
\cite[Ex.~9.1]{Joyc5}. Identify $K(\A)=\Z^2$, such that if $X$ is a
vector bundle over $P$ with rank $n$ and degree $d$ then $[X]=(n,d)$
in $K(\A)=\Z^2$, and then
\e
C(\A)=\bigl\{(n,d)\in\Z^2:\text{$n\ge 0$, and $d>0$ if
$n=0$}\bigr\}.
\label{an6eq22}
\e
Also $\Ext^i(X,Y)=0$ for all $X,Y\in\A$ and $i>1$, so \eq{an6eq8}
holds, and using the Riemann--Roch Theorem we find that
$\chi:K(\A)\t K(\A)\ra\Z^2$ is given by
\e
\chi\bigl((n_1,d_1),(n_2,d_2)\bigr)=n_1d_2-d_1n_2-(g-1)n_1n_2.
\label{an6eq23}
\e

There is an ample line bundle $E$ on $P$ with $[E]=(1,1)$ in
$K(\A)$, such that if $X\in\A$ with $[X]=(n,d)$ in $K(\A)$ then the
Hilbert polynomial of $X$ w.r.t.\ $E$ is $p_X(t)=nt+d+n(1-g)$. Let
$(\ga,G_1,\le)$ be the Gieseker stability condition on $\A=\coh(P)$
defined in \cite[Ex.~4.16]{Joyc7} using this $E$; since $\dim P=1$
this coincides with $\mu$-stability in \cite[Ex.~4.17]{Joyc7}. It is
{\it permissible} by \cite[Th.~4.20]{Joyc7}. Let $(\de,D_1,\le)$ be
the {\it purity} weak stability condition on $\A$ defined in
\cite[Ex.~4.18]{Joyc7}, so that $X\in\A$ is $\de$-semistable if and
only if it is pure. It is {\it not\/} permissible.

As $(\de,D_1,\le)$ dominates $(\ga,G_1,\le)$ Theorem \ref{an5thm4}
applies with $(\ga,G_1,\le)$, $(\de,D_1,\le)$ in place of
$(\tau,T,\le),(\ti\tau,\ti T,\le)$. We shall rewrite \eq{an5eq15}
with $\al=(n,d)$ for $n>0$ and $d\in\Z$. For $n,\ka$ as in
\eq{an5eq15}, replace $n$ by $k$ and write $\ka(i)=(n_i,d_i)$ for
$i=1,\ldots,k$, and use \eq{an6eq22}. Since $n>0$ we have
$\de(\al)=t$, and from \cite[Ex.~4.18]{Joyc7} we see that
$\de\ci\ka(i)$ is $t$ if $n_i>0$ and 1 if $n_i=0$, so the condition
$\ti\tau\ci\ka\equiv\ti\tau (\al)$ in \eq{an5eq15} is equivalent to
$n_i>0$ for all $i$. Then $\ga\ci\ka(i)=t+d_i/n_i+1-g$, so
$\tau\ci\ka(1)>\cdots >\tau\ci\ka(n)$ in \eq{an5eq15} is equivalent
to $d_1/n_1>d_2/n_2>\cdots>d_k/n_k$. Putting all this together,
Theorem \ref{an5thm4} shows that
\e
\begin{gathered}
\sum_{k=1}^n\,\sum_{\begin{subarray}{l} n_1,\ldots,n_k>0:\\
n_1+\cdots+n_k=n\end{subarray}}\,\,
\sum_{\begin{subarray}{l} d_1,\ldots,d_k\in\Z:\;
d_1+\cdots+d_k=d,\\
d_1/n_1>\cdots>d_k/n_k \end{subarray}
\!\!\!\!\!\!\!\!\!\!\!\!\!\!\!\!\!\!\!\!\!\!\!\!\!\!\!\!
\!\!\!\!\!\!\!\!\!\!\!\! }
\begin{aligned}[t]
\bdss^{(n_1,d_1)}(\ga)*\bdss^{(n_2,d_2)}(\ga)*\cdots*
\bdss^{(n_k,d_k)}(\ga)&\\
=\bdss^{(n,d)}(\de)&,
\end{aligned}
\end{gathered}
\label{an6eq24}
\e
which holds as an infinite sum in $\LSF(\fObj_\A)$ converging as in
Definition~\ref{an2def10}.

We shall show that \eq{an6eq24} is {\it strongly convergent}. To do
this we will need some dimension calculations. Suppose $\al\in
C(\A)$ and $[X]\in\Oss^\al(\ga)$. Then the Zariski tangent space to
$\fObj_\A^\al$ at $[X]$ is $\Ext^1(X,X)$, and $\Iso_\K(X)=\Aut(X)$
which is open in $\Hom(X,X)$. Since $\Ext^2(X,X)=0$ and
$\Oss^\al(\ga)$ is open in $\fObj_\A^\al(\K)$ we see by deformation
theory (see for instance \cite[\S 2.A, \S 4.5]{HuLe}) that the
dimension of $\Oss^\al(\ga)$ near $X$ is $\dim\Ext^1(X,X)-
\dim\Hom(X,X)$, which is $-\chi(\al,\al)$ by \eq{an6eq8}. As this is
independent of $[X]$ we see that $\dim\Oss^\al(\ga)=\dim\Oss^\al
(\de)=-\chi(\al,\al)$ provided $\Oss^\al(\ga),\Oss^\al(\de)\ne
\emptyset$. Therefore
\e
\bdss^\al(\ga)\in\SF(\fObj_\A)_{-\chi(\al,\al)} \quad\text{and}\quad
\bdss^\al(\de)\in\LSF(\fObj_\A)_{-\chi(\al,\al)}.
\label{an6eq25}
\e

\begin{lem} In the situation above, suppose $a,b\in\Z$, $\al,\be\in
C(\A)$ and\/ $f\in\SF(\fObj_\A)_a$ and\/ $g\in\SF(\fObj_\A)_b$ are
supported on $\fObj_\A^\al(\K)$ and\/ $\fObj_\A^\be(\K)$
respectively. Then~$f*g\in\SF(\fObj_\A)_{a+b-\chi(\be,\al)}$.
\label{an6lem1}
\end{lem}

\begin{proof} Clearly $f\ot g\in\SF(\fObj_\A\t\fObj_\A)_{a+b}$. Now
\cite[Cor.~5.15]{Joyc6} gives explicit expressions for $f\ot g$ and
$f*g$, involving vector spaces $E^0_m,E^1_m$ isomorphic to
$\Hom(Y,X)$ and $\Ext^1(Y,X)$ for $([X],[Y])$ in the support of
$f\ot g$. As $([X],[Y])\in\fObj_\A^\al(\K)\t\fObj_\A^\be(\K)$ we see
from \eq{an6eq8} that $\dim E^1_m-\dim E^0_m=-\chi(\be,\al)$. The
lemma then easily follows from~\cite[Cor.~5.15]{Joyc6}.
\end{proof}

If $(n_i,d_i)\in C(\A)$ for $i=1,\ldots,k$ with $n=n_1+\cdots+n_k$
and $d=d_1+\cdots+d_k$ then using \eq{an6eq23}, \eq{an6eq25},
induction on $k$ and biadditivity of $\chi$ we see that
\e
\begin{split}
&\bdss^{(n_1,d_1)}(\ga)*\bdss^{(n_2,d_2)}(\ga)*\cdots*
\bdss^{(n_k,d_k)}(\ga)\in\La_N,\quad\text{where}\\
N&=\ts-\sum_{i=1}^k\chi\bigl((n_i,d_i),(n_i,d_i)\bigr)\!-\!
\sum_{1\le i<j\le k}\chi\bigl((n_j,d_j),(n_i,d_i)\bigr)\\
&=\ts-\chi\bigl((n,d),(n,d)\bigr)+\sum_{1\le
i<j\le k}\chi\bigl((n_i,d_i),(n_j,d_j)\bigr)\\
&=\ts(g\!-\!1)n^2\!-\!(g\!-\!1)\sum_{1\le i<j\le k}n_in_j
\!+\!\sum_{1\le i<j\le k}(n_id_j\!-\!d_in_j)\\
&\le\ts\max\bigl(0,(g-1)\bigr)n^2+\sum_{1\le i<j\le
k}(n_id_j-d_in_j).
\end{split}
\label{an6eq26}
\e
For $n_i,d_i$ as in \eq{an6eq24} we have $d_1/n_1>\cdots>d_k/n_k$,
which implies that all terms $n_id_j-d_in_j$ in the bottom line of
\eq{an6eq26} are {\it negative}. Using this we prove:

\begin{prop} Equation \eq{an6eq24} is strongly convergent, in the
sense of Definition \ref{an2def10}. Hence $\bdss^{(n,d)}(\de)\in
\ESF(\fObj_\A)$ for all\/ $n>0$ and\/~$d\in\Z$.
\label{an6prop2}
\end{prop}

\begin{proof} We already know \eq{an6eq24} is convergent. To prove it
is strongly convergent, it is enough to show that for each $m\in\Z$
there are only finitely many choices of $k$ and $n_i,d_i$ in
\eq{an6eq24} with $N>m$ in \eq{an6eq26}. Let $k,n_i,d_i$ be some
such choice. As $n_id_j-d_in_j<0$ for all $1\le i<j\le k$ from
above, we see from \eq{an6eq26} that for all $1\le i<j\le k$ we have
\e
m-\max\bigl(0,(g-1)\bigr)n^2<n_id_j-d_in_j<0.
\label{an6eq27}
\e
Clearly there are only finitely many choices for $k$ and
$n_1,\ldots,n_k$. For each such choice \eq{an6eq27} shows there are
only finitely many possibilities for the integers $n_id_j-d_in_j$
for $1\le i<j\le k$. But these and $d_1+\cdots+d_k=d$ determine
$d_1,\ldots,d_k$, so there are only finitely many choices for
$d_1,\ldots,d_k$. This proves \eq{an6eq24} is strongly convergent.
The final part follows as $\ESF(\fObj_A)$ is closed under strongly
convergent limits, from Definition~\ref{an2def10}.
\end{proof}

If $d>0$, so that $(0,d)\in C(\A)$, it is easy to show that every
$[X]\in\fObj_\A(\K)$ is $\ga$- and $\de$-semistable, so that
$\Oss^{(0,d)}(\ga)=\Oss^{(0,d)}(\de)=\fObj_\A^{(0,d)}(\K)$.
Therefore $\bdss^{(0,d)}(\de)=\bdss^{(0,d)}(\ga)\in\SF(\fObj_\A)$,
as $(\ga,G_1,\le)$ is permissible. Together with Proposition
\ref{an6prop2} this shows $\bdss^\al(\de)\in\ESF(\fObj_\A)$ for
all~$\al\in C(\A)$.

\begin{dfn} Let Assumption \ref{an3ass} hold and $(\tau,T,\le)$ be
a weak stability condition on $\A$. Generalizing Definition
\ref{an3def12}, we call $(\tau,T,\le)$ {\it essentially permissible}
if (i) $\A$ is $\tau$-artinian and (ii) $\bdss^\al(\tau)\in
\ESF(\fObj_\A)$ for all $\al\in C(\A)$. Clearly, $(\tau,T,\le)$
permissible implies $(\tau,T,\le)$ essentially permissible.

When $\A=\coh(P)$ for $P$ a smooth curve as above, part (i) holds
for $(\de,D_1,\le)$ by \cite[Lem.~4.19]{Joyc7}, and (ii) holds as
above. Thus $(\de,D_1,\le)$ is essentially permissible, but not
permissible.

Suppose Assumptions \ref{an2ass2} and \ref{an3ass} hold and
$(\tau,T,\le)$ is an essentially permissible weak stability
condition on $\A$. For $\al\in C(\A)$, define {\it invariants}
$\Iss^\al(\tau)^{\sst\La}$ in $\La$ by $\Iss^\al(\tau)^{\sst\La}
=\Pi_\La(\bdss^\al(\tau))$, where $\Pi_\La$ is as in Definition
\ref{an2def12}. If $(\tau,T,\le)$ is permissible then
$\bdss^\al(\tau)\in\SF(\fObj_\A)$, and so the $\Iss^\al(\tau)^{
\sst\La}$ agree with those defined in \eq{an6eq10}, as $\Pi_\La$
coincides with $\Up'\ci\Pi_*$ on $\SF(\fObj_\A)$ by
Definition~\ref{an2def12}.
\label{an6def5}
\end{dfn}

Combining Propositions \ref{an2prop} and \ref{an6prop2} shows that
applying $\Pi_\La$ to \eq{an6eq24} yields a convergent identity in
$\La$. Since \eq{an6eq8} holds we are in the situation of
\S\ref{an62}, and as $\bdss^{(n_i,d_i)}(\ga)\in\SF(\fObj_\A)$ and
$\Pi_\La$ coincides with $\Up'\ci\Pi_*$ on $\SF(\fObj_\A)$, as in
the proof of \eq{an6eq16} we can rewrite $\Pi_\La\bigl(\bdss^{(n_1,
d_1)}(\ga)*\cdots*\bdss^{(n_k,d_k)}(\ga)\bigr)$ as a power of $\el$
times the product of the $\Pi_\La\bigl(\bdss^{(n_i,d_i)}(\ga)
\bigr)$. This proves:

\begin{cor} Suppose Assumption \ref{an2ass2} holds, $\A=\coh(P)$ for
$P$ a smooth projective curve over $\K$ of genus $g$, and\/
$(\ga,G_1,\le),(\de,D_1,\le)$ are as above. Let\/
$\Iss^\al(\ga)^{\sst\La},\Iss^\al(\de)^{\sst\La}$ be as in
Definition \ref{an6def5}. Then for all\/ $n>0$ and\/ $d\in\Z$ we
have
\e
\begin{gathered}
\sum_{k=1}^n\,\sum_{\begin{subarray}{l} n_1,\ldots,n_k>0:\\
n_1+\cdots+n_k=n\end{subarray}}\,\,
\sum_{\begin{subarray}{l} d_1,\ldots,d_k\in\Z:\;
d_1+\cdots+d_k=d,\\
d_1/n_1>\cdots>d_k/n_k \end{subarray}
\!\!\!\!\!\!\!\!\!\!\!\!\!\!\!\!\!\!\!\!\!\!\!\!\!\!\!\!\!\!\!\!\!
\!\! }
\begin{aligned}[t]
\el^{\sum_{1\le i<j\le n}
(n_id_j-d_in_j+(g-1)n_in_j)}\cdot\\[-2pt]
\ts\prod_{i=1}^k\Iss^{(n_i,d_i)}(\ga)^{\sst\La}
=&\,\Iss^{(n,d)}(\de)^{\sst\La},
\end{aligned}
\end{gathered}
\label{an6eq28}
\e
which holds as an infinite convergent sum in $\La,$ as in
Definition~\ref{an2def11}.
\label{an6cor2}
\end{cor}

In the same way, \eq{an5eq20} implies the inverse identity to
\eq{an6eq24}. Proposition \ref{an6prop2} and \eq{an6eq25} imply that
$\bdss^\al (\de)\in\ESF(\fObj_\A)_{-\chi(\al,\al)}$ for all $\al\in
C(\A)$, and using this and a similar argument to Proposition
\ref{an6prop2} we can show this inverse identity is strongly
convergent. Applying $\Pi_\La$ as above thus shows:

\begin{prop} In Corollary \ref{an6cor2}, for $n>0$ and\/ $d\in\Z$
we have
\e
\begin{gathered}
\sum_{k=1}^n\,\sum_{\begin{subarray}{l} n_1,\ldots,n_k>0:\\
n_1+\cdots+n_k=n\end{subarray}}\,\,
\sum_{\begin{subarray}{l} d_1,\ldots,d_k\in\Z:\;
d_1+\cdots+d_k=d,\\
(d_1+\cdots+d_i)/(n_1+\cdots+n_i)>d/n,\; 1\le i<k
\end{subarray}
\!\!\!\!\!\!\!\!\!\!\!\!\!\!\!\!\!\!\!\!\!\!\!\!\!\!\!\!\!\!\!\!\!
\!\!\!\!\!\!\!\!\!\!\!\!\!\!\!\!\!\!\!\!\!\!\!\!\!\!\!\!\!\!\!\!\!
\!\!\!\! }
\begin{aligned}[t]
&\\[-20pt]
(-1)^{k-1}\el^{\sum_{1\le i<j\le n}
(n_id_j-d_in_j+(g-1)n_in_j)}\cdot&\\[-8pt]
\ts\prod_{i=1}^k\Iss^{(n_i,d_i)}(\de)^{\sst\La}&\\[-3pt]
=\Iss^{(n,d)}(\ga)^{\sst\La}&,
\end{aligned}
\end{gathered}
\label{an6eq29}
\e
which holds as an infinite convergent sum in $\La,$ as in
Definition~\ref{an2def11}.
\label{an6prop3}
\end{prop}

We shall give an explicit expression for $\Iss^{(n,d)}(\de)^{\sst
\La}$, using the following notation. For $C$ as above, write ${\rm
Jac}(C)$ for the {\it Jacobian} of $C$, an abelian variety
parametrizing degree 0 line bundles on $C$ which is topologically a
torus $T^{2g}$ when $\K=\C$. For $m\ge 0$ write $C^{(m)}$ for the
$m^{\rm th}$ symmetric power of $C$, which is a nonsingular
projective $\K$-variety of dimension $m$. Then

\begin{thm} In the situation above, for all\/ $n>0$ and\/ $d\in\Z$
we have
\e
\Iss^{(n,d)}(\de)^{\sst\La}=\frac{\Up\bigl([{\rm
Jac}(C)]\bigr)}{\el-1}\sum_{m_2,\ldots,m_n\ge
0\!\!\!\!\!\!\!\!\!\!\!\!\!\!\!}\el^{(n^2-1)(g-1)
-\sum_{a=2}^nam_a}\prod_{a=2}^n\Up\bigl([C^{(m_a)}]\bigr).
\label{an6eq30}
\e
When $\Up,\La$ are virtual Poincar\'e series as in Example
\ref{an2ex2} this simplifies to
\e
\Iss^{(n,d)}(\de)^{\sst\La}=\frac{(z+1)^{2g}(z^3+1)^{2g}\cdots
(z^{2n-1}+1)^{2g}}{(z^2-1)^2(z^4-1)^2\cdots(z^{2n-2}-1)^2(z^{2n}-1)}\,.
\label{an6eq31}
\e
Here \eq{an6eq30} converges in $\La$, and\/ the rational functions
in $\el,z$ in \eq{an6eq30}--\eq{an6eq31} may be interpreted as
elements of\/ $\La$ by writing them as power series
in~$\el^{-1},z^{-1}$.
\label{an6thm3}
\end{thm}

Equation \eq{an6eq30} may be deduced from Behrend and Dhillon
\cite[\S 6]{BeDh}, who using important results of Bifet, Ghione and
Letizia \cite{BGL} and Bialynicki-Birula \cite{BiBi} perform
essentially the same calculation, except that they fix the
determinants of their vector bundles; allowing determinants to vary
gives our extra factor $\Up\bigl([{\rm Jac}(C)]\bigr)$. Equation
\eq{an6eq31} comes from \eq{an6eq30} using the following Poincar\'e
polynomial formulae, the second due to MacDonald:
\begin{equation*}
P\bigl({\rm Jac}(C);z\bigr)=(1+z)^{2g}\quad\text{and}\quad
\sum_{m=0}^\iy P\bigl(S^{(m)}C;z\bigr)t^m=
\frac{(1+tz)^{2g}}{(1-t)(1-tz^2)}\,.
\end{equation*}
Substituting \eq{an6eq30} or \eq{an6eq31} into \eq{an6eq29} then
gives an explicit expression for $\Iss^{(n,d)}(\ga)^{\sst\La}$ which
`counts' semistable vector bundles in class~$(n,d)$.

We now explain how the above relates to work by other authors. Let
$n\ge 1$ and $d\in\Z$ be coprime, and fix a line bundle $L$ over $P$
of degree $d$. Then there exists a moduli space ${\mathcal
N}_{g,n,d}$ of semistable rank $n$ vector bundles over $P$ with
determinant $L$, which is a smooth projective $\K$-variety, so its
Poincar\'e polynomial $P({\mathcal N}_{g,n,d};z)$ is well-defined.
By counting semistable bundles over finite fields ${\mathbb F}_p$
and applying Deligne's solution of the Weil conjectures, Harder and
Narasimhan \cite{HaNa} proved results on Betti numbers of ${\mathcal
N}_{g,n,d}$, that were strengthened by Desale and Ramanan
\cite{DeRa} to a recursive formula for $P({\mathcal N}_{g,n,d};z)$,
when $\K$ is an algebraic closure of~${\mathbb F}_p$.

These relate to our invariants for $\Up,\La$ as in Example
\ref{an2ex2} by
\begin{equation*}
\Iss^{(n,d)}(\ga)^{\sst\La}=\frac{(1+z)^{2g}}{z^2-1}\,P({\mathcal
N}_{g,n,d};z),
\end{equation*}
where $(1+z)^{2g}$ compensates for us not fixing determinants, and
$(z^2-1)^{-1}$ for the stabilizer groups $\K^\t$ of points in
${\mathcal N}_{g,n,d}$, which appear in our stack framework but not
in the moduli schemes context. Then Desale and Ramanan's recursive
formula is equivalent to \eq{an6eq28} above, with \eq{an6eq31}
substituted in for~$\Iss^{(n,d)}(\de)^{\sst\La}$.

Atiyah and Bott \cite{AtBo} also derive recursive formulae for
$P({\mathcal N}_{g,n,d};z)$ when $\K=\C$, by completely different
methods involving the topology of infinite-dimensional spaces. Their
formula is equivalent to \eq{an6eq28}, but in a surprising way. As
they explain on \cite[p.~596]{AtBo}, their formula is an infinite
sum in {\it positive} powers of $z$. When $n,d$ are coprime
${\mathcal N}_{g,n,d}$ is a nonsingular projective $\K$-variety of
dimension $(g-1)(n^2-1)$, and so Poincar\'e duality implies that
\e
P({\mathcal N}_{g,n,d};z)=z^{2(g-1)(n^2-1)}P({\mathcal
N}_{g,n,d};z^{-1}).
\label{an6eq32}
\e
If we replace $z$ by $z^{-1}$ throughout and use \eq{an6eq32} to
transform $P({\mathcal N}_{g,n,d};z^{-1})$ into $P({\mathcal
N}_{g,n,d};z)$ then the Atiyah--Bott formulae again become
equivalent to \eq{an6eq28} above, with $\el=z^2$ and \eq{an6eq31}
substituted in for~$\Iss^{(n,d)}(\de)^{\sst\La}$.

This suggests an interesting explanation for why we are dealing with
infinite series in negative powers of $\el$ and $z$. If $V$ is a
smooth $\C$-variety of dimension $d$ then we have
(compactly-supported) cohomology groups $H^k(V,\C)$ and $H^k_{\rm
cs}(V,\C)$ for $0\le k\le 2d$, with $H^k_{\rm cs}(V,\C)\cong
H^{2d-k}(V,\C)^*$ by Poincar\'e duality. If $\fF$ is a smooth
$\C$-stack of pure dimension $d$ we also have cohomology groups
$H^k(\fF,\C)$ for all $k\ge 0$, which can be nonzero for~$k>2d$.

There may be a notion of {\it compactly-supported cohomology
$H^k_{\rm cs}(\fF,\C)$ of stacks} with $H^k_{\rm cs}(\fF,\C)\cong
H^{2d-k}(\fF,\C)^*$ for $\fF$ smooth of pure dimension $d$, which
could be nonzero for integers $k\le 2d$, and in particular for {\it
negative} $k$. Compactly-supported cohomology is the right kind for
our purposes, as it behaves in a motivic way and is used to define
virtual Poincar\'e polynomials. So we can think of the negative
powers of $\el,z$ in our formulae as measuring some kind of
compactly-supported cohomology of stacks which exists in negative
dimensions.

\subsection{Counting sheaves on surfaces $P$ with $K_P^{-1}$ nef}
\label{an64}

Let $\K$ be an algebraically closed field and $P$ a smooth
projective surface over $\K$. Take $\A=\coh(P)$ with data
$K(\A),\fF_\A$ satisfying Assumption \ref{an3ass} as in
\cite[Ex.~9.1]{Joyc5}. Then $\Ext^i(X,Y)=0$ for all $i>2$ and
$X,Y\in\A$, so from \eq{an6eq7} there is a biadditive $\chi:K(\A)\t
K(\A)\ra\Z$ such that for all $X,Y\in\A$ we have
\e
\dim_\K\Hom(X,Y)-\dim_\K\Ext^1(X,Y)+\dim_\K\Ext^2(X,Y)=
\chi\bigl([X],[Y]\bigr).
\label{an6eq33}
\e
Also, by Serre duality, writing $K_P$ for the canonical bundle of
$P$ we have
\e
\Ext^2(X,Y)\cong\Hom(Y,X\ot K_P)^*\quad\text{for all $X,Y\in\A$.}
\label{an6eq34}
\e
Let $(\ga,G_2,\le)$ be the Gieseker stability condition on $\A$
defined in \cite[Ex.~4.16]{Joyc7} using an ample line bundle $E$ on
$P$. It is permissible by \cite[Th.~4.20]{Joyc7}. Let
$(\de,D_2,\le)$ be the purity weak stability condition on $\A$
defined in~\cite[Ex.~4.18]{Joyc7}.

First we consider how much of \S\ref{an63} can be generalized to the
surface case. Unfortunately the answer is very little. The next
example shows Proposition \ref{an6prop2} does not generalize to
$\KP^2$; a similar argument works for all surfaces~$P$.

\begin{ex} Let $P=\KP^2$. Computations with Chern classes show we
may identify $K(\A)$ with $\big\{(a,b,c)\in\Q^3:a,b,c-\ha
b\in\Z\bigr\}$ so that $[E]=\bigl(\rk(E),c_1(E),c_2(E)-\ha c_1(E)^2)
\bigr)$ for $E\in\A$. Then $\chi$ is given explicitly by
\e
\chi\bigl((a_1,b_1,c_1),(a_2,b_2,c_2)\bigr)=a_1a_2+
\ts\frac{3}{2}(a_1b_2-b_1a_2)-b_1b_2+a_1c_2+c_1a_2.
\label{an6eq35}
\e

Suppose $V_n\in\A$ is a rank 2 vector bundle with $[V_n]=(2,0,0)$ in
$C(\A)$, such that $V_n\ot\O(-n)$ has a section $s$ vanishing
transversely at finitely many points for some $n\ge 1$. Calculation
with Chern classes show there must be $n^2$ such points
$x_1,\ldots,x_{n^2}$, and as $[\O(n)]=(1,n,\ha n^2)$ in $C(\A)$
there is a pure sheaf $X_n$ with $[X_n]=(1,-n,-\ha n^2)$ in $C(\A)$
fitting into short exact sequences
\begin{equation*}
0\ra\O(n)\,{\buildrel s\over\longra}\,V_n\ra X_n\ra
0\quad\text{and}\quad \ts 0\ra
X_n\ra\O(-n)\ra\bigop_{i=1}^{n^2}\O_{x_i}.
\end{equation*}
Here $\O(n)\subset V_n$ is the $\ga$ Harder--Narasimhan filtration
of $V_n$, so $V_n$ is $\ga$-unstable. Clearly $X_n$ is determined up
to isomorphism by $x_1,\ldots,x_{n^2}$. We shall describe the family
of such vector bundles~$V_n$.

We have $\Hom(X_n,\O(n))\cong\Hom(\O(-n),\O(n))\cong H^0(\O(2n))
\cong\K^{2n^2+3n+1}$, and $\Ext^2(X_n,\O(n))^*\cong\Hom(\O(n),X_n
\ot\O(-3))\subseteq \Hom(\O(n),\O(-n-3))\cong H^0(\O(-2n-3))=0$ by
\eq{an6eq34}, so $\Ext^2(X_n,\O(n))=0$. As $\chi([X_n],[\O(n)])=
n^2+3n+1$ by \eq{an6eq35} we see that $\Ext^1(X_n,\O(n))\cong
\K^{n^2}$. One can show $\Ext^1(X_n,\O(n))$ is the direct sum of a
copy of $\K$ located at each $x_i$, and the extension of $\O(n)$ by
$X_n$ corresponding to an element of $\Ext^1(X_n,\O(n))$ is a vector
bundle $V_n$ if and only if the components in each of these $n^2$
copies of $\K$ are nonzero. We also find that $\Aut(V_n)\cong\K^\t
\ltimes\Hom(X_n,\O(n))$, so that $\dim\Aut(V_n)=2n^2+3n+2$, and the
family of all such vector bundles $V_n$ has dimension $3n^2-1$, that
is, $2n^2$ parameters for the choice of $x_1,\ldots,x_{n^2}$ in
$\KP^2$, plus $n^2-1$ for the extensions~$P(\Ext^1(X_n,\O(n)))$.

Thus the family of all such points $[V_n]$ is a $\K$-substack of
$\Oss^{(2,0,0)}(\de)\subset\fObj_\A^{(2,0,0)}$, the pure sheaves in
class $(2,0,0)$, with dimension $n^2-3n-3$, that is, the na\"\i ve
dimension $3n^2-1$ of the family minus the dimension $2n^2+3n+2$ of
the stabilizer groups. Since $n^2-3n-3\ra\iy$ as $n\ra\iy$ we have
$\dim\Oss^{(2,0,0)}(\de)=\iy$. This implies that $\bdss^{(2,0,0)}
(\de)\notin\ESF(\fObj_\A)$, so $(\de,D_2,\le)$ is {\it not\/}
essentially permissible in the sense of Definition~\ref{an6def5}.
\label{an6ex2}
\end{ex}

Because of this, for surfaces the invariants
$\Iss^\al(\de)^{\sst\La}$ counting pure sheaves on $P$ are
undefined, and the analogues of \eq{an6eq28}--\eq{an6eq31} do not
make sense. The example also shows that it will not help to work
with vector bundles rather than pure sheaves, as the $\K$-substack
of $\fObj_\A$ of vector bundles in class $(2,0,0)$ on $\KP^2$ has
dimension $\iy$ and is not essentially permissible either. So we
cannot hope to define invariants counting all vector bundles in
class $\al$ on a surface $P$ in any meaningful sense.

However, note that if $P$ is a {\it ruled surface} with ruling
$\pi:P\ra C$, Yoshioka \cite{Yosh1} computes the Betti numbers of
moduli spaces of stable rank 2 coherent sheaves on $P$ in a similar
way to Harder and Narasimhan for curves. His method involves
counting over finite fields the number of sheaves in class $\al$ on
$P$ whose restriction to generic fibres of $\pi$ is semistable. By
analogy with \S\ref{an63}, the author expects that semistability on
generic fibres of $\pi$ comes from an essentially permissible
stability condition $(\ze,Z,\le)$ on $\coh(P)$, and one could hope
to prove an analogue of Theorem \ref{an6thm3} evaluating the
invariants~$\Iss^\al(\ze)^{\sst\La}$.

In the rest of the section we will show that if the anticanonical
bundle $K_P^{-1}$ of $P$ is {\it numerically effective} ({\it
nef\/}) then the invariants $\Iss^\al(\ga)^{\sst\La}$ of
\eq{an6eq10} transform according to \eq{an6eq16} under change of
Gieseker stability condition, even though \eq{an6eq9} is {\it not\/}
an algebra morphism. The basic idea is that for $K_P^{-1}$ nef we
use \eq{an6eq34} to force $\Ext^2(X,Y)=0$ for $X,Y$ satisfying some
conditions, so that \eq{an6eq33} reduces to \eq{an6eq8} for
these~$X,Y$.

The classification of algebraic surfaces, described for instance in
Iskovskikh and Shafarevich \cite{IsSh}, determines the possible
surfaces $P$ with $K_P^{-1}$ nef very explicitly. As $K_P^{-1}$ is
nef either (a) $K_P^n$ admits no sections for all $n>0$, or (b)
$K_P^n$ is trivial for some $n>0$. In case (a) $P$ has {\it Kodaira
dimension} $-\iy$, and so is either rational or ruled. If also
$K_P^{-1}$ is {\it ample} $P$ is called a {\it Del Pezzo surface}
(or a {\it Fano $2$-fold\/}), and is $\KP^2$, or $\KP^1\t\KP^1$, or
the blow-up of $\KP^2$ in $d$ points, $1\le d\le 9$. In case (b),
$P$ has Kodaira dimension 0, and is a $K3$ surface, an Enriques
surface, an abelian surface, or a bielliptic surface.

The next three results show how we will use $K_P^{-1}$ nef.

\begin{lem} In the situation above, with\/ $K_P^{-1}$ nef, for
all\/ $0\not\cong W\in\coh(P)$ we have $\ga([W])\le\ga([W\ot
K_P^{-1}])$ in~$G_2$.
\label{an6lem2}
\end{lem}

\begin{proof} Using Chern classes and the Riemann--Roch Theorem
as in \S\ref{an56} we compute the Hilbert polynomials of $W$ and
$W\ot K_P^{-1}$ with respect to $E$ as
\begin{small}
\begin{align*}
p_W(t)&=\bigl(\ha\rk(W)c_1(E)^2\bigr)t^2+
\bigl(c_1(W)c_1(E)+\ha\rk(W)c_1(E)(c_1(K_P)-c_1(E))\bigr)t\\
&+\bigl(\ts\frac{1}{12}\rk(W)(c_1(K_P)^2-2c_2(TP))+\ha
c_1(W)c_1(TP)+\ha c_1(W)^2-c_2(W)\bigr),\\
p_{W\!\ot\!K_P^{-1}}(t)\!&=p_W(t)-\bigl(\rk(W)c_1(E)c_1(K_P)
\bigr)t-\bigl(c_1(W)c_1(K_P)-\ha\rk(W)c_1(K_P)^2\bigr).
\end{align*}
\end{small}

As $E$ is ample and $K_P^{-1}$ nef we can show that
$c_1(E)c_1(K_P)\le 0$ with equality if and only if $c_1(K_P)=0$.
Also $\rk(W)\ge 0$, and if $\rk(W)=0$ then $c_1(W)c_1(K_P)\le 0$ as
$K_P^{-1}$ is nef. Using this we can show that $p_{W\ot
K_P^{-1}}(t)-p_W(t)$ has smaller degree than $p_W(t)$, and has
positive leading coefficient if it is nonzero. Since
$\ga([W]),\ga([W\ot K_P^{-1}])$ are these Hilbert polynomials
divided by their leading coefficients, it follows that
$\ga([W])\le\ga([W\ot K_P^{-1}])$ in the total order `$\le$' on
$G_2$ defined in~\cite[\S 4.4]{Joyc7}.
\end{proof}

\begin{prop} In the situation above, with\/ $K_P^{-1}$ nef, suppose
$X,Y\in\coh(P)$ are $\ga$-semistable with\/ $\ga([X])<\ga([Y])$.
Then~$\Ext^2(X,Y)=0$.
\label{an6prop4}
\end{prop}

\begin{proof} By \eq{an6eq34} it is enough to show that
$\Hom(Y,X\ot K_P)=0$. Suppose for a contradiction that $\phi:Y\ra
X\ot K_P$ is a nonzero morphism in $\A$, and let $W\in\A$ be the
image of $\phi$. Then $W\not\cong 0$ is a quotient object of $Y$, so
$\ga([Y])\le\ga([W])$ as $Y$ is $\ga$-semistable. Also $W$ is a
subobject of $X\ot K_P$, so $W\ot K_P^{-1}$ is a subobject of $X$,
giving $\ga([W\ot K_P^{-1}])\le\ga([X])$ as $X$ is $\ga$-semistable.
Putting this together with Lemma \ref{an6lem2} gives
$\ga([Y])\le\ga([W])\le\ga([W\ot K_P^{-1}])\le\ga([X])$, which
contradicts~$\ga([X])<\ga([Y])$.
\end{proof}

\begin{prop} In the situation above with\/ $K_P^{-1}$ nef, let
Assumption \ref{an2ass1} hold, and define $\Iss^\al(\ga)^{\sst\La}
\in\La$ for $\al\in C(\A)$ by $\Iss^\al(\ga)^{\sst\La}\!=\!\Up'
\ci\Pi_*\,\bdss^\al(\ga)$, as in \eq{an6eq10}. Suppose
$(\{1,\ldots,n\},\le,\ka)$ is $\A$-data with\/
$\ga\ci\ka(1)>\cdots>\ga\ci\ka(n)$. Then
\e
\Up'\ci\Pi_*\bigl(\bdss^{\ka(1)}(\ga)*\cdots*
\bdss^{\ka(n)}(\ga)\bigr)=\el^{-\sum_{1\le i<j\le
n}\chi(\ka(j),\ka(i))}\ts\prod_{i=1}^n\Iss^{\ka(i)}(\tau)^{\sst\La}.
\label{an6eq36}
\e
\label{an6prop5}
\end{prop}

\begin{proof} Let $[Y_i]\in\Oss^{\ka(i)}(\ga)$ for $i=1,\ldots,n$.
Then Proposition \ref{an6prop4} and $\ga\ci\ka(1)\!>\!\cdots\!>\!
\ga\ci\ka(n)$ give $\Ext^2(Y_j,Y_i)=0$ for $1\le i<j\le n$. Let
$[X_i]$ lie in the support of $\bdss^{\ka(1)}(\ga)*\cdots*
\bdss^{\ka(i)}(\ga)$. Using exact sequences and induction on $i$ we
can show that $\Ext^2(Y_j,X_i)=0$ whenever $1\le i<j\le n$. For
example, when $i=2$ there is an exact sequence $0\ra Y_1\ra X_2\ra
Y_2\ra 0$ for $[Y_1],[Y_2]$ in the supports of $\bdss^{\ka(1)}(\ga),
\bdss^{\ka(2)}(\ga)$. This induces an exact sequence
$\cdots\ra\Ext^2(Y_j,Y_1)\ra\Ext^2(Y_j,X_2)\ra\Ext^2(Y_j,Y_2)\ra 0$.
As $\Ext^2(Y_j,Y_1)=\Ext^2(Y_j,Y_2)=0$ from above this gives
$\Ext^2(Y_j,X_2)=0$. More generally $\Ext^2(Y_j,X_i)$ is built from
terms $\Ext^2(Y_j,Y_a)=0$ for $a=1,\ldots,i$, and so is zero.

Apply \cite[Cor.~5.15]{Joyc6} with $f=\bdss^{\ka(1)}(\ga)*\cdots
*\bdss^{\ka(j-1)}(\ga)$ and $g=\bdss^{\ka(j)}(\ga)$. This gives
expressions for $f\ot g$ and $f*g$, involving vector spaces
$E^0_m,E^1_m$ isomorphic to $\Hom(Y_j,X_{j-1}),\Ext^1(Y_j,X_{j-1})$
for $([X_{j-1}],[Y_j])$ in the support of $f\ot g$. As
$\Ext^2(Y_j,X_{j-1})=0$ and $[X_{j-1}]=\ka(\{1,\ldots,j-1\})$,
$[Y_j]=\ka(j)$ in $C(\A)$ for such $X,Y$ we see from \eq{an6eq33}
that $\dim E^1_m-\dim E^0_m=-\sum_{i=1}^{j-1}\chi(\ka(j),\ka(i))$.
Hence
\begin{align*}
&\Up'\ci\Pi_*\bigl(\bdss^{\ka(1)}(\ga)*\cdots*\bdss^{\ka(j)}(\ga)
\bigr)=\\
&\el^{-\sum_{i=1}^{j-1}\chi(\ka(j),\ka(i))}
\Up'\ci\Pi_*\bigl(\bdss^{\ka(1)}(\ga)*\cdots*\bdss^{\ka(j-1)}(\ga)
\bigr)\Up'\ci\Pi_*\bigl(\bdss^{\ka(j)}(\ga)\bigr),
\end{align*}
by the proof of \cite[Th.~6.1]{Joyc6}. Equation \eq{an6eq36} follows
by induction on~$j$.
\end{proof}

We can now prove our main result on invariants counting sheaves on
surfaces. First we give a sketch of why it is true, which may be
more helpful than the actual proof. In the situation of the theorem,
suppose $(\ti\ga,G_2,\le)$ dominates $(\ga,G_2,\le)$. Then Theorem
\ref{an5thm4} applies, and \eq{an5eq15} writes $\bdss^\al(\ti\ga)$
as a sum of $\bdss^{\ka(1)}(\ga)*\cdots*\bdss^{\ka(n)}(\ga)$ with
$\ga\ci\ka(1)>\cdots>\ga\ci\ka(n)$. So applying $\Up'\ci\Pi_*$ and
using Proposition \ref{an6prop5} writes $\Iss^\al(\ti\ga)^{\sst\La}$
as a sum of $\el^{\cdots}\prod_i\Iss^{\ka(i)}(\ga)^{\sst\La}$, a
special case of \eq{an6eq37} and an analogue of~\eq{an6eq28}.

We then invert these identities in $\La$ to write the
$\Iss^\al(\ga)^{\sst\La}$ in terms of the $\Iss^{\be}(\ti\ga)^{\sst
\La}$. Exchanging $\ga,\ti\ga$ proves another special case of
\eq{an6eq37}, when $(\ga,G_2,\le)$ dominates $(\ti\ga,G_2,\le)$, an
analogue of \eq{an6eq29}. It is important that this inversion is
done in $\La$ rather than trying to apply $\Up'\ci\Pi_*$ to
\eq{an5eq20}, as in \eq{an5eq20} we do not have
$\ga\ci\ka(1)>\cdots>\ga\ci\ka(n)$ and so cannot use
Proposition~\ref{an6prop5}.

For the general case, suppose we could find a permissible weak
stability condition $(\hat\ga,G_2,\le)$ dominating both
$(\ga,G_2,\le)$ and $(\ti\ga,G_2,\le)$. Then the two special cases
above allow us to write $\Iss^\al(\ti\ga)^{\sst\La}$ in terms of the
$\Iss^\be(\hat\ga)^{\sst\La}$ and $\Iss^\be(\hat\ga)^{\sst\La}$ in
terms of the $\Iss^{\ka(i)}(\ga)^{\sst\La}$, and substituting one in
the other yields~\eq{an6eq37}.

The problem in proving Theorem \ref{an6thm4} is that we have no
suitable $(\hat\ga,G_2,\le)$. The purity weak stability condition
$(\de,D_2,\le)$ will not do as it is not permissible, nor even
essentially permissible, so the $\Iss^\be(\de)^{\sst\La}$ are
undefined. So instead we introduce a 1-parameter family of stability
conditions $(\ga_s,G_2,\le)$ for $s\in[0,1]$ with $\ga_0=\ga$ and
$\ga_1=\ti\ga$, in a similar way to the proof of Theorem
\ref{an5thm7}. There are finitely many `walls' in $[0,1]$ where the
expression for $\bdss^\al(\ga_s)$ in terms of the $\bdss^\be(\ga)$
changes. Jumping onto a wall $s$ from a nearby point $s'$ is like
transforming to a dominant weak stability condition, so we can use
the ideas above.

\begin{thm} Suppose $\K$ is an algebraically closed field and\/
$P$ a smooth projective surface over $\K$ with\/ $K_P^{-1}$
numerically effective. Take $\A=\coh(P)$ with data $K(\A),\fF_\A$
satisfying Assumption \ref{an3ass} as in {\rm\cite[Ex.~9.1]{Joyc5},}
with biadditive $\chi:K(\A)\t K(\A)\ra\Z$ satisfying \eq{an6eq33}
for all\/ $X,Y\in\A$. Let\/ $(\ga,G_2,\le),(\ti\ga,G_2,\le)$ be
Gieseker stability conditions on $\A$ defined using ample line
bundles $E,\ti E$ on $P$ as in {\rm\cite[Ex.~4.16]{Joyc7}}. Suppose
Assumption \ref{an2ass1} holds, and define invariants
$\Iss^\al(\ga)^{\sst \La},\Iss^\al(\ti\ga)^{\sst\La}$ for $\al\in
C(\A)$ as in \eq{an6eq10}. Then for all\/ $\al\in C(\A)$ the
following holds, with only finitely many nonzero terms:
\e
\begin{gathered}
\Iss^\al(\ti\ga)^{\sst\La}=
\sum_{\begin{subarray}{l}
\text{$\A$-data $(\{1,\ldots,n\},\le,\ka):$}\\
\text{$\ka(\{1,\ldots,n\})=\al$}\end{subarray}
\!\!\!\!\!\!\!\!\!\!\!\!\!\!\!\!\!\!\!\!\!\!\!\!
\!\!\!\!\!\!\!\!\!\!\!\!\!\!\!\!\! }
\begin{aligned}[t]
S(\{1,\ldots,n\},\le,\ka,\ga,\ti\ga)\cdot
\el^{-\sum_{1\le i<j\le n}\chi(\ka(j),\ka(i))}\cdot&\\
\ts\prod_{i=1}^n\Iss^{\ka(i)}(\ga)^{\sst\La}&.
\end{aligned}
\end{gathered}
\label{an6eq37}
\e
\label{an6thm4}
\end{thm}

\begin{proof} For $X\in\coh(P)$ define the {\it joint Hilbert
polynomial\/} $q_X$ w.r.t.\ $E,\ti E$ by
\begin{equation*}
q_X(k,l)=\ts\sum_{i=0}^{\dim P}(-1)^i\dim_\K H^i\bigl(P,X\ot E^k\ot
\ti E^l\bigr) \quad\text{for $k,l\in\Z$.}
\end{equation*}
As for conventional Hilbert polynomials, one can show that
\begin{equation*}
q_X(k,l)=\ts\sum_{a,b=0}^{\dim P}c_{a,b}k^al^b/a!b! \quad\text{for
$c_{a,b}\in\Z$,}
\end{equation*}
where $c_{a,b}=0$ if $a+b>\dim X$ and $c_{a,b}>0$ if $a+b=\dim X$.
Also $q_X$ is additive in $[X]\in K(\A)$, so there is a unique group
homomorphism $\Pi_{E,\ti E}:K(\A)\ra\Q[t,u]$ with $\Pi_{E,\ti
E}:[X]\mapsto q_X(t,u)$ for all~$X\in\coh(P)$.

Define $\ga_s:C(\A)\ra G_2$ by $\ga_s(\al)= L_{\al,s}^{-1}\Pi_{E,\ti
E}\bigl((1-s)t,st\bigr)$ for each $s\in[0,1]$, where $L_{\al,s}$ is
the leading coefficient of $\Pi_{E,\ti E}\bigl((1-s)t,st\bigr)$,
which is a polynomial in $t$ with positive leading coefficient. Then
$\ga_0=\ga$ and $\ga_1=\ti\ga$, and the proof in \cite[\S
4.4]{Joyc7} shows that $(\ga_s,G_2,\le)$ is a stability condition on
$\A$. As $q_X(t,0),q_X(0,t)$ are the Hilbert polynomials of $X$
w.r.t.\ $E,\ti E$ respectively, we see that $\ga_0=\ga$ and
$\ga_1=\ti\ga$, so the $\ga_s$ for $s\in[0,1]$ interpolate between
$\ga$ and $\ti\ga$. When $s=p/(p+q)$ for integers $p,q\ge 0$ with
$p+q>0$ it is equivalent (after rescaling $t$) to Gieseker stability
with respect to the ample line bundle $E^q\ot\ti E^p$. Thus
$(\ga_s,G_2,\le)$ is permissible for all $s\in[0,1]\cap\Q$
by~\cite[Th.~4.20]{Joyc7}.

Each $\al\in C(\A)$ has a Hilbert polynomial $p_\al(t)=a_dt^d/d!
+\cdots+a_0$ w.r.t.\ $E$, with $a_i\in\Z$ and $a_d>0$, where
$d=\dim\al$. We shall prove \eq{an6eq37} by induction on $a_d$. Here
is our inductive hypothesis, for~$r\ge 1$:
\begin{itemize}
\setlength{\itemsep}{0pt}
\setlength{\parsep}{0pt}
\item[$(*\kern .05em{}_r)$] Suppose that for all $\al\in C(\A)$ whose
Hilbert polynomial $p_\al(t)=a_dt^d/d!+\cdots+a_0$ w.r.t.\ $E$ has
$0<a_d\le r$, equation \eq{an6eq37} holds with $\ga_s,\ga_{s'}$ in
place of $\ga,\ti\ga$, with only finitely many nonzero terms, for
all~$s,s'\in[0,1]\cap\Q$.
\end{itemize}
From Corollary \ref{an5cor3} and Theorems \ref{an5thm1} and
\ref{an5thm3} we see that \eq{an5eq3} holds in $\SF(\fObj_\A)$ with
$\ga_s,\ga_{s'}$ for $s,s'\in[0,1]\cap\Q$ in place of $\tau,\ti\tau$
with finitely many nonzero terms. This implies there are only
finitely many nonzero terms in \eq{an6eq37} with $\ga_s,\ga_{s'}$ in
place of $\ga,\ti\ga$, proving part of~$(*\kern .05em{}_r)$.

Let $\al,d,a_d,s,s'$ be as in $(*\kern .05em{}_r)$, and $n,\ka$ as
in \eq{an5eq3} or \eq{an6eq37} with $S(\{1,\ldots,n\},\ab\le,\ab
\ka,\ab\ga_s,\ga_{s'})\ne 0$. Then the Hilbert polynomials
$p_{\ka(i)}(t)$ are also of the form $a_d^it^d/d!+\cdots+a_0^i$ for
$a_d^i\ge 1$, and sum to $p_\al(t)$. Thus $\sum_{i=1}^na_d^i=a_d\le
r$, which forces $n\le r$. In particular, when $r=1$ the only
nonzero terms in \eq{an5eq3} and \eq{an6eq37} are $n=1$ and
$\ka(1)=\al$. Thus \eq{an5eq3} reduces to
$\bdss^\al(\ga_{s'})=\bdss^\al(\ga_s)$, so
$\Iss^\al(\ga_{s'})^{\sst\La}=\Iss^\al(\ga_s)^{\sst\La}$, which is
\eq{an6eq37}. This proves $(*_1)$, giving the first step.

Suppose by induction that $(*\kern .05em{}_r)$ holds for some $r\ge
1$, and let $\al\in C(\A)$ with $p_\al(t)=a_dt^d/d!+\cdots+a_0$ and
$a_d=r+1$. Using the methods of \S\ref{an56} we can show that there
are only {\it finitely many} sets of $\A$-data $(\{1,\ldots,n\},\le,
\ka)$ with $\ka(\{1,\ldots,n\})=\al$, such that for some
$s,s'\in[0,1]\cap\Q$ we have $S(\{1,\ldots,n\},\ab\le,\ab\ka,
\ab\ga_s,\ab \ga_{s'})\ne 0$ and $\bdss^{\ka(i)}(\ga_s)\ne 0$
for~$i=1,\ldots,n$.

Let $U$ be the finite set of all such pairs $(n,\ka)$. Let $V$ be
the finite set of elements of $C(\A)$ of the form $\ka(i)$ or
$\ka(\{1,\ldots,i\})$ or $\ka(\{i,\ldots,n\})$ for $(n,\ka)\in U$
and $1\le i\le n$. Let $W$ be the set of $w\in[0,1]$ such that for
some $v,v'\in V$ we have $\ga_w(v)=\ga_w(v')$ but $\ga_s(v)\ne\ga_s
(v')$ for generic $s\in[0,1]$. In fact $\ga_w(v)=\ga_w(v')$ is
equivalent to $aw+b=0$ for $a,b\in\Z$, so there is at most one $w$
for each pair $v,v'$, which is rational. Therefore $W$ is a finite
subset of~$[0,1]\cap\Q$.

The point of this is that in \eq{an6eq37} with $\ga_s,\ga_{s'}$ in
place of $\ga,\ti\ga$ only $(n,\ka)$ in $U$ can give nonzero terms,
and the $S(\{1,\ldots,n\},\le,\ka,\ga_s,\ga_{s'})$ depend only on
whether or not $\ga_s(v)\le\ga_s(v')$ and
$\ga_{s'}(v)\le\ga_{s'}(v')$ hold for pairs $v,v'$ in $V$, by
Definition \ref{an4def1}. But these inequalities only change when
$s$ or $s'$ pass through a point of $W$. Thus, we have a finite set
of `walls' $W$ such that $S(\{1,\ldots,n\},\le,\ka,\ga_s,\ga_{s'})$
is locally constant for $s,s'\in[0,1]\sm W$, for all $n,\ka$ that
could contribute nonzero terms in~\eq{an6eq37}.

Order the elements of $W\cup\{0,1\}$ as $0=w_0<w_1<\cdots<w_{k-1}
<w_k=1$. We divide the remainder of the proof into three steps:
\begin{list}{}{
\setlength{\itemsep}{1pt}
\setlength{\parsep}{1pt}
\setlength{\labelwidth}{40pt}
\setlength{\leftmargin}{40pt}
}
\item[{\bf Step 1.}] Show \eq{an6eq37} holds with $\ga_s,\ga_{s'}$ in
place of $\ga,\ti\ga$ when $s'=w_i$ and $s\in(w_{i-1},w_i)\cap\Q$ or
$s\in(w_i,w_{i+1})\cap\Q$;
\item[{\bf Step 2.}] Show \eq{an6eq37} holds with $\ga_s,\ga_{s'}$ in
place of $\ga,\ti\ga$ when $s=w_i$ and $s'\in(w_{i-1},w_i)\cap\Q$ or
$s'\in(w_i,w_{i+1})\cap\Q$;
\item[{\bf Step 3.}] Show \eq{an6eq37} holds with $\ga_s,\ga_{s'}$ in
place of $\ga,\ti\ga$ for any~$s,s'\in[0,1]\cap\Q$.
\end{list}

\noindent{\bf Step 1.} Since $s'\in W$ and no elements of $W$ lie
between $s$ and $s'$, we see by definition of $W$ that $\ga_s(v)\le
\ga_s(v')$ implies $\ga_{s'}(v)\le\ga_{s'}(v')$ for all pairs $v,v'$
in $V$. If $(n,\ka)\in U$ we see using the proof of \eq{an5eq13}
that $S(\{1,\ldots,n\},\le,\ka,\ga_s,\ga_{s'})$ is 1 if
$\ga_{s'}\ci\ka\equiv\ga_{s'}(\al)$ and
$\ga_s\ci\ka(1)>\cdots>\ga_s\ci\ka(n)$ and zero otherwise. Applying
$\Up'\ci\Pi_*$ to \eq{an5eq3} with $\ga_s,\ga_{s'}$ in place of
$\tau,\ti\tau$ and using Proposition \ref{an6prop5} then gives
\eq{an6eq37} with $\ga_s,\ga_{s'}$ in place of $\ga,\ti\ga$.
\medskip

\noindent{\bf Step 2.} For $s,s'$ as in Step 2, in Step 1 we proved
\eq{an6eq37} with $\ga_{s'},\ga_s$ in place of $\ga,\ti\ga$. The
only term on the right with $n=1$ is
$\Iss^{\al}(\ga_{s'})^{\sst\La}$, so we may rewrite it as
\begin{equation*}
\Iss^\al(\ga_{s'})^{\sst\La}\!=\!\Iss^\al(\ga_s)^{\sst\La}-\!\!\!\!\!
\sum_{\begin{subarray}{l}
\text{$\A$-data $(\{1,\ldots,n\},\le,\ka):$}\\
\text{$n\ge 2$, $\ka(\{1,\ldots,n\})=\al$}\end{subarray}
\!\!\!\!\!\!\!\!\!\!\!\!\!\!\!\!\!\!\!\!\!\!\!\!
\!\!\!\!\!\!\!\!\!\!\!\!\!\!\!\!\!\!\!\!\!\! }
\begin{aligned}[t]
S(\{1,\ldots,n\},\le,\ka,\ga_{s'},\ga_s)\cdot
\el^{-\sum_{1\le i<j\le n}\chi(\ka(j),\ka(i))}\cdot&\\
\ts\prod_{i=1}^n\Iss^{\ka(i)}(\ga_{s'})^{\sst\La}.&
\end{aligned}
\end{equation*}
For $n,\ka$ on the right hand side we have Hilbert polynomials
$p_{\ka(i)}=a_d^it^d/d!+\cdots+a_0^i$ for $i=1,\ldots,n$ which sum
to $p_\al$, so $\sum_{i=1}^na_d^i=a_d=r+1$ by choice of $\al$. As
$n\ge 2$ and $a_d^i\ge 1$ we have $0<a_d^i\le r$ so the inductive
hypothesis $(*\kern .05em{}_r)$ applies with $\ka(i)$ in place of
$\al$, yielding expressions for $\Iss^{\ka(i)}(\ga_{s'})^{\sst\La}$
in terms of the $\Iss^\be(\ga_s)^{\sst\La}$. Substituting these into
the right hand side above gives an expression for
$\Iss^\al(\ga_{s'})^{\sst\La}$ in terms of the
$\Iss^\be(\ga_s)^{\sst\La}$ and powers of $\el$. But this is exactly
\eq{an6eq37} with $\ga_s,\ga_{s'}$ in place of $\ga,\ti\ga$, since
\eq{an6eq37} with $\ga,\ti\ga$ exchanged is the combinatorial
inverse of \eq{an6eq37} by properties of the transformation
coefficients~$S(\cdots)$.
\medskip

\noindent{\bf Step 3.} Let $s,s'\in[0,1]\cap\Q$, and suppose for
simplicity that $s'<s$. Then there exist unique $1\le i\le j\le k$
with $w_{i-1}\le s'<w_i$ and $w_{j-1}<s\le w_j$. Choose
$x_i,\ldots,x_j$ in $\bigl([0,1]\cap\Q\bigr)\sm W$ with $s'\le
x_i<w_i<x_{i+1}<w_{i+2}<\cdots<w_{j-1}<x_j\le s$, where we take
$x_i=s'$ if $w_{i-1}<s'$ and $x_i>s'$ if $w_{i-1}=s'$, and $x_j=s$
if $s<w_j$ and $x_j<s$ if~$s=w_j$.

Then we write $\Iss^\al(\ga_{s'})^{\sst\La}$ in terms of the
$\Iss^\be(\ga_s)^{\sst\La}$ by $2(j-i+1)$ substitutions, as follows.
For the first substitution, if $x_i<s'$ then $s'=w_{i-1}$, and Step
2 writes $\Iss^\al(\ga_{s'})^{\sst\La}$ in terms of the
$\Iss^\be(\ga_{x_i})$, where either $\be=\al$ or
$p_\be(t)=b_dt^d/d!+\cdots+b_0$ with $0<b_d\le r$, so that $(*\kern
.05em{}_r)$ applies with $\be$ in place of $\al$. If $x_i=s'$ then
$\Iss^\al(\ga_{s'})^{\sst\La}=\Iss^\al(\ga_{x_i})^{\sst\La}$ and the
first substitution is trivial.

For substitution number $2k$ for $k=1,\ldots,j-i$ we have already
written $\Iss^\al(\ga_{s'})^{\sst\La}$ in terms of
$\Iss^\be(\ga_{x_{i+k-1}})^{\sst\La}$ for $\be=\al$ or $\be$ to
which $(*\kern .05em{}_r)$ applies. We then use Step 1 for $\be=\al$
and $(*\kern .05em{}_r)$ otherwise to write
$\Iss^\be(\ga_{x_{i+k-1}})^{\sst\La}$ in terms of the
$\Iss^{\be'}(\ga_{w_{i+k-1}})^{\sst\La}$ for $\be'=\al$ or $\be'$ to
which $(*\kern .05em{}_r)$ applies, and substitute this into the
previous expression to write $\Iss^\al(\ga_{s'})^{\sst\La}$ in terms
of these $\Iss^{\be'}(\ga_{w_{i+k-1}})^{\sst\La}$. For substitution
number $2k+1$ for $k=1,\ldots,j-i$ we use Step 2 for $\be'=\al$ and
$(*\kern .05em{}_r)$ otherwise to write these
$\Iss^{\be'}(\ga_{w_{i+k-1}})^{\sst\La}$ in terms of
$\Iss^{\be''}(\ga_{x_{i+k}})^{\sst\La}$ for $\be''=\al$ or $\be''$
to which $(*\kern .05em{}_r)$ applies. Finally, for substitution
number $2(j-i+1)$ we use Step 1 and $(*\kern .05em{}_r)$ if $x_j<s$
and do nothing if~$x_j=s$.

As all we are doing at each stage is substituting in finitely many
copies of \eq{an6eq37} with different values for $\al,\ga,\ti\ga$,
by repeated use of \eq{an4eq7} we see that the expressions we get
for $\Iss^\al(\ga_{s'})^{\sst\La}$ in terms of
$\Iss^\be(\ga_t)^{\sst\La}$ for $t=x_i,w_i,x_{i+1},\ldots,x_j,s$
respectively are just \eq{an6eq37} with $\ga_t,\ga_{s'}$ in place of
$\ga,\ti\ga$. So at the last substitution we prove \eq{an6eq37} with
$\ga_s,\ga_{s'}$ in place of $\ga,\ti\ga$, as we want. The case
$s<s'$ is similar, and $s=s'$ is trivial.
\medskip

In Step 3 we have proved $(*\kern .05em{}_{r+1})$ supposing $(*\kern
.05em{}_r)$. Thus by induction $(*\kern .05em{}_r)$ holds for all
$r\ge 1$. The theorem follows by setting $s=0$ and~$s'=1$.
\end{proof}

An interpretation of the theorem will be proposed after Problem
\ref{an7prob3} below. When $P$ is a ruled surface with $K_P^{-1}$
nef, Yoshioka \cite[Cor.~3.3]{Yosh2} proves a result related to
Theorem \ref{an6thm4}, a wall-crossing formula for Poincar\'e
polynomials of moduli spaces $\Oss^\al(\ga)$. It is valid only when
$\ga,\ti\ga$ are in the interiors of chambers separated by a single
wall, in contrast to \eq{an6eq37} which holds for
arbitrary~$\ga,\ti\ga$.

We now define invariants $\bar J^\al(\ga)^{\sst\La}$ related to the
$J^\al(\ga)^{\sst\La^\ci}$ of~\eq{an6eq11}.

\begin{dfn} In the situation above, motivated by \eq{an6eq13}
define invariants $\bar J^\al(\ga)^{\sst\La}\in\La$ for all $\al\in
C(\A)$ by
\e
\bar J^\al(\ga)^{\sst\La}=
\sum_{\begin{subarray}{l}\text{$\A$-data $(\{1,\ldots,n\},\le,\ka):$}\\
\text{$\ka(\{1,\ldots,n\})=\al$, $\ga\ci\ka\equiv\ga(\al)$}
\end{subarray}
\!\!\!\!\!\!\!\!\!\!\!\!\!\!\!\!\!\!\!\!\!\!\!\!\!\!\!\!\!\!\!\!
\!\!\!\!\!\!\!\!\!\!\!\!\!\!\!\!\!\!\!\!\!\!\! }
\el^{-\sum_{1\le i<j\le n}\chi(\ka(j),\ka(i))}
\frac{(-1)^{n-1}(\el-1)}{n}\,
\prod_{i=1}^n\Iss^{\ka(i)}(\ga)^{\sst\La}.
\label{an6eq38}
\e
Using the proof of \cite[Th.s 6.4 \& 7.7]{Joyc7} in the algebra
$A(\A,\La,\chi)$ rather than $\SF(\fObj_\A)$ we deduce an analogue
of~\eq{an6eq14}:
\e
\Iss^\al(\tau)^{\sst\La}=
\sum_{\begin{subarray}{l}\text{$\A$-data $(\{1,\ldots,n\},\le,\ka):$}\\
\text{$\ka(\{1,\ldots,n\})=\al$, $\tau\ci\ka\equiv\tau(\al)$}
\end{subarray}
\!\!\!\!\!\!\!\!\!\!\!\!\!\!\!\!\!\!\!\!\!\!\!\!\!\!\!\!\!\!\!\!
\!\!\!\!\!\!\!\!\!\!\!\!\!\!\!\!\!\!\!\!\!\!\! }
\el^{-\sum_{1\le i<j\le n}\chi(\ka(j),\ka(i))}\frac{(\el-1)^{-n}}{n!}\,
\prod_{i=1}^n\bar J^{\ka(i)}(\tau)^{\sst\La}.
\label{an6eq39}
\e
There are finitely many nonzero terms in both equations, as for
\eq{an3eq14} and~\eq{an3eq16}.
\label{an6def6}
\end{dfn}

If \eq{an6eq8} held in $\A=\coh(P)$ then equation \eq{an6eq13} of
Theorem \ref{an6thm2} would give $\bar J^\al(\ga)^{\sst\La}=
J^\al(\ga)^{\sst\La^\ci}$, so that $\bar J^\al(\ga)^{\sst\La}
\in\La^\ci$. We can prove this in a special case, but as \eq{an6eq8}
does not hold, in general it is likely that $J^\al(\ga)^{\sst
\La^\ci}\ne\bar J^\al(\ga)^{\sst\La}\notin\La^\ci$. Thus we do not
define invariants $\bar J^\al(\ga)^{\sst\Om}$, as projecting to
$\Om$ may not be possible.

\begin{prop} In the situation above, suppose $K_P^{-1}$ is ample
and\/ $\al\in C(\coh(P))$ with\/ $\dim\al>0$. Then $\bar
J^\al(\ga)^{\sst\La}=J^\al(\ga)^{\sst\La^\ci}$ and lies
in~$\La^\ci$.
\label{an6prop6}
\end{prop}

\begin{proof} If $K_P^{-1}$ is ample and $0\not\cong W\in\coh(P)$
with $\dim W>0$, so that $\rk(W)\ne 0$ or $c_1(W)\ne 0$, the proof
of Lemma \ref{an6lem2} also shows that $\ga([W])<\ga([W\ot
K_P^{-1}])$ in $G_2$. Then we can modify Proposition \ref{an6prop4}
to show that if $X,Y\in\coh(P)$ are $\ga$-semistable with $\dim
X,\dim Y>0$ and $\ga([X])=\ga([Y])$ then $\Ext^2(X,Y)=0$, and
Proposition \ref{an6prop5} to show that if $(\{1,\ldots,n\},\le,
\ka)$ is $\A$-data with $\ka(\{1,\ldots,n\})=\al$ and $\ga\ci\ka
\equiv\ga(\al)$ and $\dim\ka(i)>0$ then \eq{an6eq36} holds. So
applying $\Up'\ci\Pi_*$ to \eq{an3eq14} gives \eq{an6eq13}, and the
proposition follows.
\end{proof}

Combining \eq{an6eq37}--\eq{an6eq39} and the definition \eq{an4eq5}
of the coefficients $U(\cdots)$ we see that the $\bar J^\al(\ga)^{
\sst\La}$ transform according to \eq{an6eq17} under change of
Gieseker stability condition. Then the proof of Corollary
\ref{an6cor1} shows that $\bar J^\al(\ga)^{\sst\La}$ is independent
of $\ga$ if $\chi$ is symmetric. Now using Chern classes and the
Riemann--Roch theorem as in Hartshorne \cite[App.~A]{Hart} we find
that for all $X,Y\in\A$ we have
\begin{equation*}
\chi\bigl([X],[Y]\bigr)-\chi\bigl([Y],[X]\bigr)=
\bigl(\rk(X)c_1(Y)-\rk(Y)c_1(X)\bigr)c_1(K_P).
\end{equation*}
Thus $\chi$ is symmetric if $c_1(K_P)=0$ in $H^2(P,\Q)$ or
$H^2(P,\Q_l)$, though not necessarily in $H^2(P,\Z)$. This holds if
$K_P^n$ is trivial for some $n>0$. So from the classification of
surfaces \cite{IsSh} we deduce:

\begin{thm} In the situation of Theorem \ref{an6thm4}, for all\/
$\al\in C(\A)$ we have
\e
\begin{gathered}
\bar J^\al(\ti\ga)^{\sst\La}=
\sum_{\begin{subarray}{l}
\text{$\A$-data $(\{1,\ldots,n\},\le,\ka):$}\\
\text{$\ka(\{1,\ldots,n\})=\al$}\end{subarray}
\!\!\!\!\!\!\!\!\!\!\!\!\!\!\!\!\!\!\!\!\!\!\!\!
\!\!\!\!\!\!\!\!\!\!\!\!\!\!\!\!\! }
\begin{aligned}[t]
U(\{1,\ldots,n\},\le,\ka,\tau,\ti\tau)\cdot
\el^{-\sum_{1\le i<j\le n}\chi(\ka(j),\ka(i))}\cdot&\\
(\el-1)^{1-n}\ts\prod_{i=1}^n\bar J^{\ka(i)}(\ga)^{\sst\La}&,
\end{aligned}
\end{gathered}
\label{an6eq40}
\e
with only finitely many nonzero terms. If\/ $c_1(K_P)=0$ then the
$\bar J^\al(\ga)^{\sst\La}$ are independent of the choice of
Gieseker stability condition $(\ga,G_2,\le)$ on $\A=\coh(P)$, that
is, independent of the ample line bundle $E$. This holds if\/ $P$ is
a $K3$ surface, an Enriques surface, an abelian surface, or a
bielliptic surface.
\label{an6thm5}
\end{thm}

Yoshioka \cite[Rem.~3.2]{Yosh2} proves a related result, that if $P$
is a $K3$ or abelian surface for which the Bogomolov--Gieseker
inequality holds and $(\ga,G_2,\le)$ is a suitably generic Gieseker
stability condition on $\coh(P)$, certain weighted counts of
$\ga$-semistable sheaves in class $\al$ over finite fields are
independent of~$\ga$.

We now propose a conjectural means to compute many of the invariants
$\bar J^\al(\ga)^{\sst\La}$ when $c_1(K_P)=0$, motivated by ideas of
Bridgeland \cite{Brid1,Brid2,Brid3}. It involves the extension of
our whole programme to {\it triangulated categories}, and the {\it
bounded derived category} $D^b(\coh(P))$ of coherent sheaves on $P$.
This extension will be discussed in \S\ref{an7}. Suppose for the
moment that:
\begin{itemize}
\setlength{\itemsep}{0pt}
\setlength{\parsep}{0pt}
\item There is a good notion of permissible stability condition $\tau$
on $D^b(\coh(P))$, based on Bridgeland stability \cite{Brid1}, which
includes the extension of Gieseker stability on $\coh(P)$ to
$D^b(\coh(P))$, perhaps as a limit. These stability conditions form
a moduli space, with a topology. Write $\Stab(P)$ for the connected
component of the moduli space including Gieseker stability.
\item One can define invariants $\hat J^\al(\tau)^{\sst\La}\in\La$
`counting' $\tau$-semistable objects in class $\al\in K(\A)$. They
transform according to a generalization of \eq{an6eq40} under change
of stability conditions, at least for `nearby' stability conditions.
When $\tau$ is the extension of $\ga$ to $D^b(\coh(P))$ and $\al\in
C(\A)$ we have~$\hat J^\al(\tau)^{\sst\La}=\bar J^\al(\ga)^{\sst
\La}$.
\item When $c_1(K_P)=0$ the $\hat J^\al(\tau)^{\sst\La}$ are
independent of choice of $\tau$ in the connected component
$\Stab(P)$.
\item The whole framework is preserved by autoequivalences
$\Phi:D^b(\coh(P))\ra D^b(\coh(P))$ of the derived category. Write
$\Aut^+\bigl(D^b(\coh(P))\bigr)$ for the group of autoequivalences
preserving the connected component~$\Stab(P)$.
\end{itemize}

Let $c_1(P)=0$, $\tau\in\Stab(P)$ and $\Phi\in\Aut^+\bigl(D^b(\coh
(P))\bigr)$. Then for $\al$ in $K(\A)$ we have $\hat J^\al(\tau)^{
\sst\La}=\hat J^{\Phi_*(\al)}(\Phi_*(\tau))^{\sst\La}=\hat
J^{\Phi_*(\al)}(\tau)^{\sst \La}$, so the $\hat J^\al(\tau)^{\sst
\La}$ are unchanged by the action of $\Phi_*$ on $K(\A)$. This
suggests the following:

\begin{conj} Let Assumption \ref{an2ass1} hold and\/ $P$ be a
smooth projective surface over $\K$ with\/ $c_1(K_P)=0$, and define
$\A=\coh(P),K(\A)$ and invariants $\bar J^\al(\ga)^{\sst\La}$ as
above. Then there exist $\hat J^\al\in\La$ for $\al\in K(\A)$
satisfying:
\begin{itemize}
\setlength{\itemsep}{0pt}
\setlength{\parsep}{0pt}
\item[{\rm(a)}] If\/ $\al\in C(\A)$ then~$\bar J^\al(\ga)^{\sst
\La}=\hat J^\al$.
\item[{\rm(b)}] If\/ $\Phi\in\Aut^+\bigl(D^b(\coh(P))\bigr)$ then
$\hat J^{\Phi_*(\al)}=\hat J^\al$ for all\/~$\al\in K(\A)$.
\end{itemize}
\label{an6conj1}
\end{conj}

The conjecture could be applied in the following way. We first
compute the invariants $\bar J^\al(\ga)^{\sst\La}$ for some small
subset $S$ of $\al\in C(\A)$ for which the moduli spaces
$\Oss^\al(\ga)$ can be explicitly understood. For example, if
$\rk(\al)=1$ and $\chi(\al,\al)=\chi(\O,\O)-n$ for $n\in\Z$ one can
show that $\Oss^\al(\ga)$ is empty if $n<0$ and otherwise is
1-isomorphic to ${\rm Jac}(P)\t{\rm Hilb}^n(P)\t[\Spec\K/\K^\t]$,
where ${\rm Jac}(P)$ is the {\it Jacobian variety} of line bundles
$L$ on $P$ with $[L]=[\O]\in K(\A)$, and ${\rm Hilb}^n(P)$ is the
{\it Hilbert scheme} of $n$ points on $P$. So $\bar
J^\al(\ga)^{\sst\La}=0$ for $n<0$ and
\begin{equation*}
\bar J^\al(\ga)^{\sst\La}=(\el-1)\Up'\bigl(\bigl[\Oss^\al(\ga)
\bigr]\bigr)=\Up\bigl(\bigl[{\rm Jac}(P)\bigr]\bigr)
\Up\bigl(\bigl[{\rm Hilb}^n(P)\bigr]\bigr) \quad\text{for $n\ge 0$.}
\end{equation*}

This gives $\hat J^\al$ for $\al\in S$, so (b) determines $\hat
J^\al$ for $\al\in\Aut^+\bigl(D^b(\coh(P))\bigr)\cdot S$, and then
(a) gives $\bar J^\al(\ga)^{\sst\La}$ for $\al\in\bigl(\Aut^+(D^b
(\coh(P)))\cdot S\big)\cap C(\A)$. Thus, provided we understand the
action of $\Aut^+\bigl(D^b(\coh(P))\bigr)$ on $K(\A)$ reasonably
well, we can compute the invariants $\bar J^\al(\ga)^{\sst\La}$ on a
much larger subset of $C(\A)$, perhaps even the whole of it.

Let $P$ be an algebraic $K3$ surface over $\C$. Using his notion of
stability condition on triangulated categories \cite{Brid1},
Bridgeland \cite{Brid2} (surveyed in \cite[\S 6]{Brid3})
parametrizes a connected component $\Stab(P)$ of the moduli space of
stability conditions on $D^b(\coh(P))$. Bridgeland's definition does
not include the extension of Gieseker stability on $\coh(P)$ to
$D^b(\coh(P))$, but can probably be generalized so that it does,
perhaps following Gorodentscev et al.~\cite{GKR}.

Using results of Orlov \cite{Orlo} on autoequivalences of
$D^b(\coh(P))$, Bridgeland \cite[Conj.~1.2]{Brid2} conjecturally
describes $\Aut^+(D^b(\coh(P)))$. His description implies that
$\Aut^+(D^b(\coh(P)))$ acts on $H^*(P,\Z)$ as a certain index 2
subgroup of the group of automorphisms of $H^*(P,\Z)$ preserving the
Mukai form and the subspace $H^{2,0}(P)\subset H^2(P,\Z)\ot_\Z\C$.
Combining this with Conjecture \ref{an6conj1} should make an
effective tool for conjecturally computing invariants $\bar
J^\al(\ga)^{\sst\La}$ when $P$ is a complex algebraic $K3$ surface,
which might lead to new formulae for Poincar\'e and Hodge
polynomials for moduli spaces of $\ga$-semistable sheaves on~$P$.

As evidence for Conjecture \ref{an6conj1} we note some results of
Yoshioka \cite{Yosh3}. Let $P$ be an abelian surface, $\al\in C(\A)$
{\it primitive} with $\chi(\al,\al)\le -2$ and $(\ga,G_2,\le)$ a
Gieseker stability condition with $\Ost^\al(\ga)=\Oss^\al(\ga)$, so
that $\Oss^\al(\ga)\cong\Mss^\al(\ga)\t[\Spec\K/\K^\t]$ for a
nonsingular projective symplectic variety $\Mss^\al(\ga)$ of
dimension $2-\chi(\al,\al)$. Then \cite[Th.~0.1]{Yosh3} shows
$\Mss^\al(\ga)$ is deformation equivalent to ${\rm Jac}(P)\t{\rm
Hilb}^{-\chi(\al, \al)/2}(P)$. If instead $P$ is a $K3$ surface then
\cite[Th.~8.1]{Yosh3} shows $\Mss^\al(\ga)$ is deformation
equivalent to~${\rm Hilb}^{1-\chi(\al,\al)/2}(P)$.

Now if our choice of $\Up$ in Assumption \ref{an2ass1} is unchanged
by deformations of smooth projective $\K$-varieties, which holds for
virtual Poincar\'e polynomials and virtual Hodge polynomials when
$\K=\C$, these results allow us to evaluate $\bar J^\al(\ga)^{\sst
\La}$ for most primitive $\al\in C(\A)$. The answer depends only on
$P$ and $\chi(\al,\al)$, and so is invariant under
$\Aut^+\bigl(D^b(\coh(P))\bigr)$ as in Conjecture~\ref{an6conj1}(b).

We quote Bridgeland~\cite[\S 6]{Brid3}:
\begin{quotation}
`As a final remark in this section note that Borcherd's work on
modular forms \cite{Borc} allows one to write down product
expansions for holomorphic functions on $\Stab(P)$ that are
invariant under the group $\Aut(D^b(\coh(P)))$. It would be
interesting to connect these formulae with counting invariants for
stable objects in $D^b(\coh(P))$.'
\end{quotation}
If Conjecture \ref{an6conj1} and the reasoning behind it are
correct, then the $\hat J^\al$ we propose are invariants `counting'
$\tau$-semistable objects in class $\al$, which are independent of
choice of $\tau$. Following Bridgeland's suggestion, the author
wonders if the $\hat J^\al$ can be combined in a generating function
on $\Stab(P)$ to give one of Borcherd's automorphic forms of weight
$k\ge 1$. One might try something like:
\e
f(\tau)=\ts\sum_{\al\in K(\A)\sm\{0\}}\hat J^\al Z(\al)^{-k},
\label{an6eq41}
\e
where $Z:K(\A)\ra\C$ is the `central charge' associated to $\tau$.
If the sum converges absolutely then $f:\Stab(P)\ra\La\ot_\Q\C$ is a
holomorphic function. See \cite{Joyc8} for results on how to encode
the invariants of this paper into holomorphic generating functions
on complex manifolds of stability conditions.

{\it Donaldson invariants} are invariants of smooth 4-manifolds $P$
described in Donaldson and Kronheimer \cite{DoKr}, which when $P$ is
a K\"ahler surface effectively `count' stable rank 2 holomorphic
vector bundles in class $\al$ on $P$. When $b^2_+(P)>1$ these
invariants depend only on $P$ as a smooth 4-manifold, and so are
unchanged under deformations of $P$. When $b^2_+(P)=1$ they depend
on a little more: they are defined using a metric $g$ on $P$, which
determines a splitting $H^2(P,\R)=\H^2_+\op\H^2_-$, and the
invariants depend on this splitting and have wall-crossing behaviour
when~$\H^2_-\cap H^2(P,\Z)\ne\{0\}$.

If $K_P^{-1}$ is nef but not trivial then $H^{2,0}(P)=0$, so
$b^2_+(P)=1$, and $\H^2_+=\langle[\om]\rangle$ is spanned by the
cohomology class of the K\"ahler form $\om$. This depends on the
ample line bundle $E$ used to embed $P$ into projective space, and
so on the stability condition $(\ga,G_2,\le)$. Thus, the
wall-crossing behaviour for Donaldson invariants when $b^2_+(P)=1$
is directly analogous to the transformation laws \eq{an6eq37},
\eq{an6eq40} for our invariants under change of $(\ga,G_2,\le)$.
When $K_P$ is trivial we have $b^2_+(P)=3$, so Donaldson invariants
are independent of $\H^2_\pm$; this is analogous to $\bar
J^\al(\ga)^{\sst\La}$ being independent of $\ga$ in this case.

The author wonders whether there exist invariants related to
Donaldson invariants and similar to $\Iss^\al(\ga)^{\sst\La},\bar
J^\al(\ga)^{\sst\La}$ above, for which the analogues of Theorems
\ref{an6thm4} and \ref{an6thm5} hold, but which are {\it independent
of deformations of\/ $P$ which do not change} $K(\coh(P))$. When
$\K=\C$ this means not changing $H^{1,1}(P)\cap H^2(P,\Q)$ in
$H^2(P,\C)$. We discuss a similar question for Calabi--Yau 3-folds
and Donaldson--Thomas invariants in Conjecture \ref{an6conj2} below.

\subsection{Invariants of Calabi--Yau 3-folds}
\label{an65}

We now prove that when $\A=\coh(P)$ for $P$ a Calabi--Yau 3-fold the
invariants $J^\al(\tau)^{\sst\Om},\Jsib(I,\pr,\ka,\tau)^{\sst \Om}$
of \eq{an6eq12} have special properties. Suppose Assumptions
\ref{an2ass1} and \ref{an3ass} hold, and $\bar\chi:K(\A)\t
K(\A)\ra\Z$ is biadditive and satisfies
\e
\begin{split}
&\bigl(\dim_\K\Hom(X,Y)-\dim_\K\Ext^1(X,Y)\bigr)-\\
&\bigl(\dim_\K\Hom(Y,X)-\dim_\K\Ext^1(Y,X)\bigr)=
\bar\chi\bigl([X],[Y]\bigr)\quad\text{for all $X,Y\in\A$.}
\end{split}
\label{an6eq42}
\e
We showed in \cite[\S 6.6]{Joyc6} using Serre duality that this
holds if $\A=\coh(P)$ for $P$ a Calabi--Yau 3-fold. Also,
\eq{an6eq8} implies \eq{an6eq42} with $\bar\chi(\al,\be)=
\chi(\al,\be)-\chi(\be,\al)$, so \eq{an6eq42} is a weakening of
\eq{an6eq8}, and thus holds for $\A=\coh(P)$ with $P$ a smooth
projective curve and for $\A=\modKQ$, as in~\S\ref{an62}.

Under these assumptions, in \cite[\S 6.6]{Joyc6} we defined a {\it
Lie algebra morphism}
\e
\Psi^{\sst\Om}\ci\bar\Pi^{\Th,\Om}_{\fObj_\A}:\SFai(\fObj_\A)\ra
C^\ind(\A,\Om,\ha\bar\chi),
\label{an6eq43}
\e
where $C^\ind(\A,\Om,\ha\bar\chi)$ is an explicit Lie algebra
contained in an explicit algebra $C(\A,\Om,\ha\bar\chi)$. Now
$C^\ind(\A,\Om,\ha\bar\chi)$ is an $\Om$-module with $\Om$-basis
$c^\al$ for $\al\in C(\A)$, and comparing Definition \ref{an6def4}
with \cite[Def.~6.10]{Joyc6} shows that
\e
\begin{split}
\Psi^{\sst\Om}\ci\bar\Pi^{\Th,\Om}_{\fObj_\A}\bigl(\bdsib(I,\pr,\ka,\tau)
\bigr)&=\Jsib(I,\pr,\ka,\tau)^{\sst\Om}\,c^{\ka(I)}\\
\text{and}\quad
\Psi^{\sst\Om}\ci\bar\Pi^{\Th,\Om}_{\fObj_\A}\bigl(\bep^\al(\tau)
\bigr)&=J^\al(\tau)^{\sst\Om}\,c^\al.
\end{split}
\label{an6eq44}
\e

We use this to prove {\it multiplicative identities\/} on
the~$\Jsib(I,\pr,\ka,\tau)^{\sst\Om}$.

\begin{thm} Let Assumptions \ref{an2ass1} and \ref{an3ass} hold,
$\bar\chi:K(\A)\t K(\A)\ra\Z$ be biadditive and satisfy
\eq{an6eq42}, and\/ $(\tau,T,\le)$ be a permissible weak stability
condition on $\A$. Suppose $(I,\pr,\ka),(J,\ls,\la)$ are $\A$-data
with\/ $(I,\pr),(J,\ls)$ connected and\/ $I\!\cap\!J\!=\!\emptyset$.
Define $K\!=\!I\!\amalg\!J$ and\/ $\mu:K\!\ra\!C(\A)$ by
$\mu\vert_I\!=\!\ka$ and\/ $\mu\vert_J\!=\!\la$. Then
\e
\begin{gathered}
\sum_{\substack{\text{p.o.s $\tl$ on $K$: $(K,\tl)$ connected,}\\
\text{$\tl\vert_I=\pr$, $\tl\vert_J=\ls$, and}\\
\text{$i\in I$, $j\in J$ implies $j\ntl i$}}
\!\!\!\!\!\!\!\!\!\!\!\!\!\!\!\!\!\!\!\!\!\!\!\!\!\!\!\!\!\!\!}
\Jsib(K,\tl,\mu,\tau)^{\sst\Om}-
\sum_{\substack{\text{p.o.s $\tl$ on $K$: $(K,\tl)$ connected,}\\
\text{$\tl\vert_I=\pr$, $\tl\vert_J=\ls$, and}\\
\text{$i\in I$, $j\in J$ implies $i\ntl j$}}
\!\!\!\!\!\!\!\!\!\!\!\!\!\!\!\!\!\!\!\!\!\!\!\!\!\!\!\!\!\!\!}
\Jsib(K,\tl,\mu,\tau)^{\sst\Om}\\
=\bar\chi\bigl(\ka(I),\la(J)\bigr)
\Jsib(I,\pr,\ka,\tau)^{\sst\Om}\Jsib(J,\ls,\la,\tau)^{\sst\Om}.
\end{gathered}
\label{an6eq45}
\e
\label{an6thm6}
\end{thm}

\begin{proof} As \eq{an6eq43} is a Lie algebra morphism, by
\eq{an6eq44} we have
\begin{align*}
\Psi^{\sst\Om}\ci\bar\Pi^{\Th,\Om}_{\fObj_\A}\bigl(
\bigl[\bdsib(I,\pr,\ka,\tau),\bdsib(J,\ls,\la,\tau)\bigr]\bigr)=&\\
\Jsib(I,\pr,\ka,\tau)^{\sst\Om}\Jsib(J,\ls,\la,\tau)^{\sst\Om}
&[c^{\ka(I)},c^{\la(J)}].
\end{align*}
Also $[c^{\ka(I)},c^{\la(J)}]=\bar\chi\bigl(\ka(I),\la(J)\bigr)
c^{\ka(I)+\la(J)}$ by definition of $C^\ind(\A,\Om,\ha\bar\chi)$.
The multiplicative relations between the $\bdsib(I,\pr,\ka,\tau)$ in
$\bHp_\tau$ are given explicitly in \cite[\S 7.1 \& \S 8]{Joyc7}.
Using these we can write $\bigl[\bdsib(I,\pr,\ka,\tau),
\bdsib(J,\ls,\la,\tau)\bigr]$ as a linear combination of
$\bdsib(K,\tl,\mu,\tau)$ over different $\tl$. Combining all this
and \eq{an6eq44} gives \eq{an6eq45}, but without the conditions
$(K,\tl)$ connected. However, as $(I,\pr),(J,\ls)$ are connected the
only possibility for $\tl$ with $(K,\tl)$ disconnected is given by
$a\tl b$ if and only if $a,b\in I$ with $a\pr b$ or $a,b\in J$ with
$a\ls b$. As this appears in both sums with opposite signs, the
terms in disconnected $(K,\tl)$ cancel out, so we may restrict to
connected $(K,\tl)$ in~\eq{an6eq45}.
\end{proof}

Next we discuss how the invariants depend on $(\tau,T,\le)$. The
$\Jsib(I,\pr,\ka,\tau)^{\sst\La^\ci}$ and $\Jsib(I,\pr,\ka,\tau)^{
\sst\Om}$ satisfy {\it additive} transformation laws similar to
\eq{an6eq1}, which may be deduced from \eq{an6eq1} and identities in
\cite[\S 8]{Joyc7}. The $J^\al(\tau)^{\sst\Om}$ satisfy a {\it
multiplicative} transformation law, \eq{an6eq48} below. To deduce it
we need some facts about multiplication in~$C(\A,\Om,\ha\bar\chi)$.

Suppose $\al_1,\ldots,\al_n\in C(\A)$, so that $c^{\al_1},\ldots,
c^{\al_n}\in C^\ind(\A,\Om,\ha\bar\chi)$. Then using the
multiplication relations in \cite[Def.~6.10]{Joyc6} we can compute
$c^{\al_1}\star c^{\al_2}\star\cdots\star c^{\al_n}$ in the algebra
$C(\A,\Om,\ha\bar\chi)$. This is rather complicated, but we will
only need to know the component in $C^\ind(\A,\Om,\ha\bar\chi)$,
that is, the coefficient of $c^{\al_1+\cdots+\al_n}$ in this sum.
Calculation shows this is given by:
\begin{gather}
c^{\al_1}\star\cdots\star c^{\al_n}=\text{ terms in $c_{[I,\ka]}$,
$\md{I}>1$, }
\label{an6eq46}\\
+\raisebox{-6pt}{\begin{Large}$\displaystyle\biggl[$\end{Large}}
\frac{1}{2^{n-1}}\!\!\!\!\!
\sum_{\substack{\text{connected, simply-connected digraphs
$\Ga$:}\\
\text{vertices $\{1,\ldots,n\}$, edge $\mathop{\bu} \limits^{\sst
i}\ra\mathop{\bu}\limits^{\sst j}$ implies $i<j$}}} \,\,\,
\prod_{\substack{\text{edges}\\
\text{$\mathop{\bu}\limits^{\sst i}\ra\mathop{\bu}\limits^{\sst
j}$}\\ \text{in $\Ga$}}}\bar\chi(\al_i,\al_j)
\raisebox{-6pt}{\begin{Large}$\displaystyle\biggr]$\end{Large}}
c^{\al_1+\cdots+\al_n}.
\nonumber
\end{gather}
Here a {\it digraph\/} is a directed graph.

Since \eq{an5eq7} can be regarded as an identity in the Lie algebra
$\SFai(\fObj_\A)$ by Theorem \ref{an5thm2}, we can apply the Lie
algebra morphism \eq{an6eq43} to \eq{an5eq7} to get an identity
writing $J^\al(\ti\tau)^{\sst\Om}$ in terms of the
$J^{\ka(i)}(\tau)^{\sst\Om}$. This is expressed as a sum over all
$\A$-data $(\{1,\ldots,n\},\le,\ka)$ followed by a sum over digraphs
$\Ga$ with vertices $\{1,\ldots,n\}$. To simplify this formula we
define some new notation.

\begin{dfn} Suppose Condition \ref{an4cond} holds for $(\tau,T,\le)$
and $(\ti\tau,\ti T,\le)$, and let $\Ga$ be a connected,
simply-connected digraph with finite vertex set $I$, and $\ka:I\ra
C(\A)$. Define $V(I,\Ga,\ka,\tau,\ti\tau)\in\Q$ by
\e
V(I,\Ga,\ka,\tau,\ti\tau)=\frac{1}{2^{\md{I}-1}\md{I}!}
\!\!\sum_{\substack{\text{total orders $\pr$ on $I$:}\\ \text{edge
$\mathop{\bu} \limits^{\sst i}\ra\mathop{\bu}\limits^{\sst j}$ in
$\Ga$ implies $i\pr j$}}} \!\!\!\!\! U(I,\pr,\ka,\tau,\ti\tau).
\label{an6eq47}
\e
\label{an6def7}
\end{dfn}

With this we prove a transformation law for
the~$J^\al(\tau)^{\sst\Om}$.

\begin{thm} Let Assumptions \ref{an2ass1} and \ref{an3ass} hold
and\/ $\bar\chi:K(\A)\t K(\A)\ra\Z$ be biadditive and satisfy
\eq{an6eq42}. Suppose $(\tau,T,\le),(\ti\tau,\ti T,\le)$ are
permissible weak stability conditions on $\A$, the change from
$(\tau,T,\le)$ to $(\ti\tau,\ti T,\le)$ is globally finite, and
there exists a weak stability condition $(\hat\tau,\hat T,\le)$ on
$\A$ dominating $(\tau,T,\le),(\ti\tau,\ti T,\le)$ with the change
from $(\hat\tau,\hat T,\le)$ to $(\ti\tau,\ti T,\le)$ locally
finite. Then for all\/ $\al\in C(\A)$ the following holds in $\Om$,
with only finitely many nonzero terms:
\e
J^\al(\ti\tau)^{\sst\Om}\!=\!\!\!\!
\sum_{\substack{\text{iso.}\\ \text{classes}\\
\text{of finite}\\ \text{sets $I$}}}\,\,
\sum_{\substack{\ka:I\ra C(\A):\\ \ka(I)=\al}}\,\,
\sum_{\begin{subarray}{l} \text{connected,}\\
\text{simply-connected}\\ \text{digraphs $\Ga$,}\\
\text{vertices $I$}\end{subarray}\!\!\!\!\!\!\!\!\!\!\!\!\!\!\!\!\!
\!\!\!\!\!\!\!\!\!\!} V(I,\Ga,\ka,\tau,\ti\tau)
\begin{aligned}[t]
&\cdot\prod\limits_{\text{edges \smash{$\mathop{\bu}\limits^{\sst
i}\ra\mathop{\bu}\limits^{\sst j}$} in
$\Ga$}\!\!\!\!\!\!\!\!\!\!\!\!\!\!\!\!\!\!\!\!\!\!\!
\!\!\!\!\!\!\!\!\!\!\!\!\!\!} \bar\chi(\ka(i),\ka(j))\\
&\cdot\prod\nolimits_{i\in I}J^{\ka(i)}(\tau)^{\sst\Om}.
\end{aligned}
\label{an6eq48}
\e
\label{an6thm7}
\end{thm}

\begin{proof} Theorem \ref{an5thm1} implies that \eq{an5eq7} holds,
with only finitely many nonzero terms. Theorem \ref{an5thm2} shows
\eq{an5eq7} can be rewritten as a Lie algebra identity, a linear
combination of multiple commutators in $\SFai(\fObj_\A)$. So
applying the Lie algebra morphism \eq{an6eq43} yields a Lie algebra
identity in $C^\ind(\A,\Om,\ha\bar\chi)$, which by \eq{an6eq44}
writes $J^\al(\ti\tau)^{\sst\Om}c^\al$ as a linear combination of
multiple commutators of $J^{\ka(i)}(\tau)^{\sst\Om}c^{\ka(i)}$. But
$C^\ind(\A,\Om,\ha\bar\chi)$ is a Lie subalgebra of the algebra
$C(\A,\Om,\ha\bar\chi)$, so we can regard this as an identity in
$C(\A,\Om,\ha\bar\chi)$, avoiding the need to rewrite \eq{an5eq7} in
terms of multiple commutators. This proves:
\begin{equation*}
J^\al(\ti\tau)^{\sst\Om}c^\al=\sum_{\begin{subarray}{l}
\text{$\A$-data $(\{1,\ldots,n\},\le,\ka):$}\\
\text{$\ka(\{1,\ldots,n\})=\al$}\end{subarray}
\!\!\!\!\!\!\!\!\!\!\!\!\!\!\!\!\!\!\!\!\!\!\!\!\!\!\!\!\!
\!\!\!\!\!\!\!\!\!\!\!\!\!\! }
U(\{1,\ldots,n\},\le,\ka,\tau,\ti\tau)\bigl[\ts\prod_{i=1}^n
J^{\ka(i)}(\tau)^{\sst\Om}\bigr]c^{\ka(1)}\star\cdots\star
c^{\ka(n)}.
\end{equation*}
Equating coefficients of $c^\al$ on both sides of this equation, and
noting that the coefficient in $c^{\ka(1)}\star\cdots\star
c^{\ka(n)}$ is given by \eq{an6eq46}, yields:
\e
\begin{gathered}
J^\al(\ti\tau)^{\sst\Om}=\!\!\!\!\!\!
\sum_{\begin{subarray}{l}
\text{$\A$-data $(\{1,\ldots,n\},\le,\ka):$}\\
\text{$\ka(\{1,\ldots,n\})=\al$}\end{subarray} }
\sum_{\substack{\text{connected, simply-connected digraphs $\Ga$:}\\
\text{vertices $\{1,\ldots,n\}$, edge $\mathop{\bu} \limits^{\sst
i}\ra\mathop{\bu}\limits^{\sst j}$ implies $i<j$}}}\\
\frac{1}{2^{n-1}}\, U(\{1,\ldots,n\},\le,\ka,\tau,\ti\tau)
\!\!\!\!\! \prod_{\text{edges $\mathop{\bu}\limits^{\sst
i}\ra\mathop{\bu}\limits^{\sst j}$ in $\Ga$}}\!\!\!\!\!
\bar\chi(\ka(i),\ka(j))\prod_{i=1}^n J^{\ka(i)}(\tau)^{\sst\Om}.
\end{gathered}
\label{an6eq49}
\e

We claim that substituting \eq{an6eq47} into \eq{an6eq48} gives a
sum equivalent to \eq{an6eq49}. For substituting \eq{an6eq47} into
\eq{an6eq48} gives sums over $I$, $\ka$, $\Ga$, and total orders
$\pr$ on $I$. The relationship with \eq{an6eq49} is that having
chosen $I,\Ga,\pr$, there is a unique isomorphism
$(I,\pr)\cong(\{1,\ldots,n\},\le)$, where $n=\md{I}$. So identifying
$I$ with $\{1,\ldots,n\}$ relates the sums over $\Ga,\pr$ in
\eq{an6eq47}--\eq{an6eq48} to the sums over $n,\Ga$ in \eq{an6eq49}.
But this is not a 1-1 correspondence: rather, to each choice of
$n,\Ga$ in \eq{an6eq49} there are $\md{I}!$ choices of $I,\Ga,\pr$
in the substitution of \eq{an6eq47} in \eq{an6eq48}, as there are
$\md{I}!$ total orders $\pr$ on $I$. This is cancelled by the factor
$1/\md{I}!$ in \eq{an6eq47}. All the other terms in
\eq{an6eq47}--\eq{an6eq48} and \eq{an6eq49} immediately agree.
\end{proof}

\begin{rem} The author can prove that in Definition \ref{an6def7},
if $\ti\Ga$ is a directed graph obtained from $\Ga$ by reversing the
directions of $k$ edges in $\Ga$ then $V(I,\ti\Ga,\ka,\tau,\ti\tau)=
(-1)^kV(I,\Ga,\ka,\tau,\ti\tau)$. Since $\bar\chi$ is antisymmetric,
replacing $\Ga$ by $\ti\Ga$ also multiplies the product
$\prod_{\text{\smash{$\mathop{\bu}\limits^{\sst i}\ra \mathop{\bu}
\limits^{\sst j}$} in $\Ga$}}\bar\chi(\ka(i),\ka(j))$ in
\eq{an6eq48} by~$(-1)^k$.

Thus the product of terms on the r.h.s.\ of \eq{an6eq48} is actually
{\it independent\/} of the orientation of $\Ga$, and depends only on
the underlying undirected graph. As there are $2^{\md{I}-1}$
orientations on this graph, we could omit the factor
$1/2^{\md{I}-1}$ in \eq{an6eq47} and write \eq{an6eq48} as a sum
over undirected graphs rather than digraphs.

Equation \eq{an6eq48} will play an important r\^ole in a sequel
\cite{Joyc8}. Neglecting issues to do with convergence of infinite
sums, we encode the invariants $J^\al(\tau)^{\sst\Om}$ in {\it
holomorphic generating functions\/} on the complex manifold of
stability conditions. Because the $J^\al(\tau)^{\sst\Om}$ satisfy
the multiplicative transformation law \eq{an6eq48}, these generating
functions satisfy a {\it nonlinear p.d.e.}, which can be interpreted
as the {\it flatness of a connection\/} over the complex manifold of
stability conditions.
\label{an6rem2}
\end{rem}

In the remainder of the section we discuss how our results should
relate to other proposed invariants of Calabi--Yau 3-folds, and the
whole Mirror Symmetry picture. Motivated by Donaldson and Thomas
\cite[p.~33--34]{DoTh}, Thomas \cite{Thom} defines invariants
$DT^\al(\tau)$ `counting' $\tau$-stable coherent sheaves in class
$\al\in C(\coh(P))$ on a Calabi--Yau 3-fold $P$, which are now known
as {\it Donaldson--Thomas invariants}. We compare these with our
invariants $J^\al(\tau)^{\sst\Om}$ above, which also `count'
$\tau$-semistable coherent sheaves in class $\al$ on~$P$:
\begin{itemize}
\setlength{\itemsep}{0pt}
\setlength{\parsep}{0pt}
\item The main good property of Donaldson--Thomas invariants
$DT^\al(\tau)$ is that they are {\it unchanged\/} under deformations
of the complex structure of $P$. Thomas achieves this by using a
{\it virtual moduli cycle} to cut the moduli schemes down to the
expected dimension (zero) and then counting points.

In contrast, our invariants $J^\al(\tau)^{\sst\Om}$ are {\it not\/}
expected to be unchanged under deformations of $P$. This is because
rather than using virtual moduli cycles we just take a motivic
invariant, such as an Euler characteristic, of the moduli scheme as
it stands. This is quite a crude thing to do.
\item Donaldson--Thomas invariants are (so far) defined only for
$\al\in C(\A)$ for which $\Oss^\al(\tau)=\Ost^\al(\tau)$, that is,
$\tau$-semistability and $\tau$-stability coincide. This is because
strictly $\tau$-semistable objects would give singular points of the
moduli space which the virtual moduli cycle technology is presently
unable to cope with.

In contrast, our invariants $J^\al(\tau)^{\sst\Om}$ are defined for
all $\al\in C(\A)$. Moreover, a great deal of the work in this paper
and \cite{Joyc3,Joyc4,Joyc5,Joyc6,Joyc7} is really about how to deal
with strictly $\tau$-semistables -- for instance,
$\ep^\al(\tau),\bep^\al(\tau)$ coincide with
$\dss^\al(\tau),\bdss^\al(\tau)$ except over strictly
$\tau$-semistables, and much of \cite{Joyc3,Joyc4,Joyc6} concerns
how best to include stabilizer groups when forming invariants of
subsets of stacks, and this is not relevant for $\tau$-stable
objects as they all have stabilizer group~$\K^\t$.
\item The transformation laws for Donaldson--Thomas invariants under
change of stability condition are not known. In contrast, our
invariants $J^\al(\tau)^{\sst\Om}$ transform according
to~\eq{an6eq48}.
\item Donaldson--Thomas invariants are defined uniquely and take
values in $\Z$. In contrast, our invariants depend on a choice of
motivic invariant $\Th$ in Assumption \ref{an2ass1}, with values in
a $\Q$-algebra $\Om$. As in \cite[Ex.~6.3]{Joyc4}, possibilities for
$\Th$ include the Euler characteristic $\chi$, and the sum of the
virtual Betti numbers.
\item Donaldson and Thomas worked hard to produce invariants whose
behaviour under deformation of $P$ is understood. In contrast, we
have worked hard to produce invariants whose behaviour under change
of weak stability condition $(\tau,T,\le)$ is understood.
\end{itemize}

We would now like to conjecture that there exist invariants
`counting' $\tau$-semistable coherent sheaves on $P$ which combine
the good features of both Donaldson--Thomas invariants and our own.

\begin{conj} Fix $\K=\C$, let\/ $P$ be a Calabi--Yau $3$-fold, and
define $\A=\coh(P)$, $\fF_\A$ and $K(\A)$ as in {\rm\cite[\S
9.1]{Joyc5}}. Let\/ $(\tau,T,\le)$ be a permissible stability
condition of Gieseker type on $\A$, as in {\rm\cite[\S 4.4]{Joyc7}}.
Then there should exist extended Donaldson--Thomas invariants
$\bar{DT}^\al(\tau)\in\Q$ defined for all\/ $\al\in C(\A)$, which
should be unchanged under deformations of\/ $P$, and should
transform according to \eq{an6eq48} under change of stability
condition. If\/ $\al\in C(\A)$ with\/
$\Oss^\al(\tau)=\Ost^\al(\tau)$, so that the Donaldson--Thomas
invariant\/ $DT^\al(\tau)$ is defined as in {\rm\cite{Thom}}, we
have~$\bar{DT}^\al(\tau)=DT^\al(\tau)\in\Z$.

There may also exist more complicated systems of invariants
analogous to the $\Iss(I,\pr,\ka,\tau)$ or
$\Jsib(I,\pr,\ka,\tau)^{\sst\Om}$ above, which take values in $\Q$,
are unchanged under deformations of\/ $P$, and transform in the
appropriate way under change of stability condition.
\label{an6conj2}
\end{conj}

I believe that proving this conjecture is feasible, although
difficult, and that proving it will probably involve an extension to
virtual moduli cycle technology. I have plans to attempt this in the
next few years, in collaboration with others.

The invariants of Conjecture \ref{an6conj2} should be defined using
0-dimensional virtual moduli cycles, which are basically finite sets
of points with integer multiplicities. The only information in them
independent of choices is the number of points $[X]$, counted with
multiplicity and weighted by the correct function of the stabilizer
group $\Aut(X)$, probably $\chi(\Aut(X)/T_X)^{-1}$, where $T_X$ is a
maximal torus in $\Aut(X)$. This is why we say the invariants should
take values in $\Q$, in contrast to the rest of the section where
our invariants take values in more general
algebras~$\La,\La^\ci,\Om$.

Now under the {\it Homological Mirror Symmetry programme} of
Kontsevich \cite{Kont}, (semi)stable coherent sheaves on $P$ are
supposed to be mirror to {\it special Lagrangian $3$-folds} ({\it
SL\/ $3$-folds\/}) in the mirror Calabi--Yau 3-fold $M$. So we
expect invariants counting (semi)stable (complexes of) coherent
sheaves on $P$ to be equal to other invariants counting SL 3-folds
in $M$, at least if these invariants count anything meaningful in
String Theory.

Motivated by a detailed study \cite{Joyc2} of the singularities of
SL $m$-folds, in \cite{Joyc1} I made a conjecture that there should
exist interesting invariants $K^\al(J)$ of $M$ that `count' SL
homology 3-spheres in a given class $\al\in H_3(M,\Z)$. I expected
these invariants to be independent of the choice of K\"ahler form on
$M$, and to transform according to some wall-crossing formulae under
deformation of complex structure $J$, but at the time I could only
determine these transformation laws in the simplest cases. I can now
expand this conjecture to specify these transformation laws.

\begin{conj} Let\/ $(M,J,\om,\Om)$ be an (almost) Calabi--Yau
$3$-fold, with compact\/ $6$-manifold\/ $M$, complex structure $J$,
K\"ahler form $\om$ and holomorphic volume form $\Om$. Then there
should exist invariants $K^\al(J)\in\Q$ for $\al\in H_3(M,\Z)$ which
`count' special Lagrangian homology $3$-spheres $N$ in $M$ with
$[N]=\al$, and probably other kinds of immersed or singular SL\/
$3$-folds as well. These invariants $K^\al(J)$ are independent of
the K\"ahler form $\om$ and of complex rescalings of\/ $\Om$, and so
depend only on the complex structure~$J$.

The holomorphic volume form $\Om$ determines a stability condition
$\tau$ of Bridgeland type {\rm\cite{Brid1}} on the derived Fukaya
category $D^b(F(M,\om))$ of the symplectic manifold $(M,\om)$, a
triangulated category, and SL\/ $3$-folds $N$ whose Floer homology
is unobstructed correspond to $\tau$-semistable objects in
$D^b(F(M,\om))$. Thus, the $K^\al(J)$ are invariants counting
$\tau$-semistable objects in class $\al$ in
$H_3(M,\Z)=K(D^b(F(M,\om)))$. Under deformation of complex structure
of\/ $M$, the $K^\al(J)$ transform according to the triangulated
category extension of\/~\eq{an6eq48}.

If\/ $P$ is a mirror Calabi--Yau $3$-fold to $M$, then the
$K^\al(J)$ should be equal to invariants counting
$\ti\tau$-semistable objects in the derived category $D^b(\coh(P))$
of coherent sheaves on $P$, with respect to a stability condition
$\ti\tau$ on $D^b(\coh(P))$ of Bridgeland type. These invariants are
a triangulated category version of the extended Donaldson--Thomas
invariants $\bar{DT}^\al(\tau)$ of Conjecture \ref{an6conj2}, and
coincide with them if\/ $\ti\tau$ is induced by Gieseker stability
on~$\coh(P)$.

There may also exist more complicated systems of invariants
analogous to the $\Iss(I,\pr,\ka,\tau)$ or
$\Jsib(I,\pr,\ka,\tau)^{\sst\Om}$ above, that count configurations
of SL\/ $3$-folds, take values in $\Q$, are independent of\/ $\om$,
and transform in the appropriate way under change of complex
structure. This set-up may generalize to (almost) Calabi--Yau
$m$-folds for all\/~$m\ge 2$.
\label{an6conj3}
\end{conj}

The extension of our programme to triangulated categories will be
discussed at length in \S\ref{an7}, so we will say no more about it
here. I expect Conjecture \ref{an6conj3} to be extremely difficult
to prove, much more so than Conjecture \ref{an6conj2}, even just the
part concerning the definition of $K^\al(J)$, independence of $\om$,
and transformation laws under deformation of $J$. This is because we
cannot use the machinery of algebraic geometry, and instead need to
know a lot about the singular behaviour of SL 3-folds, which is at
the moment only partially understood.

In fact clarifying Conjecture \ref{an6conj3} was the beginning of
this whole project, which eventually grew to \cite{Joyc3,Joyc4,
Joyc5,Joyc6,Joyc7} and this paper, and may continue to expand. I
wanted to know more about the proposed invariants $K^\al(J)$, in
particular their transformation laws under change of $J$, so I
decided to study the mirror problem of counting semistable coherent
sheaves on Calabi--Yau 3-folds using algebraic geometry. I didn't
realize at the time how large an undertaking this would be.

The necessity of introducing, and counting, {\it configurations} in
order to understand transformation laws under change of stability
condition came directly out of my work \cite[\S 9]{Joyc2} on
creation of new SL $m$-folds by multiple connected sums under
deformation of the underlying Calabi--Yau $m$-fold.

I believe that the invariants $\bar{DT}^\al(\tau)$ and $K^\al(J)$
proposed in Conjectures \ref{an6conj2} and \ref{an6conj3} should
play an important part in a chapter of the Mirror Symmetry story
that has not yet been understood. They should encode a lot about the
structure of the `stringy K\"ahler moduli space', teach us about
{\it branes} and $\Pi$-{\it stability} in String Theory, and perhaps
have other applications.

Note that Conjecture \ref{an6conj3} has predictive power for the
kind of wall-crossing phenomena in special Lagrangian geometry
considered in \cite[\S 9]{Joyc2}. For example, we can consider a
smooth family of Calabi--Yau 3-folds $M_t:t\in(-\ep,\ep)$ with
compact nonsingular SL homology 3-spheres $N_1,N_2$ in $M_0$ of the
same phase, intersecting transversely in a single point. Then
\cite[Th.~9.7]{Joyc2} gives necessary and sufficient criteria for
the existence of a connected sum SL 3-fold $N_1\#\,N_2$ in $M_t$ for
small $t>0$ or~$t<0$.

Using Conjecture \ref{an6conj3} and a little geometric and algebraic
intuition, we can make predictions for the number of SL 3-folds in
$M_t$ of the general form $a_1N_1\#\,a_2N_2$ for $a_1,a_2>0$. This
is not at all clear using the techniques of \cite{Joyc2}, as for
$a_1>1$ the connected sum `necks' might occur at any point in $N_1$,
and we might have to consider branching or covering phenomena over
$N_1$ too.

\section{Questions for further research}
\label{an7}

We have already discussed a number of research problems in
\S\ref{an63}--\S\ref{an65}, in particular Conjectures
\ref{an6conj1}, \ref{an6conj2} and \ref{an6conj3}. We now pose some
more, and in the process sketch out how the author would like to see
the subject develop in future. Perhaps the most important is:

\begin{prob} Develop an extension of the work of\/
{\rm\cite{Joyc5,Joyc6,Joyc7}} and this paper from abelian categories
to triangulated categories, and apply it to derived categories
$D^b(\A)$ when $\A$ is a category of quiver representations
$\modKQ,\modKQI$ or coherent sheaves~$\coh(P)$.
\label{an7prob1}
\end{prob}

Here are some issues and difficulties involved in this:
\begin{itemize}
\setlength{\itemsep}{0pt}
\setlength{\parsep}{0pt}
\item {\it Defining configurations in triangulated categories}
$\T$. This looks straightforward: in Definition \ref{an3def2} we
replace the exact sequences \eq{an3eq1} by distinguished triangles.
This involves introducing a third family of morphisms
$\partial(L,J):\si(L)\ra\si(J)[1]$ in $\T$, so that a configuration
is a quadruple $(\si,\io,\pi,\partial)$, and we must add some extra
axioms involving the $\partial(L,J)$. For example, a
$(\{1,\ldots,n\},\le)$-configuration will be equivalent to a {\it
distinguished $n$-dimensional hypersimplex} in $\T$,
\cite[p.~260-1]{GeMa}.
\item {\it Constructing moduli stacks of objects and configurations
in} $\T$. This seems difficult, and neither triangulated categories
nor Artin stacks may be the right frameworks. To\"en \cite{Toen1,
Toen2,ToVa} works instead with {\it dg-categories} \cite{Toen1}, a
richer structure from which one can often recover the triangulated
categories we are interested in as the homotopy category. When $\T$
is a dg-category To\"en and Vaqui\'e \cite{ToVa} construct moduli
$\iy$-{\it stacks} and $D^-$-{\it stacks} of objects in $\T$, some
of which can be truncated to Artin stacks.

For our applications we may not need the moduli `stack' of all
objects in $\T$, but only those which could be $\tau$-semistable for
some stability condition $\tau$, or appear in the heart of some
t-structure. So we can probably restrict to $X\in\T$ with
$\Ext^i(X,X)=0$ for $i<0$, which may improve the problem.
\item {\it Non-representable morphisms}. For abelian categories
$\bs\si(I):\fM(I,\pr)_\A\ra\fObj_\A$ is a {\it representable}
1-morphism, but in the triangulated category case this will no
longer hold. This is because for abelian categories the automorphism
group of an exact sequence injects into the automorphism group of
its middle term, but the analogue fails for distinguished triangles
in triangulated categories.

Because of this, the multiplication $*$ and operations $P_\sIp$ on
$\CF(\fObj_\A)$, $\SF(\fObj_\A)$ in \cite{Joyc6} cannot be defined
in the triangulated category case. We can still work with the spaces
$\uSF(\fObj_\A)$, though, as they do not require representability.
\item {\it Non-Cartesian squares}. In \cite[Th.~7.10]{Joyc5} we
proved some commutative diagrams of 1-morphisms of configuration
moduli stacks are Cartesian. In the triangulated category case the
corresponding diagrams may be commutative, but not Cartesian. So the
multiplication $*$ on $\CF,\SF(\fObj_\A)$ may no longer be
associative, which spoils most of \cite{Joyc6,Joyc7} and this paper.

The reason is the `nonfunctoriality of the cone' in triangulated
categories, discussed by Gelfand and Manin \cite[p.~245]{GeMa}. If
$\phi:X\ra Y$ is a morphism in an abelian category $\A$, the kernel
and cokernel of $\phi$ are determined up to canonical isomorphism.
But in a triangulated category $\T$ the cone on $\phi$ is determined
up to isomorphism, but not canonically. This means that the
triangulated category versions of operations on configurations such
as substitution \cite[Def.~5.7]{Joyc5} will be defined up to
isomorphism, but not canonically, which will invalidate the proof
of~\cite[Th.~7.10]{Joyc5}.

I expect the right solution is to change the definition of
triangulated category, perhaps following Neeman \cite{Neem} or
Bondal and Kapranov \cite{BoKa}. Neeman's proposal includes as data
a category $\mathcal S$ of triangles in $\T$, but the morphisms in
$\mathcal S$ are different from morphisms of triangles in $\T$,
which restores the functoriality of the cone. If in defining
configurations we took the distinguished triangles to be objects in
$\mathcal S$, and their morphisms to include morphisms in $\mathcal
S$, I hope the relevant squares will be Cartesian.

In the dg-category approach of To\"en and Vaqui\'e \cite{ToVa} this
problem may not occur. In particular, To\"en \cite{Toen2} defines
{\it derived Hall algebras} for certain dg-categories $\T$, which
are associative. This is an analogue of associativity $*$ on
$\uSF(\fObj_\A)$, and suggests that the programme ought to work.
\item {\it Defining stability conditions on triangulated categories}.
Here we can use Bridgeland's wonderful paper \cite{Brid1}, and its
extension by Gorodentscev et al.\ \cite{GKR}, which combines
Bridgeland's idea with Rudakov's definition for abelian categories
\cite{Ruda}. Since our stability conditions are based on
\cite{Ruda}, the modifications for our framework are
straightforward.
\item {\it Invariants and changing stability conditions}. Given a
triangulated category $\T$ and a Bridgeland stability condition
$\tau$ on $\T$, there is a t-structure on $\T$ whose heart is an
abelian category $\A$, and $\tau$ determines a slope function $\mu$
on $\A$, and $X\in\T$ is $\tau$-semistable if and only if $X=Y[n]$
for some $\mu$-semistable $Y\in\A$ and $n\in\Z$. Here $\A$ is not
unique, but two $\A,\A'$ coming from $\T,\tau$ are related by
tilting. Thus invariants counting $\tau$-semistable objects in $\T$
are essentially the same as invariants counting $\mu$-semistable
objects in $\A$, which we have studied in this paper.

We can also reduce changing stability conditions in triangulated
categories to the abelian category case, in the following way. Let
$\tau,\ti\tau$ be stability conditions on $\T$ which are
`sufficiently close' in some sense. Then the author expects that we
can find a third stability condition $\hat\tau$ on $\T$ and two
t-structures on $\T$ with hearts $\A_1,\A_2$, such that the first
t-structure is compatible with $\tau,\hat\tau$, which arise from
slope functions $\mu_1,\hat\mu_1$ on $\A_1$, and the second
t-structure is compatible with $\hat\tau,\ti\tau$, which arise from
slope functions $\hat\mu_2,\ti\mu_2$ on~$\A_2$.

As both t-structures are compatible with $\hat\tau$, $\A_1,\A_2$ are
related by tilting, and $\hat\mu_1$-semistable objects in $\A_1$ are
essentially the same as $\hat\mu_2$-semistable objects in $\A_2$.
Thus transforming between invariants counting $\tau$- and
$\ti\tau$-semistable objects in $\T$ is equivalent to transforming
between invariants counting $\mu_1$- and $\hat\mu_1$-semistable
objects in the abelian category $\A_1$, and then transforming
between invariants counting $\hat\mu_2$- and $\ti\mu_2$-semistable
objects in the abelian category $\A_2$, which we already know how to
do.

Because of this reduction to the abelian category case, the author
is confident that there should be a well-behaved theory of
invariants counting $\tau$-semistable objects in triangulated
categories, despite all the problems described above.
\end{itemize}

In \S\ref{an64} we proposed that invariants counting
$\ga$-semistable sheaves on a $K3$ surface should be combined in a
holomorphic {\it generating function} \eq{an6eq41}. This had a
simple form because for $K3$ surfaces the invariants $\bar
J^\al(\ga)$ and conjectured invariants $\hat J^\al$ were independent
of the choice of weak stability condition. Our next problem concerns
the `right' way to form generating functions when the invariants are
not independent of weak stability condition.

\begin{prob}{\rm(a)} Let\/ $\A,K(\A),\fF_\A$ be one of the examples
defined using quivers in {\rm\cite[Ex.s 10.5--10.9]{Joyc5}}, such as
$\A=\modKQ$. Let\/ $c,r:K(\A)\ra\R$ be group homomorphisms with
$r(\al)>0$ for all\/ $\al\in C(\A)$. Define a slope function
$\mu:C(\A)\ra\R$ by $\mu(\al)=c(\al)/r(\al)$. Then $(\mu,\R,\le)$ is
a permissible stability condition by {\rm\cite[Ex.~4.14]{Joyc7}}.
Define $Z:K(\A)\ra\C$ by $Z(\al)=-c(\al)+ir(\al)$. Then $Z$
determines $\mu$, and the family of\/ $Z$ arising from slope
functions $\mu$ in this way is an open subset\/ $U$ in the complex
vector space~$\Hom(K(\A),\C)$.

Form invariants $\Iss(I,\pr,\ka,\mu)$ or $\Iss^\al(\mu)^{\sst\La}$
or $J^\al(\mu)^{\sst\La^\ci}$ or $J^\al(\mu)^{\sst\Om}$ counting
$\mu$-semistable objects or configurations in $\A$, as in
\S\ref{an6}. Find natural ways to combine these invariants in
generating functions $F:U\ra V$, where $V$ is $\C,\La\ot_\Q\C,
\La^\ci\ot_\Q\C$ or $\Om\ot_\Q\C$, and\/ $F(Z)$ depends on $Z$ and
the values of the invariants for the stability condition
$(\mu,\R,\le)$ depending on $Z$. Is $F$ holomorphic? If it is
defined by an infinite sum, does it converge? Note that the
$\Iss(I,\pr,\ka,\mu),\ldots$ change discontinuously with\/ $\mu$
and\/ $Z$ according to the transformation laws \eq{an6eq1},
\eq{an6eq16}, \eq{an6eq17} or \eq{an6eq48}. Can we make $F$
continuous under these changes?
\medskip

\noindent{\rm(b)} Let $\T$ be a triangulated category and\/
$\Stab(\T)$ the moduli space of Bridgeland stability conditions
{\rm\cite{Brid1}} on $\T$. Each\/ $\tau\in\Stab(\T)$ has an
associated central charge $Z:K(\A)\ra\C$, and the map $\tau\mapsto
Z$ induces a map $\Stab(\T)\ra\Hom(K(\T),\C)$ which is a local
diffeomorphism, and gives $\Stab(\T)$ the structure of a complex
manifold. When $\T=D^b(\A)$ for $\A$ as in {\rm(a)}, each\/ $\mu,Z$
in {\rm(a)} determines a unique $\tau\in\Stab(\T)$ with the same
$Z$, such that\/ $X\in\A\subset D^b(\A)$ is $\mu$-semistable in $\A$
if and only if it is $\tau$-semistable in~$\T$.

Suppose we can solve Problem \ref{an7prob1} and have a good theory
of invariants $\Iss^\al(\tau),\ldots$ counting $\tau$-semistable
objects in $\T$, for $\tau\in\Stab(\T)$. As in {\rm(a)}, try to
encode these in holomorphic generating functions $F:\Stab(\T)\ra V$.
Can we make $F$ continuous when the $\Iss^\al(\tau)$ change
discontinuously with~$\tau$?
\medskip

\noindent{\rm(c)} Now let\/ $\T=D^b(\coh(P))$ for $P$ a Calabi--Yau
$3$-fold. Can we form holomorphic generating functions $F$ on
$\Stab(\T)$ encoding derived category versions of the $J^\al(\tau)$
of\/ \S\ref{an65}, or the $\bar{DT}^\al(\tau)$ of Conjecture
\ref{an6conj2}? If so, does $F$ encode some structure on $\Stab(\T)$
important in String Theory and Mirror Symmetry; for instance, can we
recover the stringy K\"ahler moduli space, the subset of\/
$\Stab(\T)$ corresponding to complex structures on the mirror
Calabi--Yau $3$-fold, from~$F$?
\label{an7prob2}
\end{prob}

In a sequel \cite{Joyc8}, we present the author's attempt at solving
this problem. Rather than a single generating function $F$, we
define a family of generating functions $F^\al$ for all $\al\in
C(\A)$, or $\al\in K(\T)\sm\{0\}$ in the triangulated case. For
example, the invariants $J^\al(\tau)^{\sst\Om}$ of \S\ref{an65} are
encoded in maps $F^\al:U$ or $\Stab(\T)\ra\Om\ot_\Q\C$ of the form
\e
\begin{split}
F^\al(Z)=\!\!\!\!\!\!\sum_{\substack{\text{$n\ge
1$, $\al_1,\ldots,\al_n\in C(\A)$ or $K(\T)\sm\{0\}$:}\\
\text{$\al_1+\cdots+\al_n=\al$, $Z(\al_k)\ne 0$ all
$k$}}}\!\!\!\!\!\!\!\!\!\!\!\!\!\!\! F_n\bigl(Z(\al_1),
\ldots,Z(\al_n)\bigr) \prod_{i=1}^nJ^{\al_i}(\mu)^{\sst\Om}\cdot
\\
\raisebox{-6pt}{\begin{Large}$\displaystyle\biggl[$\end{Large}}
\frac{1}{2^{n-1}}\!\!\!\!\!
\sum_{\substack{\text{connected, simply-connected digraphs
$\Ga$:}\\
\text{vertices $\{1,\ldots,n\}$, edge $\mathop{\bu} \limits^{\sst
i}\ra\mathop{\bu}\limits^{\sst j}$ implies $i<j$}}} \,\,\,
\prod_{\substack{\text{edges}\\
\text{$\mathop{\bu}\limits^{\sst i}\ra\mathop{\bu}\limits^{\sst
j}$}\\ \text{in $\Ga$}}}\bar\chi(\al_i,\al_j)
\raisebox{-6pt}{\begin{Large}$\displaystyle\biggr]$\end{Large}},
\end{split}
\label{an7eq1}
\e
where $F_n:(\C^\t)^n\ra\C$ are some functions to be determined,
and~$\C^\t=\C\sm\{0\}$.

Supposing that the $J^\al(\tau)^{\sst\Om}$ transform according to
\eq{an6eq48}, and neglecting issues to do with convergence of
infinite sums, in \cite{Joyc8} we show that there is an essentially
{\it unique\/} family of functions $F_n$ such that \eq{an7eq1}
yields $F^\al$ which are both continuous and holomorphic. These
$F_n$ are special functions of polylogarithm type, and are
holomorphic at $(z_1,\ldots,z_n)\in(\C^\t)^n$ except along the real
hypersurfaces $\arg(z_k)=\arg(z_{k+1})$ for $1\le k<n$, where $F_n$
is discontinuous.

The point is that as we cross a real hypersurface in $U$ or
$\Stab(\T)$ where $\arg(Z(\be))=\arg(Z(\ga))$ for some $\be,\ga\in
C(\A)$ or $K(\T)\sm\{0\}$, both the invariants
$J^{\al_i}(\mu)^{\sst\Om}$ and the functions
$F_n\bigl(Z(\al_1),\ldots,Z(\al_n)\bigr)$ in \eq{an7eq1} can jump
discontinuously. We arrange that these jumps exactly cancel out, so
that $F^\al$ remains continuous. Very surprisingly, it turns out
that with these choices of functions $F_n$ the $F^\al$ satisfy the
p.d.e.
\e
{\rm d}F^\al(Z)=-\!\!\!\!\!\!\sum_{\text{$\be,\ga\in C(\A)$ or
$K(\T)\sm\{0\}:\al=\be+\ga$}}\!\!\!\!\!\!\bar\chi(\be,\ga)
F^\be(Z)F^\ga(Z)\frac{{\rm d}(Z(\be))}{Z(\be)}\,,
\label{an7eq2}
\e
which can be interpreted as the flatness of an infinite-dimensional
connection over $U$ or $\Stab(\T)$. This still leaves many questions
unanswered, for example, mathematical questions on the convergence
of the infinite sums \eq{an7eq1} and \eq{an7eq2}, and physical
questions on the interpretation of \eq{an7eq1} and \eq{an7eq2} in
String Theory.

Here is our final problem.

\begin{prob} Redo {\rm\cite{Joyc4,Joyc5,Joyc6,Joyc7}} and this paper
using instead of Artin stacks some other kind of stack, which at
points $[X]$ includes the data $\Ext^i(X,X)$ for $i>1$ as part of
its structure. In this framework try to generalize results in
{\rm\cite{Joyc6,Joyc7}} and this paper which work when
$\Ext^i(X,Y)=0$ for all\/ $i>1$ and\/ $X,Y\in\A$ to the general
case, in particular the morphisms to explicit algebras
in~{\rm\cite[\S 6]{Joyc6}}.
\label{an7prob3}
\end{prob}

Some good candidates for the appropriate notion of stack are the
{\it D-stacks} of Toen and Vessozi \cite{ToVe} and the {\it
dg-stacks} of Ciocan-Fontanine and Kapranov \cite[\S 5]{CiKa}. Both
of these papers explain relevant ideas in {\it derived algebraic
geometry}.

The motivation behind this problem is that in constructing the
algebra morphism $\Phi^{\sst\La}$ in \cite[\S 6.2]{Joyc5}, if we
could include the groups $\Ext^i(X,X)$ for $i>1$ appropriately then
$\Phi^{\sst\La}$ would be an algebra morphism without assuming
$\Ext^i(X,Y)=0$ for all $X,Y$ and $i>1$. So we must work with a kind
of stack including the data $\Ext^i(X,X)$ for $i>1$, and new
`derived stack function' spaces~$\SF^{\rm der}(\fObj_\A)$.

There is a cost to this, though: the author expects Theorem
\ref{an5thm1} to fail for the spaces $\SF^{\rm der}(\fObj_\A)$,
since in the proof of Theorem \ref{an5thm4}, the 1-morphism
$\bs\si(\{1,\ldots,n\})$ is a 1-isomorphism of Artin stacks with a
substack of $\fObj_\A$, but it is probably {\it not\/} a
1-isomorphism of D-stacks or dg-stacks, as it may not preserve the
data $\Ext^i(X,X)$ for $i>1$. So our material on change of stability
condition will not work in $\SF^{\rm der}(\fObj_\A)$ in general.

However, in the situation of \S\ref{an64} when $\A=\coh(P)$ for
$K_P^{-1}$ nef, the author expects that at least equation
\eq{an5eq3} of Theorem \ref{an5thm1} should hold in some suitable
class of spaces $\SF^{\rm der}(\fObj_\A)$, as \S\ref{an64} works
precisely by forcing $\Ext^i(X,Y)=0$ for $i>1$ and the relevant
$X,Y$. Then we could interpret Theorem \ref{an6thm4} as saying that
we have an algebra morphism $\SF^{\rm der}(\fObj_\A)\ra
A(\A,\La,\chi)$, and the invariants $\Iss^\al(\ga)$ encode the
restriction of this to a derived version $\bar\H^{\rm to,der}_\ga$
of $\bHt_\ga$, which is independent of stability condition
$(\ga,G_2,\le)$. This would clear up a mystery about Theorem
\ref{an6thm4}, that is, why we have invariants with multiplicative
transformation laws but no underlying (Lie) algebra morphism.

\medskip

\noindent{\small\sc The Mathematical Institute, 24-29 St. Giles,
Oxford, OX1 3LB, U.K.}

\noindent{\small\sc E-mail: \tt joyce@maths.ox.ac.uk}

\end{document}